\documentclass[a4paper,11pt]{book}
\usepackage{amssymb}
\usepackage{amsmath} 
\usepackage[dvips]{graphicx}

\begin{document}

\def\V{\mathbb V}
\def\M{\mathbb M}
\def\R{\mathbb R}

\def\Aut{\mathop{\rm Aut}}
\newcommand{\ZZ}{\mathbb Z}
\newcommand{\F}{\mathbb F}
\newcommand{\pa}{\partial}
\newcommand{\ri}{\rightarrow}
\newcommand{\ov}{\overline}
\newcommand{\wh}{\widetilde}
\newcommand{\bg}{\Delta}
\newcommand{\st}{\stackrel}
\renewcommand{\bigtriangledown}{\nabla}
\renewcommand{\tilde}{\widetilde}


\begin{center}
{\bf \Large JORDAN STRUCTURES IN MATHEMATICS AND PHYSICS} \\[6pt]
Radu IORD\u ANESCU \footnote{The author was partially supported from the contract PN-II-ID-PCE 1188 517/2009.} \\[6pt]
Institute of Mathematics \\[6pt]
of the Romanian Academy \\[6pt]
P.O.Box 1-764 \\[6pt]
014700 Bucharest, Romania \\[6pt]
E-mail: Radu.Iordanescu@imar.ro
\end{center}
\thispagestyle{empty}




\centerline{\bf FOREWORD}

\vskip6pt

The aim of this paper is to offer an overview of the most important applications of Jordan structures inside mathematics and also to physics, up-dated references being included.

For a more detailed treatment of this topic see - especially - the recent book Iord\u anescu [364w], where sugestions for further developments are given through many open problems, comments and remarks pointed out throughout the text.

Nowadays, mathematics becomes more and more nonassociative (see \S 1 below), and my prediction is that in few years nonassociativity will govern mathematics and applied sciences.

{\bf MSC} 2010: 16T25, 17B60, 17C40, 17C50, 17C65, 17C90, 17D92, 35Q51, 35Q53, 44A12, 51A35, 51C05, 53C35, 81T05, 81T30, 92D10.

{\bf Keywords}: Jordan algebra, Jordan triple system, Jordan pair, JB-,\break JB$^{\ast}$-, JBW-, JBW$^{\ast}$-, JH$^{\ast}$-algebra, Ricatti equation, Riemann space, symmetric space, R-space, octonion plane, projective plane, Barbilian space, Tzitzeica equation, quantum group, B\" acklund-Darboux transformation, Hopf algebra, Yang-Baxter equation, KP equation, Sato Grassmann manifold, genetic algebra, random quadratic form.

\pagebreak

\centerline{\bf CONTENTS}

\vspace{48pt}

\noindent \S 1. Jordan structures \dotfill 2

\noindent \S 2. Algebraic varieties (or manifolds) defined by Jordan pairs \dotfill 11

\noindent \S 3. Jordan structures in analysis \dotfill 19

\noindent \S 4. Jordan structures in differential geometry \dotfill 39

\noindent \S 5. Jordan algebras in ring geometries \dotfill 59

\noindent \S 6. Jordan algebras in mathematical biology
and mathematical statistics \dotfill 66

\noindent \S 7. Jordan structures in physics \dotfill 77

\noindent Bibliography \dotfill 102

\vspace{1cm}

\centerline{\bf\Large \S 1. JORDAN STRUCTURES}

\vspace{48pt}

In modern mathematics, an important notion is that of {\it nonassociative structure}. This kind of structures is characterized by the fact the product of elements verifies a more general law than the associativity law.

\vskip6pt

There are two important classes of nonassociative structures: {\it Lie structures} (introduced in 1870 by the Norwegian mathematician Sophus Lie in his study of the groups of transformations) and {\it Jordan structures} (introduced in 1932-1933 by the German physicist Pasqual Jordan (1902-1980)\footnote{For details on the life and scientific career of Pascual Jordan see the impressive big paper: "{\it Pascual Jordan, his contributions to quantum mechanics and his legacy in contemporary local quantum physics}" by Bert SCHROER, arXiv:hep-th/0303241v2, 15 May 2003, 36 pp.} in his algebraic formulation of quantum mechanics [379a,b,c]). These two kinds of structures are interconnected, as it was remarked - for instance - by Kevin McCrimmon [480c, p. 622]: $\ldots$ "We are saying that if you open up a Lie algebra and look inside, 9 times out of 10 there is a Jordan algebra (of pair) which makes it work" $\ldots$, as well as by Efim Zelmanov [729a]: $\ldots$ "Lie algebra with finite grading may be by right included into the Jordan theory" $\ldots$

\vskip6pt

There exist three kinds of {\it Jordan structures}, namely, {\it algebras, triple systems}, and {\it pairs} (see the definitions below). Since the creation of Jordan algebras (at the beginning of '30) by Pasqual Jordan, an improvement of the mathematical foundation of quantum mechanics was made, but the problem was not definitively solved. Fifty years later (in '80), another tentative was made using a kind of mixed structures, namely {\it Jordan-Banach algebras} (see the definition in \S 3), but the problem of the mathematical foundation of quantum mechanics is still an open problem!

\vskip6pt

Anyway, is the meantime, Jordan structures have been intensively studied by mathematicians, and a big number of important results have been obtained. At the same time, an impressive variety of applications have been explored with several surprising connections. An explanation of this fact could be the following: at the beginning, mathematics was {\it associative} and {\it commutative}, then (after the invention of matrices) it became {\it associative} and {\it noncommutative}, and now (after the invention of nonassociative structures) it becomes {\it nonassociative} and {\it noncommutative}. The study of Jordan structures and their applications is at present a wide-ranging field of mathematical research. \footnote{see IORD\u ANESCU, R., {\it Romanian contributions to the study of Jordan structures and their applications}, Mitteilungen des Humboldt-Clubs Rum\" anien, No.8-9 (2004-2005), Bukarest, 29-35.}

\vskip6pt

I shall give here only algebraic basic definitions. For an {\bf intrinsic} treatment of Jordan structures the reader is referred to the excellent books by Braun and Koecher [131] (see also Koecher [408c]), Jacobson [371a, c], Loos [448e, g], Meyberg [483d], Zhevlakov, Slinko, Shestakov, Shirshov [733]. See also {\it Jordan algebras} -- Proceedings of the Oberwolfach Conference (August, 1992), Kaup, McCrimmon, Petersson (eds.), Walter de Gruyter--Berlin--New York, 1994, and {\it Proceedings of the International Conference on Jordan Structures} (Malaga, June, 1997), Castellon Serrano, Cuenca Mira, Fern\'andez L\'opez, Mart\'\i n Gonzalez (eds.), Malaga, 1999, as well as the WEB site:

\vskip6pt

\centerline{\footnotesize Jordan Theory preprints (\tt http://homepage.uibk.ac.at/$^\sim$c70202/jordan/index.html).}

\vskip6pt

Jordan algebras emerged in the early thirties with Jordan's papers\break [379a,~b,~c] on the algebraic formulation of quantum mechanics. The name ``Jordan algebras" was given by A.A. Albert in 1946.

\vspace{6pt}

{\bf Definition 1.} Let $J$ be a vector space over a field
${\bf F}$ with
characteristic different from two. Let $\varphi : J \times J \rightarrow J$
be an ${\bf F}$-bilinear map, denoted by $\varphi : (x,y) \rightarrow xy$,
satisfying the following conditions:
$$xy = yx \quad  \mbox{and} \quad  x^{2}(yx) = (x^{2}y)x
\quad\mbox{ for all } x,y \in J.$$
Then $J$ together with the product defined by $\varphi$ is called a
{\it linear Jordan algebra over} ${\bf F}$.

\vskip6pt

{\bf Comments.} The second condition came from quantum mechanics: observables in physics (temperature, pressure, etc.) satisfy it. This condition is less restrictive than the associativity and, from this fact it follows the impressive variety of applications of Jordan structures.

\vskip6pt

{\bf Example.} If $A$ is an associative algebra over ${\bf F}$, and we
define a new product $\varphi (x, y) : = \frac{1}{2}(x \cdot y + y \cdot x)$,
where the dot denotes the associative product of $A$, then
we obtain a Jordan
algebra. It is denoted by $A^{(+)}$.

\vspace{6pt}

{\bf Remark 1.} Jordan was the first who studied the properties
of the
product $xy$ from the above example in the case when ${\bf F}$ is the
field of reals. He proved a number of properties of this product, and
showed that these were all consequences of the two identities above.

\vspace{6pt}

{\bf Remark 2.} Zelmanov and others obtained important and
interesting
results for {\it infinite-dimensional} Jordan algebras. A good account of
Zelmanov work as well as of McCrimmon's extension to quadratic Jordan
algebras, can be found in [480d] (see also [364g]).

\vspace{6pt}

{\bf Notation.} On a Jordan algebra $J$ consider the left
multiplication {\it L} given by
$$L(x) y : = x y, \quad  x, y \in J.$$

\vspace{6pt}

{\bf Remark 3.} In general, $L(x y) \not = L(x) L(y)$ for $x,
y$
from a
Jordan algebra (which is not associative). This hold for, e.g.,
$J = A^{(+)}$ for $A$ (commutative or not).

\vspace{6pt}

{\bf Definition 2.} The map $P$ defined by $P(x) : = 2 L^{2}(x)
-
L(x^{2})$, $x \in J$, is called the {\it quadratic representation} of
${\it J}$.
When $J = A^{(+)}$ it assumes the form $P(x) y = x \cdot y \cdot x$.

\vskip6pt

{\bf Proposition 1.} {\it For any $x, y \in J$ the following
fundamental formula holds}: $$P(P(x) y) = P(x) P(y) P(x).$$

\vspace{4pt}

{\bf Remark 4.} For $P(x, y)$ given by $P(x, y) : = 2(L(x) L(y)
+ L(y) L(x) - L(xy))$ we have $P(x + y) = P(x) + P(x, y) + P(y)$,
$x,y \in J$.

\vspace{6pt}

{\bf Remark 5.} In general, $P(x, y) \not = P(x) P(y)$, $x, y \in J$
(as can easily be seen for $J = A^{(+)})$.

\vspace{6pt}

{\bf Proposition  2.} {\it Suppose that $J$ has a unit element~$e$
and let $x$
be an element of $J$. Then $P(x)$ is an automorphism of $J$ if and only if
$x^{2} = e$. If $P(x)$ is an automorphism of $J$, then it is involutive}.

\vspace{6pt}

{\bf Definition 3.} An element $x \in J$ is called {\it
invertible} if the
map $P(x)$ is bijective. In this case the inverse of $x$ is given by
$x^{-1} : = (P(x))^{-1} x.$ (In case $J = A^{(+)}$, this is the usual inverse
in the associative algebra $A$.)

\vspace{6pt}

{\bf Remark 6.} We have $(P(x))^{-1} = P(x^{-1})$, $x \in J$.

\vspace{6pt}

{\bf Remark 7.} An element $x$ is invertible with the inverse
$y$
if and only $x y = e$, $x^{2} y = x$.

\vspace{6pt}

{\bf Definition 4.} Let $f$ be an element of $J$. Define a new
product on the vector
space $J$ by
$$\{x f y\}: = x(y f) + y(x f) - (x y) f.$$
The vector space $J$ together with this product is called the
{\it mutation $($homotope$)$ of $J$ with respect to} $f$ and is denoted by
$J_{f}$.

\vspace{6pt}

{\bf Proposition 3.} {\it Any mutation $J_{f}$ of $J$, $f \in J$ is a
Jordan algebra and its quadratic representation $P_{f}$ is given by}
$P_{f} (x) = P(x) P(f)$.

\vspace{6pt}

{\bf Proposition 4.} {\it The algebra $J_{f}$, $f \in J$, has a unit
element if
and only if $f$ is invertible in $J$; in this case the unit element of
$J_{f}$ is $f^{-1}$. In this situation we call $J_{f}$
the $f$-homotope of} $J$.

\vspace{6pt}

{\bf Remark 8.} If $f$ is invertible in $J$ then the set of
invertible elements of $J$ coincides with the set of invertible elements of $J_{f}$.

\vspace{6pt}

{\bf Note.} From this point on, many of the results require that $J$ be
finite-dimensional.

\vspace{6pt}

{\bf Notation.} Denote by
$${\rm Invol}(J) := \{ w \mid w \in J,\, w^{2} = e \}, \quad
{\rm Idemp}(J) : =  \{ c \mid c \in J, \, c^{2} = c \}$$
the set of involutive, respectively idempotent (zero included), elements of
$J$. Here $J$ is supposed to contain a unit element~$e$. Note that
${\rm Invol} (J) = \{w \in J$, $w^{-1} = w\}.$

\vspace{6pt}

{\bf Remark 9.} The map ${\rm Idemp}(J) \rightarrow {\rm Invol}(J)$
given by $c \rightarrow 2 c - e$ is a bijection.

\vspace{6pt}

For an element $c$ of ${\rm Idemp}(J)$ we have
$$L(c) (L(c) - {\rm Id}) (2 L(c) - {\rm Id}) = 0.$$
This leads to the {\it Peirce decomposition of $J$ with respect to the
idempotent}~$c$:
$$J = J_0(c) \oplus J_{1/2} (c) \oplus J_{1} (c),$$
where $J_{\alpha} (c) : = \{x \mid x \in J, \, c x = \alpha x\}$, for
$\alpha = 0, 1/2 ,1$.

\vspace{6pt}

{\bf Theorem 5.} $J_0(c)$ {\it and $J_{1} (c)$ are
subalgebras of $J$, and we have
$$J_{0} (c) J_{1} (c) = \{ 0\}, \quad  J_{\nu} (c) J_{1/2} (c)
\subset J_{1/2} (c), \quad  \mbox{for }  \nu = 0,1,$$
and}
$$J_{1/2}(c) J_{1/2} (c) \subset J_0 (c) \oplus J_{1} (c).$$

\vspace{6pt}

{\bf Definition 5.} Let $c$ be an indempotent of $J$,
$c \not = e$.
Then the map $P(2 c - e)$, which by virtue of Proposition 1, is an
automorphism of $J$, is called the {\it Peirce reflection with respect to
the idempotent} $c$ of $J$.

\vspace{6pt}

{\bf Notation.} ${\rm Idemp}_1 (J):=\{c \in {\rm Idemp} (J),
\, \dim J_1(c)=1\}$.

\vspace{6pt}

{\bf Definition 6.} The dimension of $J_1(c)$ is called the
{\it rank} of the idempotent~$c$.

\vspace{6pt}

{\bf Definition 7.} An idempotent $c$ of $J$ is called {\it
primitive} if it cannot be decomposed as sum $c_1+c_2$ of two
orthogonal (i.e., $c_1c_2=0$) idempotents $c_1$ and $c_2$,
$c_i\neq 0$ $(i=1,2)$.

\vspace{6pt}

{\bf Remark 10.} Every element of ${\rm Idemp}_1 (J)$ is
primitive. The converse is not true in general.

\vspace{6pt}

{\bf Definition 8.} A Jordan algebra over the real numbers is called
{\it formally real} if, for any two of its elements
$x$ and $y$, $x^{2} + y^{2} = 0$ implies that $x = y = 0$.

\vspace{6pt}

{\bf Proposition 6.} {\it A primitive idempotent of a formally real
finite-dimen\-sional Jordan algebra is of rank one.}

\vspace{6pt}

{\bf Proposition 7.} {\it The set $\Gamma(J)$ of bijective linear
maps $W$ on $J$ for which there exists a bijective map $W^{*}$ on
$J$ such that $P(Wx) = W P(x) W^{*}$ for all $x \in J$ is a
linear algebraic group.}

\vspace{6pt}

{\bf Note.} The notation $W^{\ast}$ is justified by the fact that if $J$
is real semi-simple and $\lambda$ denote the {\it trace form} on
$J$ (i.e., $\lambda (x, y) : = {\rm Tr}\, L(xy),\, x,y \in J)$ then for
$W \in \Gamma (J)$, $W^{*}$  coincides with the adjoint of $W$ with respect
to $\lambda$.

\vspace{6pt}

{\bf Definition 9} (Koecher [408b]). The (linear algebraic)
group
$\Gamma(J)$ from Proposition 7 is called the {\it structure
group of} $J$.

\vspace{6pt}

{\bf Remark 11.} The fundamental formula (see Proposition 1)
implies that $P(x) \in \Gamma(J)$ where $x$ is an invertible element of $J$. Also, every
automorphism of $J$ belongs to $\Gamma(J)$. Indeed, the automorphism group
${\rm Aut}(J)$ is just the set of elements $W \in \Gamma(J)$ fixing the unit element
$e$ of $J$, $We = e$.

\vspace{6pt}

Formally real Jordan algebras have been clasified (in the
finite-dimen\-sional case) by Jordan, von Neuman and Wigner [380]:

\vspace{6pt}

{\bf Theorem 8.} {\it Every formally real finite-dimensional Jordan
algebra
is a direct sum of the following algebras
$$H_{p}({\bf R})^{(+)}, \quad H_{p}({\bf C})^{(+)},
\quad H_{p}({\bf H})^{(+)}, \quad
H_{3}({\bf O})^{(+)}, \quad J(B),$$
here $H_{p}({\bf F)}^{(+)}$ denotes the algebra of Hermitian
$(p \times p)$-matrices with entries in ${\bf F}$ $({\bf F} = {\bf R, C, H}$
or ${\bf O})$\footnote{Throughout this paper ${\bf R, C, H}$,
and ${\bf O}$ denotes the set of reals, complex numbers, quaternions, and octonions
(Cayley numbers), respectively.}, the multiplication in
$H_{p} ({\bf F})^{(+)}$ is given by
$xy : = \frac{1}{2} (x \cdot y + y \cdot x)$
$(x \cdot y  \mbox{ denotes the usual matrix product})$, and
$J(B) = {\bf R} 1 \oplus X$, where $B$ is a real-valued symmetric bilinear
positive definite form on the real vector space $X$, equipped with the product}
$$(\lambda, u)(\mu, v) : = (\lambda \mu + B(u,v), \lambda v + \mu u).$$

\vspace{4pt}

{\bf Definition 10.} $A$ Jordan algebra $J$ is called {\it special}
if it is isomorphic
to a (Jordan) subalgebra of some $A^{(+)}$ for $A$ associative.

\vspace{6pt}

{\bf Remark 12.} The first three and the fourth algebras in
Theorem 8
are special, while $H_{3} ({\bf O})^{(+)}$ is not, that is why it is called
{\it exceptional}.

\vspace{6pt}

{\bf Proposition 9.} {\it A Jordan algebra is formally real if and
only if its trace form is positive definite}.

\vspace{6pt}

{\bf Notation.} Suppose that $J$ has a unit element $e$. Then we set
$x^0 : = e$,
${\rm exp}\, x : = \sum\limits_{n \geq 0} \frac{x^{n}}{n!}$,
and exp $J : =
\{{\rm exp}\, x \mid x \in J\}$.

\vskip6pt

{\bf Proposition 10.} {\it If $J$ is a formally  real Jordan
algebra, then
it possesses a unit element and}
${\rm exp} \,J = \{ x^{2} \mid x  \mbox{ invertible in } J\}$.

\vspace{6pt}

{\bf Theorem 11.} {\it Suppose that $J$ is a formally real Jordan
algebra
endowed with the natural topology of ${\bf R}^{n}$. Then the identity
component of the set of invertible elements of $J$ coincides with $\exp J$
and the map $x \rightarrow \exp x$ is bijective}.

\vspace{6pt}

If in a given Jordan algebra we define a triple product by
$$\{x y z \}: = (x y) z + (z y) x - y (x z),$$
then it satisfies the following two identities:

\vspace{6pt}

({\rm JTS 1}) $\{x y z\} = \{z y x\}$;

({\rm JTS 2}) $\{uv \{x y z\}\} =   \{\{u v x\}y z\} -
\{x \{v u y\} z\} + \{x y \{u v z\}\}$.

\vspace{6pt}

In general, a module with a trilinear composition $\{x y z\}$ satisfying
(JTS~1) and (JTS~2) is called a {\it Jordan triple system.} Meyberg [483a,~b] has used Jordan triple systems in the study of the extension of
Koecher's construction [408e] of a Lie algebra from a given Jordan algebra.

\vspace{6pt}

{\bf Theorem 12} (Meyberg [483a,~b]). {\it If $T$ is a Jordan triple
system and $a$
an element of $T$, then $T$ together with the product
$(x, y) \rightarrow \frac{1}{2} \{ x a y \}$ becomes a Jordan algebra,
denoted
$T_{a}$. Conversely, a Jordan algebra induces a Jordan triple system in the same
vector space by setting} $\{x y z\} : = P(x, z)y$.

\vskip6pt

There exist generalizations of Jordan triple systems, studied - especially - by Kamiya and his collaborators (see [383a,b], [385c,d], [233a, b]). I take this opportunity to point out here the existence of a new class off nonassociative algebras with involution (including the class of structurable algebras) recently defined by Kamiya and Mondoc [384].

\vskip6pt

In 1969, Meyberg introduced [483a] the "verbundene Paare" (connected pairs),
which correspond in Loos terminology to the linear Jordan pairs; see the
definition below. Such connected pairs first arose in Koecher's work on Lie
algebras.

\vspace{6pt}

{\bf Definition 11} (Loos [448g]). Let ${\bf K}$ be a unital commutative ring
such that 2 is
invertible in ${\bf K}$. Assume all ${\bf K}$-modules to be unital and to possess
no 3-torsion (i.e. no nonzero elements $x$ such that $3x = 0)$. A pair
$V = (V^{+}, V^{-})$ of ${\bf K}$-modules endowed with two trilinear maps
$V^{\sigma} \times V^{- \sigma} \times V^{\sigma} \rightarrow V^{\sigma}$,
written as $(x, y, z) \rightarrow \{x y z\}_{\sigma}$,
$\sigma = \pm$, satisfying the identities
\begin{eqnarray*}
\{x y z\}_{\sigma}\!\!\! & = &\!\!\!\{z y x\}_{\sigma}, \\
\{x y \{u v z\}_{\sigma}\}_{\sigma} - \{u v\{x y z\}_{\sigma}\}_{\sigma}
\!\!\!& = &\!\!\!
\{ \{x y u\}_{\sigma} v z\}_{\sigma} - \{ u \{y x v\}_{- \sigma} z \}_{\sigma }
\end{eqnarray*}
for $\sigma = \pm$, is called a {\it $($linear$)$ Jordan pair} over {\bf K}.

\vspace{6pt}

{\bf Remark 13.} Jordan algebras can be regarded as a
generalization of
{\it symmetric} matrices, while the (linear) Jordan triple structures
(systems or pairs) can be regarded as a generalization of
{\it rectangular} matrices.

\vskip6pt

From the recent papers focussed on Lie algebras, but from a Jordan point of view, I would like to mention these written by Fernandez Lopez and his collaborators (see 
[257 a,b], [258], [80]), while for these focussed on Lie superalgebras, see the paper by Cunha \& Elduque [187].

\vskip6pt

{\bf Remark 14.} Jordan triple systems are equivalent with Jordan
pairs with involution.

\vspace{6pt}

McCrimmon [480a] extended the theory of Jordan algebras to the case of an arbitrary commutative unital underlying ring. So, instead of considering a Jordan algebra as a vector space with a nonassociative bilinear composition satisfying certain identities, McCrimmon considered it as a module over a ring together with a quadratic representation satisfying a number of conditions. In this way, one gets unital {\it quadratic} Jordan algebras, which for the case that the underlying ring is a field of characteristic different from two (or any underlying ring containing 1/2) turn out to be the well-known Jordan algebras considered from a different point of view.

By analogy with McCrimmon's concept of quadratic Jordan algebras, Meyberg
[483d] defined {\it quadratic} Jordan triple systems, while the notion of
{\it quadratic} Jordan pairs was introduced by Loos [448f].

\vskip6pt

The Tits-Kantor-Koecher construction and its generalizations to Jordan pairs,
structurable algebras, Kantor pairs, etc., are the basis of the 
recent application of Jordan structure in the theory of Lie algebras, mainly the 
so-called root-graded Lie algebras. For details on this topic, we refer the reader to the
recent excellent survey paper [512j] by Neher. Let us mention here
from it the recent (or very recent) papers by Allison, Azam, Berman, Gao, and Pianzola
[16], Allison, Benkart, and Gao [17a,~b], Allison and Faulkner [18b], Allison and Gao [19],
Benkart [79], Benkart and Smirnov [83], Benkart and Zelmanov [84],
Berman, Gao, Krylyuk, and Neher [94], Neher [512f], and Yoshii [722d]. An important
paper in this context is also the paper [722c] by Yoshii. In this paper the author determines
the coordinate algebra of extended affine Lie algebras of type $A_1$. It turns
out that such an algebra is a unital $\mathbb{Z}^n$-graded Jordan algebra of 
a certain type, called a {\it Jordan torus}. The author also gives the classification
of Jordan tori.

For Jordan structures in the {\it super} setting (i.e., Jordan superpairs, quadra\-tic
Jordan superpairs, etc.), we refer the reader to Garc\'\i a \& Neher [282a], and Neher [512i],
while for the applications of Jordan techniques in the theory of Lie algebras, we refer
the reader to the recent paper by Fern\'andez L\'opez, Garc\'\i a, and G\'omez Lozano [257b],
as well as to the papers referred inside [257b].

Very recently, Vel\'asquez and Felipe [688a,~b] introduced and studied a new algebraic
structure, called {\it quasi-Jordan algebra}.

\vskip6pt

The mixed Jordan and Lie structures become more and more important nowadays.
{\it Jordan-Lie algebras} were defined in 1984 by Emch in his book [236b],
although the concept seems to have appeared first in the paper [306] by 
Grgin \& Petersen. {\it Lie-Jordan algebras} were defined in 2001 by
Grishkov and Shestakov in [308].

\vskip6pt

{\bf Note.} Every Jordan-Lie algebra gives rise to a Lie-Jordan algebra, but the converse is
false.

\vskip6pt

In September 2009, Makhlouf defined {\it Hom-Jordan algebras}, and in\break January 2010 he gave a survey on nonassociative Hom-algebras and\break Hom-superalgebras, see Makhlouf [463a,b], where his papers from 2007 and 2008 in co-operation with Silvestrov are also refered, see Makhlouf \& Silvetrov [464a,b,c,d]. I take this opportunity to mention here also {\it Hom-Nambu-Lie algebras}, induced by Hom-Lie algebras, see Arnlind, Makhlouf and Silvestrov [44], as well as Arnlind [43].

\vskip6pt

A lot of contributions are now devoted to {\it Moufang loops} (see Moufang [500b]), which extend groups, and the analog of Lie algebra for local Moufang loops are {\it Malcev algebras} (see Malcev [465]). There exists today a rich bibliography on Moufang loops, from which I would refer the reader to Paal [528], who wrote many papers on Moufang symmetry. Very recent related contributions are the papers by Benkart, Madariaga and Perez-Izquierdo [81], and Madariaga \& Perez-Izquierdo [457].

For a big paper on algebras, hyper-algebras, nonassociative bialgebras and loops, see Perez-Izquierdo [537b].

I would like to mention here also the papers [254a, b, c] by Faybusovich, as well as a very recent paper by Elduque [232].

\vskip6pt

It is worth to be remarked here that the study of nonassociative structures is a tradition in the Romanian universitary center of Ia\c si. I would like to give here only two examples: a paper by Climescu [180] mentioned in the famous book [131] by Braun \& Koecher, and - recently - Burdujan [149a,b] studied the deviation from associativity. As Burdujan remarked\footnote{Private communication in July 2010.}, the construction from his paper [149a] could be applied also to Jordan algebras, and - also - some his LT-algebras are Jordan algebras (see [149c]).

\vskip6pt

Let us mention here some {\bf recent} Ph.D. Theses in the field of Jordan structures and their
applications defended at universities from USA: {\it Orbits of exceptional groups and Jordan
systems} by Sergei Krutelevich (Yale University), {\it The ring of fractions of a quadratic Jordan
algebra} by James Bowling (Virginia University), {\it Affine remoteness planes} by Karen
Klintworth (Virginia University), {\it Centroids of quadratic Jordan superalgebras} by
Pamela Richardson (Virginia University), {\it Laguerre functions associated to Euclidean Jordan
algebras} by Michael Aristidon (Louisiana State University), {\it Derivations of $8$ simple
Jordan superalgebras} by Michael Smith (Virginia University).

\vskip6pt

From the {\bf recent} or {\bf very recent} papers of interest here, I would like to mention [49] by Ashihara \& Miyamoto, [136] by Bremner \& Peresi, [494] by Montaner, [453] by MacDonald, [706] by Wilson, and [33] by Anquela, Cort\'ez, and McCrimmon.

\vskip6pt

I was {\bf very recently} informed by Prof. Dr. Harald UPMEIER (Philips-Universit\" at Marburg, Germany) through e-mail (May 20, 2011) on his progress made in a {\bf new field}, concerning the conformal compactification of Jordan algebras and applications to harmonic analysis. It seems that the conformal geometry of Jordan algebras can explain deep facts about the representation theory of semisimple Lie groups and also about the geometry of non-convex cones and the higher-dimensional Radon transform.

\vspace{1cm}

\centerline{\bf \Large \S 2. ALGEBRAIC VARIETIES (OR MANIFOLDS)}
\centerline{\bf\Large DEFINED BY JORDAN PAIRS\footnote{The topics presented in this paragraph were only mentioned in the book Iord\u anescu [364w].}}

\vspace{48pt}

Loos [448i] showed that every (quadratic) Jordan pair defines an affine algebraic group, the projective group of the Jordan pair. On the other hand, he remarked\footnote{LOOS, O., {\it Homogeneous algebraic varieties defined by Jordan pairs}, Monatsh. Math. {\bf 86} (1978/79), {\it 2}, 107-129.} that every Jordan pair defines an algebraic variety, related with the projective group of the Jordan pair in a natural manner.

\vskip6pt

Let $V = (V^{+}, V^{-})$ be a finite-dimensional Jordan pair over an algebraically closed field ${\bf F}$ and let $X = X(V)$ be the quotient of $V^{+}\times V^{-}$ by the equivalence relation $(x,y) \sim (x', y')$ if and only if $(x, y-y')$ is quasi-invertible and $x' = x^{y-y'}$. (As usual, $x^y$ denotes the quasi-inverse in the Jordan pair $V$). Loss [448h] proved that $X$ is a quasi-projective variety containing $V^{+}$ as an open dense subset, and that $X$ is projective if $V$ is semisimple. Moreover, Loos showed that under the projective group of $V$ the space $X$ is homogeneous in a natural way, and that this projective group is isomorphic to the group of automorphisms of $X$ if $V$ is semisimple. This is essentially due to Chow [165] in four classical cases. For a geometric characterization of the projective group see Faulkner [248a], Freudenthal [272a], Springer [632a], as well as Chow [165], Faulkner\footnote{FAULKNER, J.R., {\it A geometry for E$_{7}$}, Trans. Amer. Math. Soc. {\bf 167} (1972), 49-58.}, Hua [352b].

\vskip6pt

For each of the Jordan pairs from the classification of simple finite-dimensional Jordan pairs over an algebraically closed field ${\bf F}$ given by Loos [448g, p. 201] the variety $X$ is as follows:\footnote{These Jordan pairs arise in a fairly obvious way from separable Jordan pairs over $\mathbb{Z}$ by reduction {\it modulo} char ${\bf F}$ and extending the prime field to ${\bf F}$. Thus, we are dealing with the fibres of a smooth projective $\mathbb{Z}$-scheme.}

Type I$_{p,q}$. $V^{+} = V^{-} = M_{p,q}({\bf F})$, $(p \times q)$-matrices over ${\bf F}(p \leq q)$, with $Q(x)y := xy'x$, where $y'$ denotes the transpose of $y$. In this case $X$ is isomorphic to the Grassmannian $G_p({\bf F}^{p+q})$ of $p$-dimensional subspaces of ${\bf F}^{p+q}$.

Type II$_n$. $V^{+} = V^{-} = A_n({\bf F})$, alternating $(n \times n)$-matrices with $Q(x)y := xy'x$. In this case $X$ is isomorphic to the subvariety of $G_n({\bf F}^{2n})$ consisting of all totally isotropic subspaces of ${\bf F}^{2n}$ of fixed parity with respect to the quadratic form $q(x_1, \ldots, x_n) := \mathop{\sum}\limits_{i=1}^{n}x_ix_{n+i}$.

Type III$_n$. $V^{+} = V^{-} = H_n({\bf F})$, symmetric $(n \times n)$-matrices with $Q(x)y:= xy'x$. In this case $X$ is isomorphic to the subvariety of $G_n({\bf F}^{2n})$ consisting of all maximal isotropic subspaces of ${\bf F}^{2n}$ with respect to alternating form $\alpha(x,y) := \mathop{\sum}\limits_{i=1}^{n}(x_iy_{n+i} - x_{n+i}y_i)$.

Type IV$_n$. $V^{+} = V^{-} = {\bf F}^n$ with $Q(x)y := q(x,y)x - q(x)y$, where $q$ is the standard quadratic form of ${\bf F}^n$, given by
$$
q(x_1, \ldots, x_{2m}) := \mathop{\sum}\limits_{i=1}^{m}x_ix_{m+i} \;\; \mbox{if} \;\; n = 2m,
$$
$$
q(x_0, \ldots, x_{2m}) := x_0^2 + \mathop{\sum}\limits_{i=1}^{m}x_ix_{m+1} \;\; \mbox{if} \;\; n = 2m + 1.
$$
In this case $X$ is isomorphic to the quadric of all isotropic lines through the origin in ${\bf F}^{2n}$.

Type V. $V^{+} = V^{-} = M_{1,2}({\bf O})$, $(1\times 2)$-matrices over the octonion (Cayley) algebra ${\bf O}$ over ${\bf F}$. Here $X$ is isomorphic to the projective octonion plane $\mathcal{P}({\bf O})$ defined by the exceptional Jordan algebra $H_3({\bf O})^{+}$ (see Faulkner [248a]).

Type VI. $V^{+} = V^{-} = H_3({\bf O})^{+}$, Hermitian $(3 \times 3)$-matrices over the octonion (Cayley) algebra ${\bf O}$ over ${\bf F}$. In this case $X$ is isomorphic to the space of lines ${\bf F}u$, where $u$ is an element of rank one in the 56-dimensional space of all matrices 
$
\left(
\begin{array}{cc}
{\alpha} & a \\ [5pt]
b & {\beta}
\end{array}\right),
$ $\alpha, \beta \in {\bf F}$, $a,b \in H_3({\bf O}, {\bf F})$.

\vskip6pt

In 1984 the fundamental paper by Petersson appeared\footnote{PETERSSON, H.-P., {\it Generic reducing fields for Jordan pairs}, Trans. Amer. Math. Soc. {\bf 285} (1984), {\it 2}, 825-843.}. It provides a uniform framework for generic splitting fields of associative algebras studied by Amitsur [24], Heuser [342], Kovacs [417], Saltman [588] and Roquette [567a,b], and generic zero fields of quadratic forms investigated by Knebusch [403]. This is achieved by first observing that both, central simple associative algebras over a field ${\bf K}$ of arbitrary characteristic and absolutely nondegenerate quadratic forms over ${\bf K}$, are examples of absolutely simple Jordan pairs over ${\bf K}$, and then constructing generic reducing fields of Jordan pairs in a classification-free treatment.

\vskip6pt

Let $V = (V^{+}, V^{-})$ be an absolutely simple Jordan pair over ${\bf K}$. A nonzero idempotent $c$ of $V$ is called {\it reduced}, if the Peirce space $V_2(c)^{\varepsilon} = {\bf K}c^{\varepsilon}$, $\varepsilon = \pm$, and $V$ itself is called {\it reduced} if it contains a reduced idempotent.

\vskip6pt

{\bf Remark 1.} For Jordan structures "split" is a more restrictive concept than "reduced", however "split" = "reduced" if $V$ is the Jordan pair associated to a central simple associative algebra.

\vskip6pt

A field extension $\mathcal{K}/{\bf K}$ is called a {\it reducing field} of $V$ in case the extended Jordan pair $V_{{\bf K}}$ is reduced, and such an extension is called a {\it generic reducing field} if an arbitrary field extension $\mathcal{L}/{\bf K}$ reduces $V$ if and only if there is a ${\bf K}$-place from $\mathcal{K}$ to $\mathcal{L}$.

Petersson has constructed two, in general not ${\bf K}$-isomorphic, generic reducing fields. Both are ${\bf K}$-rational function fields of irreducible projective ${\bf K}$-variables associated to $V$. A detailed discussion of examples, most natably Brauer-Severi varieties, and with an application to exceptional simple Jordan algebras arising from the first Tits construction is also given by Petersson.

In 1985, Jacobson [371d] considered two classes of projective varieties, namely norm hypersurfaces and varieties of reduced elements defined by\break finite-dimensional central simple Jordan algebras. Jacobson's paper [371d] overlaps substantially with the paper by Petersson. However, the methods and points of view are different, supplementing each other.

\vskip6pt

In Ch.II of his Ph.D. Thesis [701], Watson defined a class of manifolds, called {\it local Jordan manifolds} (LJM), the ground algebraic structure being (linear) Jordan pairs. For the sake of simplicity we shall omit the word "linear". As one can see from the examples given below, many familiar manifolds are local Jordan manifolds.

\vskip6pt

{\bf Note.} All Jordan pairs considered in the following will be {\it finite-dimensional} vector spaces over ${\bf R}$ and the vector spaces $V^{+}$ and $V^{-}$ of a Jordan pair $V$ will be assumed to have the {\it ordinary real topology}.

\vskip6pt

{\bf Comment.} Very recently (in 2009), Stach\' o \& Werner [636b] defined the notion of {\it Jordan manifold} in a very different setting (see the next \S 3).

\vskip6pt

{\bf Definition 1.} A {\it Jordan chart} on a topological space $M$ is a pair $(\Phi, V)$ consisting of a Jordan pair $V = (V^{+}, V^{-})$ and a map $\Phi : U \to V^{+}$, where $U \neq \emptyset$ is open in $M$, and $\Phi$ is a homeomorphism of $U$ onto an open subset of $V^{+}$.

\vskip6pt

{\bf Definition 2.} A {\it Jordan atlas} on a topological space $M$ is a set $\mathcal{A}$ of Jordan charts on $M$ satisfying the following two conditions:

(a) the domain of the Jordan charts in $\mathcal{A}$ conver $M$;

(b) if $\Phi_1 : U_1 \to V_1^{+}$ and $\Phi_2 : U_2 \to V_2^{+}$ (where $U_i$ is open in $M$ and $V_i$ is a Jordan pair) are elements of $\mathcal{A}$ and if $U_1 \cap U_2 \neq \emptyset$, then in a neighbourhood of each point in $\Phi_1 (U_1 \cap U_2)$, $\Phi_2\Phi_1^{-1} : \Phi_1(U_1 \cap U_2) \subset V_1^{+} \to V_2^{+}$ is the restriction of a map from $LF_{+}(V_1, V_2)$.

\vskip6pt

Let $\Phi_i : U_i \to V_i^{+}$, $i=1,2$, be two Jordan charts in a Jordan atlas $\mathcal{A}$ on $M$. Suppose that $U_1 \cap U_2 \neq \emptyset$. Then for every $u \in U_1 \cap U_2$, in a neighbourhood of $\Phi_1(u), \Phi_2\Phi_1^{-1}$ is the restriction of a unique linear fractional map, called the {\it coordinate transition map} and denoted by $C_{\Phi_1\Phi_2}^{u}$.

\vskip6pt

{\bf Definition 3.} Let $\mathcal{A}$ be a Jordan atlas on a topological space $M$. Then a set $\mathcal{B}$ of Jordan charts on $M$ is said to be {\it compatible} with $\mathcal{A}$ if $\mathcal{A} \cup \mathcal{B}$ is a Jordan atlas on $M$.

Now let $\hat{\mathcal{A}}$ be the class of all Jordan charts on $M$ compatible with $\mathcal{A}$. It can easily be checked that if $\mathcal{A}_1$ is a set, if $\mathcal{A}_1 \subset \hat{\mathcal{A}}$, and if the domains of the charts in $\mathcal{A}_1$ cover $M$, then $\mathcal{A}_1$ is a Jordan atlas on $M$, and $\hat{\mathcal{A}}_1 = \hat{\mathcal{A}}$. The class $\hat{\mathcal{A}}$ is called the {\it Jordan structure} of $M$.

\vskip6pt

{\bf Definition 4.} A {\it local Jordan manifold} is a pair $(M, \hat{\mathcal{A}})$ where $M$ is a topological Hausdorff space and $\hat{\mathcal{A}}$ is a Jordan structure on $M$ determined by a Jordan atlas $\mathcal{A}$ on $M$.

\vskip6pt

{\bf Remark 2.} Note that $M$ is not assumed to be connected, nor must all components of $M$ have the same dimension. However, every component of $M$ is a real analytic manifold.

\vskip6pt

{\bf Definition 5.} Let $(M_1, \hat{\mathcal{A}}_1)$ and $(M_2, \hat{\mathcal{A}}_2)$ be local Jordan manifolds. A {\it morphism} from $M_1$ to $M_2$ is a continuous map $f : M_1 \to M_2$ such that, for every pair of charts $\Phi_1 : U_1 \to V_1^{+}$ in $\mathcal{A}_1$ and $\Phi_2 : U_2 \to V_2^{+}$ in $\hat{\mathcal{A}}_2$ with $U_1 \cap f^{-1}(U_2) \neq \emptyset$, in a neighbourhood of each point in its domain the map $\Phi_2 f\Phi_1^{-1} : \Phi_1(U_1 \cap f^{-1}(U_2)) \to V_2^{+}$ is the restriction of a map from $LF_{+}(V_1, V_2)$.

\vskip6pt

{\bf Remark 3.} The composite of two morphisms (as maps) between local Jordan manifolds is again a morphism, and in this way the local Jordan manifolds and the morphisms between them form a category, denoted by LJM. Every morphism in LJM between connected local Jordan manifolds is a morphism in the category of real analytic manifolds.

\vskip6pt

{\bf Definition 6.} Let $V$ be a Jordan pair and let $M$ be a topological space. A Jordan chart on $M$ of the form $(\Phi, V)$ is called a {\it Jordan $V$-chart}. A {\it Jordan $V$-atlas} on $M$ is a Jordan atlas $\mathcal{A}$ on $M$ consisting of Jordan $V$-charts.

\vskip6pt

{\bf Definition 7.} Let $V$ be a Jordan pair. A {\it loocal Jordan $V$-manifold} is a topological Hausdorff space $M$ with a maximal Jordan $V$-atlas $\mathcal{A}$.

\vskip6pt

The standard local Jordan $V$-manifold is $V^{+}$.

\vskip6pt

{\bf Remark 4.} Every local Jordan $V$-manifold is a local Jordan manifold.

\vskip6pt

{\bf Remark 5.} Let $(M, \hat{\mathcal{A}})$ be a connected local Jordan manifold. One can see that on $M$ we can specify a Jordan $V$-atlas contained in $\hat{\mathcal{A}}$, where $V$ is any Jordan pair occurring in a Jordan chart in $\mathcal{A}$. The various Jordan pairs $V$ generated by a Jordan atlas $\mathcal{A}$ on $M$ need not be isomorphic.

\vskip6pt

{\bf Remark 6.} A morphism $f : M_1 \to M_2$ between two local Jordan manifolds $M_1$ and $M_2$ cam be a diffeomorphism but can never be an isomorphism in LJM.

\vskip6pt

{\bf Theorem 1.} Let $f : M_1 \to M_2$ be a morphism between two local manifolds $M_1$ and $M_2$. Suppose that $M_1$ is a local Jordan $V$-manifold, where Rad$V^{+} = 0$. Then $f$ is an isomorphism in LJM if and only if it is a homeomorphism.

\vskip6pt

Let $(M_1, \hat{\mathcal{A}}_1)$ and $(M_2, \hat{\mathcal{A}}_2)$ be local Jordan manifolds. Since any\break $(f_1, f_2) \in LF_{+}(V_1, W_1) \times LF_{+}(V_2, W_2)$ can be identified in a natural way with an element of $LF_{+}(V_1 \oplus V_2, W_1 \oplus W_2)$, it is clear that $M_1 \times M_2$ possesses a natural Jordan strcture. The ({\it direct product}) local Jordan manifold $M_1 \times M_2$ has the expected universal property requiring that whenever $M$ is a local Jordan manifold and $p_i : M \to M_i$, $i = 1,2$, are morphisms, there exists a unique morphism $p : M \to M_1 \times M_2$ such that $\pi_ip = p_i$, where $\pi_i : M_1\times M_2 \to M_i$ is a projection morphism, for $i = 1,2$. Moreover, any other local Jordan manifold $M'$ having this property (for a pair of morphism $\pi'_i : M' \to M_i, i = 1,2$) is isomorphic in LJM to $M_1 \times M_2$.

\vskip6pt

{\bf Remark 7.} If $M_i$ is a local Jordan $V_i$-manifold, $i = 1,2$, then $M_1 \times M_2$ is a local Jordan $V_1 \oplus V_2$-manifold.

\vskip6pt

{\bf Definition 7.} Let $(M, \hat{\mathcal{A}})$ be a local Jordan manifold. A nonempty subset $N$ of $M$ is called a {\it local Jordan submanifold} of $M$ if there exists a set $\mathcal{A}_1 \subset \hat{\mathcal{A}}$ such that for every $n \in N$ there exists a Jordan chart $\Phi : U \to V^{+}$ in $\mathcal{A}_1$ with $n \in U$ and such that $\Phi(U \cap N)$ is an open subset of a vector subspace $W^{+}$ of $V^{+}$.

\vskip6pt

Open subsets of local Jordan manifolds and affine subspaces of $V^{+}$, where $V = (V^{+}, V^{-})$ is a Jordan pair, are simple examples of local Jordan submanifolds.

\vskip6pt

{\bf Proposition 2.} Let $N_1$ and $N_2$ be local Jordan submanifolds of $M_1$ and $M_2$, respectively, and let $f : M_1 \to M_2$ be a morphism. Then $f \mid_{N_1} : N_1 \to M_2$ is a morphism. If $f(N_1) \subset N_2$ and $f : N_1 \to N_2$ is an open mapping, then $f : N_1 \to N_2$ is a morphism.

\vskip6pt

{\bf Definition 8.} Let $f : M_1 \to M_2$ be a morphism between local Jordan manifolds. Then $f$ is called a {\it local isomorphisms} if every point $p \in M_1$ has a neighbourhood $U$ such that $f(U)$ is open in $M_2$ and $f : U \to f(U)$ is an isomorphism in LJM (where $U$ and $f(U)$ are local Jordan submanifolds of $M_1$ and $M_2$, respectively).

\vskip6pt

{\bf Proposition 3.} Let $(M, \hat{\mathcal{A}})$ be a local Jordan manifold and let $X$ be a topological Hausdorff space. Let $p : X \to M$ be a local homeomorphism. Then $X$ has a unique Jordan structure such that $p$ is a local isomorphism.

\vskip6pt

Let $\tilde{M}$ be a covering space of a local Jordan manifold $M$ and let $p : \tilde{M} \to M$ be the covering map. By Proposition 3 it follows that $\tilde{M}$ has a unique Jordan structure such that $p$ is a local isomorphism. If $M$ is a local Jordan $V$-manifold, then $\tilde{M}$ is also a local Jordan $V$-manifold. Let $\tilde{f} : \tilde{M} \to \tilde{M}$ be a continuos map which lies over a continuos map $f : M \to M \; (p\tilde{f} = fp)$; then $\tilde{f}$ is a LJM morphism if and only if $f$ is a LJM morphism.

If $\tilde{M}$ is a connected local Jordan manifold and $G$ is a group of LJM automorphisms acting freely and properly discontinuosly on $\tilde{M}$, then the orbit space $M = \tilde{M}/G$ possesses a unique Jordan structure such that the natural projection $p : \tilde{M} \to M$ is a local isomorphism. A Jordan atlas for $M$ is constructed by taking charts of the form $\hat{\Phi} = \Phi p^{-1} : p(U) \to V^{+}$, where $p : U \to p(U)$ is a homeomorphism and $\Phi : U \to V^{+}$ is a Jordan chart on $\tilde{M}$. If $U_i$ and $\Phi_i, i=1,2$, determine two such charts on $M$ and if $x \in p(U_1) \cap p(U_2)$, then there exists $x_1 \in U_1, x_2 \in U_2$ and $g \in G$ such that $p(x_1) = p(x_2) = x$ and $g(x_1) = x_2$. Hence $\hat{\Phi}_2\hat{\Phi}_1^{-1} = \Phi_2g\Phi_1^{-1}$ in a neighbourhood of $\hat{\Phi}_1(x)$, and hence locally it is the restriction of a linear fractional map.

As an example of an orbit space, let $\tilde{M} = V^{+}$, where $V = (V^{+}, V^{-})$ is a Jordan pair, and let $G$ be a group of translations of $V^{+}$ acting freely and properly discontinuosly. Then $M = \tilde{M}/G$ is a local Jordan $V$-manifold. This shows that a torus, for instance, possesses many local Jordan manifold structures.

\vskip6pt

Examples of local Jordan manifolds are: affine quadrics, projective quadrics, projective spaces, Grassmann manifolds.

\vskip6pt

{\bf Comment.} It would be interesting to reconsider Martinelli's results [471a,b] on quaternionic projective planes in this more general settings of Jordan pairs. An open problem is to find an algebraic characterization of integrability of a geometric structure.

\vskip6pt

Examples of two-dimensional local Jordan manifolds are: the torus, the sphere, the real projective plane, the cylinder.

\vskip6pt

In September 1980, at a conference organized in Romania, Gelfand gave a lecture on integral geometry mentioning the results obtained by himself and some of his coworkers on transformations between Grassmann manifolds. Then, after his lecture, I have suggested him to reconsider these results in a Jordan algebra (pair) setting. I have made the same suggestion also to MacPherson, concerning his results (in cooperating with Gelfand) on polyhedra in real Grassmann manifolds (GELFAND, I.M., \& MacPHERSON, R.D., {\it Geometry in Grassmannians and a generalization of the dilogarithm}, Preprint IHES, Paris, June 1981).

\vskip6pt

In Ch.III of his Ph.D. Thesis [701], Watson has considered the important subcategory of LJM, consisting of global Jordan manifolds (or simply, the Jordan manifolds (JM)).

\vskip6pt

{\bf Important remark.} It is obvious that the Jordan manifolds defined by Watson have much in common with differentiable manifolds. (The main examples are differentiable manifolds). It follows the following

\vskip6pt

{\bf Open problem.} Find an exact relationship between the category of JM and that of differential manifolds.

\vskip6pt

In 1951 Vagner began (see [684b,c]) a series of studies which led to a mathematical tool necessary to formulate and solve the following problem: To find the geometrical properties of differentiable manifolds which are derived {\it only} from algebraic properties of the pseudogroup of local homeomorphisms and the atlas.

By means of $[a,b,c] := c\cdot b^{-1}\cdot a$, where $a,b,c$ are elements of an atlas, a ternary operation is defined and this led to the algebraic notion of {\it heap (heath)}, a structure defined by a weakened set of a group postulates - see papers by Certain [159], Clifford [179] and Sushkevich [645].

This algebraic notion was used in differential geometry by Vagner (see [684b,d,e]). In the paper [684e] by Vagner, the following question has been raised: Do there exist properties of differentiable manifolds which depend {\it only} on the atlas?

Br\^ anzei [134a,b] has answered this question affirmatively and has defined a kind of generalized manifolds, called $\mathcal{SH}$-{\it manifolds}, which he studies in detail. A large number of examples have been given by him.

\vskip6pt

{\bf Open problem.} To compare the category of Jordan manifolds with category of $\mathcal{SH}$-manifolds.

\vskip6pt

{\bf Remark 8.} It would be of interest to re-consider the results on the above-mentioned ternary operation, as well as the results of Br\^ anzei [134a,b] in the theory of Jordan manifolds.

\vskip6pt

Let us mention now other categories whose objects are constructed in a similar way as differentiable manifolds. Among them we note topological manifolds, Vagner's compound manifolds [684a] Abraham's $(C^r, C^s)$-manifolds [3], Sikorski's differentiable spaces [621], Smith's differentiable spaces [628], Aronszajn's paper [45], Marshall's $C^{\infty}$-subcartesian spaces [470] (for the differential topology of these spaces see the paper [499b] by Motreanu), and Spallek's $N$-differentiable spaces [630].

Motreanu [499a,b,c] has constructed a category of so-called {\it preringed manifolds}, which contains as particular cases all categories mentioned above. Roughly speaking a preringed manifold is a topological space $M$ which is locally determined by a triple $(E, \mathcal{F}, V)$, where $E$ is a topological space, $\mathcal{F}$ is a presheaf of real-valued functions on $E$, and $V$ is a vector space, two such triples being compatible with respect to change of charts. Here $E$ describes locally the topology of $M$ and $V$ plays the role of tangent space to $M$. By a suitable choice of $E, \mathcal{F}$, and $V$ one obtains the particular theories recalled above. So, the following problem is immediate:

\vskip6pt

{\bf Open problem.} Describe the Jordan manifolds in terms of preringed manifolds.

\vskip6pt

Finally, I would like to mention four interesting papers which are very recent: two papers by Pumpl\" un [551a,b], and two papers by Pirio \& Russo [544a,b].

\vspace{1cm}

\centerline{\bf \Large \S 3. JORDAN STRUCTURES IN ANALYSIS}

\vspace{48pt}

The relationship between formally real Jordan algebras, self-dual homogeneous cones
and  symmetric upper half-planes in finite dimensions due to Koecher [408a,~c,~f] is the background
for the study of the infinite-dimensional case. The objects here are $JB^{\ast}$-algebras and their
real analogues (the so-called $JB$-algebras).

A generalization of formally real Jordan algebras to the infinite-dimen\-sional case was introduced
and studied by Alfsen, Shultz and St\o rmer [13] as follows:

\vspace{6pt}

{\bf Definition 1.} A (linear) real Jordan algebra ${\cal J}$ with unit element $e$ which
is also a Banach space and in which
the product and the norm satisfy:

\hskip6pt (i) $\|xy\|\leq \|x\|\, \|y\|$;

\hskip3pt (ii) $\|x^2\| = \|x\|^2$;

(iii) $\|x^2\| \leq \|x^2+y^2\|$;

\noindent for all $x,y \in {\cal J}$ is called a {\it Jordan Banach algebra} (or, briefly, a $JB$-{\it algebra}).

\vspace{6pt}

{\bf Remark 1.} In the finite-dimensional case, condition (iii) is equivalent to the fact that ${\cal J}$ is
a formally real Jordan algebra.

\vspace{6pt}

{\bf Note.} The term $JB$-algebra arose as the Jordan analogue of $B^{\ast}$-algebra, much the same as
$JC$-algebras and $JW$-algebras
were termed after $C^{\ast}$- and $W^{\ast}$-algebras, respectively.

\vskip6pt

{\bf Comments.} Hanche-Olsen and St\o rmer introduced [330] the concept of $JB$-algebra as follows:
A {\it Jordan Banach algebra} is a real
Jordan algebra $A$ (not necessarily unital) equipped with a complete norm satisfying
 $\|ab\| \leq \|a\|\,\|b\|$, $a,b \in A$. A $JB$-{\it algebra}
is a Jordan Banach algebra $A$ in which the norm satisfies the following two additional
conditions for $a, b \in A$,
$$(1^{\circ}) \qquad \|a^2\| =\|a\|^2, \qquad \mbox{and} \qquad (2^{\circ})
\qquad \|a^2\| \leq \|a^2 + b^2\|.$$

\vspace{6pt}

{\bf Theorem 1.} {\it A real Jordan algebra of finite-dimension is a
$JB$-algebra if and only if it is formally real}.

\vspace{6pt}

The algebras $H_p ({\bf F})^{(+)}$ (see Theorem 8 from \S 1) can be extended to arbitrary
cardinality $p$ as follows. Let
${\bf F} = {\bf R}, {\bf C}$ or ${\bf H}$, and let $H$ be a (right) ${\bf F}$-Hilbert space of
dimension $p$ over ${\bf F}$.
Denote by ${\cal L} (H)$ the algebra of all bounded ${\bf F}$-linear operators on $H$.
Then there
exists a natural involution (the adjoint ${\ast}$) on ${\cal L} (H)$ and ${\cal H}_p ({\bf F})
:={\cal H}(H) : = \{\lambda \in {\cal L} (H)\mid \lambda^{\ast} =\lambda \}$ is a
$JB$-algebra with respect to the operator norm.

\vspace{6pt}

{\bf Remark 2.} For every compact topological space $S$ and every $JB$-algebra
${\cal J}$, the algebra ${\cal C} (S, {\cal J})$ of all continuous
functions $S \rightarrow {\cal J}$ is also a $JB$-algebra. In particular, ${\cal C}(S, {\bf R})$
is an associative $JB$-algebra.

\vspace{6pt}

{\bf Proposition 2.} {\it Every associative $JB$-algebra is isometrically
isomorphic to ${\cal C} (S, {\bf R})$ for some compact topological
space} $S$.

\vskip6pt

{\bf Notation.} Let ${\cal J}$ be a $JB$-algebra and denote by ${\cal J}^2: =
\{x^2 \mid x \in {\cal J} \}$ the positive cone in ${\cal J}$.

\vskip6pt

{\bf Definition 2.} The elements of $\Omega $ are called {\it positive definite},
and an ordering $<$ on ${\cal J}$ is defined by: $x < y$ for $x,y \in {\cal J}$ if and
only if $y-x \in \Omega $.

\vspace{6pt}

{\bf Definition 3.} A subcone of ${\cal J}^2$ is said be a {\it face} of ${\cal J}^2$ if it
contains all elements $a$ of ${\cal J}^2$ such that $a \leq b$ for some $b$ of it.

\vspace{6pt}

{\bf Definition 4.} A $JB$-subalgebra $B$ of ${\cal J}$ is said to be
{\it hereditary} if its positive cone $B^2$ is a face
of ${\cal J}^2$.

\vspace{6pt}

Edwards [226a] proved the following results (see also [231a]):

\vspace{6pt}

{\bf Theorem 3.} {\it The norm-closed quadratic ideals of a $JB$-algebra ${\cal J}$
 coincide with the hereditary
$JB$-algebras $B$ of ${\cal J}$, and the norm-closed faces of ${\cal J}^2$ are the positive
cones $B^2$ of such subalgebras} $B$.

\vspace{6pt}

{\bf Proposition 4.} {\it The norm closure of a face of ${\cal J}^2$ is a face of} ${\cal J}^2$.

\vspace{6pt}

Putter and Yood [552] generalized a number of well-known Banach algebra results to the
Jordan algebra situation by appropriately modifying the proofs.
They confined themselves to special $JB$-algebras.

\vskip6pt

Consider now the complexification ${\cal J}^{{\bf C}}:={\cal J} \oplus {\rm i} {\cal J}$ of ${\cal J}$.
${\cal J}^{{\bf C}}$ is a complex Jordan algebra with involution $(x+{\rm i}y)^{\ast}: = x-{\rm i}y$.

\vspace{6pt}

{\bf Definition 5.} $D:=D(\Omega): =\{z \in {\cal J}^{\bf C} \mid {\rm Im} (z) \in \Omega \}$
is called the {\it tube domain $($generalized upper half-plane$)$} associated
with the cone $\Omega $.

\vspace{6pt}

Let us mention here that Tsao [675] proved that under certain conditions the Fourier
coefficients of the Eisenstein series for an
arithmetic group acting on a tube domain are rational numbers. The proof involves a
mixture of Lie groups, Jordan algebras, Fourier analysis,
exponential sums, and $L$-functions.

From the results reported by Alfsen, Hanche-Olsen, Shultz and St\o rmer [11], [13],
it follows that for each $JB$-algebra ${\cal J}$
there exist a canonical $C^{\ast}$-algebra ${\frak A}$ and a homomorphism
$\Psi : {\cal J} \rightarrow {\frak A}$ such that $\Psi ({\cal J})$ generates
${\frak A}$. The kernel of $\Psi $ is the exceptional ideal ${\cal I}$ in ${\frak A}$. Using Takesaki and Tomiyama's methods, Behncke and B\"{o}s showed
[75] that ${\cal I}$ may be described as an $H_3({\bf O})$-fibre bundle over its primitive ideal space.

\vspace{6pt}

{\bf Comments.} As it was remarked by Upmeier [682b], a promising application of $JB$-algebras is to be found in complex analysis, based
on the one-to-one correspondence between $JB^{\ast}$-algebras and bounded symmetric domains in complex Banach spaces with tube
realization (Koecher [408f], and Braun, Kaup and Upmeier [133b]).

\vskip6pt

{\bf Definition 6.} A $JB^{\ast}$-{\it algebra} is a complex Jordan
algebra ${\cal J}$ with unit element $e$, (conjugate linear) involution
$\ast $,
and complete norm such that

\hskip3pt(i) $ \|xy\|\leq \|x\|\,\|y\|$;

(ii) $\|P(z)z^{\ast} \| = \|z\|^3$,

\noindent for all $x,y,z \in {\cal J}$.

\vspace{6pt}

{\bf Note.} The concept of $JB^{\ast}$-algebra was formulated by
Kaplansky (lecture at the 1976 St. Andrews Colloquim of the Edinburgh
Math. Soc. -- see [712]) and introduced as ``Jordan $C^{\ast}$-algebra".

\vspace{6pt}

{\bf Comments.} Youngson [723d] studied $JB^{\ast}$-algebra in the nonunital case. He stated, among other results, that
nonunital $JB^{\ast}$-algebras are $C^{\ast}$-triple systems in the sense of Kaup [392b].

\vskip6pt

{\bf Proposition 5.} {\it The selfadjoint part of a $JB^{\ast}$-algebra
is a $JB$-algebra}.

\vspace{6pt}

Wright [712] proved the converse:

\vspace{6pt}

{\bf Theorem 6.} {\it For every $JB$-algebra ${\cal J}$ there exists
a unique complex norm on ${\cal J}^{{\bf C}}$ such that
${\cal J}^{{\bf C}}$ is a $JB^{\ast}$-algebra with selfadjoint part
${\cal J}$. The correspondence ${\cal J} \leftrightarrow
{\cal J}^{\bf C}$ defines an equivalance of the category of $JB$-algebras
onto the category of $JB^{\ast}$-algebras}.

\vspace{6pt}

Russo and Dye [576] proved that the closed  unit ball of a
$C^{\ast}$-algebra with identity is the convex hull of its unitary
elements. The same result was proved by Wright and Youngson [713a] for $JB^{\ast}$-algebras.

Using the fact that the extreme points of the positive ball in a
$JB$-algebra are projections, Wright and Youngson first showed [713b]
that a surjective unital linear isometry between two $JB$-algebras is a Jordan isomorphism, and
then used this to obtain the same result for $JB^{\ast}$-algebras.

Bonsall [119] showed that if $B$ is a real closed Jordan subalgebra of a
complex unital Banach algebra $A$, containing the unit and such that
$B \cap {\rm i}B=\{0\}$ and $B \subset H(B) \oplus {\rm i}H(A)$, where
$H(A)$ denotes the set of Hermitian elements of $A$, then $B \oplus {\rm i}B$
is homeomorphically $\ast $-isomorphic to a $JB^{\ast}$-algebra. Using
Wright's and Youngson's results
[712], [713a], [723a,~b], Mingo [489] gave a $JB^{\ast}$-analogue of a
$C^{\ast}$-algebra result St\o rmer [641a], as follows:

\vspace{6pt}

{\bf Proposition 7.} {\it Suppose $A$ is a $JB^{\ast}$-algebra and $B$
is a real selfadjoint subalgebra with unit such that
$B \cap {\rm i}B=\{0\}$. Then $B \oplus {\rm i}B$ is a $JB^{\ast}$-algebra}.

\vspace{6pt}

Mingo used Proposition 7 to prove the above-mentioned result of Bonsall,
dispensing with the assumption $B \cap {\rm i}B=\{0\}$, and
also to prove that the isomorphism is an isometry.

\vspace{6pt}

{\bf Definition 7.} A bounded domain $B$ in a complex Banach space is
called {\it symmetric} if for every $a$ of $B$
there exists a holomorphic map $s_a:B \rightarrow B$ with $s_a^2 =
{\rm Id}_B$ and $a$ an isolated fixed point ($s_a$ is uniquely determined if it exists and is called the symmetry at $a$.)

\vskip6pt

{\bf Definition 8.} For every open cone $C$ in a real Banach space $X$,
the domain $T:=\{z \in X \oplus {\rm i} X \mid  {\rm Im} (z) \in C\}$ is called a
{\it symmetric tube domain} if $T$ is
biholomorphically equivalent to a bounded symmetric domain.

\vspace{6pt}

{\bf Theorem 8.} {\it Let ${\cal J}$ be a $JB$-algebra and let
${\cal J}^C ={\cal J} \oplus {\rm i}{\cal J}$ be the corresponding
$JB^{\ast}$-algebra. Then $D := \{z \in {\cal J}^C \mid  {\rm Im} (z)
\in \Omega \}$ is
a symmetric tube domain. The symmetry at the point ${\rm i} e \in D$ is given
by $s(z) =-z^{-1}$, and $z \rightarrow {\rm i}(z-{\rm i} e)^{-1}$ maps $D$
biholomorphically on the open unit ball $\bg $ of ${\cal J}^{{\bf C}}$.
In particular, $\bg $ is a homogeneous domain}.

\vspace{6pt}

Braun, Kaup and Upmeier [133a] proved

\vspace{6pt}

{\bf Theorem 9.} {\it If $B$ is a real Banach space and $D$ is the
symmetric tube domain
for ${\cal B}^{{\bf C}}:={\cal B} \oplus {\rm i}{\cal B}$, then for every
$e \in \Omega $ there exists a unique
Jordan product on ${\cal B}$ such that ${\cal B}$ is a $JB$-algebra
with unit $e$, and $D$ is the upper half-plane}.

\vspace{6pt}

{\bf Remark 3.} It follows that $JB$-algebras, as well as
$JB^{\ast}$-algebras, are in one-to-one
correspondence with symmetric tube domains.

\vspace{6pt}

In the theory of formally real Jordan algebras of finite dimension an
important fact is the minimal
decomposition of elements of such an algebra with respect to a complete
orthogonal system of primitive idempotents $\{e_1,
\ldots , e_k\}$. The importance of the minimal decomposition follows from
the fact that $\{e_1, \ldots , e_k\}$ determine a Peirce decomposition of the
algebra which, for instance, diagonalizes the operator $L(x)$, and hence
also $P(x)$.

The analogue for an arbitrary $JB$-algebra ${\cal J}$ is the fact that for
every $\alpha \in {\cal J}$ the
unital closed subalgebra $C(\alpha)$ generated by $\alpha $ is isomorphic
to some ${\cal C}(S, {\bf R})$, where $S$ is a compact topological space.
However, in case $S$
is connected, $e$ is the only nontrivial idempotent in $C(\alpha )$ and the
Peirce decomposition cannot be applied.

Shultz [618] proved that the bidual of $JB$-algebra with the Arens product
is also a $JB$-algebra. Hence, every $JB$-algebra
is a norm-closed subalgebra of a $JB$-algebra which is a dual Banach space.
Algebras of this type admit not only a continuous but also
an $L^{\infty}$-functional calculus.

\vspace{6pt}

{\bf Remark 4.} Edwards [226c] showed how some of the results on
multipliers and
quasi-multipliers of $C^{\ast}$-algebras can be extended to $JB$-algebras.

\vspace{6pt}

{\bf Definition 9.} A $JB$-algebra ${\cal J}$ is called a
$JBW$-{\it algebra} if ${\cal J}$ is a dual
Banach space (i.e., there exists a Banach space ${}^{\prime} {\cal J}$
with ${\cal J} = {}^{\prime}{\cal J}^{\prime}$ as dual
Banach space; ${}^{\prime} {\cal J}$ is uniquely determined by ${\cal J}$
(see Sakai [586]) and is called the {\it predual} of ${\cal J}$).

\vspace{6pt} 

{\bf Example.} The selfadjoint part of a von Neumann algebra is a 
$JBW$-algebra.

\vspace{6pt}

{\bf Remark 5.} For every $\alpha $ in the $JBW$-algebra ${\cal J}$ the
$w^{\ast}$-closed unital subalgebra $W(\alpha)$ of ${\cal J}$ generated
by $\alpha $ is a commutative von Neumann algebra, i.e., $W(\alpha )
\approx {\cal C} (S, R)$ for $S$ hyperstonian or,
equivalently, $W(\alpha ) \approx L^{\infty} (\mu)$, where $\mu $ is a
localizable measure (see Sakai [586]).

\vspace{6pt}

{\bf Remark 6.} $JBW$-algebras (weakly closed analogues of
$JB$-algebras) are the abstract analogues of von Neumann algebras in the
Jordan case.

\vskip6pt

{\bf Definition 10.} Let ${\cal J}^{\prime}$ be the dual Banach space of a $JB$-
algebra ${\cal J}$ and denote by $({\cal J}^2)^{\prime}:=\{\lambda \in {\cal J}^{\prime}
\mid \lambda ({\cal J}^2)\geq 0\}$ the dual cone of ${\cal J}^2$. Then
$K:=\{\lambda \in ({\cal J}^2)^{\prime} \mid
\lambda (e) =1\}$ is called the {\it state space of} ${\cal J}$, the
elements of $K$ being called {\it states on}
${\cal J}$.

\vspace{6pt}

{\bf Remark 7.} $K$ is a $w^{\ast}$-compact, convex subset and ${\cal J}$ can be
identified (as a Banach space) with the
space of all $w^{\ast}$-continuous affine functions on $K$. The bidual of ${\cal J}$
coincide with the set of all
bounded affine functions on $K$.

\vspace{6pt}

In a comprehensive study of state spaces of a $JB$-algebra, Alfsen and Shultz [12a]
 gave necessary and sufficient conditions for a compact convex set to be a state space of
a $JB$-algebra. Araki [35] improved the characterization  of state spaces
of $JB$-algebra given in [12a] to a form with more physical appeal in the simplified finite-dimensional case.

Alfsen, Hanche-Olsen and Shultz [11] characterized the state spaces of
$C^{\ast}$-algebras among the state spaces of all $JB$-algebras. Together, [12a]
 and [11] give a complete characterization of the state spaces of $C^{\ast}$-algebras. As is shown in [11],
a $JB$-algebra ${\cal J}$ is the selfadjoint part of a $C^{\ast}$-algebra if and only
if ${\cal J}$ is of complex type and the state space of ${\cal J}$ is orientable.
Stacey [634c] showed that the state space of a $JBW$-algebra of complex type is
orientable if and only if it is locally orientable. For local and global splittings
in the state space of a $JB$-algebra, see Stacey [634b].

Every state $\lambda \in K$ defines by $(x|y)_{\lambda}: = \lambda (xy)$ a positive
inner product on ${\cal J}$ and, in particular, by $|x|_{\lambda} :=\lambda (x^2)^{\frac{1}{2}}$ a seminorm
on ${\cal J}$.

\vspace{6pt}

{\bf Definition 11.} A state $\lambda $ on a $JB$-algebra ${\cal J}$ is called {\it faithful}
if $|\cdot |_{\lambda}$ actually is norm
on ${\cal J}$, i.e. if $\lambda (x^2)=0$ implies that $x=0$.

\vspace{6pt}

{\bf Definition 12.} A state $\lambda $ on a $JBW$-algebra ${\cal J}$ is called {\it normal}
if $\lim \lambda (x_{\alpha})= \lambda (x)$ for every increasing net $x_{\alpha}$ in
${\cal J}$ with $x=\sup x_{\alpha}
\in {\cal J}$.

\vspace{6pt}

{\bf Definition 13.} A normal state $\lambda $ on a $JBW$-algebra ${\cal J}$ is called
a {\it finite trace} if it is associative in the sense that
$\lambda ((xy)z) =\lambda (x(yz))$ for all $x,y,z \in {\cal J}$.

\vspace{6pt}

{\bf Remark 8.} The condition from Definition 13 states that every $L(y)$, $y \in {\cal J}$,
is selfadjoint with respect to the inner
product $(\cdot \mid \cdot )_{\lambda }$.

\vspace{6pt}

A complete study of $JBW$-algebras with a faithful finite trace was undertook by Janssen [373b]. On the
basis of this paper, Janssen [373c] studied the properties of the lattice of idempotents in a
finite weakly closed Jordan algebra. He proved that such an algebra admits a unique decomposition into a direct sum of a discrete
Jordan algebra and a continuous Jordan algebra. Janssen [373c, II] gave
a completely description of the discrete finite weakly closed Jordan algebras by
finite-dimensional simple formally real Jordan algebras and by
simple formally real Jordan algebras of quadratic forms of real Hilbert spaces.

Pedersen and St\o rmer [531] showed that the different definitions of trace on a
Jordan algebra are all equivalent for $JBW$-algebras, and that conditions that do not
involve projections are equivalent for $JB$-algebras. They have considered only
finite traces. Iochum [358a] extended the results to semifinite traces. By a suitable definition
of semifiniteness, he showed that for any $JBW$-algebra we have a unique central
decomposition in finite (semifinite) and proper-infinite (pure-infinite) parts
exactly as in the case of von Neumann algebras (see [358a, Theorem V.1.6]). Iochum proved also (for the
semifinite case) the equivalence between the category
of facially homogeneous self-dual cones and the category of $JBW$-algebras of
selfadjoint derivations (see [358a, Theorem V.5.1]), and ([358a, Ch. VII]) his
main theorem, which establishes the equivalence between the category of facially
homogeneous self-dual cones in Hilbert spaces and the category of $JBW$-algebras
(see also [77d]).

Assume now that ${\cal J}$ is a $JBW$-algebra with a faithful finite trace $\lambda $. Then
$\lambda $ is essentially uniquely determined (every other faithful finite trace is of the form
$\lambda \circ P(h) =\lambda \circ L(h^2)$ for some $h>0$ in the centre of ${\cal J}$), and ${\cal J}^{{\bf C}}$ is a
complex pre-Hilbert space with respect to the inner product
$(z|w):=(z|w)_{\lambda} := \lambda (zw^{\ast})$, where $\lambda $ is
extended ${\bf C}$-linearly to ${\cal J}^{{\bf C}}$.

\vspace{6pt}

{\bf Notation.} Denote by $\widehat{{\cal J}}^{{\bf C}}$ the completion of
${\cal J}^{{\bf C}}$ with respect to the norm $\|z\|_2:=\|z\|_{\lambda}: =
\lambda(zz^{\ast})^{\frac{1}{2}}$, and consider the closures $\widehat{{\cal J}}$
and $\widehat{{\cal J}}^{2}$ of ${\cal J}^2$~in~${\cal J}^{{\bf C}}$.

\vspace{6pt}

The operators $L(z)$ and $P(z)$, $z \in {\cal J}^{{\bf C}}$, can be continuously extended to
 bounded operators on $\wh{{\cal J}^{{\bf C}}}$
satisfying $L(z)^{\ast}=L(z^{\ast})$ and $P(z)^{\ast} =P(z^{\ast})$. The cone $\wh{{\cal J}^2}$
is self-dual in $\wh{{\cal J}}$, satisfies
a certain geometrical homogeneity condition, and has $e$ as trace vector (i.e., as quasi-interior
point of $\wh{J^2}$ fixed by every connected set
of isometries in $GL (\wh{{\cal J}})$). On the other hand, every cone of this type in a real Hilbert
space is obtained in this way from a $JB$-algebra with faithful finite trace (see Bellissard and
Iochum [77a]). This result can be viewed as a generalization to the
infinite-dimensional case of the following theorem of Koecher: {\it The self-dual cones with
homogeneous interior in real Hilbert spaces of finite dimension are precisely $($up to linear equivalence$)$
the cones of squares in formally real Jordan algebras}.

\vskip6pt

{\bf Proposition 10.} {\it The $JB^{\ast}$-norm $\|\cdot \|_{\infty}$ on ${\cal J}^{{\bf C}}$
satisfies $\|\cdot \|_2 \leq \|\cdot \|_{\infty} $ on ${\cal J}^{{\bf C}}$, and
$\bg = \{ z \in {\cal J}^{{\bf C}} \mid  1-P(z)P(z)^{\ast} >0\}$, $\Sigma =
 \exp ({\rm i}{\cal J}) =\{z \in {\cal J}^{{\bf C}} \mid
P(z)$ unitary on} ${\cal J}^{{\bf C}}\} =\{z \in \ov{\bg}  \mid  \|z\|_2 =1 \}$.

\vspace{6pt}

{\bf Proposition 11.} {\it If ${\cal J}$ is a $JB$-algebra, then the following conditions are equivalent:

\hskip6pt{\rm (i)} there exists a maximal associative subalgebra of finite dimension in~${\cal J}$;

\hskip3pt{\rm (ii)} ${\cal J}$ is locally finite $($i.e., every finitely generated subalgebra has finite dimension$)$;

{\rm (iii)} for every $a\in {\cal J}$ the operator $L(a) \in {\cal L} ({\cal J})$ satisfies a polynomial equation over ${\bf R}$;

\hskip0.5pt{\rm (iv)} there exists a natural number $r$ such that every $a \in {\cal J}$ admits a representation $a= \alpha_1 e_1 +
\cdots + \alpha_r e_r$, where $\{e_1,\ldots , e_r\}$ is a set of orthogonal idempotents and $\alpha_1, \ldots , \alpha_r \in {\bf R}$;

\hskip3.3pt{\rm (v)} there exists a faithful finite trace $\lambda $ on ${\cal J}$ such that the corresponding Hilbert norm $\|x\|_2 = \lambda(x^2)^{\frac{1}{2}}$ on
${\cal J}$ is equivalent to the $JB$-norm $\|x\|_{\infty}$;

\hskip0.5pt{\rm (vi)} ${\cal J}$ is reflexive.}

\vspace{6pt}

Chu [166a] studied the Radon-Nikodym property (for definition see below) in the context
of $JBW$-algebras.

\vspace{6pt}

{\bf Definition 14.} A (real or complex) Banach space $X$ is said to possess the {\it Radon-Nikodym property}
if for any finite measure space $(\Omega, \Sigma, \mu )$ and $\mu$-continuous vector measure $L: \Sigma
\rightarrow X$ of bounded total variation, there exists a Bochner integrable function $g: \Omega \rightarrow
X$ such that $L(E) = \int_E g\,{\rm d}\mu $ for\break all $E$ in $\Sigma$.

\vskip6pt

Using a result of Shultz, Chu [166a] established the following result:

\vskip6pt

{\bf Theorem 12.} {\it Let ${\cal J}$ be a $JBW$-algebra. Then its dual ${\cal J}^{\prime}$ has
the Radon-Nikodym property if and only if ${\cal J}$ is a finite direct sum of Jordan algebras, each of which is one of the following algebras:}

\hskip6pt(i) {\it Jordan $(n \times n)$-matrix algebras over ${\bf R}, {\bf C}$, or} ${\bf H}$;

\hskip3pt(ii) {\it spin factors;}

(iii) {\it the exceptional Jordan algebra of Hermitian $(3 \times 3)$-matrices over}~${\bf O}$.

\vspace{6pt}

Chu [166b] proved that the dual of a $JB$-algebra ${\cal J}$ possesses the Radon-Nikodym
property if and only if the state space of ${\cal J}$ is the $\sigma $-convex hull of its pure states. Namely, he proved:

\vspace{6pt}

{\bf Theorem 13.} {\it If ${\cal J}$ is a $JB$-algebra with state space $K$, then the following conditions are equivalent:

\hskip6pt{\rm (i)} $K$ is $\sigma $-convex hull of the pure states, i.e.,
$$K=\Big\{  \sum\limits_{n=1}^{\infty} \lambda_n k_n \, \Big|\,  \sum\limits_{n=1}^{\infty} \lambda_n=1, \,
\lambda_n \geq 0,\,  k_n \mbox{ being pure states} \Big\};$$

\hskip3pt{\rm (ii)} ${\cal J}^{\prime}$ has the Radon-Nikodym property;

{\rm (iii)} ${\cal J}^{\prime \prime}$ is a direct sum of type $I$ $JBW$-algebras
$($i.e.,\ $JBW$-algebras which contains a $($non-zero$)$
minimal idempotent$)$}.

\vspace{6pt}

In 1991--1992, Bunce and Chu [144a,~b] studied the Radon-Nikodym property in $JB^{\ast}$-triples
(see Definition 19 below).

\vskip6pt

{\bf Definition 15.} An idempotent $p$ in a $JBW$-algebra ${\cal J}$ is called an
{\it atom} if it is minimal (i.e., if $0 \leq q \leq p$ with $q^2 =q$ implies that $q=0$ or $q=p$).
${\cal J}$ is said to be {\it atomic} if every idempotent is the least upper bounded of orthogonal atoms.

\vspace{6pt}

{\bf Notation.} If $\lambda $ is a state on the $JBW$-algebra ${\cal J}$, then
$$V_{\lambda} : = \{f \mid f \in {\cal J}^{\prime},\, \exists\, a \in {\bf R}_+
\mbox{ with } -a\lambda \leq f  \leq a \lambda\}.$$

\vspace{6pt}

{\bf Theorem 14.} {\it Let ${\cal J}$ be an atomic $JBW$-algebra and let $\lambda $
be a faithful normal state on ${\cal J}$.
Then there exists an order isomorphism
$\varphi : V_{\lambda} \rightarrow {\cal J}$ with} $\varphi (\lambda ) =e$.

\vspace{6pt}

{\bf Theorem 15.} {\it Suppose that a $JBW$-algebra ${\cal J}$ admits a faithful normal trace $\lambda $
$($i.e., for all idempotents $p,q$ we have $\lambda (U_pq-U_qp)=0$, where $U_pq:=(pq)q-(qp)q+q^2p)$. Then there exists an
order isomorphism $\varphi : V_{\lambda} \rightarrow {\cal J}$ with $\varphi (\lambda)=e$. Moreover, for every
positive $\mu$ $($i.e., $\mu \in ({\cal J}^2)^{\prime})$ from $V_{\lambda}$, there exists
a positive element $y$ in ${\cal J}$ such that} $\mu (x) =\lambda (U_y x)$.

\vspace{6pt}

{\bf Theorem 16.} {\it Let ${\cal J}$ be a $JBW$-algebra satisfying the quadratic Radon-Nikodym property $($i.e., for any $f,g \in {}^{\prime}{\cal J}$
with $0 \leq f(x^2) \leq g (x^2)$ for every $x \in {\cal J}$, there exists a positive $y $ in ${\cal J}$ such that $f(x)=g(U_yx)$ for every
$x \in {\cal J})$ and let $\mu $ and $\nu $ be faithful normal states on ${\cal J}$. Then $V_{\mu}$ and $V_{\nu} $ are
order isomorphic}.

\vspace{6pt}

{\bf Corollary.} {\it Let ${\cal J}$ be as in Theorem 16 and suppose that ${\cal J}$
admits a faithful normal trace. Let $\lambda $ be a faithful normal state on ${\cal J}$. Then
$V_{\lambda }$ is order isomorphic to} ${\cal J}$.

\vskip6pt

{\bf Definition 16.} A family $(v_t)_{t \in {\bf R}}$ of linear operators on a linear space $M$
satisfying $v_0 ={\rm Id}$ and the cosine identity
$$2v_sv_t=v_{s+t}+v_{s-t}, $$
is called a (one-parameter) {\it cosine family on} $M$.

\vspace{6pt}

{\bf Remark 9.} If $(u_t)$ is a one-parameter group, then $(u_t+u_{-t})/2$ is
a cosine family.

\vspace{6pt}

{\bf Definition 17.} Let ${\cal J}$ be a $JBW$-algebra and $\lambda $ a normal state
on ${\cal J}$. A bilinear, symmetric, positive semidefinite form $s$ on ${\cal J}$ satisfying

\hskip6pt(i) $a(a,b) \geq 0$, $a \geq 0$, $b \geq 0$;

\hskip3pt(ii) $s(1,a)=\lambda (a)$, $a \in {\cal J}$;

(iii) if $0\leq \mu \leq \lambda$, there is $0 \leq b \leq 1$
so that  $\mu(a) =s(a,b)$, $a \in {\cal J}$,

\noindent is called a {\it self-polar form associated with} $\lambda $.

\vspace{6pt}

{\bf Remark 10.} There exists at most one self-polar form associated~with~$\lambda$.

\vspace{6pt}

{\bf Theorem 17.} {\it Let ${\cal J}$ be a $JBW$-algebra and let
$\lambda $ be a faithful normal state on
${\cal J}$. Then there exists a unique cosine family $(\theta_t)$ of
positive, unital linear mappings of ${\cal J}$ into
itself, having the following properties:}

\hskip6pt(i) {\it for each} $a\in{\cal J}$, $t \rightarrow \theta_t(a)$ {\it is
weakly continuous};

\hskip3pt(ii) $\lambda (\theta_t(a) \circ b) = \lambda (a \circ \theta_t(b)); $

(iii) $s(a,b):= \int \lambda (a \circ \theta_t(b)) \cos h (\pi t)^{-t}
{\rm d} t$ {\it  defines a
self-polar form asso\-ciated with} $\lambda $.

\vspace{6pt}

Let us mention the following result of St\o rmer [641e], related to $JW$-algebras.

Let $M$ be a von Neumann algebra and let $\alpha $ be a central involution of $M$, i.e., $\alpha $ is $\ast$-antiautomorphism
of order $2$ leaving the centre of $M$ elementwise fixed. Then the set $M^{\alpha}:= \{x \in M \mid x=
x^{\ast}=\alpha (x)\}$ is a $JW$-algebra with Jordan product $xy:=1/2(x \circ y + y \circ x)$. St\o rmer
studied the relationship between $M^{\alpha}$ and $M^{\beta}$ for two central involutions
$\alpha $ and $\beta $. The main result states
that $\alpha $ and $\beta $ are (centrally) conjugate, i.e., there exists a
$\ast $-automorphism $\phi $ of $M$ leaving the centre elementwise
fixed, such that $\beta = \phi \alpha \phi^{-1}$ if and only if $M^{\alpha}$
and $M^{\beta}$ are isomorphic as Jordan algebras
{\it via} an isomorphism which leaves the centre elementwise fixed. Now
$M^{\alpha}$ generates $M$ as a von Neumann algebra (except in
a few simple cases) and there are von Neumann algebras with many conjugate
classes of central involutions.

Thus there may be many, even an uncountable number, of non-isomorphic
$JW$-algebras which generate the same von Neumann algebra. Such examples may be
found in [641e, Section 5].

\vspace{6pt}

{\bf Definition 18.} A $JH^{\ast}$-{\it algebra} is a complex Hilbert space $H$ together with a
complex Jordan algebra structure and a continuous involution $x \rightarrow x^{\ast}$ such that the Jordan
product is continuous and $L(a)^{\ast}= L(a^{\ast})$ for every $a \in H$, $L$ being the left multiplication.

\vspace{6pt}

{\bf Definition 19.} A complex Banach space $Z$ with a continuous mapping
$(a,b,c) \rightarrow \{abc\}$ from $Z \times Z \times Z$
to $Z$ is called a $JB^{\ast}$-{\it triple} if the following conditions are
satisfied for all $a,b,c, d \in Z$, where the
operator $a \Box b$ from the Banach algebra ${\cal L} (Z)$ of all bounded
linear operators on $Z$ is defined by $z \rightarrow \{abz\}$ and $[\cdot\,,\cdot]$ is the commutator
product:

1. $\{abc\}$ is symmetric complex linear in $a,c$ and conjugate linear in $b$;

2. $[a \Box b, c \Box d] =\{abc\} \Box d-c \Box \{dab\}$;

3. $a \Box a$ is Hermitean and has spectrum $\geq 0$;

4. $\|\{aaa\}\|=\|a\|^3$.

\vspace{6pt}

{\bf Definition 20.} $A\; JB^{\ast}$-triple is called a $JBW^{\ast}$-{\it triple} if it is the dual of a Banach
space.

\vspace{6pt}

{\bf Definition 21.} An element $z \in Z$ is called {\it tripotent} if $\{eee\} =e$.

\vspace{6pt}

{\bf Remark 11.} The set ${\rm Tri} (Z)$ of tripotent elements is endowed with the induced topology of $Z$.
It has been showed by J. Sauter in his Ph.D. Dissertation {\it Randstrukturen
beschr\"{a}nkter symmetrischer Gebiete} (T\"{u}bingen Univ., 1995) that
${\rm Tri}(Z)$ is a real analytic direct submanifold of $Z$.

\vspace{6pt}

If $e \in {\rm Tri} (Z)$, then $e \Box e \in {\cal L} (Z)$ has the eigenvalues $0,1/2, 1$, and we have the following
{\it Peirce decomposition of} $Z$ {\it with respect to} $e$
$$Z=Z_1(e) \oplus Z_{1/2} (e) \oplus Z_0 (e), $$
the {\it Peirce projections} being
$$P_1(e)=Q^2(e), \quad P_{1/2}=2(e \Box e - Q^2(e)), \quad P_0={\rm Id}
-2e \Box e+Q^2(e),$$
where $Q(e)z:= \{eze\}$ for $z \in Z$.

\vskip6pt

A recent paper on closed tripotents is that that by Fernandez-Polo \& Peralta [262b].

Let us mention here the research monograph [512d] by Neher, where a theory of grids (i.e., special families
of tripotents in Jordan triple systems) is presented. Among the applications there is also the structure
theory of $JBW^{\ast}$-triples.

Concerning the spectrum preserving linear maps on $JBW^{\ast}$-triples, see the paper by Neal [510].

Some applications of Jordan theory to harmonic analysis have been found by Chu in [166d] and by
Chu and Lau in [169a,b].

Very recently (in 2009), Stach\'o and Werner [636b] defined the notion of {\it Jordan manifolds} as Banach manifolds whose tangent spaces are endowed with Jordan triple products depending smoothly on the underlying points. They show examples of Jordan manifolds with various features giving rise to problems for further studies.

\vskip6pt

{\bf Comment.} In 1978, Watson [701] defined the notion of {\it Jordan manifold} in a completely different setting (see the previous \S 2).

\vskip6pt

Two interesting very recent papers are those by Arazy, Engli\v s and Kaup [38], and by Kennedy [398].

\vskip6pt

In April 2011, Werner [704] pointed out very interesting relations between Jordan $C^{\ast}$-triple systems and $K$-theory. These results were very recently obtained by him in cooperation with one of his Ph.D. students.

\vskip6pt

Let us concern now with the applications of Jordan structures to Riccati differential equations,
to Hua equations and Szeg\"{o} kernel, to the (reproducing) kernel functions, as well as to
dynamical systems, and to Shilov boundary.

\vspace{6pt}

1. The Riccati differential equation
$$\dot{x}=p(x),$$
$x \in {\bf R}^n$ and $p: {\bf R}^n \rightarrow {\bf R}^n$ homogeneous and quadratic, plays an
important role in biology, genetics, ecology, and chemistry. Koecher [408f,~g,~i] and Meyberg [483b]
studied the relations of this equation with nonassociative algebras, in particular with Jordan algebras.

We consider a commutative algebra ${\cal A}$ over ${\bf R}^n$ with product $xy:=\frac{1}{2} (p(x+y)-p(x)-p(y))$, and
let ${\cal A}_a$ be the mutation of ${\cal A}$ with respect to $a$ (see  Definition 4 from \S 1).

\vspace{6pt}

{\bf Notation.} Denote by ${\cal R}_n$ the vector space of power series in ${\bf R}^n$
converging in a neighbourhood of zero.

\vspace{6pt}

For $p,q \in {\cal R}_n$ we define $p \cdot q \in {\cal R}_n$ by
$$[(p \cdot q ) (u)]_i: = \sum\limits_{j=1}^n \frac{\partial p_i (u)}{\partial u_j} q_j (u). $$

\vskip6pt

{\bf Remark 12.} The vector space ${\cal R}_n$ with the product $(p,q) \rightarrow
p\cdot q$ becomes a nonassociative algebra over ${\bf R}$.

Now define $g_{{\cal A}} (u) \in {\cal R}_n$ by $g_{\cal A}  (u):=\sum\limits_{m=0}^{\infty}
\frac{1}{m!}g_m (u)$, where $g_0(u):=u$, $g_{m+1}:=g_m\cdot p$, and $p(u):=u^2$. (Powers in
${\cal A}$ are defined as follows: $u^1:=u$, $u^{m+1}:=uu^m$.)

The elements $f \in {\cal R}_n$ such that $f(x(\xi))$ is a solution of the above mentioned Riccati
equation whenever $x (\xi)$ is a solution, form a group $S({\cal A})$ under composition,
the solution-preserving group of the equation.

\vspace{6pt}

{\bf Notation.} Let ${\cal J} ({\cal A})$ denote the subspace of all $a \in {\bf R}^n$ satisfying
$2u(u(ua))+u^3a=2u(u^2a)+u^2 (ua)$ for all $u \in {\bf R}^n$.

\vspace{6pt}

{\bf Theorem 18.} {\it If ${\cal A}$ has a unit element, then $a \rightarrow g_{{\cal A}_a}$ is an isomorphism
of the additive group ${\cal J}({\cal A})$ to} $ S ({\cal A})$.

\vspace{6pt}

{\bf Theorem 19.} {\it If ${\cal A}$ is a commutative algebra over ${\bf R}^n$, then
${\cal J} ({\cal A})$ is a Jordan subalgebra of} ${\cal A}$.

\vspace{6pt}

Moreover, the following theorem holds:

\vspace{6pt}

{\bf Theorem 20.} {\it If ${\cal A}$ is a finite-dimensional commutative algebra
over a field of characteristic different from two or three, then ${\cal J} ({\cal A})$
is a Jordan subalgebra of} ${\cal A}$.

\vspace{6pt}

Concerning the Riccati differential equation in Jordan pairs, Braun [130] proved a result recalled below (see
Theorem 21). Linearization of the matrix Riccati differential equation derived from $(m \times n)$-matrices (see Levin
[442]) and the Riccati differential equation for operators in a Banach space (see Tartar [658]) are assumed to be known
to the reader.

Let $V$ be a Jordan pair with $V^{\sigma}$, $\sigma = \pm $, Banach spaces, and let $D$ and $Q$ be the derivation and the quadratic
representation defined as usually. Let $I$ be an ${\bf R}$-interval, let $\eta $ be an initial
point in $I$, and let $k$ be a given initial value, $k \in V^{+}$. Let $v (\xi), w(\xi)$ be given continuous functions,
$v:I \rightarrow V^{-}$, $w: I \rightarrow V^{+}$, and let $D$ and $Q$ be continuous. The Riccati differential equation
(without linear term) is defined by
$$\frac{\partial x }{\partial \xi} =Q(x) v+ w.$$
The solution $x: I \times I \rightarrow V^{+}$ with initial value $k$ at the point $\eta $
will be denoted by $x (\xi, \eta )$.

\vspace{6pt}

{\bf Notation.} $B(u, t):={\rm Id} -D(u, t) + Q(u) Q(t)$,
$u^t:= B(u, t)^{-1} (u-Q(u)t)$, for $u \in V^{+}$,
$t \in V^{-}$, if the inverse of $B(u,t)$ exists.

\vspace{6pt}

{\bf Theorem 21.} {\it Let $x_0$ be the solution of the Riccati equation with initial value $k=0 $ at $\eta =0$.
Put
$$x (\xi, \eta): = x_0 + h_+ (k)^{h_{-}(z)} $$
$h_{\sigma} : I \times I \rightarrow {\rm Aut}\, V^{\sigma}$, $z : I \times I \rightarrow V^{-}$.
Solve the linear system
$$\frac{\partial h_+}{\partial \xi} =D(x_0, v)h_+, \quad
\frac{\partial h_{-}}{\partial \xi} =-D(v, x_0)h_{-}, $$
so that $\frac{\partial z}{\partial \xi} = h^{-1}_{-} (v)$ with $h_{\sigma }(\eta, \eta) ={\rm Id}$,
$z (\eta, \eta) =0$. Then
$x(\xi, \eta )$ is the solution with initial value $x (\eta, \eta)=k$ $($in a neighborhood of} $\eta$).

\vspace{6pt}

Walcher [699a] gave a characterization of regular Jordan pairs and its application to the Riccati differential eaquation
as follows. Let $V$ be a finite-dimensional vector space over ${\bf R}$, $P:V \rightarrow  {\rm Hom} (V, V)$ a quadratic map,
$G \subset V$ open $(G \neq \emptyset)$, and $\varphi \in C^1 (G, V)$. Suppose that for all $a \in V$ one has
$({\rm d}/ {\rm d}t) \varphi (z(t))=-a$ whenever
$z(t) \subset G$ is a solution of the Riccati differential equation $\dot{x} =P(x)a$.
By differentiation, $D \varphi (x) \cdot P(x)=- {\rm Id}$. Moreover,
Walcher showed that the identity $P(x, P(x) z)y=P(x, P(x)y)z$ is satisfied for all $x,y,z \in V$. Thus there exists a Jordan pair structure $(P,Q_-)$ on
${\cal V} =(V, V)$ and by Theorem 21 the following is true: Let $a:I \rightarrow V^{-}$,
$c:I \rightarrow V^+$, $(B_+, B_-):I \rightarrow {\rm Der}\, {\cal V}$
be continuous. If $z(t)$ solves
$$\dot{x} =P(x)a+B_+x +c $$
and $P(z(t))$ is invertible, then $P(z(t))^{-1} z(t)$ solves
$$\dot{x} =-Q_- (x) c + B_-x-a.$$
Let us recall that a system of ordinary differential equations $\ddot{x} = F(t, x)$ is
said to have a fundamental system of solutions if there exist finitely many solutions that determine (almost) all
other solutions; it is called a system of {\it polynomial differential equations} if, for all values of $t$, $F(t, x)$ is a
polynomial in $x$. A theorem of Lie implies that a system of polynomial differential equations
has a fundamental system of solutions
if $F(t, x)=\sum \lambda_i (t) f_i (x)$ and the polynomials
$f_i(x)$ generate a finite-dimensional subalgebra of the Lie algebra ${\rm Pol}\, V$, where $V$
is the vector space on which the system is defined.

Walcher [699b] determined these subalgebras in the case $\dim V=1$ and
showed that they correspond to the Riccati (including linear) and the
Bernoulli equation. For $\dim V >1$, Walcher investigated the finite-dimen\-sional,
graded subalgebras $L$ of ${\rm Pol}\, V$. Denoting by ${\rm Pol}_i V$ the subspace of all
polynomials of degree $i+1$, it is shown that the semisimplicity of
$$L=L_{-1}\oplus L_0 \oplus \cdots \oplus L_m, $$
with $L_i \subseteq {\rm Pol}_i V$, implies $m=1$.

$L$ is said to be transitive if $L_{-1}=V$. By a result of Kantor, it is known
that a finite-dimensional, graded, transitive subalgebra with
$m >1$ is reducible; that is, there exists a subspace $U$ of $L_{-1}$ with $0 \neq U \neq V$
such that for all $k$ with $0 \leq k \leq m$ and all
$p \in L_k$, $p(V, \ldots , V, U) \subseteq U$. This allows one to reduce the discussion of
transitive subalgebras to those whose degree equals~$1$. The latter are shown to arise
from finite-dimensional Jordan pairs. In case $\dim V=2$, this permits a complete
enumeration of all finite-dimensional, maximal, transitive subalgebras of ${\rm Pol }\, V$ of fixed degree
$m$. Walcher also discussed how these results can be used to find all
 solutions of certain types of systems of polynomial differential equations.

Let us recall now the construction given in 2000 by Liu (see Liu [445b]).

If $\V$ is a finite-dimensional real vector space endowed with an inner pro\-duct denoted by
the dot, then let us consider the real vector space $\M$ as follows 
$$
\M := \{x_0 + g x_n \mid x_0 \in \R, \ x_n \in \V\}
$$
with $x+ y := (x_0+y_0) + g(x_n + y_n)$ for $x,y \in \M$ and define also a product
$$
xy = (x_0 + gx_n)(y_0+gy_n):=(x_0y_0 + x_n \cdot y_n) + g(x_0 y_n + y_0x_n)
$$
where $g$ (called the ``$g$-number" by Liu) satisfies
\vskip6pt

\begin{center}
\begin{tabular}{cccc}
& $\big|$ & 1 & $g$ \\
\hline 
1 & $\big|$ & 1 & $g$ \\
$g$ & $\big|$ & $g$ & 1 
\end{tabular}
\end{center}

\vskip6pt

It is easily to prove that $\mathbb{M}$ is a commutative Jordan algebra.

Liu [445b] has called $\M$ the $g$-{\it based Jordan algebra}.

\vskip6pt

{\bf Remark 13.}
The $g$-based Jordan algebra $\M$ is a particular {\it linear} Jordan algebra. It is the
underling algebraic structure of a dynamical system defined on $\V$ which possesses
one or more constraints.

\vskip6pt

Some applications of this new formulation include the {\it perfect elastoplasticity}
(see Hong and Liu [349a,b]), the {\it magnetic spin equation} (see Landau and Lifshitz [427]),
the {\it suspension particle orientation equation} (see Liu [445a]). They prove the {\it usefulness}
of this new formulation of Liu.

\vskip6pt

{\bf Open Problem.} As Liu suggested, there exists the possibility to describe the non-linear
dissipative phenomena of physical systems by using the above mentioned $g$-based Jordan algebra
$\M$ (see Liu [445b, p.~428]).

\vskip6pt

One year later, in 2001, Liu has proceeded to examin above mentioned type of dynamical
systems from the view point of Lie algebras and Lie groups. Then he has derived a new dynamical
system based on the composition of the $g$-based Jordan algebra and Lie algebras (see
Liu [445c]).

\vskip6pt

{\bf Remark 14.}
Based on the symmetry study, Liu has developed a {\it numerical scheme} which preserves
the group properties for every time increment.

\vskip6pt

{\bf Important remark.} Because the above mentioned scheme is easy to implement numerically and has
high computational efficiency and accuracy, it is highly recommended for engineering applications.

\vskip6pt

In 2002, Liu has examined previous mentioned dynamical systems from the view point of Lie
algebras and Lie groups (see Liu [445d]).

\vskip6pt

{\bf Open Problem.} Consider other {\it particular} Jordan algebras suitable to be the algebraic
foundation for various dynamical systems.
\vskip6pt

On the other hand, it could be formulated also the following

\vskip6pt

{\bf Open Problem.} Taking into account of the fact that the study of linear Jordan algebras
can be included in the more general study of Jordan triple systems, it would be interesting
to develop a {\it more general} algebraic background for (various) dynamical systems, and
finally formulate an {\it unitary} mathematical theory for {\it all} dynamical systems.

\vskip6pt

In 2004, Liu [445e] used the real $g$-based Jordan algebra (defined by
himself in 2000) to the study of the Maxwell equations without appealing
to the imaginary number $i$. In terms of the $g$-based Jordan algebra formulation, the usual
Lorentz gauge condition is found to be a necessary and sufficient condition to render
the second pair Maxwell equations, while the first pair of Maxwell equations is proved to be an
intrinsic algebraic property.

The $g$-based Jordan and Lie algebras are a suitable system to implement the Maxwell equations into
a more compact form.

Finally, Liu has studied in [445d] the problem about a single formula of the
Maxwell equations.

\vskip6pt

{\bf Remark 15.}
In 1966, Hestenes has proved in his book [341] that -- in terms
of spacetime algebra (16 components) -- the four Maxwell equations can be organized into a
single one. Similarly, Liu [445e] has achieved his goal in his algebraic formulation. See also the recent paper by Liu [445f].

\vskip6pt

{\bf Remark 16.}
It is impressive that $g$-based Jordan (and Lie) algebras of Liu can be so usefull to different topics.

\vskip6pt

2. If $M=G/K$ is a bounded symmetric domain and $S=K/L$ is its Shilov boundary, then one can define
a Poisson kernel on $M \times S$ and the Poisson integral for any hyperfunction on $S$.
An open problem, formulated more than twenty years ago by Stein, is to characterize these
integrals as solutions of a system of differential
equations, established for certain cases by Hua (see~[352c]).

In [432a], Lassalle dealt with bounded symmetric domains of {\it tube type}. Poisson integrals over the Shilov
boundary are then characterized by the system of differential equation given by Johnson and Kor\'anyi [375].
Results of Berline and Vergne [93] for the domain $(I)_{n,n}$ and that of Kor\'anyi and
Malliavin [411] for the Siegel disc of dimension two, prove that the Johnson-Kor\'anyi system [375] has,
for these particular cases, too many equations. Lassalle proved that this is a general property.
In particular, he established that among
the $\dim K$ differential equations of the Johnson-Kor\'anyi system, a subsystem of $\dim S$
equations is sufficient to characterize the bounded
functions on $M$ which are Poisson integrals of a function on $S$. This new characterization
has a very natural interpretation in terms of Jordan algebras (see Lassalle [432a, pp. 326--327]).

As Lassalle [432b] proved, such an interpretation is also possible if $M$ is a symmetric
Hermitian space of tube type with Shilov boundary $S$
and can be realized as a bounded symmetric domain. The main idea is to formulate the Hua
differential equations [432c] in terms of ``polar coordinates" with
respect~to~$S$.

Consider, again, a bounded symmetric domain $D$, $S$ its Shilov boundary and let $S(z,u)$ be
the Szeg\"{o} kernel on $D \times S$. Hua [352c] was the first to calculate explicitly the expression
of $S(z, u)$ for each of the four series of irreducible domains. Later Kor\'anyi gave a general
proof that made Hua's case-by-case calculations unnecessary. However, Kor\'anyi's proof is not direct;
it uses the unbounded realization of $D$ as ``generalized half-plane"; in unbounded realization
the Szeg\"{o} kernel is given by an integral, which could be calculated by methods of Bochner
and Gindikin. In bounded realization, on the other hand, the Szeg\"{o}
kernel is given by a Fourier series, and for this no methods of calculation had been devised.

Lassalle's result [432h] offers a solution to this problem: he manages to calculate the Fourier series
defining $S(z,u)$ directly, without going through the unbounded realization of $D$. In fact, Lassalle is in position
to solve the following much more difficult problem: For every positive real number $\lambda $, what is the Fourier series
expansion of $(S(z, u))^{\lambda}$? This difficult problem had been open for nearly thirty years.
The only known solution had been given
implicitly by Hua [352c, p.~25] in the particular case of an irreducible domain of type $I_{n,m}$.
But the solution was unknown in all the other cases, including
each of the three other series of irreducible domains.

Lassalle's goal in [432h] is to present a general answer to this problem, one independent of
any classification argument. What is noteworthy in [432h] is that the framework
and tools of Lassalle's proof are provided by Jordan algebra theory. In particular,
his central result is a ``binomial formula" in the complexification of a formally real Jordan algebra.
The solution thus obtained is particularly simple and natural.

\vskip6pt

3. In a series of papers, Clerc and \O rsted [177a,~b,~c], Clerc [174b,~c], Clerc \& Neeb [176], and Clerc and
Koufany [175], defined and studied the generalization of Maslov index by making use of formally real
Jordan algebras (and their complexifications), as well as of Hermitian Jordan triple systems,
which turn out to be very convenient for stating and proving the results. Some details on these papers are pointed out in Iord\u anescu [364w, pp. 46-47].

Another paper by Clerc, which does not belong to this series, but uses Jordan algebras as main tool is [174d]. 

\vskip6pt

4. Whereas the study of Toeplitz operators for the strongly pseudoconvex domains uses methods of partial
differential equations, their structure and Toeplitz $C^{\ast}$-algebras over symmetric domains
is closely related to the Jordan algebraic
structure underlying these domains (see Upmeier [682j, Section 2]). The relation between bounded
Toeplitz operators and Weyl operators of boson quantum
mechanics was examined by Berger and Coburn in [90].

As Upmeier pointed out in [682j, p.~42] even though finite-dimensional bounded symmetric
domains have been classified and
their geometry is totally understood, there are still many open problems concerning their analysis
(i.e., the structure of function
spaces defined over these domains). Since symmetric domains are homogeneous under a (semisimple)
Lie group, these problems are related to harmonic
analysis and the theory of group representations. On the other hand, the occurring function spaces
are often Hilbert spaces of holomorphic functions which give rise
to (reproducing) kernel functions. These kernel functions can be defined in terms of certain basic
``norm functions" derived from the Jordan algebraic structure.

Ion and Scutaru [363], and Ion [362d,~e] introduced new scattering theories {\it via} optimal states.
These states are reproducing kernels in the Hilbert spaces of scattering matrices,
just as the coherent states are reproducing kernels in the Hilbert spaces of wave functions.

Using reproducing kernels, Upmeier [682j, Section 6] outlined a quatization procedure
for certain curved phase spaces
of possibly infinite dimension, namely the ``symmetric Hilbert domains". In the
finite-dimensional setting, Berezin [88a,~b] has considered quatizations for more
general complex (K\"{a}hler) manifolds.

The general formalism for quantum fields on any reproducing kernel Hilbert space is presented
by Schroeck [601], along with a discussion of the operator and distribution properties of those
fields. The Galilean and Poincar\'{e} examples are given along with considerations of the
general relativistic cases.

\vskip6pt

Let us firstly recall that if $G$ is a semi-simple Lie group, and $\sigma : G \to U(H_\sigma)$,
$\tau : G \to U(H_\tau)$ the two (irreducible) unitary representations of $G$,
acting on Hilbert space $H_\sigma$ and, respectively, $H_\tau$, 
a unitary isomorphism $\Phi : H_\sigma \to H_\tau$ is called an {\it intertwining operator}
if the identity 
$$
\Phi\, \sigma(y)=\tau(y)\Phi
$$
holds for all $y \in G$.

\vskip6pt

{\bf Remark 17.}
The construction of intertwining operators (e.g., Poisson integral) is an important
method in harmonic analysis.

\vskip6pt

In the Jordan theoretic framework, the Jordan determinant function and the associated
differential operator give rise to intertwining operators of the group $G$
of all biholomorphic automorphisms of a tube domain, which
generalize the Capelli identity from classical invariant theory.

\vskip6pt

{\bf Remark 18.} 
The Capelli identity (see Capelli [153]) was a centerpiece of the invariant theory 
in the $19^{\rm th}$ century.

\vskip6pt

In 1991, Kostant \& Sahi [414a] have proved a large class of identities generalizing
the Capelli identity.
Their approach -- which have led to the above mentioned generalization -- was
to regard the $n^2$-dimensional space $M(n,{\bf R})$ of $n\times n$ real matrices
{\it not} as a Lie algebra, but rather a {\it Jordan} algebra (with the Jordan product
$a \circ b := (ab + ba)/2$,  $a,b \in M(n,{\bf R})$).

In 1992, Sahi published the paper [583b] devoted to Capelli identity and unitary representation,
essentially based on Kostant \& Sahi [414a].

In 1993, Kostant \& Sahi [414b] -- making use of their previous paper [414a] and Sahi's paper
[583b] -- have established a connection between semisimple Jordan algebras
and certain invariant differential operators on symme\-tric spaces. They also
proved an identity for such operators which generalizes the classical Capelli
identity.

In 1995, \O rsted \& Zhang [525] have studied the 
composition series of certain generalized principal series representations
of the automorphism group of a bounded symmetric space of tube type. As
applications they obtained new proofs of the Capelli identity
of Kostant \& Sahi [414a] and some results of Faraut \& Kor\'anyi [246a,b], and
they gave the full decomposition of the $L^2$-space on the Shilov boundary.
The Shilov boundary of a tube domain can also be viewed as a compactification
of a formally real Jordan algebra, the automorphism group is then the
conformal group of the formally real Jordan algebra. Zhang has studied in [730] the case
of non-formally real Jordan algebra. (N.B. After the main part of [730] was completed,
Zhang was informed by Sahi that he has also obtained in [583d] most of
the results contained in~[730].)

A related construction using the quasi-inverse in Jordan triples leads to the transvectants
introduced by Peetre in his study of Hankel forms of arbitrary weight over symmetric
domains [533]. The word {\it transvectant} (in German: {\it \"Uberschiebung}) comes
from classical invariant theory, where objects called transvectants were defined by Gordan 
(see P. Gordan, {\it Invariantentheorie}, Teubner, Leipzig, 1887).
Janson \& Peetre [372] have ``rediscovered'' it hundred years later! Let us recall that
classical invariant theory is mainly about $SL(2,{\bf C})$. The problem of Peetre in [533]
was thus to generalize the transvectant to the case of an arbitrary semi-simple Lie group.
There are three types of approaches to symmetric domains, namely: 1) the case by case study
(see Hua [352c]),  2) the Lie approach (see Helgason [335a,b]), and 3) the Jordan approach
(see Loos [448h] and Upmeier [682j]). Peetre [533] used the approach of type 3), taking some 
advantage of the Jordan triple system structure.

In 2004, Peng \& Zhang [534], studying the tensor products of holomorphic representations
and bilinear differential operators, gave the irreducible decomposition of the 
tensor product of the representations for any two unitary weights and they have found
the highest weight vectors of the irreducible components. Peng \& Zhang have also
found -- by using some refinements of the ideas of Peetre [533] -- certain bilinear
differential intertwining operators realizing the decomposition,
and they generalize the classical transvectants in invariant theory of $SL(2,{\bf C})$.

The Kirillov orbit method gives a realization of irreducible representations of a Lie group
$G$ in terms of the orbits of $G$ in the dual space $\frak{g}^\sharp$ of the Lie algebra
$\frak{g}$. 

Representations of $G$ corresponding to nilpotent orbits are called
{\it unipotent}. The explicit realization of unipotent representations (in the semisimple
case) is an important, and not completely solved, problem in harmonic analysis. For details, we
refer the reader to Vogan [695].

In the Jordan theoretic framework, decisive progress in this direction has been made by 
Sahi (in collaboration with Dvorsky) (see [583a,c,d], [225a,b]).

\vskip6pt

{\bf Remark 19.}
The Jordan algebra determinant and its powers is -- again -- the starting point
of this construction.

\vskip6pt

The paper [225f] by Dvorsky and Sahi is a culmination of a series of papers dedicated
to the problem of constructing explicit analytic models for small unitary representations of 
certain semisimple Lie groups. Let us remark here that -- for instance -- Sahi [583a], studying
explicit Hilbert spaces for certain unipotent representations, has made use of some 
Jordan algebra results from Braun \& Koecher [131], Koecher [408c], Kor\'anyi \& Wolf [412], and
Kostant \& Sahi [414a], while the paper Dvorsky \& Sahi [225a] is devoted to {\it non-formally real}
Jordan algebra case.

In the paper Dvorsky \& Sahi [225b], the authors construct a family of small unitary representations
for {\it real} semi-simple Lie groups associated with Jordan algebras. These representations
are realized on $L^2$-spaces of certain orbits in the Jordan algebra. The representations are
spherical and one of the authors' key results is a precise $L^2$-estimate for the Fourier transform
of the spherical vector. The authors also consider the tensor products of these representations
and describe their decomposition.

\vskip6pt

{\bf Notations.} Let $G$ be a simple Lie group with Lie algebra $\frak{g}$, and let $K$
be the maximal compact subgroup corresponding to a Cartan involution $\theta$.
Suppose $G$ has a parabolic subgroup $P=LN$ such that:

(i) the nilradical $N$ is abelian, and

(ii) $P$ is conjugate to $\overline{P} =\theta(P)$.

\vskip6pt

As Sahi has mentioned [583d, p.~1], the spherical (degenerate) $\overline{P}$-prin\-cipal
series representations of $G$ are obtained by starting with a positive real character
of $L$, extending trivially to $\overline{P}$, and inducing up to $G$. For such a representation,
(i) implies that the $K$-types have multiplicity 1 and (ii) implies that each irreducible constituent
has an invariant Hermitian form.

In the paper [583d], Sahi has provided a rather detailed analysis of these representations. He has
explicitly  determined the $K$-types of their irreducible constituents and the signature of 
the Hermitian form on each $K$-type.

\vskip6pt

{\bf Remark 20.}
The results contained in the paper [583d] extend and generalize those of the previous
papers by Sahi [583b] and [583c], as well as those from several other papers.

\vskip6pt

{\bf Remark 21.}
The unitary representations described in Theorems 4B and 5A from [583d] are
of particular interest since they are all {\it unipotent}
and correspond to some of the smallest nilpotent coadjoint orbits of the
group $G$. As Sahi has mentioned [583d, p.~2], it is expected that these
representations (and their analogs over other fields) will have some interesting
applications.

\vskip6pt

Concerning the applications of Jordan structures to analysis I like to refer the reader to the {\bf very recent} book [166i] by Chu, not yet published, but available from January 2012.

\vspace{1cm}

\centerline{\bf \Large \S 4. JORDAN STRUCTURES}
\vskip4pt
\centerline{\bf \Large IN DIFFERENTIAL GEOMETRY}

\vspace{48pt}

Let ${\cal A}$ be a formally real Jordan algebra of
dimension $n$. Theorem 3.4 from Braun \& Koecher's book
[131, Ch. XI] implies that ${\cal A}$ has a unit
element, which we shall denote by $e$.

In this case, by Proposition 6 from \S 1, we have
$${\rm Idemp}_1 ({\cal A})= \{c \mid c \in {\rm Idemp}
({\cal A}),\, c  \mbox{ primitive}\}.$$

\vspace{6pt}

{\bf Definition 1.} A system of idempotents
$c_1,\ldots, c_s \in {\cal A}$ is called a {\it complete
orthogonal system of idempotents} of ${\cal A}$ {\it if}
$\sum\limits_{i=1}^{s} c_i=e
\mbox{ and } c_ic_j=\delta_{ij}c_i$ $(i,j=1, \ldots,s)$.

\vspace{6pt}

{\bf Proposition 1.} {\it A formally real Jordan
algebra contains a complete orthogonal system of
primitive idempotents.}

\vspace{6pt}

{\bf Comments}. Tillier [662] gave a geometric
characterization of primitive idempotents in a
formally real Jordan algebra, namely: every primitive
idempotent belongs to an extremal ray of the domain
of positivity of the algebra, and, conversely, such a
ray always contains a primitive idempotent.

\vspace{6pt}

{\bf Propostion 2.} {\it All complete orthogonal systems
of primitive idempotents of a formally real Jordan
algebra have the same number of elements.}

\vspace{6pt}

{\bf Definition 2.} The number of elements of a
complete orthogonal system of primitive
idempotents of a formally real Jordan algebra ${\cal A}$
is called the {\it degree of} ${\cal A}$.

\vspace{6pt}

In 1965, Hirzebruch [346] showed that the set of
primitive idempotents in a {\it finite-dimensional}
simple formally real Jordan algebra is a compact
Riemannian symmetric space of rank one and that any such
space arises in this~way.

Ten years later, Neher undertook in [512a] a detailed
differential-geometric study of idempotents in a {\it
real} Jordan algebra.

\vskip6pt

Let us recall from Hirzebruch [346] some of his
important results.

Suppose that ${\cal A}$ is {\it simple} and denote
its degree by $s$. Then the form
$$\mu(u):=\frac{n}{s} {\rm Tr}\, L(u), \quad u \in
{\cal A}, $$
is an {\it associative} (i.e., $\mu(x(yz))=\mu((xy)z)$
for any $x,y,z \in {\cal A})$ {\it linear form} on
${\cal A}$ with $\mu (c) =1$ for every $c \in {\rm
Idemp}_1 ({\cal A})$.

\vspace{6pt}

{\bf Remark.} Suppose that a formally real Jordan
algebra is not simple. Then it is semisimple
(therefore it is a sum of simple ideals), and the
associative linear form $\mu $ with value $1$ on the
primitive idempotents is constructed by means on the
forms $\mu_i$ on the components.

\vspace{6pt}

{\bf Notation.} For every $c \in {\rm Idemp}_1({\cal
A})$ define $S_c$ by
$$S_c: =\{x \mid x \in {\cal A}_{1/2} (c), \,
\mu(x^2)=2\}.$$

{\bf Theorem 3.} {\it Let ${\cal A}$ be a simple
formally real Jordan algebra and let $c \in {\rm
Idemp}_1 ({\cal A})$. For every $d \in {\rm Idemp}_1
({\cal A})$ there exists a unique real number $t$, $0 \leq t
\leq \pi/2$, and a unique element $x$ in $S_c$ such that
$d={\rm d}(t)$, where
$${\rm d}(t)=(\cos 2t)c+\left(\frac{1}{2} \sin 2t\right)x+
\frac{1}{2}(1-\cos 2t)x^2.$$
Conversely, for each such $t$, ${\rm d}(t)$ is an element of
${\rm Idemp}_1({\cal A}) $. The primitive idempotents
which are orthogonal to $c$ are exactly those of the
form $x^2-c$ with $x \in S_c$. For $x \in S_c$, $x^2
=c+d$ if and only if} $x \in {\cal A}_{1/2} (c) \cap
{\cal A}_{1/2} (d)$.

\vspace{6pt}

{\bf Corollary.} {\it A formally real Jordan algebra is
simple if and only if the set of its primitive
idempotents is connected.}

\vspace{6pt}

{\bf Theorem 4.} {\it Let ${\cal A}$ be a simple
formally real Jordan algebra and let $c$ be an
element of ${\rm Idemp}_1({\cal A})$. For $x, y \in
S_c$ there exists a product of Peirce reflections
with respect to idempotents of ${\rm Idemp}_1({\cal
A})$ that fixes $c$ and maps $x$ to} $y$.

\vspace{6pt}

{\bf Theorem 5.} {\it Let ${\cal A}$ be a simple
formally real Jordan algebra and let $c_1, c_2, d_1,
d_2 \in {\rm Idemp}_1({\cal A})$ such that $\mu (c_1
c_2)=\mu(d_1 d_2)$. Then there exists a product of
Peirce reflections with respect to idempotents of
${\rm Idemp}_1({\cal A})$ which maps $c_1$ to $d_1$
and $c_2$ to} $d_2$.

\vspace{6pt}

Consider now on ${\cal A}$ the {\it symmetric
bilinear form} $\mu (x, y):=\mu (xy)$ and the {\it
Euclidean metric} $\rho_E(x,y):=(\mu ((x-y)^2))^{1/2}$
it determines.

\vspace{6pt}

{\bf Remark.} The automorphisms of ${\cal A}$ are
isometries of the metric space $({\rm Idemp}_1({\cal
A}), \rho_E)$.

\vspace{6pt}

{\bf Definition 3.} A metric space $(M, \rho )$ is
called {\it two-point homogeneous} if there exists an
isometry $\iota $ of $M$ such that $\iota (c_1)=d_1$
and
$\iota (c_2)=d_2$ for any $c_1, c_2, d_1, d_2 \in M$
with
$\rho (c_1, c_2)=\rho (d_1, d_2)$.

\vspace{6pt}

{\bf Corollary of Theorem 5.} {\it If ${\cal A}$ is a
simple formally real Jordan algebra, then $({\rm
Idemp}_1({\cal A}), \rho_E)$ is a connected, compact
and two-point homogeneous metric space}.

\vspace{6pt}

{\bf Theorem 6.} {\it Let ${\cal A}$ be a simple
formally real Jordan algebra and let $T$ be a
one-to-one map of ${\rm Idemp}_1({\cal A})$ onto
itself such that $\mu (T(c), T(d))=$ $=\mu(c, d)$ for all
$c, d \in {\rm Idemp}_1({\cal A})$. Then $T$ can be
extended uniquely to an automorphism of} ${\cal A}$.

\vspace{6pt}

{\bf Notation.} Consider the Riemannian structure
induced on ${\rm Idemp}_1({\cal A})$ by $\mu (xy)$.\
The Riemannian manifold thus obtained will be
denoted by\break $({\rm Idemp}_1({\cal A}), R)$.

\vspace{6pt}

{\bf Remarks.} The automorphisms of ${\cal A}$ are
isometries of the Riemannian manifold $({\rm
Idemp}_1({\cal A}), R)$. $({\rm Idemp}_1({\cal A}), R)$ is a
compact symmetric Riemannian space.

\vspace{6pt}

{\bf Notation.} The Riemannian distance between two
elements $c$ and $d$ of $({\rm Idemp}_1({\cal A}),
R)$ will be denoted by $\rho_R(c, d)$.
\vspace{6pt}

Since the relations
$$0 \leq \rho_R (c, d) \leq \pi/\sqrt{2} \quad
\mbox{and} \quad \rho_E(c, d)=\sqrt{2}
\sin\left(\frac{1}{\sqrt{2}} \rho_R(c,d)\right)$$
hold, it follows that $\rho_R(c_1, c_2)=\rho_R(d_1,
d_2)$ is equivalent to
$\mu(c_1, c_2)=\mu(d_1,d_2)$ for $c_1,c_2,d_1,d_2 \in
{\rm Idemp}_1({\cal A})$.

\vspace{6pt}

{\bf Remark.} Consequently, $({\rm Idemp}_1({\cal
A}), R)$ is a two-point homogeneous symmetric
Riemannian space and hence (see Helgason [335a, p.~355])
of rank~one.

\vspace{6pt}

Let $c$ be an element of ${\rm Idemp}_1({\cal A})$.
Clearly, $\rho_R(c, d)$, $d \in {\rm Idemp}_1({\cal
A})$, is maximal only when $\rho_E(c,d)$ is maximal.
Because $\mu (c)=\mu(d)=1$, we have $\rho_E(c,
d)=\sqrt{2} \sqrt{1-\mu(cd)}$, which is maximal only
when $\mu(cd)=0$, i.e., when $cd=0$.

\vspace{6pt}

{\bf Notation.} For every $c \in {\rm Idemp}_1({\cal
A})$ we define
$$A_c := \{d \mid d \in ({\rm Idemp}_1({\cal A}), R), \,
cd=0\}.$$

{\bf Remark.} By Theorem 3 we have $A_c =
\{x^2-c \mid x \in S_c\}$.

\vspace{6pt}

{\bf Remark.} $A_c$ is a submanifold of $({\rm
Idemp}_1({\cal A}), R)$ and is called the {\it
antipodal manifold of} $c$.

\vspace{6pt}

{\bf Notation.} For any simple formally real Jordan
algebra ${\cal A}$ of dimension $n>1$ there exists a
natural number $q({\cal A})$ such that, for every
pair of orthogonal primitive idempotents $c_1,c_2 \in
{\cal A}$, the relation $\dim ({\cal
A}_{1/2}(c_1)\cap {\cal A}_{1/2}(c_2))=$ $= q({\cal A})$
holds. If $s=s({\cal A})$ denotes the degree of
${\cal A}$, then ${\cal A}$ is said to be of {\it
type} $(s, q({\cal A}))$.

\vspace{6pt}

{\bf Remark.} If ${\cal A}_1$ and ${\cal A}_2$
are simple formally real Jordan algebras with
$s({\cal A}_1) = s({\cal A}_2)$ and $q({\cal A}_1) =
q({\cal A}_2)$, then ${\cal A}_1$ and ${\cal A}_2$
are isomorphic.

\vspace{6pt}

{\bf Comments.} It would be interesting to extend
Hirzebruch's results [346, pp. 350--351] on Betti
numbers of ${\rm Idemp}_1({\cal A})$, critical points
of differentiable functions on ${\rm Idemp}_1({\cal
A})$, etc., to other kinds of Jordan algebras.

\vspace{6pt}

Using the well-known classification of compact
symmetric Riemannian spaces of rank one, Hirzebruch
[346] proved that each of these spaces can be
described in terms of a suitable formally real Jordan
algebra, namely:

a) Type $(1,0)$; ${\cal A} = {\bf R}$, and ${\rm
Idemp}_1({\cal A})$ consists of point alone.

b) Type $(2,q)$; $q \geq 1$. Let $V'$ be a
$(q+1)$-dimensional vector space over ${\bf R}$ and let
$\sigma $ be a positive definite bilinear form on
$V'$. Define $V:={\bf R}e \oplus V'$ a bilinear
product by $uv:=\sigma (u,v) e$ for $u,v \in V'$, $e$
being the unit element. It is immediate, that $V$
endowed with this product is a Jordan algebra $J(Q)$
(as in Theorem 8 from \S 1). ${\rm Idemp}_1(V)$ is
homeomorphic to the $q$-dimensional sphere $S^q$.

c) Type $(s,1)$; $s \geq 3$. Let $V$ be the vector
space of symmetric $(s \times s)$- matrices over
${\bf R}$. For $A, B \in V$, let $AB: = \frac{1}{2}
(A\cdot B + B \cdot A)$, where $A \cdot B$ denotes the
usual matrix product in $V$. We have $\mu (A)= {\rm
Tr}\, A$ and ${\rm Idemp}_1(V)= \{A \mid A \in V,\, A^2 =A,\,
{\rm Tr}\, A=1\}$. It follows that ${\rm Idemp}_1(V)$
is homeomorphic with the real $(s-1)$-dimensional
projective space $P_{s-1} ({\bf R})$. For $c \in {\rm
Idemp}_1(V)$, the antipodal manifold $A_c$ is
$P_{s-2} ({\bf R})$.

d) Type $(s,2)$; $s \geq 3$. Let $V$ be the ordinary
real vector space of complex Hermitian $(s \times
s)$-matrices. Define the product as in c). Then
${\rm Idemp}_1(V)$ is homeomorphic with $P_{s-1}
({\bf C})$. For $c \in {\rm Idemp}_1(V)$, $A_c$ is $P_{s-2}
({\bf C})$.

e) Type $(s, 4)$; $s \geq 3$. Let $V$ be the real
vector space of Hermitian $(s \times s)$-matrices
over ${\bf H}$. Define the product as in c). Then ${\rm
Idemp}_1(V)$ is homeomorphic with $P_{s-1} ({\bf H})$. For $c
\in {\rm Idemp}_1(V)$, $A_c$ is $P_{s-2} ({\bf H})$.

f) Type $(3,8)$. Let $V $ be the real vector space of
Hermitian $(3 \times 3)$-matrices over ${\bf O}$. Define
the product as in c). Then ${\rm Idemp}_1(V)$ is
homeomorphic with the projective octonion plane. For
$c \in {\rm Idemp}_1(V)$, $A_c$ is an eight-dimensional sphere.

\vspace{6pt}

{\bf Remark.} The geodesics through a point $c
\in {\rm Idemp}_1({\cal A})$, as a set of points, are
${\rm Idemp}_1({\cal A}) \cap {\cal A}(c, x)$, $x \in
{\cal A}_{1/2} (c)$, ${\cal A} (c, x)$ being a simple
three-dimensional subalgebra of ${\cal A}$
containing $c$ (see Hirzebruch [346, p.~348]). For a
detailed discussion on geodesics in a more general
case see Neher [512a].

\vspace{6pt}

{\bf Comments.}
I must point out that the long series of applications of Jordan structures to
differential geometry has its roots in 1962 lecture notes by Koecher [408c], and
in the papers [387a,~b] from 1964 and 1966 by Kantor (but published -- at that time --
{\it only} in Russian!).

\vspace{6pt}

In 1993, as a natural step of generalization of the
just recalled Hirzebruch's results, Nomura [518a]
considered an {\it infinite-dimensional} analogue of
formally real Jordan algebras and treated the set of
primitive idempotents. As it is the associative inner
product (instead of the algebraic formally real property what
plays an important role in the study of
Hirzebruch [346]), Nomura [518a] based his study on
Jordan-Hilbert algebras (i.e., real Hilbert spaces
with inner products $\langle \,\cdot \mid \cdot \,\rangle$, which are
also real Jordan algebras and $\langle xy\mid z\rangle =
\langle y\mid xz\rangle$ for
all their elements $x,y,z$). He deals with {\it the set
$J_1$ of primitive indempotents} as a Riemann Hilbert
manifold. He proved that $J_1$ is two-point
homgeneous and derived a unified formula for the
sectional curvature of $J_1$ (see Nomura [518a, Th.
5.3 and, resp., Th. 6.4]).

Let us point out the matters which do not occur in
the finite-dimen\-sional~cases.

Firstly, Jordan-Hilbert algebras do not necessarily
have a unit element. On the other hand, the
adjunction of unit element does not in general agree
with the Hilbert space structure. However, this lack
of unit element is compensated to some extent by (a
version of) McCrimmon's result (see McCrimmon, K.,
{\it Peirce ideals in Jordan
algebras}, Pacific J. Math. {\bf 78}(1978), 379--414,
Th.~1.11 and Th. 1.12).

Nomura's version of McCrimmon's result states that
in a topologically simple (non-trivial) Jordan-Hilbert
algebra, the Peirce $1$-spaces is also topologically
simple (see Nomura [518a, Prop. 1.6]). This enables
Nomura to carry out computations concerning
idempotents as in the finite-dimensional cases.

Secondly, infinite-dimensional connected complete
(in the sense of Riemann distance) Riemannian
manifolds may carry two points which cannot be
joined by a minimal geodesic (see Grossman [310] and
Klingenberg [402b, p.~127]). Because of this
possibility of missing minimal geodesic, Nomura
computed the Riemannian distance on the set $J_1$
along the standard line of textbooks on Riemannian
geometry: the inclusion map of $J_1$ into the
ambient Hilbert space being an embedding, Nomura
defined the canonical Levi-Civita connection,
exhibited a geodesic, examined the diffeomorphism
domain of the exponential mapping and made use of
Gauss lemma to derive the minimality of geodesic.

One year later, in 1994, Nomura [518b] extended his
above mentioned study to the space $J_p$ of rank $p$
idempotents in a topologically simple Jordan-Hilbert
algebra $V$. To provide the necessary
differential-geometric structure on $J_p$,
subalgebras with two idempotents generators are
studied with the aid of the Peirce decomposition and
associated Jordan-Hilbert systems. This allows the
construction of an atlas and permits to Nomura to
identify the tangent space at a point with the
corresponding Peirce $\frac{1}{2}$-space. Geodesics
in $J_p$ turn out to be orbits of one-parameter
subgroups of ${\rm Aut} \, V$, which can be
represented explicitly, and thus permit the
computation of the Riemann distance and sectional
curvature.

In December 1997, Chu and Isidro [168]
extended Nomura's results to the manifold of
extremal (minimal or maximal) projections in the
complex $C^{\ast}$-algebra ${\cal L}(H)$ of bounded
linear operators in a Hilbert space $H$. In 2002, Isidro and Mackey [368] have extended the results from [168] on the manifold
of minimal projections in ${\cal L} (H)$ to the manifolds of finite rank projections in $\mathcal{L}(H)$.

\vskip6pt

A symmetric Riemannian space, defined as usual, is a Riemannian manifold
such that the geodesic symmetry $S_{x}$ around every point $x$ is an isometry.
By writting $x \cdot y : = S_{x}(y)$, Loos has obtained a (nonassociative)
multiplication on the manifold satisfying certain algebraic identities, which
in turn suffice to characterize symmetric spaces. In this way, one obtains an
elementary ``algebraic" definition of a symmetric space not involving the
manifold structure of the underlying topological space. This definition was
firstly given by Loos in his remarkable Ph.D. Thesis [448a].

\vspace{6pt}

{\bf Definition 4.} (see Loos [448a]). A manifold
${\cal M}$ with a differentiable
multiplication ${\cal M} \times {\cal M} \rightarrow {\cal M}$, denoted
$(x, y) \rightarrow x \cdot y$, and having the properties

(1) $x \cdot x = x$;

(2) $x \cdot (x \cdot y) = y$;

(3) $x \cdot (y \cdot z) = (x \cdot y) \cdot (x \cdot z)$;

(4) every $x$ has a neighbourhood $U$ such that $x \cdot y = y$ implies
$y = x$ for\linebreak 
\phantom{aaaaaaa} all $y$ in $U$,

\noindent is called a {\it symmetric space.}

\vspace{6pt}

{\bf Note.} In Loos' terminology, a manifold is a differentiable manifold of
class $C^{\infty}$ which is Hausdorff and paracompact as a topological space.
It may have several connected components, which may be of different (yet
finite) dimensions.

\vspace{6pt}

{\bf Comments.} In contemporary mathematics there exists the notion of {\bf quandle} (related to knot theory), introduced in 1982 by Joyce [381]. Symmetric spaces in the sense of Loos are important examples of quandles (see -  for instance - the paper by Buliga [143]).

\vskip6pt

{\bf Remark.} Spaces satisfying only (1), (2), and (3)
(``reflection spaces") have been studied by Loos in
[372b]. They turn out
to be fibre bundles over symmetric spaces
(see, for instance, Neher [512a]).

\vspace{6pt}

{\bf Note.} As Vanhecke has pointed out to me\footnote{Private
communication on February $21^{\rm st}$, 1994.}, it seems that
all the globally KTS-spaces studied by himself together with
J.C. and M.C. Gonz\'{a}lez-D\'{a}vila are examples of reflection spaces. In fact, two years later, they have mentioned: ``It is worthwhile to note that
the globally KTS-spaces provide a large class of examples of reflection
spaces.'' (see [299, p.~322]).

\vspace{6pt}

{\bf Definition 5.} Left multiplication by $x$ in ${\cal M}$ is
denoted by $S_{x}$, i.e., $S_{x}(y) = x \cdot y$ for all
$x, y \in {\cal M}$, and is called {\it symmetry around} $x$.

\vspace{6pt}

{\bf Remark.} The following properties are immediate:

\hskip3pt{\rm (i)} $x$ is an isolated fixed point of $S_{x}$;

\rm{(ii)} $S_{x}$ is an {\rm Inv}olutive automorphism of ${\cal M}$.

\vspace{6pt}

{\bf Examples of symmetric spaces} (in the sense of the above {\bf Definition 4}): Lie groups, spheres, Grassmann manifolds, Jordan algebras, homogeneous spaces and spaces of symmetric elements.

\vskip6pt

In 1981, Lutz [451] introduced the notion of $\Gamma$-{\it symmetric spaces} which is a generalization of the classical notion of symmetric space. In 2008, Bahturin and Goze [64] defined $\mathbb{Z}_2 \times \mathbb{Z}_2$-{\it symmetric spaces}, and very recently Kollroos [410] gave a classification of the $\mathbb{Z}_2 \times \mathbb{Z}_2$-symmetric spaces $G/K$, where $G$ is an exceptional compact Lie group or Spin(8), complementing recent results of Bahturin and Goze.

\vskip6pt

Let ${\cal A}$ be a real Jordan algebra of dimension $n$ and with unit
element~$e$.

\vspace{6pt}

{\bf Notation.} The mutation of ${\cal A}$ with respect to
$q^{-1}$, $q \in {\rm Inv}({\cal A})$, will be denoted by ${\cal A}^{q}$,
${\rm Inv}({\cal A})$ denoting the set of invertible elements of ${\cal A}$.

\vspace{6pt}

{\bf Remarks.} The product of two elements $a, b \in {\cal A}^{q}$
is given by $a \perp b =  \\
= a(bq^{-1}) + b(aq^{-1}) - (a b) q^{-1}$. Propositions 3 and 4 from \S 1 imply that the mutation
${\cal A}^{q}$ is a Jordan algebra with unit element $q$, and that the
quadratic representation $P_{q}$ of ${\cal A}^{q}$ is given by
$P_{q}(a) = P(a) P^{-1}(q)$, $a \in {\cal A}^{q}$. One can also see that ${\rm Inv} ({\cal A}^{q}) = {\rm Inv}
({\cal A})$ and $\Gamma({\cal A}^{q}) =$ $= \Gamma({\cal A})$. The connected component of ${\rm Inv}({\cal A})$ containing $e$,
denoted as usually by ${\rm Inv}_{\circ}({\cal A})$,
with the multiplication
$q \cdot p : = P(q) p^{-1}$ becomes a symmetric space in the sense of
Definition 4.

\vspace{6pt}

Suppose now that ${\cal A}$ is endowed with an involution
${\cal J}$ (i.e.,
${\cal J} \in {\rm Aut} ({\cal A})$, ${\cal J}^{2} = {\rm Id}$).

\vspace{6pt}

{\bf Notation.} Write ${\rm Inv} ({\cal A, J}) : = \{a \mid a \in {\rm Inv} ({\cal A}),\,
a^{-1} = {\cal J} a\}$ and denote by ${\rm Inv}_{\circ} ({\cal A, J})$ the component
of ${\rm Inv} ({\cal A, J})$ containing $e$. For every $q$ of ${\rm Inv} ({\cal A, J})$
define ${\cal J}_{q} : = P(q) {\cal J}$.

\vskip6pt

{\bf Definition 6.} In the real vector space ${\cal A}$ we define a
{\it new product} of any two elements $a,b \in {\cal A}$, denoted by $a * b$,
as follows
$$2(a*b) : = ab + a ({\cal J}(b)) + b({\cal J}(a))  - {\cal J}(ab)$$
and the real algebra defined by means of $*$ in the vector space
${\cal A}$ will be denoted by
${\cal A}_{\cal J}$.

\vspace{6pt}

{\bf Remark.} The fact that ${\cal A}$ is a Jordan algebra implies that
${\cal A}_{\cal J}$ is also a Jordan algebra,
${\cal J} \in {\rm Aut} ({\cal A}_{\cal J})$, $({\cal A}_{\cal J})_{\cal J} =
{\cal A}$ and $\lambda_{\cal J}$ (defined by
$\lambda_{\cal J}(a,b) : = \lambda (a, {\cal J}(b))$ is the trace form of
${\cal A}_{\cal J}$.

\vspace{6pt}

{\bf Theorem 7.}(Helwig [337b]) {\it All Grassmann manifolds, as well as all compact
symmetric spaces of rank one, are contained in the form
${\rm Inv}_{\circ}({\cal A}_{\cal J}, {\cal J})$, where ${\cal A}$ is a
simple formally real Jordan algebra and ${\cal J}$ is a Peirce reflection of
${\cal A}$. The noncompact spaces associated with the above mentioned spaces have
the form} ${\rm Inv}_{\circ}({\cal A, J})$.

\vspace{6pt}

The description of important symmetric spaces due to Helwig [337b] embraces many other earlier descriptions of symmetric spaces using Jordan algebras and triple systems as those of Braun \& Koecher [131], Hirzebruch [346], Koecher [408h]. For a review of them see Iord\u anescu [364w, Ch. 3].

\vskip6pt

Let us recall here the simple but very ingenious description
given by Koecher [408h]: Let ${\bf A}$ be a formally real Jordan algebra, and suppose
its trace form $\lambda$ nondegenerate. Then the (not neccesarily positive
definite) line element ${\rm d}s^{2} : = \lambda (\dot x, P(x^{-1})\dot x) {\rm d}t^{2}$,
where $x = x(t)$ is a curve in ${\it {\rm Inv}}({\bf A})$, is invariant under the maps
$x \rightarrow Wx$,
$W \in \Gamma({\bf A})$, and $x \rightarrow x^{-1}$.
In order to discuss the induced (pseudo-) Riemannian structure, it suffices
to consider ${\rm Inv}_{\circ} ({\bf A})$ (indeed, if $C$ is a connected component of
${\rm Inv}({\bf A})$, then there exists an $f \in C$ such that $f^{2} = e$ and
$C = {\rm Inv}_{\circ} ({\bf A}_{f}))$. Then ${\rm Inv}_{\circ} ({\bf A})$ is a
symmetric Riemannian space and:

a) at the point $e$ the geodesic symmetry is the inversion
$x \rightarrow x^{-1}$, and the exponential map is ${\rm exp}_e
(x) = \exp  x$ and

b) the coefficients of the affine connection coincide with the structure
constants of ${\bf A}$.

With the ``metric" ${\rm d}s^{2}$, ${\rm Inv}_{\circ}({\cal A, J})$ is a regular analytic
submanifold of ${\rm Inv} ({\cal A,J})$ and pseudo-Riemannian and symmetric.

\vspace{6pt}

{\bf Comments.} Since 1966, independently of Koecher's research,
Iord\u a\-nescu, Popovici, and Turtoi, following a suggestion of their professor Gheor\-ghe Vr\v{a}nceanu, studied the spaces (Riemannian or pseudo-Riemannian) associated with various kinds of real Jordan algebras
(see [364a,~b,~c,~d], [546a],~[547],~[548], [676a]). The above recalled construction of Koecher [408h] from
1970 gave a subsequent proof of {\it the existence} of the spaces studied by Iord\u anescu-Popovici-Turtoi.
For details and also for open problems, see Iord\u anescu [364i,~t].

In 1999, Shima (see SHIMA, H.,
{\it Homogeneous spaces with invariant projectively flat affine
connections}, Trans. Amer. Math. Soc. {\bf 351} (1999), {\it 12},
4713--4726) showed that semi-simple symmetric spaces with invariant
projectively flat affine connections correspond to central-simple
Jordan algebras, and are described as centro-affine hypersurfaces in
the algebras. He also proved that Riemannian semisimple symmetric
spaces with invariant projectively flat affine connections correspond
to simple formally real Jordan algebras (see also~[490]).

\vspace{6pt}

{\bf Comments.}\ It would be interesting to connect the above
recalled\break Shima's results with Iord\u anescu-Popovici-Turtoi results
(mentioned in the previous {\bf Comments}).

\vspace{6pt}

Making use of Koecher's conference [408h], let us summarize
Helwig's construction [337b]:

\vspace{6pt}

{\bf Theorem 8.} ${\rm Inv}_{\circ} ({\cal A, J})$
{\it is a totally geodesic submanifold
of ${\rm Inv}_{\circ} ({\cal A})$. In case the pair $({\cal A, J})$ is simple
$($i.e., ${\cal A}$ contains no proper ${\cal J}$-invariant ideals$)$, 
${\rm Inv}_{\circ} ({\cal A, J})$ is an Einstein space if and only if
${\cal A}$ is central simple}.

\vskip6pt

Forty years ago (during the universitary year 1971-1972), following the suggestion of the late Prof.Dr. Enzo MARTINELLI (University "La Sapienza" of Rome, Italy), I have studied quaternionic Grassmann structures\footnote{IORD\u ANESCU, R., {\it On Grassmann quaternionic structures} (in italian), Boll. Un. Mat. It. (4) {\bf 10}(1974), 406-411.} by using Hangan's previous results (see HANGAN, Th., {\it Tensor product tangent bundles}, Arch. Math. (Basel), {\bf 19} (1968), {\it 4}, 436-440), and Pontryagin's local coordinates (see PONTRYAGIN, L.S., {\it Characteristic cycles on differentiable manifolds}, Mat. Sb. {\bf 21} (63) (1947), 233-284; Amer. Math. Soc. Transl. {\bf 32} (1950)). Making use of my results, Marchiafava extended his previous results from 1970 (see Marchiafava [469a]) - see Marchiafava [469b].  Let us point out -- from Alekseevsky \& Marchiafava [10b, p. 15] -- that `` ... the quaternionic Grassmannian carries a geometrical structure that is a natural generalization of the quaternionic structure (...); it would be interesting to examine more deeply such structures, which are integrable in case of Grassmannians, and find other examples''. (For a survey on this topic see Alekseevsky \& Marchiafava [10a,~b,~c,~d,~e], Marchiafava [469c], and the references therein.). In order to be able to remember here two interesting open problems (pointed out by me - in January 1994 - in my lectures given at Katholieke Universiteit Leuven, Belgium), let us briefly recall some facts from quaternionic geometry.

\vskip6pt

Let $L_n^{\bf H}$ (resp. $U_n^{\bf H}$) be the {\it homogeneous linear} (resp. {\it unitary}) quaternionic group acting on the left in the right quaternionic vector space ${\bf H}^n$ of dimension $n$.

Denote by $\Theta := (\theta^{\alpha})$ the matrix of a vector $\theta$ of ${\bf H}^n$ and by $A := (a_{\beta}^{\alpha})$ an $(n \times n)$-matrix over ${\bf H}$. Then  a transformation $T$ of $L_n^{\bf H}$ (resp. $U_n^{\bf H}$) is
$$
T : \Theta \to A\Theta,
$$
where $A$ is an invertible (resp. unitary) matrix.

\vskip6pt

{\bf Definition 7.} A real differentiable manifold $V_{4n}$ is called endowed with (a right) {\it almost} (resp. {\it almost Hermitian}) {\it quaternionic} structure if its structure group is $L_n^{\bf H}$ (resp. $U_n^{\bf H}$).

\vskip6pt

Recall that a $G$-{\it structure} on a real differentiable $m$-dimensional manifold $V_m$ is defined by a subbundle with the structure group $G$ of the tangent bundle $T(V_m, {\bf R}^m, L_m)$, where $G$ is a certain subgroup of $L_m$ (one assumes - of course - ${\bf R}^{4n} \equiv {\bf H}^{n}$, canonically).

\vskip6pt

{\bf Definition 8.} (see Martinelli [471a]) A real differentiable manifold $V_{4n}$ is called endowed with a {\it generalized} (right) almost (resp. almost Hermitian) quaternionic structure if its structure group is $\tilde{L}_n^{\bf H}$ (resp. $\tilde{U}_n^{\bf H}$), which consists of all tramsformations $T : \Theta \to A\Theta b$, where $A$ is an invertible matrix and $b \in {\bf H}-\{0\}$ (resp. $A$ is unitary matrix and $b\overline{b} = 1$).

\vskip6pt

{\bf Example.} The quaternionic projective plane (which suggested to Martinelli the above definition).

\vskip6pt

{\bf Question:} Do the quaternionic Grassmann manifolds have the same stucture? {\bf Answer:} No.

\vskip6pt

I have proved that their structure group $\mathcal{G}$ consists of all transformations
$$
T : \Theta \to (A\otimes Id_q)\Theta(B\otimes Id_p),
$$
where $A$ and $B$ are invertible quaternionic matrices of order $p$, resp. $q$. So, the notion of {\it locally Grassmann quaternionic manifold} arised and Marchiafava studied them later.

\vskip6pt

{\bf Geometrical open problem.} Define differentiable manifolds endowed with structure whose structure groups be similar to $\mathcal{G}$, but for which the factors $(A\otimes Id_q), (B\otimes Id_p)$ be replaced by unitary or invertible matrices.

{\bf Algebraic open problem.} What kind of algebras could describe the new differentiable manifolds from the above geometrical open problem.

\vskip6pt

It is worth to be mentioned here that recently Dubois-Violette - a specialist in noncommutative geometry - discovered the importance of (real) Jordan algebras for his field of research\footnote{Private communication on March 2010.}. For instance, in his lecture on noncommutative differential geometry [223] he has indicated that... "instead of taking Hermitian elements of $\ast$-algebras as the analogues real functions it would be more general (and radical) to take elements of real Jordan algebras".

\vskip6pt

There are two main topics of differential geometry related to Jordan triple
systems: symmetric $R$-spaces and hypersurfaces in spheres. Let us mention in
this respect that some of Cartan's results are now-being successfully
re-examined against the background of the theory of Jordan triple systems.

$R$-spaces constitue an important class of homogeneous submanifolds in
the Euclidean spheres. This class includes many examples appearing in
differential geometry of submanifolds. For example, all homogeneous
hypersurfaces and all parallel submanifolds in spheres are realized as
$R$-spaces.

Ferus [267] has characterized the $R$-spaces as compact symmetric
submanifolds of Euclidean spaces.

The connection between Jordan algebras and symmetric $R$-spaces was first
illuminated by Kantor, Sirota and Solodovnikov [389a], Koecher [408e,~II],
and Loos [448d]. Then, significant results were
obtained by Makarevich
(see [462a,~b,~c]) and Rivilis (see [558]).

In 1996, starting from the paper [462a], Bertram [98c] determined the causal
transformations of a class of {\it causal symmetric spaces} (see [98c, Th. 2.4.1]). As a basic tool he used {\it
causal imbeddings} of these spaces as open orbits in the conformal
compactification of formally real Jordan algebras. Firstly, he gave
elementary constructions of such imbeddings for the classical
matrix-algebras. Then he generalized these constructions for arbitrary
semisimple Jordan algebras: he introduced {\it Makarevich spaces} (which
are open symmetric orbits in the conformal compactification of a
semisimple Jordan algebra) and described examples and some general
properties of them which are the starting point of an
algebraic and geometric theory developed by Bertram in [98c].

In fact, in the above mentioned paper,
Bertram generalizes
features of bounded symmetric domains to a bigger class of symmetric
spaces (i.e., the above mentioned Makarevich spaces): he associates a
generalized {\it Bergman operator} to such a space and describes the
invariant pseudo-metric and the invariant measure on the space by means of
this family of operators. The space itself can be characterized
essentially as the domain where the generalized Bergman operator is
nondegenerate. These results are then applied to the theory of {\it
compact causal symmetric spaces}.

Generalizing Hermitian and pseudo-Hermitian spaces, Bertram defined in [98d, I]
{\it twisted complex symmetric spaces}, and showed that they
correspond to an algebraic object called {\it Hermitian Jordan triple
products}. He investigated the class of {\it real forms} of twisted
complex symmetric spaces called the category of {\it symmetric
spaces with twist}. Then he showed that this category is
equivalent to the category of
all real Jordan triple systems, and, using Makarevich [462a],
classified the
irreducible spaces. The classification shows that most irreducible
symmetric spaces have exactly one twisted complexification. This leads to
open problems concerning the relation of Jordan and Lie triple systems.
Using a geometric approach, Bertram [98d, II] defined and investigated the
{\it conformal group} of a symmetric space with twist. In the
non-degenerate
case he characterized this group by a theorem generalizing
the Fundamental Theorem of Projective Geometry.

\vskip6pt

{\bf Definition 9.} (see Takeuchi [655]). A {\it symmetric $R$-space}
is a compact symmetric space on which there exists a group of
transformations containing the group of motions as a proper
subgroup (see also Nagano [504]).

\vspace{6pt}

There exists a one-to-one correspondence between compact Jordan triple systems
and symmetric $R$-spaces, as was established by Loos [448d]; see Theorem 9 below.

The noncompact dual of a symmetric $R$-space can be realized as a bounded
domain $D$ in a real vector space. Loos [448d] proved that there is a one-to-one correspondence between boundary component of $D$ and idempotents of
the corresponding Jordan triple system.

A compact Jordan triple system $T$ becomes a Euclidean vector space with the
scalar product $(x, y) : = {\rm Tr}\, L(x, y)$,
where $L(x, y) z : = P(x, z) y : =
P(x + z)y - P(x) y - P(z) y$ and, by Theorem 12 from \S 1, is equal to
$\{x y z\}$. By the second equality in the definition of a Jordan triple system,
the vector space ${\cal H}$ spanned by $\{ L (x, y) \mid x, y \in T\}$ is a Lie
algebra of linear transformations of $T$, which is closed under taking
transposes with respect to $(\cdot\,,\cdot)$. The contragradient ${\cal H}$-module
$T'$ of $T$ can thus be identified with $T$ as a vector space, and
$$X\cdot v' = -\, {}^{t}X(v')\quad  \mbox{ for } X \in {\cal H}
\mbox{ and } v' \in T'.$$
The map $\tau : X \rightarrow {}^{t}X$, $X \in {\cal H}$, $v \rightarrow v'$, is a Cartan involution of the Lie algebra
${\cal L} : = T \oplus {\cal H} \oplus T'$ and $\sigma|_{\cal H} =
\pm 1$, $\sigma|_{T \oplus T'} = -1$ defines an involutive automorphism $\sigma$
of ${\cal L}$ commuting with $\tau$.

Recall that by a result of Koecher, ${\cal L} = T \oplus  {\cal H}
\oplus T'$ becomes a semisimple Lie algebra, adopting
$$[X,Y] : = XY - YX, \quad [X, v] = - [v, X] : = X \cdot v$$
for $X,Y \in {\cal H}$ and $v \in T \cup T'$,
$$[T, T] : = [T', T'] = 0, \quad [u, v'] : = - 2 L (u, v)$$
for $u \in T$ and $v' \in T'$. Also, Koecher's result states that
$Z = - {\rm Id}_{T}$ is an element
of ${\cal H}$, $({\rm ad}\, Z)^{3} = {\rm ad}\,Z$ and the $-1$-, $0$-,
$+1$-eigenspaces of ${\rm ad}\,Z$ are\break $T$, ${\cal H}$, $T'$.

Let $L$ be the centre free connected Lie group with Lie algebra ${\cal L}$,
let $H$ be the centralizer of $Z$ in $L$, let $U$ be a maximal compact
subgroup of $L$ determined by $\tau$, and let $K: = U \cap H$. Then $K$ lies
between the full set of fixed points of $\sigma$ in $U$ and the identity
component of $U$. If we denote by $P$ the normalizer of $T$ in $L$,
then $P$ is parabolic and $U/K \cong L/P$. It follows that
$M : = U/K$ is a symmetric $R$-space.

\vspace{6pt}

{\bf Theorem 9.} {\it The map $T \rightarrow M$ establishes a one-to-one
corespondence between isomorphism classes of compact Jordan triple systems
and symmetric $R$-spaces.}

\vspace{6pt}

In his paper [387b], Kantor generalized the notions of Jordan triple
system and symmetric space to cover more general cases of Riemannian
manifolds. To this end he introduced (generalized) Jordan triple systems of
second order and constructed an associated graded Lie algebra analogous to
the Kantor-Koecher-Tits construction. The corresponding Lie triple system
gives rise to ``bisymmetric'' spaces: Riemannian homogeneous fibre spaces
with symmetric base and (locally) symmetric fibre. One also has duality
theory for such spaces and can generalize the embedding of non-compact type
into the compact dual. (For bisymmetric Riemannian spaces see Kantor,
Sirota, Solodovnikov~[389b].)

\vskip6pt

Several papers of Dorfmeister and Neher [217a,~b,~c] deal principally with
isoparametric hypersurfaces in spheres and show that homogeneous examples with
four distinct principal curvatures are closely related to certain Jordan triple
systems. Then isoparametric triple systems of certain type are studied. For a detailed presentation of this topic see Iord\u anescu [364w, pp. 70-73].

\vskip6pt

Hulett and Sanchez [355] studied an algebraic structure, called
``Euclidean double-triple systems'' (because of their
analogies with Euclidean Jordan triple systems), associated
with a standard imbedding of an $R$-space. This structure
determines completely the geometry of an $R$-space and reduces
to a Jordan triple system if the $R$-space is symmetric.

Grassmann and flag manifolds associated with a Hermitian Jordan triple system
were defined by Arazy \& Upmeier in [40], and these differential-geome\-tric
objects were used to give a new proof for the intertwining formula
genera\-lizing ``Bol's Lemma" for symmetric domains which are not of tube type.

Tripotents are natural generalizations of partial isometries in $C^*$-algebras to the context
of $JB^*$-triples that is complex Banach spaces with symmetric unit ball.\
A survey on the main results contained in Chu \& Isidro [168], and Isidro \& Stach\'o
[369a,\,b,\,c] concerning the structure of the tripotents as a direct real-analytic
submanifold in a $JB^*$-triple, as well as some recent achievements are
presented in Stach\'o [635a] (see also Stach\'o [635b]).

Recently, Di Scala and Loi [210] have studied the symplectic geometry of Hermitian symmetric spaces
of noncompact type and their compact dual. Using the theory of Jordan triple systems,
they constructed an explicit {\it symplectic duality} (see Definition 10 below).

\vskip6pt

{\bf Definition 10.} Let $M \subset {\bf C}^n$ be a complex $n$-dimensional
Hermitian symmetric space endowed with the hyperbolic form $\omega_{hyp}$. Denote by
$(M^*,\omega_{FS})$ the compact dual of $(M,\omega_{hyp})$, where
$\omega_{FS}$ is the Fubini-Study form on $M^*$.
A {\it symplectic duality} is a diffeomorphism $\Psi_M : M \to {\bf R}^{2n} = {\bf C}^n \subset M^*$ satisfying $\Psi_M^* \omega_0 = \omega_{hyp}$ and
$\Psi_M^* \omega_{FS} = \omega_0$ for the pull-back of $\Psi_M$, where $\omega_0$ is the restriction
to $M$ of the flat
K\"ahler form of the Hermitian positive Jordan triple system associated to $M$.

\vskip6pt

Di Scala and Loi proved that the map $\Psi_M$ takes (complete) complex and totally
geodesic submanifolds of $M$ through the origin to complex linear subspaces
of $\mathbb{C}^n$. They also get an interesting characterization of the Bergman form
of a Hermitian symmetric space in terms of its restriction to classical
complex and totally geodesic submanifolds passing through the origin.

More recently, Di Scala, Loi, and Roos [211] have determined the group of diffeomorphisms
of a bounded symmetric domain, which preserve simultaneously the hyperbolic and the flat symplectic form.

In a very recent paper, Di Scala, Loi, and Zuddas [212] after extending the definition of
symplectic duality (given by Di Scala and Loi in [210] for bounded symmetric domains) to arbitrary
complex domains of $\mathbb{C}^n$ centered at the origin, generalize some of the results proved
in [210] and [163] to those domains.

\vskip6pt

{\bf Comments.} It would be possible that the symplectic duality be useful in the 
next future to physics. At least, Antonio Di Scala [209] expects to be proved this fact.

\vskip6pt

I like to point out here recent papers of Roos and his collaborators (see Yin, Lu, Roos [719]; Roos [565c], Wang, Yin, Zhang, Roos [700]), as well as a more recent paper by Roos [565d] on bounded symmetric domains.

In a recent (accepted 4 August 2008) and very extended (29 pages) paper, C.-H.Chu [166h]
introduced a class of real Jordan triple systems, called {\it $JH$-triples},
and showed (via the Tits-Kantor-Koecher construction of Lie algebras) that they correspond
to a class of Riemannian symmetric spaces including the Hermitian symmetric spaces and the 
symmetric $R$-spaces.

Let us recall here that there are two essentially equivalent ways to study Hermitian
symmetric spaces (via Tits-Kantor-Koecher construction, see Satake [592]), namely:

1) by using the semisimple Lie groups (it was in this way that the basic facts of the
theory were established by \'Elie Cartan in the 1930's and by Harish-Chandra in the 1950's), and

2) by using Jordan structures (algebras and triple systems), this way being essentially due to
Max Koecher and his school in the 1980's.

In his lecture notes from 2005, Koufany [416a] used the Jordan theory to point out some 
development in the geometry and analysis on Hermitian symmetric spaces. Let us point out here 
that the first part of these lecture notes in a nice survey about the geometry and the topology
of some homogeneous spaces associated with formally real Jordan algebras: Hermitian symmetric spaces
of tube type, their Shilov boundaries, and causal symmetric spaces of Cayley type. In particular,
Koufany reviews recent results by Clerc, \O rsted and himself about Maslov, Souriau, and
Arnold-Leray indices.

From recent applications of Jordan triple systems to differential geometry I must mention also the contributions of Kaneyuki [386] and Naitoh [506c].

\vskip6pt

In \S 4 of Chapter 3 from the book Iord\u anescu [364w], a discussion of Jordan structures that occur in some infinite-dimensional manifolds, for example, in symmetric Banach manifolds, is given. These are infinite-dimensional generalization of the
Hermitian symmetric spaces classified by Cartan [154a] in the
1930's using Lie theory. A connection of this classification to
Jordan algebras was later pioneered by Koecher [408h]. Jordan
triple systems first entered the scene when Loos [448h] showed
the correspondence, in finite dimensions, between bounded circular
domains and Jordan triples. The full generalization of this
correspondence to that between infinite-dimensional symmetric
Banach manifolds and Jordan triples was developed by Kaup in his outstanding papers
[392b,~d,~e].

\vskip6pt

The classical Grassmann manifolds can be regarded as manifolds of
projections in spaces of matrices. The infinte-dimensional
analogues are the manifolds of projections in $C^*$-algebras.
Since projections are tripotents, the concept of Grassmann
manifolds can be extended to include the manifolds of tripotents
in $JB^{\ast}$-triples (cf. [168], [347], [392a], [518a,~b], [587]).

\vskip6pt

Let us remember the Bertram's work up to 1999 (see Bertram [98h]): the
framework is the one of classical, finite-dimensional real differential
geometry, especially the theory of symmetric spaces.
The aim was to define and study geometric objects (manifolds with additional structure)
corresponding to real, finite-dimensional Jordan structures (pairs, triple systems, algebras).

As it is well known, in modern mathematics it is of fundamental importance the
correspondence between Lie algebras and Lie groups: there is a functor which
Bertram called the {\bf Lie functor for Lie groups} assigning to a Lie group its
Lie algebra, and every (real, finite-dimensional) Lie algebra belongs to some
Lie group. One may ask whether there is also a {\bf Jordan functor}:
can we find
a ``global'' or ``geometric'' object to which a given Jordan algebra is
associated in a similar way as a Lie algebra to a Lie group?

As in the case of Lie structures one may ask whether there is a
{\bf Jordan functor for Jordan
triple systems}: is a Jordan triple system associated with a geometric object, just as a
Lie triple system is associated with a symmetric space?

The problem of defining a {\bf Jordan functor} is of considerable interest because
it is related to many topics in geometry and harmonic analysis on symmetric
spaces.

Roughly speaking, a {\it twist} is an additional structure on a symmetric space.
This additional structure can be described in several ways. Geometrically, the
additional structure given by a twist on a symmetric space $M$ can be interpreted
as a {\it twisted para-complexification} or a {\it twisted para-complexification
of} $M$. In order to illustrate what this
 means, consider the example of the real projective space $M=P_n{\bf R}$ which is
 a symmetric space
 $$M=O(n+1)/(O(n)\times O(1)).$$
 It has a natural complexification given by the complex projective space
 $$P_n{\bf C}=U(n+1)/(U(n)\times U(1)).$$
This complexification is called by Bertram [98h] ``twisted" in order to distinguish
 it from the ``straight" complexification
 $$M_{{\bf C}} =O(n+1, {\bf C})/(O(n,{\bf C})\times O(1, {\bf C})).$$
 Every symmetric space $M=G/H$ has (locally) a unique straight complexification
 $M_{{\bf C}}=G_{{\bf C}} /H_{{\bf C}}$.

 In contrast, twisted complexifications are an additional structure of a
 symmetric space which in general need neither exist not be unique.

 Bertram [98h] proved that {\it all} (real, finite-dimensional) 
 Jordan structures correspond to certain
 geometric objects, namely:

-- Jordan algebras correspond to {\it quadratic prehomogeneous spaces};

-- Jordan pairs correspond to {\it twisted polarized symmetric spaces};

-- Jordan triple systems correspond to {\it symmetric spaces together with a
 twisted complexification}.

There is a forgetful functor from all of these objects into the category of symmetric
spaces which corresponds to a natural {\bf Jordan-Lie functor}.
It is a surprising
fact, obtained by classification of simple objects (due to E. Neher and M. Berger),
that this functor is quite close to being bijective.

\vskip6pt

{\bf Comment.} For a detailed presentation of Bertram's contributions see Iord\u anescu [364w, \S 5 from Ch.3].

\vskip6pt

Concerning the Bertram's work since 2000, let us point out that his aim was to
generalize correspondence between Jordan structures and geometric objects to the case of 
arbitrary dimension and of general base-fields and -rings (cover, in particular, the ``Jordan algebra
of quantum mechanics'' Herm$(\mathcal{H})$ which is infinite-dimensional, and more general
$C^*$-algebras).

In 2002, Bertram [98i] introduced {\it generalized projective geometries} which are a natural generalization
of projective geometries over a field or ring ${\bf K}$ but also of other important geometries
such as Grassmannian, Lagrangian or conformal geometry.
He also introduced the corresponding {\it generalized polar geometries} and associated to such a geometry
a {\it symmetric space over ${\bf K}$}.

\vspace{6pt}

{\bf Remark.} In the finite-dimensional case over ${\bf K}= {\bf R}$, 
all classical and many exceptional symmetric spaces are obtained in this way.

\vspace{6pt}

In the same paper [98i], Bertram proved that generalized projective and polar geometries
are essentially equivalent to {\it Jordan algebraic structures}, na\-mely to {\it Jordan
pairs}, respectively to {\it Jordan triple systems} over ${\bf K}$ which
are obtained as a linearized tangent version of the geometries in a similar way as a Lie group is
linearized by its Lie algebra. In contrast to the case of Lie theory, the construction of 
the {\bf Jordan functor} works equally well over general base rings and in arbitrary dimension.

One year later, in 2003, Bertram completed his previous results (see Bertram [98j]): he showed
that the correct generalization of the projective line in the category of generalized projective
geometry is given by spaces corresponding to unital Jordan algebras.

\vspace{6pt}

{\bf Remark.} The case of characteristic 2 is still an open problem.

\vspace{6pt}

{\bf Note.} An overview of the just above mentioned results can be found in Bertram~[98l].

\vspace{6pt}

{\bf Question.} What about the infinite-dimensional structures on the geometric objects introduced by
Bertram and briefly recalled just above?

\vspace{6pt}

Recently, Bertram and Neeb (see [105a,~b]), based on a previous their joint paper with Gl\"ockner
(see [101]), gave a new approach to differential calculus which works naturally in the framework of
very general base fields or even rings and in arbitrary dimension and characteristic, and on
consequences for differential geometry, Lie group theory and symmetric space theory (for full
details see Bertram [98p], where a general differential geometrical framework is developed).

\vspace{6pt}

{\bf Note.} An overview of these recent results can be found in Bertram~[98n].

\vskip6pt

{\bf Comment.} Related to this question there are two {\it very recent} paper by Bertram \&
L\"owe [104] and Bertram [98o].

\vspace{6pt}

In the paper [104], Bertram and L\"owe introduced the notion of {\it intrinsic subspaces} of
{\it linear} and {\it affine pair geometries}, which generalizes the one of projective
subspaces of projective spaces. They proved that, when the affine pair geometry is the
{\it projective geometry of a Lie algebra} introduced in [105a], such intrinsic subspaces correspond
to {\it inner ideals} in the associated Jordan pair, and they investigated the case of intrinsic
subspaces defined by the Peirce decomposition which is related to 5-gradings of the
projective Lie algebra. These examples, as well as the examples of general and Lagrangian
flag geometries, lead to the conjecture that geometries of intrinsic subspaces tend to
be themselves linear pair geometries.

It is known that the homotopy is an important feature of associative and Jordan
algebraic structures: such structures always come in families whose members need not be 
isomorphic among each other, but still share many important properties. One may
regard homotopy as a special kind of deformation of a given algebraic
structure. In the paper [98o], Bertram investigates the geometric counterpart
of this phenomenon on the level of the associated symmetric spaces. 
On this level, homotopy gives rise to conformal deformations
of symmetric spaces. These results are valid in arbitrary dimension and over
general base fields and rings.

\vskip6pt

Very recently, Bertram \& Bieliavsky [99a,b] investigated a special kind of construction of symmetric spaces, called {\it homotopy}, and they gave the classification of homotopes.

\vskip6pt

The winter 1980--1981 witnessed the appearance of the famous paper [490a]
by Sato, where he proved that the totality of solutions of the
{\it Kadomtsev-Petviashvili
$(KP)$ equation},
$$3u_{yy} + (-4u_t + u_{xxx} + 12 uu_{x})_x =0,$$
forms an infinite-dimensional Grassmann manifold.

\vspace{6pt}

{\bf Note.} Let us recall that the {\rm KP} equation was discovered in 1970 in an
effort to understand the propagation of long, shallow waves in plasma (see B.
KADOMTSEV and V. PETVIASHVILI, Dokl. Akad. Nauk SSSR {\bf 192} (1970), {\it 4}, 753).

\vspace{6pt}

The evolution of $u$ in the variables $x,y,t$ is interpreted as a dynamical
motion of a point on the Sato's Grassmann manifold by the action of a three-(or
more) parameter subgroup of the group of its automorphisms. Generic points of
Sato's Grassmann manifold give generic solutions to the {\rm KP} equation, whereas
points on particular its submanifolds give solutions of particular type. For instance,
{\it rational} solutions correspond to points on {\it finite-dimensional}
Grassmann submanifolds. Also, different kinds of submanifolds give rise to generic
solutions of other soliton equations, such as the
Korteweg-de Vries $(KdV)$
equation, the modified KdV equation, the Boussinesq equation, the
Sawada-Kotera equation, the non-linear Schr\"{o}dinger equation, the Toda
lattice, the equation of self-induces transparency, the Benjamin-Ono
equation, as well as to solutions of particular type of these soliton equations.

Moreover, a multicomponent generalization of the theory shows that
solutions of
other soliton equations (such as the equation for three-wave interaction, the
multi-component nonlinear Schr\"{o}dinger equation, the sine-Gordon
equation, the Lund-Regge equation, and the equation for intermediate long
wave) also constitute submanifolds of Sato's Grassmann manifold.

Sato [595a] conjectured that {\it any} soliton equation, or
completely integrable
system, is obtained in this way. It follows that the classification of soliton
equations would be reduced to the classification of submanifolds of Sato's
Grassmann manifold which are stable under the subgroup of its automorphisms
describing space-time evolution.

\vspace{6pt}

{\bf Comments.} As I have predicted in 1983 through my talks given at the
Universities of Timi\c{s}oara and Ia\c{s}i (Romania), Sato's theory is an outstanding
contribution with a deep impact on physics, but also on mathematics. In fact, on the occasion of Sato's
visit in Romania (in August 1983) I suggested him to use
Jordan structures as main tool for his theory. It is the great merit of Josef Dorfmeister and his
collaborators to give, beginning with 1989, an impressive and important work in this direction.

\vskip6pt

{\bf Note.} For a detailed presentation of the results obtained by Dorfmeister and his collaborators see Iord\u anescu [364w, \S 6 from Ch.3].

\vskip6pt

At the end of this paragraph I would like to mention a {\bf very recent} paper by Alekseevsky [9] on pseudo-K\" ahler symmetric spaces $\dots$ "which are very closely related with Jordan pairs" (as Alekseevsky himself informed me\footnote{Private communication in March 2010.}).

\vskip6pt

Concerning the applications of Jordan structures to differential geometry, I like to refer the reader to the {\bf very recent} book [166i] by Chu, not yet published, but available from January 2012.

\vspace{1cm}

\centerline{\bf \Large \S 5. JORDAN ALGEBRAS IN RING GEOMETRIES}

\vspace{48pt}

The first investigation of octonion planes dates from 1933 and is due to
Moufang [500a]. It consisted in the construction of a projective plane
coordinatized with an octonion division algebra. In this Moufang plane,
Desargues Theorem fails but the Harmonic Point Theorem is valid. In 1989,
Faulkner [248d] described a geometric construction of the Moufang
projective octonion plane from the projective quaternionic 3-space. This
geometric construction is a translation of an algebraic construction of the
27-dimensional degree 3 exceptional simple Jordan algebra as the trace 0
elements of a 28-dimensional degree 4 Jordan algebra (see Allison \& Faulkner
[18a]). The geometric construction is a particular case of a more general
construction described by Faulkner that also yields as a special case a
construction of the complex projective plane from the complex projective
line. In 1994, using the classification of the group $SL_2 (K)$ and its
natural module, Timmesfeld [663] gave a classification of Moufang
planes and the groups $E_6^K$ together with various other results on Moufang planes.

Another approach to octonion planes was given in 1945 by Jordan [379d] {\it via}
the Jordan algebra $H_3({\bf O})^{(+)}$. Recall the definition of the
exceptional Jordan algebra $H_3({\bf O})^{(+)}$: Let $H_3({\bf O})$ be the
set of all $(3 \times 3)$-matrices with entries in an octonion algebra
${\bf O}$ and which are symmetric with respect to the involution
$x \ri\ov{x}'$. The
characteristic of the underlying field is supposed to differ from two. On
$H_3({\bf O})$ we can define a Jordan algebra structure by means of the
product $xy:= \frac{1}{2}(x \cdot y+y \cdot x)$, where the dot
means the usual matrix product. The resulting Jordan algebra is denoted by
$H_3({\bf O})^{(+)}$. Jordan focussed on a real octonion division algebra
${\bf O}$ and used the primitive idempotents in $H_3({\bf O}^{(+)})$ to
represent the points and lines of a projective plane.

In 1951, Freudenthal [272a] obtained essentially the same construction\footnote{For a
  systematic treatise on applications of the real algebra of Cayley numbers,
  I refer the reader to Brada's Ph.D. Thesis [104].}. Atsuyama [52a]
  used the embedding defined by Yokota [720] to obtain new results in this
  direction. In 1997, Allcock [14a] provided the identification of the
  classical Jordan's model of octonion plane [379d] with the elegant model given by
Aslaksen [50].

The construction was extended in 1960 by Springer [632a], who considered
${\bf O}$ as an octonion division algebra over a field of characteristic
different from two or three. In this more general setting, elements of rank one
(which are either non-zero multiples of primitive idempotents or nilpotents of
index two) are used to represent the points and lines of a projective
plane. Springer proved the fundamental theorem relating collineations of the
plane and norm semisimilarities of the Jordan algebra. Jacobson [371a]
showed that the little projective group, i.e., the group generated by elations
(transvections) of these planes, is simple and isomorphic to the
norm-preserving group of the Jordan algebra {\it modulo} its centre. Suh [644]
showed that any isomorphism between the little projective groups of two
planes is induced by a correlation of the planes. Springer and Veldkamp
[633a] undertook a study of Hermitian polarities of a projective octonion
plane and the related hyperbolic and eliptic planes. The unitary group of
collineations commuting with a hyperbolic polarity was studied by Veldkamp
[689c].

Springer and Veldkamp [633b] considered planes associated with split
(i.e., not division) octonion algebras over a field of characteristic different
from two or three. These planes are not projective. (For the study of these planes, see Veldkamp
[689b,~d].)

In 1970, Faulkner [248a] extended the notion of octonion planes in another
direction by removing the restriction that the characteristic of the
underlying field to be different from two or three. After McCrimmon [480a]
had introduced the notion of quadratic Jordan algebra and verified that
$H_3^{(+)} ({\bf O})$ possesses such a structure for any characteristic,
Jacobson suggested to Faulkner (see [248a, p.~3]) that characteristic-two
octonion planes could be approached in this way. As it turned out, in the
setting of quadratic Jordan algebras, most of the results on octonion planes
can be derived in an uniform manner, without referring to the
characteristic or the type of an octonion algebra.

For collineation groups of projective planes over degenerate octaves or
antioctaves, see Persits [538a,~b].

Davies [195] studied bi-axial actions on projective planes (including
octonion planes), making use of their Jordan algebra description.

Bix [112a] has defined and studied octonion planes over local rings. He
generalized Faulkner's result on the simplicity of {\it PS} (see Iord\u anescu [364w, Ch. 4, p. 107]) of an octonion plane over a field to
octonion planes over local
rings. Those subgroups of the collineation group of an octonion plane over a
local ring which are normalized by the little projective group have been
classified. These parallel results of Klingenberg and Bass, who classified
those subgroups of the general linear group over a local ring which are normalized by the
special linear group.

\vskip6pt

{\bf Comment.} For a detailed presentation of the contributions of Faulkner [248a], and Bix [112a,b,c] see Iord\u anescu [364w, \S 1 of Ch. 4].

\vskip6pt

Using concepts from valuation theory, Carter and Vogt [155] gave a
 characterization of all collinearity-preserving functions from one affine or
 projective Desarguesian plane into another. {\it Lineations} (i.e., point
 functions $f$ from one plane into another with the property that, whenever
 $x,y$ and $z$ are collinear points, $f(x), f(y)$ and $f(z)$ are collinear
 points) whose ranges contain a quadrangle, called in [155] {\it full
 lineations}, have been algebraically characterized in various setting by
 Klingenberg ([402a], lineations from a Desarguesian plane onto another),
 Skornyakov ([625], linieations from an arbitrary plane onto another),
 Rad\'{o} ([554b], full lineations from a Desarguesian plane into another, the
 infinite-dimensional case being considered in [554a]), and Garner
 ([285], lineations from a Pappian coordinate plane into another taking the reference quadrangle of one plane to
 that of the other).

Carter's and Vogt's results allow one or both planes to be affine and include
 cases where the range contains a triangle but no quadrangle. A key theorem is
 that, with the exception of certain embeddings defined on planes of order 2
 and 3, every collinearity-preserving function from an affine Desarguesian
 plane into another can be extended to a collinearity-preserving function
 between the enveloping projective planes. Full lineations defined on
 finite-dimensional affine spaces can also be extended to the enveloping projective space (Brezuleanu
 and R\u{a}dulescu [139a,~b]).

Faulkner and Ferrar [249a] have showed that, up to conjugation by
 col\-line\-ations, there exists at most one surjective homomorphism from an
 octonion plane to a Moufang plane. They also established the existence of
 proper homomorphisms between octonion planes and of homomorphisms
 from octonion planes onto Desarguesian planes.

Ferrar and Veldkamp [264] studied neighbour-preserving homomorphisms
 between projective ring planes (i.e., mappings preserving incidence and the
 neighbour relation between points and lines). These are generalizations of the
 homomorphisms between ordinary Desarguesian projective planes which have
 first been studied by Klingenberg [402a]. On the other hand, in the
 context of projective planes over rings of stable rank 2 as studied by
 Veldkamp [689e], an obvious question to ask was what mappings between
 such planes are induced by homomorphisms between coordinatizing rings. If one
 requires that the ring homomorphisms carry 1 to 1, then they induce
 distant-preserving homomorphisms which are mappings incidence and the
 negation of the neighbour relation between a point and a line. Veldkamp
 [689f] proved that any distant-preserving homomorphism $\psi$ is induced
 by a ring homomorphism carrying 1 to 1, provided the two planes are
 coordinatized with respect to basic quadrangles which correspond under
$\psi$. Veldkamp [689f] also studied homorphisms between projective ring
 planes which only preserve incidence. They turn out to be
 products of a bijective neighbour-preserving homomorphism followed by an
 arbitrary distant-preserving homomorphisms.

In 1983, Faulkner and Ferrar [249c], utilizing methods similar to those of
Skornyakov [625], extended the results of Klingenberg [402a] from
Desarguesian planes to Moufang planes. Let us note in this respect that a
systematic study of places of octonion algebras over discrete valuation rings
has been carried out by Petersson [539a]. Bro\v{z}ikova [141], making use
of previous results of Havel [333] and Faulkner \& Ferrar [249c],
provided a Jordan-theoretic description of all homomorphisms between Moufang
planes having the property that the points identified {\it via} Springer's
isomorphism (see Springer [632a]) with $(0,0)$, $0$, $(\infty )$
and $(1,1)$ are mapped to their analogues in the image plane.

\vspace{6pt}

{\bf Open Problem.} It would be interesting to deal with topics similar to those mentioned above in the case of octonion planes.

\vskip6pt

A systematic study of projective planes over large classes of associative
 rings was initiated by Barbilian in his very general approach [68b] (see
 also Barbilian [68a]). He proved there that the rings which can be
 underlying rings for projective geometries are (with a few exceptions) rings
 with a unit element in which any one-sided inverse is a two-sided
 inverse. Barbilian [68b, I] called these rings ``$Z$-rings'' (from
 ``Zweiseitigsingul\"{a}re Ringe'') and gave a set of 11 axioms of projective
geometry over a certain type of $Z$-ring (see [68b, II]).

For more than thirty years no development of Barbilian's study
succeeded. Beginning with 1975, outstanding mathematicians like
W. Leissner in
Germany, F.D. Veldkamp in Holland, and J.R. Faulkner in USA, developed
Barbilian's ideas, and so, notions as ``Barbilian domains'',
``Barbilian planes'', ``Barbilian spaces'', or
``Barbilian geometry'' appeared.

\vskip6pt

Leissner [437a,~b,~c,~d,~e] developed a plane geometry over an arbitrary $Z$-ring $R$, in which a
 point is an element of $R \times R$ and a line is a set of the~form
$$\{(x+ra, y+rb) \mid (x,y) \in R\times R, \, r \in R, \, (a,b) \in B\},$$
where $B$ is a ``Barbilian domain'', i.e., a set of unimodular pairs from $R
\times R$ satisfying certain axioms.

\vspace{6pt}

{\bf Note.} Let us mention in this context that Lantz extended in [431]
 the results of Benz from [87a,~b] by showing that large classes of
 commutative rings admit only one Barbilian domain.

\vspace{6pt}

Rad\'{o} [554a,~b,~c] extended Leissner's results [437a] to affine Barbilian
 planes over
 an arbitrary ring with a unit element and investigated the corresponding
 affine Barbilian
structures and translation Barbilian planes. Corresponding to the algebraic
 representation of affine Barbilian spaces as affine geometries over unitary
 free modules, Leissner [437b]
characterized algebraic properties of the underlying ring $R$, respectively
 module $M_R$,
respectively Barbilian domain $B \subset M_R$ by geometric properties of the
 affine Barbilian space and viceversa.

Veldkamp [689c] gave an axiomatic description of plane geometries of the
 kind considered by Bingen [109]. A most satisfactory situation is reached
 by extending the class of ring used for coordinatization from semiprimary
 rings, which Bingen used to rings of stable rank two. These rings have played
 a role in algebraic $K$-theory, and seem to form a natural framework
 for many geometric problems.

\vspace{6pt}

{\bf Note.} For simplicity Veldkamp [689e] confines himself to the case of
 planes (called {\it projective Barbilian planes}), but a generalization to
 higher dimensions is straighforward (see below {\it projective Barbilian spaces}
defined also by Veldkamp).

\vskip6pt

Veldkamp has chosen an approach somewhat along the lines of Artin
[46] rather than to follow Barbilian.

The basic relations in the plane are the incidence and the neighbour
relation. The
axioms
consist of a number of axioms expressing elementary relations between points
and lines
such as, e.g., the existence of a unique line joining any two
non-neighbouring
points,
and a couple of axioms ensuring the existence of transvections and
dilatations.

In 1993, Veldkamp refined in [689j] the notion of Barbilian domain to
$n$-Barbilian domain in a free module of rank $n$. This leads to results
that
bear on $n$-dimensional affine ring geometry. The case of infinite rank is
also considered by Veldkamp in [689j].

\vspace{6pt}

{\bf Note.} For a good survey on the theory of projective planes over
rings of stable rank two, see Veldkamp [689h]. Such a plane is described as a
structure of points and lines together with an incidence relation and a
neighbour relation and which has to satisfy two groups of axioms. The axioms in
the first group express elementary relations between points and lines such as,
e.g., the existence of a unique line joining any two non-neighbouring points,
and define what is called a Barbilian plane. In the second group of axioms the
existence of sufficiently many transvections, dilatations, and generalizations
of the latter, the affine dilatations and their duals is required.

\vspace{6pt}

In 1987, Veldkamp [689i] extended all this above mentioned results to
arbitrary finite dimension. Basic objects in the axioms are points and
hyperplanes, by analogy with the selfdual set-up for classical projective
spaces over skew fields given by Esser [238]. As basic relations again
serve incidence and the neighbour relation. The self-dual approach is quite
natural since incidence and the neighbour relation between points and
hyperplanes have a simple algebraic description in coordinates. Homomorphisms
are more or less the same as in the plane case, things becoming a bit more
complicated because Veldkamp included homomorphism between spaces of unequal
dimension.

\vspace{6pt}

{\bf Note.} Veldkamp confine himself to full homomorphisms, which can only
increase the dimension or leave it the same. Thus he excluded
homomorphisms which lower the dimension, an example of which was given by
Fritsch and Prestel~[276].

\vskip6pt

{\bf Comment.} For a suvey of the results of Veldkamp [689e,i] concerning the fundamental properties of projective Barbilian spaces, see Iord\u anescu [364w, \S 3 of Ch. 4].

\vskip6pt

Faulkner and Ferrar [249d] surveyed the development which leads from
classical Desarguesian projective plane {\it via} Moufang planes to
Moufang-Veld\-kamp planes. They first sketched inhomogeneous and homogeneous
coordinates in the real and projective planes and in ring planes, the Jordan
algebra construction of Moufang planes, and the representation of all these
planes as homogeneous spaces for their groups of transvections. Then attention
is focussed on {\it Moufang-Veldkamp planes}, i.e., projective Barbilian planes
in which all possible transvections exist and which satisfy the little
quadrangle section condition for quadrangles in general position. As
coordinates for the affine plane one easily obtains an alternative ring of
stable rank two. Unfortunately, the Bruck and Kleinfeld theorem for alternative
division rings does not carry over to alternative rings in general, i.e., such
a ring need neither be associative nor be an octonion algebra. Therefore, to
coordinatize the whole projective plane one cannot rely on either homogeneous
coordinates (as in the associative case) or the Jordan algebra construction
(as in the octonion case). In this case, one has to follow a more complicated
way, namely: first to construct a certain Jordan pair from the given
alternative ring, then to define a group of transformations of this Jordan
pair and, finally, to represent the projective plane as a homogeneous space
for that group.

Faulkner [248c] proved that for a {\it connected Barbilian transvection
  plane}~$P$ (i.e., a plane with incidence and neighbouring generalizing
  Moufang projective planes) one can construct a connected Barbilian plane
  $T(P)$ called {\it tangent bundle plane}. This construction agrees with the
  usual tangent bundle when it exists. If $T(P)$ is also a transvection
  plane, then the set $S$ of sections of $T(P)$ is a Lie ring. The group $G$
  generated by all transvections of $P$ acts on $S$. Since $S$ is isomorphic
  to the Koecher-Tits Lie ring constructed from the Jordan pair $(M_{12}(R),
  M_{21}(R))$, where $R$ is the associated alternative ring,
one can determine
  $G$ and thereby $P$ from $R$.

In 1987, Spanicciati [631] introduced {\it near Barbilian planes} (NBP)
and {\it strong near Barbilian planes} (SNBP) as a variation of Barbilian
planes. In 1989, Hanssens and van Maldeghem [331] showed that a NBP is an
SNBP, and classified all NBP up to the classification of linear spaces (many
examples follows as a result of a universal construction). They also showed
that only NBPs that are also BPs are those mentioned in [631], namely the
projective~planes.

Allison and Faulkner [18a] gave in 1984 an algebraic construction of
degree three Jordan algebras (including the exceptional one) as trace zero
elements in a degree four Jordan algebra. Five years later, Faulkner [248d]
translated this algebraic construction to give a geometric construction of
Barbilian planes coordinatized by composition algebras (including the Moufang
plane) as skew polar line pairs and points on the quadratic surfaces
determined by a polarity of projective $3$-space over a smaller composition
algebra.

In 1989, Faulkner [248e] defined and studied the so-called $F$-{\it planes}
which generalize the projective planes. Planes considered by Barbilian in the
{\it Zusatz} to [68b] are connected $F$-planes in Faulkner's setting. Besides
extending the class of coordinate ring, Faulkner's work [248e] introduces
some new concepts, techniques, and connections with other areas. These include
a theory of covering planes and homotopy although there is no topology, a
theory of tangent bundle planes and their sections although there is no
differential or algebraic geometry, a purely geometric and coordinate-free
construction of the Lie ring of the group generated by transvections, and
connections to the $K$-theory of the coordinate ring.

Greferath and Schmidt [305b] introduced in 1992 the notion of {\it
  Barbilian space of a projective lattice geometry} in order to investigate
  the relationship between lattice-geometric properties and the properties of
  point-hyperplane structures associated with. They obtained a
  characterization of those projective lattice geometries, the Barbilian space
  of which is a Veldkamp space (see also Greferath and Schmidt [305a]).
  
\vskip6pt

In \S 5 of Ch. 4 from the book Iord\u anescu [364w] are presented coordinatization theorems similar to Moufang's [500a] for
several polygonal (projective, quadrilateral, and hexagonal, admitting all
elations) geometries. For the classification of polygonal geometries, see Tits
[664b,~d,~e]. These coordinatize a projective Moufang plane
by an alternative division algebra. The essential information needed is
contained in the group generated by all elations. These groups are groups with
Steinberg relations for which parametrization theorems were presented in the \S 4 of the same Ch.4 from the book Iord\u anescu [364w]. For these presentations the Faulkner's formulations [248b] were used.

\vskip6pt

{\bf Comments.} For a brief survey on the relations between various
exceptional notions in algebra and geometry (e.g., non-classical Lie algebras,
nonassociative alternative algebras, non-special Jordan algebras,
non-Desargue\-sian projective planes), the reader is referred to Faulkner and
Ferrar [249a], who proved that all these notions are related, one way or
another, to the octonions.

\vskip6pt

An interesting area of research is that of {\it chain geometries} and their
generalizations, {\it chain spaces}, and their relations with Jordan algebras
(see Blunck [115] and Herzer [339]).

\vspace{1cm}

\centerline{\bf \Large \S 6. JORDAN ALGEBRAS}

\centerline{\bf \Large IN MATHEMATICAL BIOLOGY}

\centerline{\bf \Large AND IN MATHEMATICAL STATISTICS\footnote{The topics presented in \S 6 were only mentioned in the book Iord\u anescu [364w].}}

\vspace{48pt}

{\bf Note.} For a comprehensive account on algebras in genetics up to 1980, the reader is reffered to W\" orz-Busekros monograph [711a].

\vskip6pt

Etherington [239a,b] showed how a nonassociative algebra can be made to correspond to a given genetic system. The fact that many of these algebras have comon properties has prompted their study from a purely abstract standpoint. Furthermore, these algebraic studies gave new ways of tackling problems in genetics.

In a study of nonassociative algebras arising in genetics, Schafer [598a] proved that the so-called gametic and zygotic algebras (see Etherington [239b]) for a single diploid locus are Jordan algebras.

Holgate [348a] proved Schafer's results by methods which do not make use of transformation algebras (employed by Schafer [598a]), which therefore accomodate the multi-allelic case more easily, and in which the main object is to maximise the interplay between the algebraic formalism and the genetic situation to which it corresponds.

The first part of the present paragraph deals with these results as treated by Holgate [348a], while the second deals with the results due to Piacentini Cattaneo [542], W\" orz-Busekros [711c,d] and Walcher [699c].

We consider the {\it gametic algebra} $G$ of a single locus with $n+1$ alleles, i.e., the algebra $G$ over ${\bf R}$ with basis $\{a_0, \ldots, a_n\}$, whose elements correspond to the actual allelic forms, its multiplication table being
$$
a_ia_j := \displaystyle\frac{1}{2}(a_i + a_j).
$$

For an element $x = \mathop{\sum}\limits_{i=0}^{n} x_ia_i$, the {\it weight} $w$ is defined by $w(x) := \mathop{\sum}\limits_{i=0}^{n} x_i$.

\vskip6pt

{\bf Remark.} It is easily seen that $x^2 = w(x)x$.

\vskip6pt

{\bf Proposition 1.} ({\it Algebraic}). {\it Every element of unit weight in $G$ is idempotent.}

{\bf Proposition 1.} ({\it Genetic}). {\it In the absence of selection, the gametic proportions remain constant from one generation to another.}

\vskip6pt

{\bf Remark 1.} The algebraic result is more comprehensive, since only those elements of unit weight for which all the $x_i$ are nonnegative correspond to populations.

\vskip6pt

{\bf Remark 2.} The nonassociativity of genetic algebras corresponds to the fact that if $P, Q$ and $R$ are populations and if $P$ and $Q$ mate and the offspring mates with $R$, the final result is, in general different from that arising from mating between $P$, and the offspring of mating between $Q$ and $R$. The two situations are shown in the diagram below:

\begin{center}
\includegraphics[width=80mm]{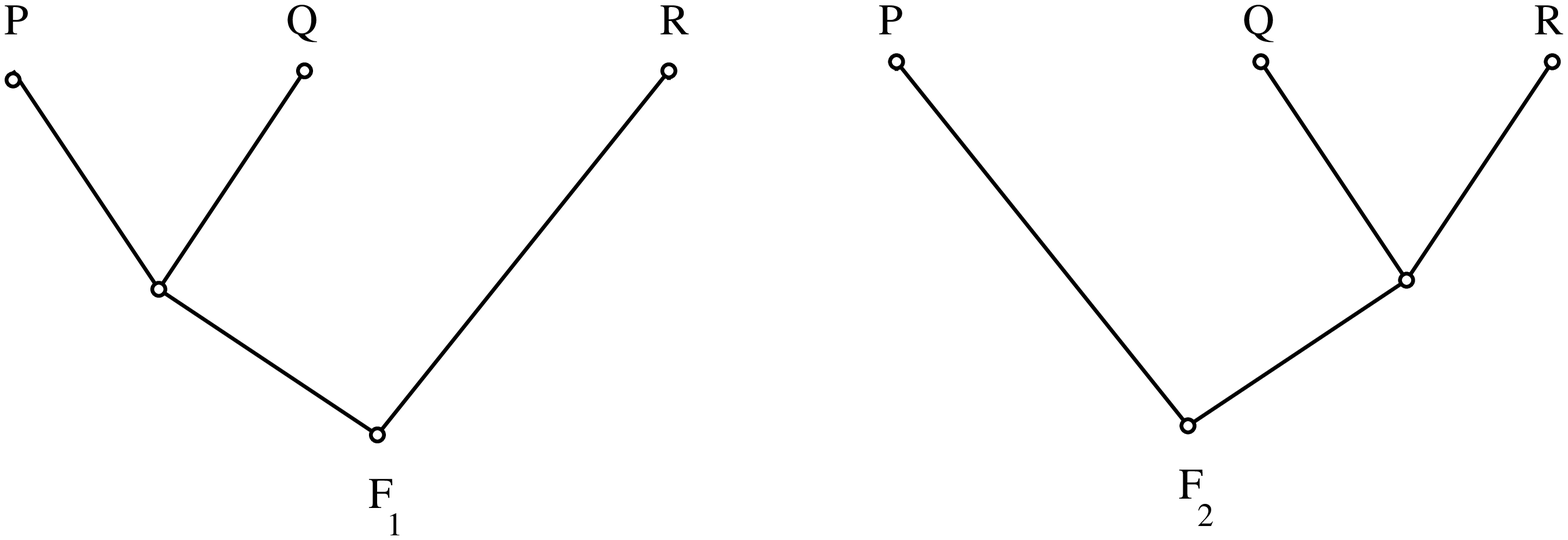}
\end{center}

\centerline{Figure 1}

\pagebreak

{\bf Proposition 2.} ({\it Algebraic}). {\it The algebra $G$ is a Jordan algebra.}

\vskip6pt

{\bf Proposition 2.} ({\it Genetic}). {\it In the mating schemes shown in Figure 1, the populations $F_1$ and $F_2$ have the same genetic proportions if $P$ is the offspring of mating of $R$ with itself.}

\vskip6pt

{\bf Notation.} The algebra $Z$, corresponding to proportions of zygotic types, is formed by duplicating $G$ (see Etherington [239a]) its basic elements are pairs $(x,y)$ of basis elements of $G$ with the multiplication rule $(x,y)(u,v) := (xy, uv)$. A canonical basis may be taken in $G$ by setting: $c_0 := a_0, c_i := a_0-a_i \; (i\neq 0)$, for which the multiplication table is
$$
c_0^2 = c_0, \; c_0c_i = \displaystyle\frac{1}{2}c_i, \; c_ic_j = 0 \quad (i,j \neq 0).
$$

Then on writing $d_{ij} := (c_i, c_j)$, the multiplication table for the duplicate $Z$ can be written as
$$
d_{00}^2 = d_{00}, \; d_{00}d_{0i} = \displaystyle\frac{1}{2}d_{0i}, \; d_{oi}d_{0j} = \displaystyle\frac{1}{4}d_{ij},
$$

other products being zero $(i,j \neq 0)$.

\vskip6pt

{\bf Remark.} The weight of an element $\mathop{\sum}\limits_{i,j=0}^{n} x_{ij}d_{ij}, \; x_{ij} = x_{ji}$ is $w(x) = d_{00}$.

\vskip6pt

{\bf Proposition 3.} ({\it Algebraic}). {\it Every element of the form $y := x^2 - w(x)x$ annihilates $Z$.}

\vskip6pt

{\bf Proposition 3.} ({\it Genetic}). {\it The extent to which the zygotic proportions in a population differ from the Hardy-Weinberg equilibrium state has no effect on the offspring distribution produced by mating between this population and any other.}

\vskip6pt

{\bf Proposition 4.} The algebra $Z$ is a Jordan algebra.

\vskip6pt

{\bf Remark.} Let $\mathcal{A}$ be the algebra over $\mathbb{\bf C}$ with basis $\{a_0, \ldots, a_n\}$ whose multiplication table is
$$
a_i . a_j := a_i. \eqno(\ast)
$$

Obviously, $\mathcal{A}$ is associative. Consider the special Jordan algebra $\mathcal{A}^{(+)}$ obtained from the vector space $\mathcal{A}$ by means of product $xy := \displaystyle\frac{1}{2}(x . y + y . x)$. It can easily be seen that $\mathcal{A}^{(+)}$ is isomorphic to $G$.

\vskip6pt

If it were possible to know in advance that the genes of only one of two given populations mating together are transmitted to the offspring, these could be written first in the product, and the system would correspond to the multiplication table $(\ast)$. The fact that $G$ is a special Jordan algebra appears as a consequence of inheritance being symmetric in the parents.

Recall from Gonshor [298b,I] the following

\vskip6pt

{\bf Definition.} A {\it special train algebra} is a commutative algebra over $\mathbb{\bf C}$ for which there exists a basis $\{a_0, \ldots, a_n\}$ with a multiplication table of the following kind: $a_ia_j := \sum x_{ijk}a_k$, where

(i)\quad\quad {$x_{000} = 1$},

(ii)\quad\quad for $k<j$, \quad $x_{0jk} = 0$,

(iii)\quad\quad for $i,j > 0$, \quad $k \leq \max(i,j)$, \quad $x_{ijk} = 0$,

\noindent and all powers of the ideal $(a_1, a_2, \ldots, a_n)$ are ideals. (The powers $I^r$ of an ideal $I$ are defined by $I^r := I^{r-1}I$)\footnote{From {\bf very recent} papers devoted to train algebras, I like to mention [71] by Bayara, Conseibo, Ouattara and Zitan.}.

\vskip6pt

{\bf Remark.} A commutative algebra over $\mathbb{\bf C}$ for which only conditions (i), (ii) and (iii) are required was called by Gonshor {\it genetic} algebra (see [298a]). Schafer's concept of genetic algebra coincides with that of Gonshor (see Gonshor [298a, Theorem 2.1]\footnote{Concerning genetic algebras, the fundamental idea has been to define a basis $\{G_1, \ldots, G_n\}$ with a one-to-one correspondence to the genotypes $g_1, \ldots, g_n$ considered, and then give a multiplication table so that the product $G_iG_j$ of two basis elements be equal to a linear combination $\sum p_{ijk}G_k$, where $p_{ijk}$ is the probability of getting genotype $g_k$ in a cross between $g_i$ and $g_j$ individuals.}). W\" orz-Busekros defined [711b] three kinds of {\it noncommutative} Gonshor genetic algebras and characterized them in terms of matrices.

\vskip6pt

{\bf Comments.} Let us mention in this respect that in the mathematical theory of algebras in genetics, whose origins are in several papers by Etherington, fundamental contributions have been made by Schafer, Gonshor, Holgate, Reiers\o l, Hech and Abraham (for a detailed account see [711a]).

\vskip6pt

{\bf Definition.} The $x_{0jj}$ are called the {\it train roots} of the algebra. (They are the characteristic roots of the operator which is multiplication by $a_0$).

\vskip6pt

{\bf Remark.} From Schafer [598a, Theorem 5], it follows that a special train algebra can only be a Jordan algebra if its train roots all have values among $1, \displaystyle\frac{1}{2}, 0$. This excludes the genetic algebras corresponding to polyploidy of several loci. Therefore, the appearance of Jordan algebra seems to be bound up with the property of attaining equilibrium after a single generation of mating

Piacentini Cattaneo considered [542] the gametic algebra $G$ (see the beginning of this paragraph) of simple Mendelian inheritance. Suppose that mutation occurs in the chromosomes, i.e., suppose that a rate of alleles $a_i$ mutate into the alleles $a_j, \; j\neq i$. If we denote this rate by $r_{ij}$ (setting $r_{kk} = 0$), we can construct a new algebra, denoted by $G_m$, called {\it a gametic algebra of mutation} (see [542, p. 180]). The new multiplication table then is
$$
a_j^2 = (1 - \mathop{\sum}\limits_{i=0}^{n} r_{ij})a_j + \mathop{\sum}\limits_{i=1}^{n} r_{ji}a_i;
$$
$$
a_ja_k = \displaystyle\frac{1}{2} (1 - \mathop{\sum}\limits_{i=0}^{n} r_{ij})a_j + \displaystyle\frac{1}{2}\mathop{\sum}\limits_{i=0}^{n} r_{ji}a_i + \displaystyle\frac{1}{2} (1 - \mathop{\sum}\limits_{i=0}^{n}r_{ki})a_k + \displaystyle\frac{1}{2}\mathop{\sum}\limits_{i=0}^{n}r_{ki}a_i \quad (j\neq k).
$$

\vskip6pt

{\bf Proposition 5.} Let $G_m$ be a gametic algebra of mutation, with mutation rates $r_{ij}$. For $G_m$ to be a Jordan algebra, it is necessary and sufficient that the following system of $n$ identities in $x_k$ ($k = 0, 1, \ldots, n$) holds:
$$
\mathop{\sum}\limits_{j=1}^{n}\beta_j(r_{0i} - r_{ji}) + \beta_i(1 + \mathop{\sum}\limits_{k=0}^{n}r_{ik}) = 0, \;\; i =1, \ldots, n,
$$
where
$$
\alpha_t = \alpha_t(x_0, \ldots, x_n) = (\mathop{\sum}\limits_{k=0}^{n}x_k)r_{0t} - \mathop{\sum}\limits_{k=1}^{n}x_kr_{kt} + \mathop{\sum}\limits_{k=0}^{n}x_tr_{tk},
$$
$$
\beta_j = \beta_j(x_0, \ldots, x_n) = (1 - \mathop{\sum}\limits_{k=0}^{n}r_{jk})\alpha_j - \mathop{\sum}\limits_{t=1}^{n}\alpha_t(r_{0j} - r_{tj}).
$$

\vskip6pt

In the same paper [542], Piacentini Cattaneo used the identities from Proposition 5 to determine the restrictions of the $r_{ij}$'s for $G$ to be a Jordan algebra in specific cases.

In 1988, Peresi [536a] proved that if $A$ is a nonassociative algebra that verifies $A^2 = A$ and has an idempotent, than $A$ and its duplicate have isomorphic automorphism groups and isomorphic derivation algebras. This result is then applied by Peresi to the gametic algebra for polyploidy with multiple alleles.

\vskip6pt

{\bf Definitions.} An algebra $A$, not necessarily associative, over a commutative field ${\bf K}$ of characteristic different from two, that admits a nontrivial homomorphism $w : A \to {\bf K}$ is said to be {\it baric}.

A baric algebra $A$ that satisfies the identity $(a^2)^2 = w^2(a)a^2$ for all $a \in A$ is called a {\it Bernstein} algebra.

\vskip6pt

{\bf Remark.} Singh and Singh [624] showed that Lie and Clifford algebras are never baric. On the other hand, starting with a baric algebra, it is possible to derive new algebras which are Lie, Jordan, alternative or associative.

\vskip6pt

Every Bernstein algebra $A$ possesses at least one idempotent $e$. It can be decomposed into the direct sum of subspaces $A = E \oplus U \oplus Z$ with $E := {\bf K}_e$, $U := \{ey \mid y \in {\rm Ker} w\}$, $Z := \{z \in A \mid ez = 0\}$.

If $A$ has finite dimension, which is at least $1$, dim $A = 1+n$, then one can associate to $A$ a pair of integers $(r+1, s)$, called {\it type} of $A$, whereby
$$
r := \rm{dim} U \quad , \;\; s := \rm{dim} Z,
$$
hence $r + s = n$.

In 1989, W\" orz-Busekros [711d] showed that for each decomposition $n = r+ s$ there exists a Bernstein algebra of type $(r+1, s)$. Thereby the so-called {\it trivial} Bernstein algebra of type $(r+1, s)$ has been introduced as Bernstein algebra of the corresponding type where $({\rm Ker} w)^2 = \{0\}$.

W\" orz-Busekros [711c] showed that the well-known decomposition of a Bernstein algebra with respect to an idempotent is nothing else but the Peirce decomposition known for finite-dimensional power-associative algebras with idempotent, especially for Jordan algebras with idempotent.

\vskip6pt

{\bf Note.} Bernstein algebras are not in general power-associative.

\vskip6pt

In terms of Peirce theory, W\" orz-Busekros [711c] showed that in a Bernstein algebra all idempotents are principal and thus primitive. Hence, the Peirce decomposition cannot be further decomposed. She deduced a necessary and sufficient condition for a Bernstein algebra to be Jordan, and obtained a number of special results from it (the principal two being Proposition 6 and Theorem 7 below).

\vskip6pt

{\bf Proposition 6.} (see W\"orz-Busekros [711c, p. 396]) {\it A trivial Bernstein algebra of type $(r+1, s)$ is a special Jordan algebra.}

\vskip6pt

{\bf Remark.} Proposition 6 is a generalization of Holgate's result [348a] (see Proposition 4 above and Remark which follows), who proved that all gametic algebras for simple Mendelian inheritance are special Jordan algebras. Thereby the gametic algebra for simple Mendelian inheritance with $n+1$ alleles is a trivial Bernstein algebra of type $(n+1, 0)$, cf. W\" orz-Busekros [711d].

\vskip6pt

{\bf Definition.} Let $A$ be an algebra over ${\bf K}$ with weight homomorphism $w : A \to {\bf K}$. Then $A$ is called a {\it normal} algebra, if the identity $x^2y = w(x)xy$ is satisfied in $A$.

\vskip6pt

{\bf Theorem 7.} (see W\" orz-Busekros [711c, p. 397]). {\it Every normal algebra is a Jordan algebra.}

\vskip6pt

In 1988, Walcher [699c] gave a characterization of Bernstein algebras which are Jordan algebras (called by him {\it Jordan Bernstein algebras}) over a field of characteristic different from $2$ or $3$, and listed some of their properties.

\vskip6pt

{\bf Theorem 8.} (see Walcher [699c, p. 219]). {\it Let $A$ be a baric algebra over a field of characteristic different from $2$ or $3$, and with the nontrivial homomorphism from the definition of $A$. The following statements are equivalent:

(i) $A$ is a Jordan Bernstein algebra

(ii) $A$ is a power-associative Bernstein algebra

(iii) $x^3 - w(x)x^2 = 0$ for all $x \in A$.}

\vskip6pt

As a corollary of Proposition 1 from Walcher [699c], it follows that every Jordan Bernstein algebra is genetic. Thus, by W\" orz-Busekros [711a, Theorem 3.18], for dim $A = m+1$, we have a chain of ideals of $A$
$$
N := {\rm Ker} w \supset N_1 \supset N_2 \supset \ldots \supset N_m \supset \{0\},
$$
such that dim $N_i = m+1 - i$ and $N_iN_j \subset N_{k+1}$, where $k := \max \{i,j\}$, for all $i$ and $j$.

\vskip6pt

{\bf Notation.} Let $c$ be an idempotent of $A$ from Theorem 8 above, and let $L(c)$ denotes, as usual, the left multiplication by $c$.

\vskip6pt

{\bf Proposition 9.} (see Walcher [699c, p. 221]). Let $A$ be a Jordan Bernstein algebra of dimension $m+1$. Then there exists a basis $\{v_1, \ldots, v_m\}$ of $N$ such that $v_i$ is an eigenvector of $L(c)$ for $1 \leq i \leq m$ and $N_i$ is spanned by $v_i, \ldots, v_m$ ($1 \leq i \leq m$).

\vskip6pt

{\bf Comments.} The Walcher's results [699c] should at least make the construction of Jordan Bernstein algebras a manageable task: Start with a basis of eigenvectors in $N$ (the eigenvalues preassigned) take into account the composition rules for the eigenspaces and note that the only thing to be checked besides this is the identity $x^3 = 0$ in $N$.

In 1989, Holgate [348e] examined conditions under which the entropic law is satisfied in genetic algebras, and the consequences of imposing it when it is not. It appears that, as with the Jordan identity (see Holgate [348a], and Micali and Quattara [485]), the entropic law only interacts inclusively with the properties of genetic algebras for small rank or dimension.

\vskip6pt

Making use of the papers [557a,b] by Resnikoff, let us mention now applications of Jordan algebras to color perception.

\vskip6pt

In order to endow the set $\mathcal{C}$ of perceived colors with a geometrical structure, various standard experimental results are taken as axioms. One can show that there exists a real vector space $\mathcal{V}$ spanned by the set $\mathcal{C}$ in which $\mathcal{C}$ is a cone of perceived colors. Denote by $GL(\mathcal{C})$ the group of orientation-preserving linear transformations of $\mathcal{V}$ which preserve the cone $\mathcal{C}$. $GL(\mathcal{C})$ is a subgroup of $GL(\mathcal{V})$, and therefore a Lie group.

Making use of standard results in the theory of homogeneous spaces, $\mathcal{C}$ can be identified with the homogeneous space $GL(\mathcal{C})/K$, where $K$ is isomorphic to the subgroup of $GL(\mathcal{C})$ which leaves some point of $\mathcal{C}$ fixed, hence to a closed subgroup of the orthogonal group, and consequently to a compact subgroup of $GL(\mathcal{C})$.

Finally it follows that $\mathcal{C}$ is a homogeneous space equivalent either to ${{\bf R}}^{+} \times {{\bf R}}^{+} \times {{\bf R}}^{+}$ or to ${\bf R}^{+} \times SL(2, {\bf R})/SO(2)$, ${\bf R}^{+}$ denoting the positive real numbers.

The $GL(\mathcal{C})$-invariant metric (see (**) below) yield in the first case Stiles's generalization [640] of Helmholtz's color metric [336], and in the second a new color metric with respect to which $\mathcal{C}$ is not isometric to a Euclidean space.

Resnikoff [557b] showed how the concept of Jordan algebra provides an unification of both cases.

Namely, let $\mathcal{J}$ be a (finite-dimensional) formally real Jordan algebra and consider $\exp \mathcal{J} := \{\exp a \mid a \in \mathcal{J}\}$. 

Consider on $\mathcal{J}$ the form given by
$$
\mu(a) := \displaystyle\frac{s}{n}Tr L(a), \quad a \in \mathcal{J}.
$$

If $\mathcal{J} = {\bf R} (= M_1({\bf R}^{(+)}))$, then $\mu(a) = a$, while if $\mathcal{J} = M_r({\bf R})^{(+)}$, $\mu(a) = Tr\;a$.

It can easily be seen that for $\alpha > 0$ the map $a \to a/\alpha$ is an isomorphism of ${\bf R}$ onto a Jordan algebra $\mathcal{J}_{(\alpha)}$ with unit element $1/\alpha$ and that
$$
\exp \mathcal{J}_{(\alpha)} = \{\exp \alpha a \mid a \in \mathcal{J}\} = \{(\exp a)^{\alpha} \mid a \in \mathcal{J}\} = \{x^{\alpha} \mid x \in \exp \mathcal{J}\}.
$$

Writing $\mathcal{J}_{(\alpha_1\alpha_2\alpha_3)} := \mathcal{J}_{(\alpha_1)} \oplus \mathcal{J}_{(\alpha_2)} \oplus \mathcal{J}_{(\alpha_3)}$, it follows that
$$
\exp \mathcal{J}_{(\alpha_1\alpha_2\alpha_3)} = \{(x_1^{\alpha_1}, x_2^{\alpha_2}, x_3^{\alpha_3}) \mid x_1 \in {\bf R}^{+}\} = {\bf R}^{+} \times {\bf R}^{+} \times {\bf R}^{+}
$$

and
$$
\exp M_2({\bf R})^{(+)} = \{ 
\left(
\begin{array}{cc}
x_1 &  x_3 \\ [5pt]
x_3 & x_2
\end{array} \right)
= x \mid x \; \mbox{is positive definite} \}.
$$

\noindent Thus,

$\exp \mathcal{J} = \mathcal{C} =$ space of perceived colors

if $\mathcal{J} = \mathcal{J}_{(\alpha_1\alpha_2\alpha_3)}$ \quad or \quad $\mathcal{J} = M_2({\bf R})^{(+)}$.

The group $GL(\exp \mathcal{J})$ is generated by the map $P(a)$ for $a \in \mathcal{J}$ ($P(a)$ being the quadratic representation of $\mathcal{J}$); $\exp \mathcal{J}$ is a homogeneous space of $GL(\exp \mathcal{J})$, and a $GL(\exp \mathcal{J})$-invariant metric on $\exp \mathcal{J}$ is given by
$$
ds^2 := \mu((P^{-1}(x)dx)dx). \leqno(**)
$$

\vskip6pt

{\bf Remark.} With the unification provided by the concept of a Jordan algebra, the arguments concerning brightness can be conceptually reversed (see Resnikoff [557b, pp. 122-123]).

\vskip6pt

As is widely accepted nowadays, proteins are the principal workhorses of the cell. They are the major organizers and manipulators of biological energy and enzymes that catalyze and maintain the life process. The proteins are responsible for the active transport of ions into and out of the cell, as well as for cellular and intracellular movement. That is why the discipline of bioenergetics, which is study of how cells generate and transfer their energy supply, is primarily the investigation of how proteins work.

At present, the composition and three-dimensional structure of about two-hundred proteins are known. However, there is no generally accepted model of how proteins operate dynamically.

The idea that the energy released in the hydrolysis of adenosine triposphate (ATP) molecules tramsforms into that of soliton excitation and is transferred with great efficiency along protein molecules was used by Davydov as early as 1973 (see [197a]) to explain the contraction machanics of transversely striated muscles of animals at the molecular level. Davydov et al. considered, in addition, the idea that $\alpha$-helical proteins may facilitate electron transport through a soliton mechanism. In this case, an extra electron causes a lattice distortion in the protein that stabilizes the electron's motion.

Thus it may be reasonable to consider charge transfer across membranes, energy coupling across membranes, and energy transport along filamentous cytoskeletal proteins in terms of a soliton mechanism, since proteins that carry out these functions contain structural units with significant $\alpha$-helical character (see Davydov [197b]). The Davydov model leads to a nonlinear Schr\" odinger equation which has solutions (see [447, p. 13]).

\vskip6pt

{\bf Note.} As was observed by Lomdahl, Layne and Bigio [447, p. 16], the soliton model is one among several concepts for protein dynamics which should attract the careful attention of biologists. Clearly, it cannot explain every aspect of protein dynamics, but it is motivating exciting questions and new experiments.

\vskip6pt

Layne [433] presented a simplified theoretical model for anesthesia activity, taking advantage of the fact that the $\alpha$ helix is an important structure in membrane and cytoskeletal proteins.

More precisely, Layne [433, p. 24] formulated the following question:

How does the binding of the anesthetic molecule to a protein modify normal protein behavior? He answered this question using the soliton model as a paradigm for normal protein functioning. The soliton model proposes that $\alpha$-helical proteins effect the transport of ATP hydrolysis energy through a coupling of vibrational excitations to displacements along the spines of the helix. This coupling leads to a self-focusing of vibrational energy that has remarkably stable qualities. Layne [433, p. 24] suggests that the binding of an anesthetic molecule to a protein interferes with soliton propagation. He suggests further that this type of interference is most important in two seperate reegions of a cell where soliton propagation is an attractive candidate: first, in the $\alpha$-helical proteins of the inner mitochondrial membrane, which appear to participate in ATP synthesis and electron transport and secondly, in the membrane proteins of neurons, which are responsible for chemical reception and signal transduction.

\vskip6pt

{\bf Remark.} [433, p. 26]. If the Davydov soliton finds experimental support in biology, then such a model may help to explain some of the molecular mechanisms behind general anesthesia.

\vskip6pt

Let us mention that Takeno [654a] studied vibron (i.e., vibrational exciton) solitons in one-dimensional molecular crystals by employing a coupled oscillator-lattice model. Takeno showed that although vibron solitons in his theory and those in the Davydov theory are both described by the non-linear Schr\" odinger equation, their nature is fairly different from each other. The nonlinear Schr\" odinger equation arises in the Takeno theory from modulations of vibrons by nonlinear coupling with acoustic phonons propagating along helics of the $\alpha$-proteins, while that in the Davydov theory follows immediately from the quantal Schr\" odinger equation for the exciton probability. In 1985, Takeno [654c] presented an exactly tractable model of an oscillator-lattice system which is capable of incorporating borh of the pictures of Fr\" ohlich (see [277], [108]) and that of Davydov in a unified way and to make a more datailed study of vibron solitons by giving a significant improvement of the theory developed in [654a].

\vskip6pt

{\bf Comments.} As was already pointed out (see end of \S 4), an open problem is to find an algebraic description of the Grassmann manifolds appearing in Sato's approach [595a,b] to soliton equations, resembling the Jordan algebra description of finite-dimensional Grassmann manifolds given by Helwig in [337b]. Taking into account the previous considerations, solving this open problem could be useful in bioenergetics.

\vskip6pt

In 1992 it was published the book "Mathematical Structures in Population Genetics" [452d] by Lyubich, which is the English version of the Russian edition published in 1983.

\vskip6pt

In 2005, Bremner [135] used computer algebra to show that a linearization of the operation of intermolecular recombination from theoretical genetics satisfies a nonassociative polynomial identity of degree 4 which implies the Jordan identity. The representation theory of the symmetric group is used to decompose this new identity into its ireducible components. Bremner showed that this new identity implies all the identities of degree $\leq 6$ satisfied by intermolecular recombination.

\vskip6pt

Concerning the applications of Jordan algebras to mathematical statistics, the year 1994 witnessed the appearance of the book [466] by Malley, where the use of Jordan algebras in mathematical statistics is presented. The kinds of such applications are presented in the above mentioned book, namely: applications to random quadratic forms (sums of squares) and applications to the algebraic simplification of maximum likehood estimation of patterned covariance matrices. In the second chapter the use of Jordan algebras in random quadratic forms and the mixed linear model are presented. Jordan algebras are used to obtain maximal extensions of the work of Cochran [181], and Rao \& Mitra [556] on the independence and chi-squared distribution of quadratic forms in multivariate normal random variables. The second chapter concludes with results that unify previous work on mixed models. In the third chapter, details on Jordan algebras are presented. In the last chapter of this book [466], Malley uses Jordan algebras, Galois field theory and the EM algorithm to obtain either closed-form or simplified solutions to the maximum likehood equations for estimation of covariance matrices for multivariate normal data.

\vspace{1cm}

\pagebreak

\centerline{\bf \Large \S 7. JORDAN STRUCTURES IN PHYSICS}

\vspace{48pt}

I shall mention here only some ideas, while for a detailed presentation the reader is referred to the book Iord\u anescu [364w, Ch. 5], which is a {\bf revised, extended} and {\bf up-dated} version of the 2003 monograph [364g]. It is worth mentioning that the above mentioned 2003 monograph [364g] - exhausted in two years - has good reviews in Zentralblatt f\" ur Mathematik and Mathematical Reviews (see Zbl 107317014 and MR 1979748), and - concerning the applications to physics - is more comprehensive than books published afterwards (in 2004 and 2005), written by famous mathematicians (K. McCrimmon and Y. Friedman), as leading experts in the field (W. Bertram and H. Upmeier) remarked in their reviews of those books (see the review [98u] by Bertram and the review [682o] by Upmeier). Jordan structures are still very much in use in modern mathematical physics, but one has to include the so-called Freudenthal triple systems, which are a little bit more complicated\footnote{UPMEIER. H., private communication (March 2010)}.

\vskip6pt

Concerning already {\bf classical} results, let us mention that - as it is well known - Jordan [379a] stressed that the most fruitful attempt at generalizing of the sstandard Hilbert space structure of quantum mechanics would be to change the algebraic structures (see also the more recent opinion expressed by Dirac [208]). Jordan [379b] formulated a quantum mechanics in terms of commutative, but nonassociative (finite-dimensional) algebras of abservables, now called (finite-dimensional) Jordan algebras. Jordan, von Neumann and Wigner [380] showed that this approach is equivalent to the construction of the standard quantum mechanics in finite-dimensional subspaces of the physical Hilbert space with the single exception of $H_3(\bf{O})^{(+)}$. The infinite-dimensional case was studied by von Neumann in [515]. A more recent tentative axiomatization was given by Emch [236a]. Exceptional quantum mechanics was investigated by G\" unaydin, Piron and Ruegg [318] and shown to be in accordance with the standard propositional formulation, with a unique probability function for the Moufang (non-Desarguesian) plane. Pedroza and Vianna's results [532] concerning the dynamical variables for contrained and unconstrained systems described by the symmetric formulation of classical mechanics can be connected with results on supersymmetry and supermanifolds of Berezin [88c]. Araki [35] improved the characterization of state spaces of JB-algebras given by Alfsen \& Shultz [12b] to a form with more physical appeal (proposed by Wittstock [708]) in the simplified case of a finite dimension. JB-algebras were fruitfully used by Guz [322b] in a tentative axiomatization for nonrelativistic quantum mechanics and Kummer [424] gave in 1987 a new approach. Results on Jordan (quantum) logics due to Morozova \& Chentsov [498a,b], and on order unit spaces arising from sum logics due to Abbati \& Mani\` a [1a,b], and also Bunce \& Wright's results [146a,b] must be mentioned.

\vskip6pt

{\bf Note.} For a comprehensive presentation of the above mentioned {\bf classical} results see the mimeographed monograph Iord\u anescu [364g, Ch. VIII]. Anyway, the bibliography of the book Iord\u anescu [364w] contains all the references, as well as herein.

\vskip6pt

Let us mention now the construction, due to
Truini and Biedenharn [672b], of a quantum mechanics for the
complexified octonion plane. This plane, denoted by $\mathcal{P}(J)$, as they showed (see
[672b, p.~1337]), has automorphism group large enough to
accomodate -- as finite-dimensional quantum-mecha\-nical charge
spaces -- a color-flavor structure which is not ruled out by
current experimental evidence. The Truini-Biedenharn
construction makes essential use of Jordan pairs. The
construction of a quantum-mechanics over a complex octonion
plane was begun by G\"{u}rsey [320a,~b], without, however, using
the concepts of linear ideals or Jordan pairs.

\vskip6pt

{\bf Remark.} The Truini-Biedenharn plane ${\cal P}
(J)$ has a {\it nonprojective} geometry (two lines may
intersect in more than one point) and, consequently, the
propositions system {\it is not a lattice}.

\vskip6pt

As it is well known, the language of quantum mechanics has
always been identified with the language of projective
geometry, the points of the geometry being identified
with the density matrices of the (pure) states, and the
lines and hyperplanes with the propositions which are not
atoms. The automorphism group of the geometry (that is,
its collineation group) is, however, larger than the
automorphism group of the quantum structure, because
collineations need not preserve the traces (which are the
canonical measure defining the quantum states) nor
orthogonality, which has no projective meaning. In
mathematical language we can say that the quantum logic
requires an automorphism group which preserves an
elliptic polarity.

Truini and Biedenharn [672b] defined the propositional
system as follows: the propositions are identified with
the geometrical objects (points and lines correspond to
the principal inner ideals of $V$). They form a partially
ordered set, with ordering given by the set inclusion of
the inner ideals. The plane itself (i.e., the principal
inner ideal generated by an invertible element) is the
trivial proposition. We have an orthocomplementation $a
\rightarrow a^{\bot}$, which is the standard polarity
$a_{\star} \rightarrow a^{\star}$. Thus we can define
orthogonality: $a \bot b$ if $a<b^{\bot}$, which is
symmetric.

\vskip6pt

{\bf Remark.} As it is well known, the lattice axiom is
the axiom least justified experimentally since it is
nonconstructive. It is the merit of the Truini-Biedenharn
construction that it provides a model in which this axiom
is denied in a natural way.

\vspace{6pt}

Because of the lack of a lattice structure, the
definition of ``state" given by Truini and Biedenharn [672b] was suitable a
``measure" with unusual properties
thereby being defined. However, this measure coincides
with the unique probability function (defined by
G\"{u}naydin, Piron and Ruegg [318] on the Moufang plane)
when restricted to the {\it real} octonion case.
Moreover, when restricted to the purely {\it complex}
case, the measure coincides with the usual modulus
(squared) of complex three-dimensional Hilbert space
quantum mechanics.

\vskip6pt

{\bf Open Problem.} (see [672b, p.~1329]). To obtain some
kind of physical understanding of the role of the
connected points which are responsible for all unusual
features of Truini-Biedenharn quantum mechanics.

\vspace{6pt}

Finally, let us mention the opinion of Truini and
Biedenharn [672b, p.~1328] that ``It is our belief (noting
the close relationship between geometries and quantum mechanics)
that the concepts of quadratic Jordan
algebras
and inner ideals will be useful in physics."

\vskip6pt

{\bf Note.} For details on Truini \& Biedenharn's paper [672] see also Iord\u anescu [364w, \S 1 of Ch. 5].

\vskip6pt

Let us refer now on the survey paper [365]
by Iord\u anescu and Truini, where an informal introduction
to quantum groups is given, and the
attention on the relationship among
quantum groups, integrable models and Jordan structures
was, in particular, called.

The {\it historical} and the {\it basic approach} of quantum group theory
presented in the second section of [365] can be completed with more
information by using, for instance, the surveys by Biedenharn [107a,~b],
Dobrev [214], Drinfeld [222a,~b], Faddeev [241],
Kundu [425], Majid [461a], Ruiz-Altaba [574], Smirnov [626], Takhtajan
[656b] -- used also by Iord\u anescu and Truini [365] -- and, for an exhaustive
information, the papers referred therein. For a deeper analysis of the concept
of {\it extended enveloping algebra} and in particular for an exhaustive discussion on the
ring of functions of the Cartan generators necessary in the construction of the extended enveloping
algebra we refer to Truini and Varadarajan [673b].

I like to point out that in 1993, Boldin, Safin and
Sharipov [117] proved a surprising connection between Tzitzeica surfaces
and the inverse scattering method (see also [582] and [622], as well as the
more recent papers [313] and [599]). The transformation that generates the
family of such surfaces
found by Tzitzeica [687] in 1907 and its slight generalizations obtained
by Jonas [376a,~b] in 1921 and 1953 are known in the modern
literature on integrable equations as Darboux or B\"acklund
tranformations\footnote{It is worthy of being mentioned with this occasion that the
title of the paper [376a] by Jonas contains {\it the definition} of
Tzitzeica surfaces (see Bibliography).}.
They are used to construct the soliton solutions starting
from some trivial solution of the equation $u_{xy} = {\rm e}^u - {\rm e}^{-2u}$.
It is worth mentioning
that the paper [678] by Tzitzeica seems to be the first in the world
where
the equation $u_{xy} = {\rm e}^u - {\rm e}^{-2u}$ (the nearest relative of the
sine-Gordon
equation $u_{xy} = \sin u$) was considered.

\vskip6pt

{\bf Comment.} For a presentation of the main ideas in the work of Gheorghe Tzitzeica, see Teleman [660c,d] and Teleman \& Teleman [659a,b,c].

\vskip6pt 

In 1999, from the existing methods of singularity analysis only,
Conte, Musette and Grundland [184] derived the two equations which define the
B\"{a}cklund transformation of the Tzitzeica equation. This is achieved by defining
a truncation in the spirit of the approach of Weiss et al., so as to preserve
the Lorentz invariance of the Tzitzeica equation. If one asumes a third-order
scattering problem, then this truncation admits a unique solution, thus leading to a
matrix Lax pair and a Darboux transformation. In order to obtain the B\"{a}cklund
transformation, which is the main new result in [184], one represents the Lax pair by an
equivalent two-component Riccati pseudopotential. This yields two different B\"{a}cklund
transformations: the first one is a B\"{a}cklund transformation for the
Hirota-Satsuma equation, while the second is a B\"{a}cklund transformation for the
Tzitzeica equation. One of the two equations defining the B\"{a}cklung transformation
is the fifth ordinary differential equation of Gambier.

Also in 1999, Grundland and Levi [312] have shown that there is a strong relationship
between Riccati equations and B\"{a}cklund transformations for integrable nonlinear
partial differential equations. As it has been established in many of the well-known cases
(see [2]), the simplest B\"{a}cklund transformation is given by the classical
first-order Riccati equation. There are, however, a few well-known cases in which
the simplest B\"{a}cklund transformation is given by a higer-order differential equation.
Grundland and Levi proved by a few examples (the Sawada-Kotera equation, the Tzitzeica
equation, and the Fitzhugh-Nagumo equation) that in such a case the  B\"{a}cklund transformation
is given by a higher Riccati equation, higher in the so-called Riccati chain.

The fact, shown in [312], that higher-order Riccati equations may also play a role in the
construction of B\"{a}cklund transformations for some nonlinear partial differential
equations indicates the possibility of introducing higher-order conditional symmetries.
This result can open the way to the construction of new classes of exact solutions for
many physical important differential equations (see [127], [313]).

\vskip6pt

{\bf Note.}\ Grundland and Levi mention at the end of their paper [312] that work on the
extension of the above mentioned results to the case of matrix\break Riccati chains and their
reduction in application to nonlinear partial differential equations is in progress.

\vskip6pt

In 1999, Ferapontov and Schief [255] reviewed some of the most important geometric properties of the Demoulin surfaces and constructed
a B\"{a}cklund transformation which may be specialized to the well-known B\"{a}cklund
transformation for the Tzitzeica equation governing affine spheres in affine geometry.

In 1997, Magri, Pedroni and Zubelli [459] tackled the problem of interpreting the
Darboux transformation (see Darboux [191]) for the {\rm KP} hierarchy and its relations
with the modified {\rm KP} hierarchy from a geometric point
of view. This is achieved by introducing the concept of a {\it Darboux covering}. They
constructed a Darboux covering of the {\rm KP} equations and obtained a new hierarchy which
they called the {\it Darboux-KP hierarchy} ({\rm DKP}). Then they used the {\rm DKP} equations
to discuss the relationships among the modified {\rm KP} equations and the discrete {\rm KP} equations.

In 1999, Fastr\'{e} [247] proposed a Grassmannian definition for the Darboux transformation.
It generalizes two other previously used version of the Darboux
transformation. This new Darboux transformation is used to build $\tau$ functions from the solutions
of the bispectral problem of Duistermaat and Gr\"{u}nbaum. These $\tau$ functions are connected to
the study of the $W$-algebra and the Virasoro algebra (highest weight vectors).

In 1996, Bakalov, Horozov and Yakimov [67] defined B\"{a}cklund-Darboux
transformation in Sato's Grassmannian (see Sato [595a]), which can be regarded as Darboux transformations on maximal
algebras of commuting ordinary differential operators. They described the
action of these transformations on related objects: wave functions, $\tau
$-functions, and spectral algebras.

In 2001, Musette, Conte, and Verhoeven [503] studied the B\"{a}cklund transformation, nonlinear super-position formula of the Kaup-Kupershmidt and Tzitzeica equations.

In 2005, at a mathematical conference in Lubbock (USA), Erxiao Wang (University of Texas at Austin)
gave a talk entitled {\it Transformations of affine spheres}. He identified the classical
Tzitzeica transformations for affine spheres as dressing actions of rational
twisted loop group element and, then he discussed the permutability
formula and the group structure of these transformations.

\vspace{6pt}

{\bf Comments.} Taking into account of the previous considerations on the importance of Tzitzeica
surfaces in the framework of recent mathematical physics researches, I think that it would be correct
to call in the future the {\it B\"acklund-Darboux} transformations as {\it B\"acklund-Darboux-Tzitzeica}
transformations. In fact, in 1998, on the occasion of the anniversary of  125 years
from the birth of Tzitzeica -- organized by the Faculty of Mathematics of the University of Bucharest --
I suggested since then the above mentioned completion in the terminology.

\vskip6pt

{\bf Remark.} Another contribution that I like to mention
here is that of Beidar, Fong, and Stolin [76] which showed that every
Frobenius algebra over a commutative ring determines a class of solutions
of the quantum Yang-Baxter equation, which forms a subbimodule of its
tensor square. Moreover, this subbimodule is free of rank one as a left
(right) submodule. An explicit form of a generator is given in terms of
the Frobenius homomorphism. It turns out that the generator is invertible
in the tensor square if and only if the algebra is Azumaya.

\vspace{6pt}

At the end of the eighties other approaches to quantum groups were given. The
objects of these approaches -- which can be called {\it quantum matrix\break
groups} -- are Hopf algebras in duality to quantum algebras.

\vspace{6pt}

{\bf Definition.} Two Hopf algebras ${\cal A}$ and ${\cal A}'$
are said to
be in {\it duality} if there exists a doubly nondegenerate bilinear
form
$$\langle  \cdot\, ,\cdot \rangle  : {\cal A} \times {\cal A}'
\to {\bf C}, \quad  \langle  \cdot\, ,\cdot \rangle  : (a, a') \to \langle  a, a'\rangle , $$
such that, for $a, b \in {\cal A}$ and $a', b' \in {\cal A}'$ the
following relations hold
$$\langle  a, a'b' \rangle  = \langle  \Delta_{\cal A} (a), a' \otimes b'\rangle, \quad
\langle  ab, a'\rangle  = \langle  a \otimes b, \Delta_{{\cal A}'} (a')\rangle ,$$
$$\langle  1_{\cal A}, a'\rangle  = \varepsilon_{{\cal A}'} (a')\ ,\
\langle  a, 1_{{\cal A}'}\rangle \ =\ \varepsilon_{\cal A} (a),
\quad \langle  S_{\cal A} (a)\ ,\ a'\rangle = \langle  a, S_{{\cal A}'} (a')\rangle .$$

\vskip6pt

{\bf Remark.} I want to emphasize that quantum groups are
also involved and studied in many different fields of mathematics and physics. Among them are:
topological quantum field theories, 2-dimensional gravity and
3-dimensional Chern-Simons theory [287, 314, 461a, 707c], rational conformal
field theory [21, 296, 497], braid and knot theory [377, 391], non-standard
quantum statistics (see GREENBERG, O.W. [in {\it Proc. Argonne Workshop on Quantum Groups}, T. Curtright, D. Fairlie and C.
Zachos (eds.), World Scientific, Singapore, 1990]), quantum Hall effect [409, 594].

\vspace{6pt}

Concerning {\it the relations of Jordan structures with quantum groups},
I would like to recall here the new topic proposed by Truini
and Varadarajan at the end of their paper [673a], namely: {\it
quantization of Jordan structures}.

Okubo [523f] considered the space $V$ endowed also with a bilinear non-degenerate
form $\langle  x|y\rangle $ satisfying
$$\langle  y|x\rangle \ = \epsilon \langle  x|y\rangle,  \quad\epsilon = \pm 1.$$

\vspace{6pt}

Let $R(\theta ) \in {\rm End} (V)\otimes {\rm End}(V)$ be the scattering matrix
with matrix elements $R^{dc}_{ab}$, defined by $R(\theta )e_a\otimes e_b =
R^{dc}_{ab} e_c \otimes e_d$, with respect to a basis $\{ e_j\} $ of $V$
and suppose that $R$ satisfies the {\it quantum Yang-Baxter equation} ({\it QYBE})
$${R}_{12} (\theta) {R}_{13} (\theta^\prime)
{R}_{23} (\theta^{\prime \prime}) =
{R}_{23} (\theta^{\prime \prime}) {R}_{13} (\theta^\prime)
R_{12}(\theta)  \leqno(7.1{\rm a})$$
with
$$\theta^\prime = \theta + \theta^{\prime \prime}. \leqno(7.1{\rm b})$$

\noindent Two $\theta$-dependent triple linear products $[x,y,z]_\theta$ and
$[x,y,z]^*_\theta$ are defined in terms of the scattering matrix
elements $R^{dc}_{ab} (\theta)$, by
$$\big[e^c,e_a,e_b \big]_\theta := e_d R^{dc}_{ab}
(\theta), \quad \big[ e^d,e_b,e_a \big]^*_\theta := R^{dc}_{ab}
(\theta)e_c
$$
or, alternatively, by
$$R^{dc}_{ab} (\theta) =\langle  e^d | \big[ e^c,e_a,e_b\big]_\theta\rangle =
\langle  e^c | [e^d,e_b,e_a ]^*_\theta\rangle,  $$
where $e^d$ is given by
$$\langle  e^d |e_c\rangle = \delta^d_c . $$
The QYBE (7.1a) can be then rewritten as a triple product
equation
$$
\sum^N_{j=1} [v,[u,e_j,z]_{\theta^\prime} ,
[e^j,x,y]_\theta ]^*_{\theta^{\prime \prime}}
\!= \! \sum^N_{j=1} \big[ u,[v,e_j,x]^*_{\theta^\prime} ,
[e^j,z,y]^*_{\theta^{\prime \prime}} \big]_\theta.
\leqno(7.2)$$

\vspace{6pt}

{\bf Proposition 1.} {\it Let $V$ be a Jordan or anti-Jordan triple system with
$\epsilon =1$ satisfying the following conditions}

\hskip6pt i) $\langle  u|xvy\rangle = \langle  v|yux\rangle$;

\hskip3pt ii) $\langle  u|xvy\rangle = \delta \langle  x|uyv\rangle = \delta \langle  y|vxu\rangle$;

iii) $(y e^j x)v e_j = a\{ \langle  x|v\rangle y + \delta \langle  y|v\rangle x\}
+ byvx $;

\hskip0.5pt iv) $(ye^jx) v(zue_j) -
(ye^jz)u(xve_j) = \alpha \{\langle  v|x\rangle zuy -
\langle  u|z\rangle xvy\}
+ \beta \{ \langle  v|y\rangle xuz -$
$- \langle  u|y\rangle  zvx+ \langle  y|uzv\rangle x
- \langle  y|vxu\rangle z\} + \gamma \{(yux)vz - (yvz) ux\}$

\noindent {\it for some constants $a,b, \alpha, \beta$, and
$\gamma$. Then
$$[x,y,z]_\theta = P(\theta) y x z + B(\theta)
\langle  x|y\rangle z + C(\theta)\langle  z|x\rangle y $$
for $P(\theta) \not= 0$ is a solution of the QYBE $(7.2)$ with
$${B (\theta) \over P (\theta)} = \delta \gamma + k \theta,
\quad {C(\theta) \over P(\theta)} =
{\beta \delta \over k \theta} $$
for an arbitrary constant $k$, provided that we have either}

\hskip3pt i) $\alpha = \beta = 0$,

\noindent {\it or}

ii) $\alpha = \beta \not= 0$, $b = -2 \gamma$, $a = 2 \beta$.

\vspace{6pt}

{\bf Remark.} The solution satisfies the unitarity condition
$$R(\theta)\, R(-\theta) = f(\theta) {\rm Id},$$
where
$$f(\theta) = P(\theta) P(-\theta) \left[ (a+\gamma^2) - (k \theta)^2
- {\beta^2 \over (k \theta)^2} \right]. $$

\vspace{6pt}

{\bf Proposition 2.} {\it Let $V$ be the Jordan triple system defined on
the vector space of the Lie-algebra $u(n)$ by means of the product
$$x y z = x \cdot y \cdot z + \delta \ z \cdot y \cdot x
$$
the dot denoting the usual associative product in $V$ and let $\langle  \cdot\, |\, \cdot \rangle $
be the trace form. Then,
$$[x,y,z]_\theta = P(\theta) x z y + A(\theta) \langle  y|z\rangle x +
C(\theta) \langle  z|x\rangle y $$
for $P(\theta) \not= 0$ offers solutions of the QYBE $(7.2)$ for
the following two cases}:
$$ {\rm (i)} \quad
\frac{A(\theta)}{P(\theta)} = -
\frac{\lambda^2
{\rm e}^{k \theta} - d }{\lambda ({\rm e}^{k \theta}-d)}, \quad \frac{C(\theta)}{P(\theta)} = -
\frac{{\rm e}^{k \theta} - \lambda^2}
{\lambda ({\rm e}^{k \theta}-1)}, \leqno(7.3{\rm a})$$
\noindent {\it where $d$ is either $\lambda^2$ or $-\lambda^4$ and $k$ is an
arbitrary constant, or}
$${\rm (ii)} \quad \frac{A(\theta)}{P(\theta)} = - \lambda,
\quad
\frac{C(\theta)}{P(\theta)} = -
\frac{1}{\lambda}.\leqno(7.3{\rm b})
$$

{\it In both cases $\lambda$ is given by}
$$\lambda = \frac{1}{2}( n \pm \sqrt{n^2 -4}). $$

\vspace{6pt}

{\bf Remark.} The first solution Eq. (7.3a) satisfies both
unitarity and crossing symmetry relations:
$$R (\theta) R (-\theta) =
C(\theta) C(-\theta) {\rm Id} \leqno(7.4{\rm a})$$
$$\frac{1}{P (\overline{\theta})}[y,x,z]_{\overline{\theta}}=
\frac{1}{P(\theta)}[x,y,z]_\theta, \leqno(7.4{\rm b})$$
where $\overline{\theta}$ in Eq. (7.4b) is related to $\theta$ by
$$\theta + \overline \theta = \frac{1}{k}\log d. $$
In view of these, the solution is likely related [603] to some exactly
solvable two-dimensional quantum field theory.

\vspace{6pt}

{\bf Remarks.} The Yang-Baxter as well as classical
Yang-Baxter equations have been recast as triple product
equation, and some solutions of these equations are obtained by
Okubo [523g]. In his paper [250], Fauser emphasized a new direction for the
application of Okubo's method.

\vspace{6pt}

I want to mention a result by Svinolupov [648b] which is interesting
in the context of this paper. He considered systems of nonlinear
equations which, in a particular case, may be reduced to the nonlinear Schr\"odinger
equation and are therefore called generalized Schr\"odinger equations. A one
to one correspondence between such integrable systems and Jordan pairs is
established. It turns out that {\it irreducible}
systems correspond to {\it simple} Jordan pairs. Later, Svinolupov [648c]
showed that to every finite-dimensional Jordan algebra with unity there corresponds
a series of integrable vector systems. Among them there are vector
analogues of the famous scalar equations as KdV, modified KdV
and sine-Gordon. Some relations between a generalization of the
Miura transformation connected with Jordan algebra is found.

In 1994, Svinolupov and Yamilov [650] built upon the results\break of
Svinolupov [648b]: multidimensional generalizations of the B\"{a}cklund\break transformations
$$\widetilde{u} =u_{xx} -u^{-1}u_x^2-u^2v, \quad \widetilde{v}
=-u^{-1}, $$
connecting two solutions $(\widetilde{u}, \widetilde{v}),
(u,v)$ of the pair of coupled nonlinear equations
$$u_t=u_{xx}-2u^2v, \quad v_t=-v_{xx}+2v^2u, $$
in $(1+1)$ dimensions, are discussed as well as
corresponding B\"{a}cklund transformations for some
integrable generalizations of the Toda chain.

In a 1993 note, Svinolupov and Sokolov [649a] deal with
{\it Jordan tops}, which are a special class of
systems of quadratic first order differential equations with
coefficients in a Jordan algebra. The general solution is
found for some particular systems.

\vskip6pt

As it is well-known, a major tool in transferring results between Jordan
theory and Lie theory is the Kantor-Koecher-Tits construction of a Lie algebra.
 This construction was first formulated for linear Jordan algebras and was modified
 as Jordan theory expanded to {\it quadratic} Jordan algebras, Jordan triple systems,
 and Jordan pairs. The reverse construction starts with a three-term $\ZZ$-graded
 Lie algebra and recovers a Jordan structure. However, since the Lie algebra involves
 only linear operations, it is impossible to recover the quadratic operations
unless the base ring contains 1/2. The main purpose of the 
 paper [248k] by Faulkner is to give constructions which substitute for the Kantor-Koecher-Tits
 construction and its reverse construction and which work for the quadratic operations for
 all base rings. The role of the Lie algebra is replaced by a certain kind of
 Hopf algebra. The primitive elements form a three-term $\ZZ$-graded Lie algebra, but
 the Hopf algebra also contains divided power sequences which capture the quadratic
 properties of the Jordan pair.

 In the paper [82], Benkart and Perez-Izquierdo have constructed a quantum analogue
 of the split octonions and have studied its properties. By its construction, the quantum
 octonion algebra is a nonassociative algebra with a Yang-Baxter operator action coming from  the $R$-matrix
 of $U_q(D_4)$.

 \vspace{6pt}

 {\bf Comments.} A quantized octonion algebra using the representation theory of
 $U_q (sl_2)$ was independently constructed by Bremner in his 1997 preprint {``Quantum octonions"}. Although both
 Bremner's quantized octonions and Benkart-Perez-Izquierdo's quantum octonions reduce
 to the octonions at $q=1$, they are different algebras.
 
\vskip6pt

 {\bf Remark.} As Perez-Izquierdo [537] remarked, since the Albert algebra is
 close related to the octonion algebra, maybe it could be interesting to have it at
 hand when quantizing Jordan algebras.

\vskip6pt

I shall present some interesting results that lead us to belive
that octonions or exceptional Jordan algebra should play
an important role in recent fundamental physical theories,
namely, in the theory of superstrings.

Let us briefly recall that the exceptional Jordan algebra
made a dramatic appearence within the framework of
supergravity theories through the work of G\"{u}naydin,
Sierra and Townsend [319a,~b,~c,~d]. In their work on the
construction and classification of $N=2$ Maxwell-Einstein
supergravity theories, they showed that there exist four
remarkable theories of this type that are uniquely
determined by simple Jordan algebras of degree three.
These are the Jordan algebras of $(3 \times 3)$-Hermitian
matrices over ${\bf R}, {\bf C}, {\bf H}$ and ${\bf O}$, and
denoted,respectivelly, by $H_3({\bf R})^{(+)}, H_3({\bf C}
)^{(+)}, H_3({\bf H})^{(+)}$, and $H_3({\bf O})^{(+)}$. Their
symmetry groups in five, four and three space-time dimensions give the
famous magic square. From this largest one, namely the
exceptional $N=2$ Maxwell-Einstein supergravity defined by
$H_3({\bf O})^{(+)}$ emerge all the remarkable features of the
maximal $N=8$ supergravity theory in the respective
space-time dimensions. In [265a,b, 271a,b,d,e] it was
speculated that a larger theory that includes the
exceptional $N=2$ theory and the $N=8$ theory may provide
us with a unique framework for a realistic unification of
all known interactions. Such a theory, if it exists, may
well turn out to be a string theory (see G\"{u}naydin and
Hyun [317a, p.~498]).

\vspace{6pt}

{\bf Remark.} The work by G\"{u}naydin, Sierra and
Townsend was reviewed by Truini in [671] which aimed at
indicating the {\it usefulness} and {\it naturalness} of implementing the Jordan pair language in such theory.

\vspace{6pt}

A crucial question in superstring theory is the following:
What mathematical structures have a large degree of
uniqueness and can be associated with strings? Foot and
Joshi suggested in [271a] that the exceptional Jordan
algebra may be such a structure. This algebra is indeed
unique as it is the only formally real Jordan algebra
whose elements cannot be expressed in terms of real
matrices. Although quantum mechanically superstring
theories appear to be consistent only in ten space-time
dimensions, classically superstring theories are
consistent in space-time dimensions of $3,4,6$ and $10$.
These dimensions are suggestive of the sequence of division
algebras ${\bf R}, {\bf C}, {\bf H}$ and ${\bf O}$ whose
respective dimensions correspond to the number of transverse degree
of freedom in $d=3,4,6$ and $10$. These remarks prompted
Foot and Joshi [271a] to look for mathematical structures
which automatically single out $d=3,4,6$ and $10$ with
$d=10$ perhaps appearing special. They investigated the
sequence of Jordan algebras $H_3({\bf K})^{(+)}$ over $
{\bf K} = {\bf R}, {\bf C}, {\bf H},{\bf O}$ and showed that variables of the
superstring can be interpreted as elements of the exceptional Jordan
algebra $H_3({\bf O})^{(+)}$.

The other algebras in this sequence correspond
to classical superstring theories. One of the motivations
for introducing the sequence of algebras $H_3({\bf K})^{(+)}$
is that it is naturally supersymmetric: for
$H_3({\bf R})^{(+)}$, the spinor corresponds to a Majorana
spinor of $SO(2,1)$, for $H_3({\bf C})^{(+)}$, the spinor
corresponds to a Weyl spinor of $SO(3,1)$, for
$H_3({\bf H})^{(+)}$, the spinor corresponds to a Weyl spinor
of $SO(5,1)$, and for $H_3(\bf O)^{(+)}$, the spinor
corresponds to a Majorana-Weyl spinor of $SO(9,1)$. In
each case the number of spinor degrees of freedom agrees
with the number of vector degrees of freedom. Thus the
sequence automatically incorporates equal Bose and Fermi
degrees of freedom.

In Foot and Joshi's approach [271a], transverse Lorentz
rotations are contained in the automorphism group of the
algebra $H_3({\bf K})^{(+)}$.

In conclusion, the classical superstring theories can be
expressed in a unified way using sequence $H_3({\bf K})^{(+)}$.
Furthermore, the $d=10$ case is especially interesting as
it corresponds to the exceptional Jordan algebra
$H_3({\bf O})^{(+)}$.

As Foot and Joshi pointed out [271a], the Green-Schwarz [304a,~b]
superstring is not the only mathematically consistent
candidate for a unified theory of all interactions.
Nevertheless, Foot and Joshi analysed the superstring
because of its central role in the other string theories.
Of particular interest is the heterotic string (see Gross,
Harvey, Martinec and Rohm [309a,~b]), which can incorporate
the exceptional gauge group $E_8 \otimes E_8$. The
appearence of the exceptional group $E_8 \otimes E_8$ is
interesting because $E_8$, like $F_4$, can be related to
octonions.

In 1986, Witten [707a] made some interesting remarks
concerning a new approach to string field theory. Witten
attempted to interpret the interactions of the open
bosonic string in terms of noncommutative differential
geometry. Furthermore he suggested that closed bosonic
strings may be connected with some kind of commutative but
nonassociative algebra.

Motivated by Witten's ideas, Foot and Joshi investigated
in [271b] the incorporation of Jordan algebras, to obtain a
manifestly {\it commutative} but {\it nonassociative}
string theory. Namely, they showed that the free
bosonic string theory can be reformulated using the special Jordan
algebra. Then they proceded to incorporate the exceptional
Jordan algebra into the bosonic string. This leads to an
exceptional group structures at the level of first
quantization, which they interpreted as the appearence of
the gauge group.

\vskip6pt

{\bf Remark.} The appearence of the transformation
group $SO(8)$ in Foot-Joshi's approach [271b] suggests that
a matrix of the exceptional Jordan algebra with fixed
eigenvalues may be related to $d=10$. It may thus be
possible to incorporate this work into the heterotic
string, which consists of closed bosonic strings in
$d=26$, and $d=10$ fermionic strings [309a,~b].

\vspace{6pt}

Let us mention now the work of Li, Peschanski, and Savoy [443] by
which a generalization of no-scale supergravity
models is presented, where scale transformations and
axion-like classical symmetries of the superstrings in
four-dimensions are explicitly realized as dilatations and
translations of the scalar fields in the K\"{a}hler
manifold. A sufficient condition is that the (dimension
one) dilation field matrices can be arranged in matrices
of a Jordan algebra. This determines four possible classes of irreducible
manifolds which are symmetric spaces.

\vskip4pt

Goddard, Nahm, Olive, Ruegg and Schwimmer [293] analysed the
algebraic structure of dependent fermions, namely ones interrelated by the
vertex operator construction. They are associated with special
sorts of lattice systems which are introduced and discussed. The
explicit evaluation of the relevant cocycles leads to the results
that the operator product expansion of the fermions is related in a
precise way to one or other of the division algebras given by
${\bf C}, {\bf H}$, or ${\bf O}$. In the paper [573], Ruegg showed that from the
fermionic operator product expansion one can define a product with the same
algebraic properties as the Jordan product.

The Goddard-Nahm-Olive-Ruegg-Schwimmer octonion result has an
important physical application in the formulation of the
superstring theory of particle interactions. The fermionic vertex
operators related to octonions are associated with short roots of
$F_4$ and fall into three orbits under the action of the Weyl group
of $D_4$, the subalgebra of $F_4$ defined by its long roots
$D_4=so(8)$ is the residue of the Lorentz invariance group of the
superstring in the light cone gauge. In superstring theory the
fermionic vertex operators are familiar and important
constructions. For points of the orbit constituting vector weights
of $D_4$ they are Ramond-Neveu-Schwartz fields. For one of the
other two orbits, comprising spinor or conjugate spinor weights,
they are the fermion emission-absorbtion vertices (see [294]).
Thus the algebra of these quantities which is essential to the
evaluation of superstring scattering amplitudes appears to be related to the
algebra of octonions or to the exceptional Jordan algebra $H_3
({\bf O})^{(+)}$.

Corrigan and Hollowood [186a,b] discussed how to represent Jordan
algebras in terms of superstring vertex operators (see [293]). The
analysis was explicitly carried out in the case of the exceptional
Jordan algebra but applied similarly to non-exceptional Jordan
algebras.

Let us mention also the works of Fairlie and Manonge [243],
Sierra [620], Chapline and G\"{u}naydin [160], and G\"{u}rsey
[320c] who speculated on the possible role that the exceptional Jordan
algebra may play in the framework of string theories.

\vskip4pt

{\bf Note.} For comments on a series of papers by Foot \& Joshi [271c,d,e,f,g], published between 1988 and 1992, see \S 3 of Ch. 5 from the book Iord\u anescu [364w].

\vskip4pt

G\"{u}rsey [320e] used the Kaluza-Klein technique of compactification
to reveal the deep correspondence between lattices generated by
discrete Jordan algebras and symmetries of superstrings, which
suggests that all known superstring theories are related and
descend from a more general theory related to the Conway-Sloane
transhyperbolic group.

\vskip6pt

G\"{u}naydin and Hyun [317a] gave a stringly construction of the
exceptional Jordan algebra $H_3 ({\bf O})^{(+)}$. Specifically,
they constructed $H_3({\bf O})^{(+)}$ using Fubini-Veneziano vertex
operators. This is a very special application of a general vertex operator construction
of nonassociative algebras and their affine extensions developed by
G\"{u}naydin [316b]. This construction gives not only $H_3({\bf O})^{(+)}$ but
also its natural affine extension in terms of the vertex operators.

G\"{u}rsey [320d] considered the discrete Jordan algebras of
$(1 \times 1)$- $(2 \times 2)$- and $(3 \times 3)$-Hermitian matrices
over integer elements of the four division algebras ${\bf R}, {\bf C}, {\bf
H}$ and ${\bf O}$. They are transformed under discrete subgroups of groups associated with the magic
square. Points corresponding to a discrete Jordan matrix belong to a
lattice generated by Weyl reflections that are expressed by means of
Jacobson's triple product. Special cases include the $O(32), E_8 \times
E_8$ and $E_{10}$ lattices that occur in superstring theories.

\vspace{6pt}

{\bf Note.} Various aspects of the connection between K\"{a}hler
manifolds and string theories are examined by Rajeev [555] (see also
Bowick and Rajeev [125]), Zanon [727], Cecotti, Ferrara, Girardello and
Porrati [157].

\vspace{6pt}

{\bf Comments.} The case of $(3 \times 3)$-Hermitian octonionic matrices is of
particular interest because it corresponds to the exceptional Jordan algebra.
As it is well known -- and we partially recalled before -- there have been
numerous attempts to use this algebra to describe quantum physics, which was
in fact Jordan's original motivation. In 1994 and 1996, Schray [600a,~b] showed how
to use the exceptional Jordan algebra to give an elegant description of the
superparticle, which Dray, Janesky, and Manogue [220] have been attempting
to extend to the superstring. Their dimensional reduction scheme extends
naturally to this case [221], and they believe it is the natural languange to
describe the fundamental particles of nature.

\vskip6pt

As it is well known, a {\it soliton} is a  nonlinear wave whose properties
are characterized as follows:

A. a localized wave propagates without changing its properties (shape,
velocity, etc.);

B. localized waves are stable against mutual collisions and each wave
conserves its individuality.

The first property has been  known in hydrodynamics since the middle of
the last century as a solitary wave condition. The second means that the
localized wave behaves like a particle. In modern physics, a suffix ``on"
implies the particle property, for instance, phonon and photon. In 1965,
emphasizing the particlelike behaviour of the solitary wave, Zabusky and
Kruskal called waves with the properties A and B ``soliton". (For more
hystorical details, see, for instance, Wadati and Akutsu [698].)

For an elementary introduction to Sato's theory we refer the reader to
the paper [522] by Ohta, Satsuma, Takahashi and Tokihiro. Starting with an
ordinary differential equation, introducing an infinite number of time
variables, and imposing a certain time dependence on the solutions, they
obtained the Sato equation which governs the time development of the
variable coeficients. It is shown that the generalized Lax equation, the
Zakharov-Shabat equation and the inverse scattering transform scheme are
generalized from the Sato equation. It is also revealed that the $\tau
$-function becomes the key function to express the solution of the Sato
equation. By using the results of the representation theory of groups,
they showed that the $\tau$-function is governed by the partial
differential equations in the bilinear forms which are closely related to
the Pl\"{u}cker relations.

Takasaki, inspired by Sato's theory for soliton equations, gave in [653a, b,~c] a new approach  to the self-dual Yang-Mills equations, which is an
alternative method also based on the view-point of a complet
integrability. It is remarkable, that the self-dual Yang-Mills equations
admit such an approach parallel to Sato's approach to soliton equations.

\vspace{6pt}

{\bf Remarks.} A close relationship with Mulase's method [501a,~b,~c,~d,~e] can be
pointed out. An application of the above-mentioned Takasaki's
approach would be expected to higher dimensional generalizations of gauge
field equations.

\vskip6pt

Another application in eight dimensions was solved by Suzuki [646a,~b]
using Grassmann manifold method. Witten's gauge fields are interpreted by
Suzuki [646c] as motions on an infinite-dimensional Grassmann manifold.
Unlike the case of self-dual Yang-Mills equations in Takasaki's work [653a,~b], 
the initial data must satisfy a system of differential equations
since Witten equations comprise a pair of spectral parameters. Solutions
corresponding to (anti-) self-dual Yang-Mills fields are characterized in
the space of initial data and in application, some Yang-Mills fields
which are not self-dual, anti-self-dual nor abelian can be constructed.

Let us also mention the Jimbo and Miwa's approach to the theory of
solition equations [374]. They considered an inifinite-dimensional Lie
algebra and its representation on a function space. The group orbit of
the highest weight vector is an infinite-dimensional Grassmann manifold.
Its defining equations on the function space, expressed in the form of
differential equations, are then exactly the soliton equations. To put it
the other way, there is a transitive action of an infinite-dimensional
group on the manifold of solutions.

Manin and Radul [468] gave a supersymmetric extension of the one-component
{\rm KP}  hierarchy as the Lax equations. The finite-dimensional version of
the {\rm KP} hierarchy was called  by Ueno the {\it Grassmann hierarchy}. In
the theory of Grassmann hierarchy the fundamental role is played by a
linear algebraic equation which is called the {\it Grassmann equation}.
Ueno and Yamada gave in [680a,~b] a supersymmetric extension of
one-component hierarchies from the viewpoint of the Grassmann equation.
Their approach is slightly different from that of Manin and Radul [468].
Yamada generalized in [714] the results of [680a,~b] to the multicomponent
case. In [680c], Ueno and Yamada revealed that the super
{\rm KP} hierarchy is equivalently transformed to the super Grassmann
equation that connects a point in the universal super Grassmann manifold
with an initial data of a solution.

As Takasaki pointed out [653d], physicists have come to recognize the
revelance of the theory of universal Grassmann manifold (sketched by Sato
[595a]) to physical new topics, such as conformal field theories and
strings (see Ishibashi, Matsuo and Ooguri [366], Vafa [683],
Alvarez-Gaum\'{e}, Gomez and Reina [20], Witten [707b], Kawamoto,
Namikawa, Tsuchiya and Yamada [397], Mickelsson [487a], see also
Arbarello, De Concini, Kats and Procesi [41] for an application to the
moduli geometry of algebraic curves which has a close relation to strings)
and anomalies (see Mickelsson [487b,~c] and Mickelsson and Rajeev [488]).

Almost all of them are based on the framework developed by Segal and
Wilson [605] and Pressley and Segal [550]. Their functional analytical
formulation have a number of advantages, and is now widley recognized as
a standard framework. Admitting this fact, Takasaki has rewrote
everything in the spirit of Sato [595a]. Their highly abstract and
algebraic standpoint is fairly distinct from common sense of most
physicists, who are much more familiar with the use of Hilbert spaces
rather than abstract vector spaces.

\vskip4pt

A particular choice of affine coordinates on Grassmann manifolds, for
both the finite- and infinite-dimensional case, made by Takasaki [653d]
turns out to be very useful for the understanding of geometric
structures therein. The so-called ``Kac-Peterson cocycle", which is
physically a kind of ``commutator anomaly", then arises as a cocycle of a
Lie-algebra of infinitesimal transformations on the universal Grassmann
manifold. These ideas are extended in [653d] to a multi-component
theory. A simple application to a nonlinear realization of current and
Virasoro algebras is also presented for illustration~in~[653d].

Let us mention here the paper [705b] by Wilson, where ``one of
the more puzzling discoveries in the theory of integrable systems", is
re-examined.

Saito [585a] (see also [585b]) showed that the vertex operator of the
three-bosonic-string interaction of Della Selva and Saito (see [199]) is
an element of the universal Grassmann manifold. The correspondence
between string theories and soliton theories is made explicite through
the transformation of evolution parameters of solitons to string
coordinates, the same transformation which relates  Fay's trisecant
formula (see [253]) to  Hirota's bilinear difference equation (see [345a,~b]).

Gilbert [289], based on the approach to infinite Grassmannian as the space
of solutions of {\rm KP} equations (see [616], [501a,~b], [224]), described in
simple terms the infinite sequence of non-linear partial differential
equations (the {\rm KP} equations) and gave possible applications to a
fundamental description of interacting strings. Gilbert also indicated in
[289] lines of research likely to prove useful in formulating a
description of non-perturbative string configurations.

An interesting connection between Witten's string field theory and the
infinite Grassmannian, and the possible characterization of the group
orbit on the Grassmannian by the bilinear identity are examined by Gao [280].

Awada and Chamseddine introduced [56a] the infinite-dimensional graded
Grassmann manifolds in terms of free field operators and studied their
properties. They showed the embedding of the graded ${\rm Diff}\, S^1/S^1$ manifold in the graded Grassmannians, and commented on the
possible supersymmetric {\rm KP} hierarchy.

Let us recall at this point that there are two attractive views of
string theory, both based on holomorphic geometry. The first is the formulation
of quantum string theory as integrable analytic geometry on the universal
moduli space of Riemann surfaces. The second is based on the concept of
loop space and formulated as a holomorphic vector bundle over the
manifold ${\rm Diff} \, S^1/S^1$. In both cases, there exists an
one-to-one embedding of the base manifold into the infinite-dimensional
Grassmannians.

\vskip4pt

In 1987, Awada and Chamseddine [56b] formulated the closed string theory as Hermitian geometry on Grassmannians.

\vskip4pt

{\bf Open Problem.} (see [56b]). Generalise the Awada-Chamseddine approach [56b] to the closed superstring and heterotic string.

\vspace{6pt}

As I already mentioned, Segal and Wilson [605], and Pressley and Segal
[550] developed a framework which is a different approach to infinite
Grassmannians. It consists of the space of choices of fermion boundary
conditions for the free fermion field theory on a disc. In the paper [705a]
is described how the modified KdV equations fit into the
Grassmannian framework, topic not tuched in the paper [605]. In 1988, Witten
[707b] clarified some aspects of the relation between quantum field
theory and infinite-dimensional Grassmannians. More precisely, he
described in physical terminology some aspects of relation, surveyed by
Segal and Wilson [605] between Riemann surfaces and infinite-dimensional
Grassmannians. This relation has been essential in the studies of the
Schottky problem (see Mulase [501b], Shiota [616]), and its relation with
quantum field theory and string theory have been subject of a discussion
from a physical point of view (see Ishibashi, Matsuo, Ooguri [366],
Alvarez-Gaum\'{e}, Gomez, Reina [20], Vafa [683]).

Mickelsson and Rajeev [488] extended the methods of Pressley and Segal [550]
for constructing cocycle representations of the restricted general linear
group in infinite dimensions to the case of a larger linear group modeled
by Schatten classes of rank $p$, $1 \leq p < \infty $ (see Simon [622]). An essential ingredient is the generalization of the determinant line
bundle over an infinite-dimensional Grassmannian to the case of an
arbitrary Schatten rank $p \geq 1$. The results are used to obtain
highest weight representations of current algebras in $(d+1)$ dimensions
when the space dimension $d$ is any odd number.

\vspace{6pt}

{\bf Conjecture} (see [608]). Similar problems to that of Mickelsson and Rajeev [488] must afflict the electric field operators constructed by Semenoff in [608].

\vspace{6pt}

In 1988, Yamagishi [715] pointed out an interesting relation between the
{\rm KP} hierarchy and the extended Virasoro algebra, namely, he showed that
the simply extended {\rm KP} equation has enough information to determine the
extended Virasoro algebra. Levi and Winternitz [441] showed that a class of
integrable nonlinear differential equations in $(2+1)$ dimensions,
including the physically important cylindrical {\rm KP} equation, has a
symmetry algebra with a specific Kac-Moody-Virasoro structure.
Kodama [407] presented a systematic method to produce a class of exact solutions
of the dispersionless {\rm KP} equation, using the conservation equations
derived from the semi-classical limit of the {\rm KP} theory. These exact
solutions include rarefaction waves (global solutions) and shock waves
(breaking solutions in finite time). Zabrodin [724] proved that the
scattering matrix for free massless fermions on a Riemann surface of
finite genus generates the quasiperiodic solutions of the {\rm KP} equation.
The operator changing the genus of the solution is constructed and the
composition law of such operators is discussed. Zabrodin's construction
extends the well-known operator approach in the case of soliton solutions
to the general case of the quasiperiodic $\tau $-functions. David, Levi
and Winternitz [196] constructed a general class of fourth order scalar
partial differential equations, invariant under the same group of local
point transformations as the {\rm KP} equation.

\vskip6pt

Recently, Dimakis \& M\" uller-Hoissen [207a] proved that on any "weakly nonassociative" algebra there is a universal family of compatible ordinary differential equations (provided that differentiability with respect to parameters can be defined), any solution of which yields a solution of the KP hierarchy with dependent variable in an associative subalgebra, the middle nucleus. More recently, they derived in [207b] a sequence of identities in the algebra of quasi-symmetric functions that are in formal correspondence with the equations of the KP hierarchy.

\vskip6pt

De Concini, Fucito and Tirozzi (see [198], [278b]) formulated conformal field
theories on the infinite-dimensional Grassmann manifold. Beside recovering the known
results for the central charge and correlation functions of the $b$-$c$
system, this formalism immediately leads itself to further
generalization. The Grassmann manifold is in fact an {\it ad hoc} model
for the geometrical interpretation of the irreducible representations of
an infinite-dimensional Kac-Moody algebra which, in turn, admit an
intertwining action of a Virasoro algebra. They also gave a proof of
bosonization from a purely Grassmann manifold point of view.

In 1992, Anagnostopoulos, Bowick, and Schwarz [26] determined the space of
all solutions to the string equation of the symmetric unitary one-matrix
model.

\vskip4pt

Aoyama and Kodama [34] generalized the Sato equation of the {\rm KP} theory on
the basis of the $W_{1+\infty}$ algebra. An infinite set of flows which
commute with the {\rm KP} hierarchy (the commuting {\rm KP} flows) are
explicitely constructed. They satisfy the $W_{\infty}$ algebra. Aoyama
and Kodama studied the generalized Sato equation by making a
$p$-truncation of the dressing pseudo-differential operator, instead of
the usual $p$-reduction. The cases of $1$- and $2$-truncations were
studied in some details. The commuting {\rm KP} flows become $W_{\infty }$
operators of the same form for both cases, when operated on the $\tau
$-function, and the central charge if the Virasoro algebra is found to
be~$-2$.

Both finite- and infinite-dimensional integrable systems can be
linearized on orbits of the infinite Abelian group $\Gamma^{+}$ on the
universal Grassmann manifold. The aim of Landi and Reina [428] is to link
these results to standard symplectic dynamics by giving an explicit
Hamiltonian formulation on the symplectic manifold
$$M= Gl_{{\rm res}} ({\cal H}) /Gl({\cal H}_+) \times Gl ({\cal H}_{-}).$$

\vskip6pt

{\bf Note.} Let us recall here that, in 1996, Bakalov, Horozov,
and Yakimov [67] defined B\"{a}cklund-Darboux transformations in Sato's
Grassmann manifold.

\vspace{6pt}

In 1997, van de Leur [439], making use of the representation theory of
infinite matrix group, showed that (in the polynomial case) the
$n$-vector $k$-constrained {\rm KP} hierarchy has a natural geometrical
interpretation on Sato's Grassmann manifold. His description generalizes
the $k$-reduced {\rm KP} or Gelfand-Dickey hierarchies.

In 1992, Bottacin [121a] showed that to each element in a large class of solutions
to the {\rm KP} hierarchy there is associated, in a natural way, a group variety. This
result was achieved by showing that these solutions are in fact functions of theta
type. In [121b] he restricted himself to a class of solutions of particular interest
in physics, namely the so-called $N$-solitons. He proved that an $N$-soliton
 to the {\rm KP} hierarchy is actually a holomorphic theta type and its associated group variety is a
product of multiplicative groups. Hence the same conclusions also hold for the KdV and all
other hierarchies which can be obtained as specializations of the {\rm KP}.
Moreover, all these constructions are made explicitly.

In the note [36], Aratyn develops an explicit construction of the constrained {\rm KP} hierarchy within the Sato Grassmannian framework. Useful relations are established
between the kernel elements of the underlying ordinary differential operator
and the eigenfunctions of the associated {\rm KP} hierarchy, as well as between the related
bilinear concomitant and the squared eigenfunction~potential.

\vskip4pt

In 1999, in Subsection 3.1 of their paper [92], Bergvelt and
Rabin introduced an infinite super Grassmannian and related constructions. The
infinite Grassmannian of Sato [595a] or of Segal $\&$ Wilson [605] consists (essentially)
of ``half-infinite-dimensional" vector subspaces $W$ of an infinite-dimensional vector
space $H$ such that the projection on a fixed subspace $H_-$ has finite-dimensional
kernel and cokernel. In the super category, Bergvelt and Rabin replace this by the
super Grassmannian of free, ``half-infinite-rank" $\Lambda $-modules of infinite-rank,
free $\Lambda $-module $H$ such that the kernel and cokernel of the projection on
$H_-$ are a submodule and a quotient module, respectively, of a free, finite-rank
$\Lambda $-module.

\vskip4pt

{\bf Open Problem.} It is well-known that the infinite-dimensional Grassmann
manifold contains moduli spaces of Riemann surfaces of all genera. Taking into
account this fact, Schwarz [602b] conjectured that non-perturbative string
theory can be formulated in terms of the Grassmannian.

\vspace{6pt}

In [602b] Schwarz presented new facts supporting the conjecture contained in the
above mentioned Open Problem. In particular, it is shown that Grassmannians can
be considered as generalized moduli spaces; this statement permits to Schwartz to define
corresponding ``string amplitudes" (at least formally).

\vskip6pt

{\bf Open Problem.} Another conjecture formulated by Schwarz [602b] is the following: It is
possible to explain the relation between non-perturbative and perturbative  string
theory by means of localization theorems for equivariant cohomology. ({\bf N.B.} This conjecture
is based on the characterization of moduli spaces, relevant to string theory, as sets consisting
of points with large stabilizers in certain groups acting on the Grassmannian.)

\vskip6pt

{\bf Comments.} In the eighties, Sasaki remarked that all
differential equations with soliton solutions (known to him)
could be
considered as immersion equations for a two-dimensional surface
into some
Euclidean space. Later, Chern and Terng have worked ont this idea for
Sine-Gordon and the modified KdV equation, but in contrast to
Sasaki's
method, they considered surfaces of constant negative Gauss curvature. It
is also possible to work with other types of surfaces (e.g., surfaces of
constant mean curvature). Gerold and Buchner [286] gave and
explicit solution of the immersion equations for the one- and two-soliton solution
of the Sine-Gordon equation. The surface corresponding to the one-soliton
solution is Dini's surface (i.e., a helicoid of Gauss curvature ($-1$)
with a tractrix as profile curve), the one corresponding to the
two-soliton solution is a Joachimsthal surface of Enneper's type that has
not yet been discussed explicitly in literature. These results stimulate
the question how the properties of the solitons are transfered to the
surface. More important are the following:

\vskip6pt

{\bf Open Problems.} (see Gerold and Buchner [286, p.~2056]).
Which immersions yield soliton solutions? Can all solitons be obtained in this~way?

\vskip6pt

In this section I shall comment on very recent applications of two
kinds of important mixed Jordan and Banach structures (namely, $JB^*$-triple and $JBW^*$-triples) to quantum mechanics, emphasizing on this occasion the decoherent 
states and the geometry of projection operators.

In 2008, Edwards and H\"ugli (see [228]) have considered pre-symmetric complex Banach spaces which have been proposed as models for state spaces of physical systems.

A complex Banach space $A_*$ is said to be {\it pre-symmetric} if the open unit
ball in its Banach dual space $A$ is a bounded symmetric domain.
In this case $A$ possesses a natural triple product with respect to which it
forms a $JBW^*$-triple with unique predual $A_*$. Consequently,
there exists a bijection between the set of pre-symmetric spaces and the set of 
$JBW^*$-triples. An important property of a $JBW^*$-triple $A$ is that
the group of linear isometries of $A$ to itself coincides with the group
of triple automorphisms $\Aut(A)$ of $A$, and because of the uniqueness of the
pre-dual $A_*$ of $A$, each such mapping is weak*-continuous, which implies that
there exists an isomorphism $\phi \to \phi_*$ from $\Aut(A)$ onto the group
$\Aut(A_*)$ of linear isometries of $A_*$.

\vskip6pt

{\bf Remark.} There exist approaches to the theory of statistical physical systems,
in which a pre-symmetric space $A_*$ represents the state space of the system,
the linear isometries representing the symmetries of the system and the
contractive projections on $A_*$ representing filters on the system (see Friedman
[273a,~b] and Friedman \& Gofman [274a,~b]).

\vskip6pt

Edwards \& H\"ugli [228] have defined structural projections as follows: A contractive
projection $R$ on the pre-symmetric space $A_*$ is said to be structural if, for
each element $x$ of $A_*$ such that $Rx$ and $x$ have equal norm, it follows that
$Rx$ and $x$ coincide, and, if $x$ is an element of $A_*$ for which, for all elements
$y$ in $A_*$,
$$
\|x + Ry\| = \|x\|+\|Ry\|,
$$
then $x$ lies in the kernel of $R$.

The range $RA_*$ of a structural projection is said to be a {\it structural subspace
of $A_*$}.

\vskip6pt

{\bf Remarks.} If pre-symmetric spaces are considered as state spaces, then the
physical aspects are more transparent when viewed as properties 
of structural projections and structural subspaces of $A_*$ (see [228, p.~220]). From the point of view of physical systems, the two conditions that
a contractive projection $R$ must satisfy in order to be structural have powerful
physical motivations (for details, see [228, p.~221]).

\vskip6pt

Two structural projections on the pre-symmetric space represent decoherent operations
when their ranges are rigidly collinear. Edwards and H\"ugli have proved in [228]
that, for decoherent elements $x$ and $y$ of $A_*$, there exists an involutive
element $\phi_x$ in $\Aut(A_*)$ which conjugates the structural projections corresponding
to $x$ and $y$, and conditions are found for $\phi_x$, to exchange $x$ and $y$.
The results are used to investigate when certain subspaces of $A_*$ are the ranges
of contractive projections and, therefore, represent systems arising from filtering
operations.

\vskip4pt

{\bf Remark.} The results of the paper [228] depend crucially upon the detailed study
of pairs of rigidly collinear tripotents in $JBW^*$-triples.

\vskip6pt

In 2007, H\"ugli has investigated in [353c] normal contractive projections
in connection with certain algebraic conditions on generalized operators.

I must mention here the paper [353b] by H\"ugli, where the set of tripotents
in a $JB^*$-triple is characterized in various ways.
This paper is one of the supporting papers for [228], another supporting paper being the paper
[354] by H\"ugli and Mackey. 

A very recent paper of H\"ugli, namely [353d], is based on the results contained in [353c]
and [354]: more exactly, it generalizes the results of [353c] (which is mainly about Hilbert
spaces as subtriples) to a more general class of subtriples (those having the Dunford-Pettis
property).

Edwards and R\"{u}ttimann presented in their survey [230] a tentative
description of certain results, concerning orthomodular partially ordered
sets of tripotents in $JBW^{\ast}$-triples, for an approach to a
characterization of latices of events of physical systems.

\vskip6pt

The expository paper [682k] by Upmeier gives a survey of results
concerning harmonic analysis and quantization of geometric phase spaces
associated with Jordan structures. As Upmeier pointed out (see
[682k], p. 301--302]), ``... recently it has become clear that quantization theory
can be carried to new directions, namely:

(1)~symmetric supermanifolds (replacing Jordan algebras by super-\break algebras),

(2)~deformation quantization (replacing the structure group of a Jordan
algebra by a so-called quantum group), and

(3)~non-convex symmetric cones and tube domains (replacing holomorphic
functions by cohomology classes).

In all these projects it is expected that the fine structure of Jordan
systems will play a crucial role.\,..."

\vskip6pt

{\bf Note.} For a survey of the Upmeier's work related to the Toeplitz-Berezin
type of quantization of symmetric (and related) domains, see [682m].

\vskip6pt

In his big survey [37], Arazy presented some basic facts and developments concerning
invariant Hilbert spaces of analytic functions on bounded symmetric domains, the
Jordan theoretic background being $JB^{\ast}$-triples.

In his paper [273a], Friedman presented in detail two examples in
theoretical physics where $JB^{\ast}$-triples appear. He showed that
the M\"{o}bius-Potapov-Harris transformations of the automorphism group of a
bounded symmetric domain occur as transformations of signals in an
ideal transmission line and as velocity transformations between
two inertial systems in special relativity.

In their 2001 paper, Friedman and Russo [275] proposed a triple product representation
of the canonical anticommutation relations which does not make use of the
associative Clifford algebra. Imposition of these commutation relations on the
natural basis of ${\bf C}^{n}$ defines a triple product making ${\bf C}^n$ into a Cartan factor
of type 4 (called also {\it spin factor}) that the authors denote by~$S^n$.

This Jordan structure (i.e., $JB^{\ast}$-triple structure) is used to represent
the Lorentz group on $S^3$ and $S^4$. The irreducible representation on $S^3$ corresponds
to the relativistic transformations of the electro-magnetic field. The irreducible
spin-1 representation on $S^4$ extends the Lorentz space-time transformation. By
taking the self-adjoint part of this representation with respect to the spin conjugation,
a reducible spin-1/2 representation on $S^4$ results. The latter is shown to induce
two spin-1 representations in the space of determinant preserving maps on $S^4$
showing that the same spin factor could be used to represent the two types of elementary
particles: bosons and fermions.

\vskip6pt

Let us mention some {\bf recent} (or {\bf very recent}) contributions of interest.

\vskip6pt

Bordemann \& Walter [120] showed that for each semi-Riemannian locally symmetric space the curvature tensor gives rise to a rational solution of the classical Yang-Baxter equation with spectral parameters. For several Riemannian globally symmetric spaces $M$ such as real, complex and quaternionic Grassmann manifolds, they explicitly compute solutions of the quantum Yang-Baxter equations (represented in the tangent space of $M$) which generalize results of Zamolodchikov \& Zamolodchikov [726].

\vskip6pt

Nichita \& Popovici [516] used $(\mathbb{G}, \theta)$-Lie algebras (which unify the Lie algebras and Lie superalgebras) to produce solutions for the quantum Yang-Baxter equation. The constant and the spectral-parameter Yang-Baxter equations and Yang-Baxter systems are also studied in [516].

\vskip6pt

Burhan \& Henrich [147], making use of geometric methods due to Polishchuk [545] and Burban \& Kreussler [148], studied unitary solutions of the associative Yang-Baxter equation with spectral parameters.

\vskip6pt

Other papers of interest are [7] by Akrami \& Majid, and [124] by Bouwknegt, Hannabuss, and Mathai. Concerning quantum projective spaces, very recent papers are [189] by D'Andrea, Dabrowski, and Landi, [188] by D'Andrea \& Dabrowski, and [190] by D'Andrea \& Landi.

\vskip6pt

{\bf Concluding comments.} There are almost 80 years since the appearence of the paper [379b] by Pasqual
Jordan, in which he has tried to use {\it for the first time} algebraic structures defined by
himself (called -- since 1946 -- Jordan algebras) for a mathematical description of quantum
mechanics. During this period, Jordan algebras have found applications in projective geometry,
algebraic geometry, differential geometry, various chapters of mathematical analysis, differential
equations, probability, genetics, statistics, and finally also returned to physics.
But, this time, as it is easily remarked, Jordan structures are involved mainly as additional
mathematical structures: one could say -- in few words -- that, concerning the applications of
Jordan structures to physics, {\it they are necessary but not sufficient\,!} A recent 
interesting paper devoted to these considerations was written by Bertram and it is entitled
``Is there a Jordan geometry underlying quantum mechanics?'' (see Bertram [98s]). This rather
big paper (it has 30 pages) is addressed ``to readers coming from physics rather than from mathematics''
(see [98s, p.~2]) and it is very interesting {\it also} from philosophical point of view.
In this paper, {\it generalized projective geometries} have been proposed as a framework for a
geometric formulation of quantum theory. The above mentioned paper [98s] by Bertram ends as follows:

\vskip6pt

`` ... it could be hoped that Jordan geometry gives some hints on what the last two items on the
following matrix might be:

\vskip6pt

{\footnotesize
\noindent\;\;\;\begin{tabular}{llll}
geometry: & linear; affine & projective & manifold \\
mechanics: & classical & special relativistic & general relativistic \\
quantum theory:& Hilbert space q.m. & projective q.m.? & Cartan geometric q.m.? ''
\end{tabular}
}

\vskip6pt

Let us mention that for {\it Cartan connection} (which generalizes the projective and conformal connnections)
one could read a modern presentation given by Sharpe in [613].

In a more recent note, Bertram [98t] refines his assertion from his previous paper [98s] (i.e.,
generalized projective geometries could be a framework for a geometric formulation of quantum theory),
by discussing further structural features of quantum theory: the link with {\it associative involutive algebras}
$\mathbb{A}$ and with {\it Jordan-Lie} and {\it Lie-Jordan algebras}.
The associated geometries are ({\it Hermitian}) {\it projective lines over} $\mathbb{A}$. Let us recall here that -- more than twenty years ago -- V.S. Varadarajan wrote `` ... quantum mechanical systems are those whose logic form some sort of projective geometry'' (see [686, p.~6]).
It is worth being mentioned here that Bertram wrote his note [98t] in September 2008, after he
attended the XXVII Workshop on Geometrical Methods in Physics from Bia\l owie$\dot{\rm z}$a (Poland), where
he had the occasion to discuss on this topic with physicists.

\vspace{1cm}

\centerline{BIBLIOGRAPHY}

\vspace{48pt}

\small

\begin{description}

\item{[1]} {\sc Abbati, M.C., Mani\`a, A.}, a) {\it A specral theory for order
unit spaces}, Ann. Inst. Henri Poincare {\bf 35} A (1981), {\it 4}, 259--285.

b) {\it Quantum logic and operational quantum mechanics}, Rep. Math.
Phys. {\bf 19} (1984), {\it 3}, 383--406.

\item{[2]} {\sc Ablowitz, M., Fokas, A.S.}, {\it Comments on the inverse scattering
transform and related nonlinear evolution equations}, in  {\it Nonlinear Phenomena},
Lecture Notes in Physics {\bf 189}, K.B. Wolf (ed.), Springer-Verlag,
Berlin, 1983.

\item{[3]} {\sc Abraham, R.}, {\it Picewise differentiable manifolds and the spacetime of general relativity}, J. Math. Mech. {\bf 11} (1962), {\it 4}, 553-592.

\item{[4]} {\sc Abraham, V.M.} a) {\it Linearizing quadratic transformations in genetic algebras}, Proc. London Math. Soc. (3) {\bf 40} (1980), 346-363.

b) {\it The induced linear transformation in a genetic algebra}, ibid. (3) {\bf 40} (1980), 364-384.

c) {\it The genetic algebra of polyploids}, ibid. (3) {\bf 40} (1980), 385-429.

\item{[5]} {\sc Achab, D.}, a) {\it Zeta functions of Jordan algebras representations}, Ann. Inst. Fourier
(Grenoble) {\bf 45} (1995), 1283--1303.

b) {\it Repr\'esentation des alg\`ebres de rang $2$ et fonctions zeta
associ\'ees}, Ann. Inst. Fourier (Grenoble) {\bf 45} (1995), {\it 2}, 437--451.

\item{[6]} {\sc Akivis, M.A.}, {\it Smooth lines on projective planes over certain
  associative algebras}, Math. Notes {\bf 41} (1987), {\it 1-2}, 131--136.
  
\item{[7]} {\sc Akrami, S.E., Majid, S.}, {\it Braided cyclic cocycles and nonassociative geometry}, J. Math. Phys. {\bf 45} (2004), 3883-3911.

\item{[8]} {\sc Alcalde, M.T., Burgue\~ no, C., Labra, A., Micali, A.}, {\it Sur les alg\'ebres de Bernstein}, Proc. London Math. Soc. (3) {\bf 58} (1989), {\it 1}, 51-68.

\item{[9]} {\sc Alekseevsky, D.V.}, {\it Pseudo-K\" ahler and para-K\" ahler symmetric spaces}, in {\it Handbook of pseudo-Riemannian Geometry and Supersymmetry}, Cortes, V. \& Batyrev, V. (eds.), IRMA-series, Publishing House of the EMS, 2007, Chapter 21 (27 pp).

\item{[10]} {\sc Alekseevsky, D.V., Marchiafava, S.}, a) {\it Quaternionic-like
structures on a manifold}:
Note I. 1 -- {\it Integrability and integrability conditions} -- Note II. {\it Automorphism
groups and their interrelations}, Rend. Mat. Acc. Lincei {\bf 4} (1993), 43--52; 53--61.

b) {\it A report on quaternionic-like structures on a manifold}, in {\it Proc. Internat.
Workshop Diff. Geom. and its Appl.} (Bucharest,
Romania, 1993), Scientific Bull. UPB {\bf 55} (1993), 9--34.

c) {\it Almost quaternionic hermitian and quasi-K\"{a}hler manifolds}, in {\it Proc.\
Internat. Workshop Complex Structures} (Sofia, Bulgaria, 1992), World Scient.,
Singapore, 1994, 150--175.

d) {\it Quaternionic structures on a manifold and subordinated structures}, Ann.
Mat. Pura Appl. {\bf 171} (1996), 205--273.

e) {\it Quaternionic geometry}, in {\it Proc. Internat. Workshop Diff. Geom. and
its Appl.} (Bra\c sov, Romania, 1999), 13--19, Bra\c{s}ov, 2000.

\item{[11]} {\sc Alfsen, E.M., Hanche-Olsen, H., Schultz, F.W.}, {\it State
spaces of $C^{\ast}$-algebras}, Acta Math. {\bf 144} (1980), 267--305.

\item{[12]} {\sc Alfsen, E.M., Shultz, F.W.}, a) {\it State spaces of Jordan
algebras}, Acta Math. {\bf 140} (1978), 155--190.

b) {\it State spaces of operator algebras}, Birkh\"auser, Boston, 2001.

\item{[13]} {\sc Alfsen, E.M, Shultz, F.W., St\o rmer, E.}, {\it A Gelfand-Neumark theorem for Jordan
algebras}, Adv. in Math. {\bf 28} (1978), 11--56.

\item{[14]} {\sc Allcock, D.}, a) {\it Identifying models of the octave projective
  plane}, Geom. Dedicata {\bf 65} (1997), 215--217.

b) {\it Reflection groups on the octave hyperbolic plane}, J. Algebra {\bf 213} (1999),
{\it 2}, 467--498.

\item{[15]} {\sc Allison, B.N.}, {\it A class of nonassociative algebras with involution containing the
class of Jordan algebras}, Math. Ann. {\bf 237} (1978), 133--156.

\item{[16]} {\sc Allison, B.N., Azam, S., Berman, S., Gao, Y., Pianzola, A.},
{\it Extended affine Lie algebras and their root systems}, Mem. Amer. Math. Soc.
{\bf 126} (1997), {\it 603}, X + 122.

\item{[17]} {\sc Allison, B.N., Benkart, G., Gao, Y.}, a) {\it Central extensions of Lie algebras
graded by finite root systems}, Math. Ann. {\bf 316} (2000), {\it 3}, 499--527.

b) {\it Lie algebras graded by the root systems ${\rm BC}_r$, $r\geq 2$}, Mem. Amer. Math. Soc.
{\bf 158} (2002), {\it 751}, X + 158.

\item{[18]} {\sc Allison, B.N., Faulkner, J.R.}, a) {\it A Cayley-Dickson process for
  a class of structurable algebras}, Trans. Amer. Math. Soc. {\bf 283} (1984),
  185--210.

b) {\it Isotopy for extended affine Lie algebras and Lie tori}, arXiv: 0709.1181 [math.RA].

\item{[19]} {\sc Allison, B.N., Gao, Y.}, {\it The root system and the core of an extended
affine Lie algebra}, Selecta Math. (N.S.) {\bf 7} (2001), {\it 2}, 149--212.

\item{[20]} {\sc Alvarez-Gaum\'e, L., Gomez, C., Reina, C.}, {\it Loop groups,
Grassmannians and string theory}, Phys. Lett. B {\bf 190} (1987),
55--62.

\item{[21]} {\sc Alvarez-Gaum\'e, L., Gomez, C., Sierra, G.}, in Nucl. Phys. {\bf
319 B} (1989), 155.

\item{[22]} {\sc Alvermann, K.}, a) {\it The multiplicative triangle inequality in non-commutative
$JB^{\ast}$-algebras}, Abh. Math. Sem. Univ. Hamburg {\bf 55} (1985), 91--96.

b) {\it Real normed Jordan Banach algebras with an involution}, Arch.\ Math. (Basel) {\bf 47} (1986),
{\it 2}, 135--150.

c) {\it The dimension lattice of a $JBW$-algebra}, Math. Z. {\bf 195} (1987), {\it 3}, 437--441.

\item{[23]} {\sc Alvermann, K., Janssen, G.}, {\it Real and complex noncommutative Jordan Banach algebras}, Math.
Z. {\bf 185} (1984), {\it 1}, 105--113.

\item{[24]} {\sc Amitsur, S.A.},. {\it Generic splitting fields of central simple algebras}, Ann. of Math. (2) {\bf 62} (1955), 8-43.

\item{[25]} {\sc Ammar, G., Mehl, Ch., Mehrmann, V.}, {\it Schur-like forms for matrix Lie groups, Lie algebras
and Jordan algebras}, Linear Algebra {\bf 287} (1999), {\it 1--3},\break 11--39.

\item{[26]} {\sc Anagnostopoulos, K.N., Bowick, M.J., Schwarz, A.}, {\it The
solution space of the Unitary Matrix Model string equation and the Sato
grassmannian}, Comm. Math. Phys. {\bf 148} (1992), 469--485.

\item{[27]} {\sc Ancochea Bermud\'ez, J.M., Campoamor-Stursberg, R., Garc\'ia Vergnolle, L., S\'anchez Hern\'andez, J.}, {\it Contractions d'alg\'ebres de Jordan en dimension 2}, J. Algebra {\bf 319} (2008), {\it 6}, 2395-2409.

\item{[28]} {\sc Andreoli, G.}, {\it Algebre non associative e sistemi differenziali di Riccati in un problema di genetica}, Ann. Mat. Pura Appl. {\bf 49} (1960), 97-116.

\item{[29]} {\sc Anquela, J.A.}, {\it Changing the element at which primitivity happens}, in
{\it Proceedings of the Intenat.\ Conference on Jordan structures} (Malaga, 1997), 9--13,
Castellon Serrano, Cuena Mira, Fern\'andez L\'opez, Mart\'\i n Gonzalez (eds.), Malaga, 1999.

\item{[30]} {\sc Anquela, J.A., Cort\'es, T.}, {\it Minimal ideals of Jordan systems},
Invent. Math. {\bf 168} (2007), {\it 1}, 83--90.

\item{[31]} {\sc Anquela, J.A., Cort\'es, T., Garc\'\i a, E.}, a) {\it Local and subquotient inheritance
of the heart of a Jordan system}, Manuscripta Math. {\bf 106} (2001), {\it 3}, 279--290.

b) {\it Herstein's theorems and simplicity of Hermitian Jordan systems}, J. Algebra {\bf 246} (2001), {\it 1}, 193-214.

\item{[32]} {\sc Anquela, J.A., Cort\'es, T., Garc\'\i a, E., McCrimmon, K.}, 
{\it Outer inheritance of simplicity in
Jordan systems}, Comm. Algebra {\bf 32} (2004), {\it 2}, 747--766.

\item{[33]} {\sc Anquela, J.A., Cort\'es, T., McCrimmon, K.}, {\it Trivial minimal ideals of Jordan systems}, J. Algebra {\bf 328} (2011), {\it 1}, 167-177.

\item{[34]} {\sc Aoyama, S., Kodama, Y.}, {\it A generalized Sato equation and the
$W_{\infty}$ algebra}, Phys. Lett. B {\bf 278} (1992), {\it 1-2},
56--62.

\item{[35]} {\sc Araki, H.}, {\it On a characterisation of the state space of
quantum mechanics}, Comm. Math. Phys. {\bf 75} (1980), {\it 1}, 1--24.

\item{[36]} {\sc Aratyn, H.}, {\it On Grassmannian description of the constained {\rm KP} hierarchy},
 J. Geom. Phys. {\bf 30} (1999), 295--312.

\item{[37]} {\sc Arazy, J.}, {\it A survey of invariant Hilbert spaces of analytic functions on
bounded symmetric domains}, in {\it Multivariable operator theory}, Curto, R.E.
et al. (eds.), Contemporary Math. {\bf 185} (1995), 7--65.

\item{[38]} {\sc Arazy, J., Engli\v s, M., Kaup, W.}, {\it Holomorphic retractions and boundary Berezin transforms}, Ann. Inst. Fourier, Grenoble, {\bf 59} (2009), {\it 2}, 641-657.

\item{[39]} {\sc Arazy, J., Kaup, W.}, a) {\it On the holomorphic rigidity of linear operators on
complex Banach spaces}, Quart. J. Math. Oxford (2) {\bf 50} (1999), 249--277.

b) {\it On continuous Peirce decompositions, Schur multipliers and the perturbation of
triple functional calculus}, Math. Ann. {\bf 320} (2001), 431--461.

\item{[40]} {\sc Arazy, J.,  Upmeier, H.}, {\it Jordan Grassmann manifolds and intertwining operators
for weighted Bergman spaces}, in {\it Contemporary Geometry and Topology and Related Topics}
(Proc. Internat. Workshop Diff. Geom. Appl., Cluj-Napoca, August 2007), Cluj University Press,
2008, 25--53.

\item{[41]} {\sc Arbarello, E., De Concini, C., Kats, V., Procesi, C.},
{\it Moduli space of curves and representation theory}, a lecture at Amer. Math.
Soc. Summer Institute on Theta Functions, Braunswick, 1987.

\item{[42]} {\sc Arnal, P.M., Laliena, J.A.}, {\it Lie and Jordan structures in associative superalgebras with
superinvolution: special cases}, in {\it Proceedings of the Internat. Conference on Jordan structures} (Malaga,
1997), 15--23, Castellon Serrano, Cuenea Mira, Fern\'andez L\'opez, Mart\'\i n Gonzalez (eds.), Malaga, 1999.

\item{[43]} {\sc Arnlind, J.}, {\it Construction of Hom-Nambu-Lie algebras from Hom-Lie algebras}, talk at AGMP-6 Workshop, Tj\"arn\"o, Oct. 2010.

\item{[44]} {\sc Arnlind, J., Makhlouf, A., Silvestrov, S.D.}, {\it Ternary Hom-Nambu-Lie algebras induced by Hom-Lie algebras}, J. Math. Phys. {\bf 51} 043515(2010), {\it 4}, 11 pp.

\item{[45]} {\sc Aronszajn, N.}, {\it Subcartesian and subriemannian spaces}, Notices Amer. Math. Soc. {\bf 14} (1967), 111.

\item{[46]} {\sc Artin, E.}, {\it Geometric Algebra}, Interscience Publishers, New York--London,
1957.

\item{[47]} {\sc Asadurian, E.}, {\it A pair of Jordan triple systems of Jordan pair type},
An. \c Stiin\c t. Univ. Ovidius Constan\c ta Ser. Math. {\bf 9} (2001), {\it 2}, 1--4.

\item{[48]} {\sc Asadurian, E., \c Stef\u anescu, M.}, {\it Jordan Algebras} (in Romanian),
Editura Nicu\-lescu s.r.l., Bucure\c sti, 2001.

\item{[49]} {\sc Ashihara, T., Miyamoto, M.}, {\it Deformation of central charges, vertex operator algebras whose Griess algebras are Jordan algebras}, J. Algebra {\bf 321} (2009), {\it 6}, 1593-1599.

\item{[50]} {\sc Aslaksen, H.}, {\it Restricted homogeneous coordinates for the Cayley projective
plane}, Geom. Dedicata {\bf 40} (1991), 245--250.

\item{[51]} {\sc Atiyah, M., Berndt, J.}, {\it Projective planes, Severi varieties and spheres}, 
in {\it Surveys in Differential Geometry}, {\bf VIII}, Internat. Press, Somerville, MA, 2003, 1--27.

\item{[52]} {\sc Atsuyama, K.}, a) {\it On the embeddings of the Cayley plane into the exceptional
Lie group of type} $F_4$, Kodai Math. Sem. Rep. {\bf 28} (1977), 129--134.

b) {\it Symmetric spaces derived from algebras}, Kodai Math.\ J. {\bf 3} (1980),
{\it 3},\break 358--363.

\item{[53]} {\sc Atsuyama, K., Tamiguchi, Y.}, {\it On generalized Jordan triple systems and their modifications},
Yokohama Math. J. {\bf 47} (2000), {\it 2}, 165--175.

\item{[54]} {\sc Aupetit, B., Youngson, M.A.}, {\it On symmetry of Banach Jordan algebras}, Proc. Amer. Math. Soc.
{\bf 91} (1984), {\it 3}, 364--366.

\item{[55]} {\sc Aupetit, B., Zra\"{\i}bi, A.}, {\it Propri\'et\'es analytiques du spectre dans les
alg\`ebres de Jordan-Banach}, Manuscripta Math.
{\bf 38} (1982), {\it 3}, 381--386.

\item{[56]} {\sc Awada, M.A., Chamseddine, A.H.}, a) {\it Superstring and graded
Grassmannians}, Phys. Lett. B {\bf 206} (1988), {\it 3}, 437--443.

b) ET-H preprint 87/3 (1987).

\item{[57]} {\sc Ayupov, Sh.A.}, a) {\it On the theory of partially ordered Jordan algebras} (in Russian),
Dokl. Akad. Nauk USSR, 1979, {\it 8}, 6--8.

b) {\it Normal states on $OJB$-algebras} (in Russian), Izv. Akad. Nauk USSR, Ser. Fiz.-Mat. Nauk, 1980, {\it 3}, 9--13.

c) {\it Extension of traces and type criterions for Jordan algebras of self-adjoint operators}, Math. Z. {\bf 181} (1982), 253--268.

d) {\it Integration on Jordan algebras} (in Russian), Izv. Akad. Nauk USSR {\bf 47} (1983), {\it 1}, 3--25.

e) {\it Classification and representations of ordered Jordan algebras} (in Russian), ``Fan", Tashkent, 1986.

f) {\it Traces on $JW$-algebras and enveloping $W^{\ast}$-algebras}, Math. Z. {\bf 194} (1987), {\it 1}, 15--23.

g) {\it On the existence of traces on modular $JW$-factors} (in Russian), Uspehi Mat. Nauk {\bf 44} (1989), 183--184.

\item{[58]} {\sc Ayupov, Sh.A., Abdullaev, R.Z.}, {\it The Radon-Nikodym theorem for weights on
semi-finite $JBW$-algebras}, Math. Z.
{\bf 188} (1985), {\it 4}, 475--484.

\item{[59]} {\sc Ayupov, Sh.A., Iochum, B., Yadgorov, N.Zh.}, {\it Geometry of state space of
finite-dimensional Jordan algebras} (in
Russian), Izv. Akad. Nauk UzSSR, Ser. Fiz.-Mat. Nauk 1990, {\it 3}, 19--22.

\item{[60]} {\sc Ayupov, Sh.A., Rakhimov, A., Usmanov, Sh.}, {\it Jordan, real and Lie structures
in operator algebras}, Mathematics and its Applications,
{\bf 418}, Kluwer, Dordrecht--Boston--London, 1997.

\item{[61]} {\sc Ayupov, Sh.A., Yadgorov, N.Zh.}, {\it The geometry of the state space of modular
Jordan algebras} (in Russian), Izv. Ross. Akad. Nauk. Ser. Mat.
{\bf 57} (1993), {\it 6}, 199--211; translation in Russian Acad. Sci. Izv. Math.
{\bf 43} (1994), {\it 3}, 581--592.

\item{[62]} {\sc Bachmann, F.}, {\it Aufbau der Geometrie aus dem Spiegelangsbegriff}, 2 Auflage, 
Springer-Verlag, Berlin, 1973.

\item{[63]} {\sc Backes, E.}, {\it Geometric applications of Euclidean Jordan triple
systems}, Manu\-scrip\-ta Math. {\bf 42} (1983), {\it 2--3}, 265--272.

\item{[64]} {\sc Bahturin, Y., Goze, M.}, {\it $\mathbb{Z}_2 \times \mathbb{Z}_2$-symmetric spaces}, Pacific J. Math. {\bf 236} (2008), {\it 1}, 1-21.

\item{[65]} {\sc Bahturin, Y., Shestakov, I.}, {\it Gradings of simple Jordan algebras and their
relation to the gradings of simple associative algebras}, Comm. Algebra {\bf 29} (2001), {\it 9},
4095--4102.

\item{[66]} {\sc Bahturin, Y., Shestakov, I.P., Zaicev, M.V.}, {\it Gradings on simple Jordan and Lie algebras}, J. Algebra {\bf 283} (2005), {\it 2}, 849-868.

\item{[67]} {\sc Bakalov, B., Horozov, E., Yakimov, M.}, {\it B\"{a}cklund
Darboux transformations in Sato's Grassmmanian}, Serdica Math. J. {\bf 22}
(1996), 571--586.

\item{[68]} {\sc Barbilian, D.}, a) {\it Galileishe Gruppen und quadratishe Algebren}, mit einen Anhang: {\it Grundiss einer Geometrischen Algebrenlehre},
Bull. Math. Soc. Roumaine Sci. {\bf 41} (1939), {\it 1}, 7-64.

b) {\it Zur Axiomatik der projektiven Ringgeometrien}, I, II,
Jahresber. Deutsch. Math. Verein. {\bf 50} (1940), {\it 1}, 179--229, {\bf 51} (1941), {\it 3}, 34--76.

\item{[69]} {\sc Barton, T.J., Dang, T., Horn, G.}, {\it Normal representations of Banach Jordan triple systems},
Proc. Amer. Math. Soc. {\bf 102} (1988), {\it 3}, 551--555.

\item{[70]} {\sc Bashkirov, E.L.}, {\it Transvection parameter sets in linear groups over
associative division rings and special Jordan algebras defined by these rings}, 
Arch. Math. {\bf 87} (2008), 295--302.

\item{[71]} {\sc Bayara, J., Conseibo, A., Ouattara, M., Zitan, F.}, {\it Power-associative algebras that are train algebras}, J. Algebra {\bf 324} (2010), {\it 6}, 1159-1176.

\item{[72]} {\sc Becerra G.J.}, $JB^{\ast}$-{\it triples and transitivity of the norm}, in
{\it Proceedings of the Internat. Conference on Jordan structures} (Malaga, 1997), 31--37,
Castellon Serrano, Cuenca Mira, Fern\'andez L\'opez, Mart\'\i n Gonzalez (eds.), Malaga, 1999.

\item{[73]} {\sc Becerra, G., Rodr\'\i guez Palacios, A.}, {\it $C^*$- and $JB^*$-algebras
generated by a nonself-adjoint idempotent}, J. Funct. Anal.
{\bf 248} (2007), {\it 1}, 107--127.

\item{[74]} {\sc Behncke, H.}, a) {\it Hermitian Jordan Banach algebras}, J. London Math. Soc.
(2) {\bf 20} (1979), {\it 2}, 327--333.

b) {\it Finite-dimensional representations of $JB$-algebras}, Proc. Amer. Math. Soc. {\bf 88} (1983),
{\it 3}, 426--428.

\item{[75]} {\sc Behncke, H., B\"{o}s, W.}, $JB$-{\it algebras with an exceptional ideal},
Math. Scand. {\bf 42} (1978), {\it 2}, 306--312.

\item{[76]} {\sc Beidar, K.I., Fong, Y., Stolin, A.}, {\it On Frobenius algebras
and the quantum Yang-Baxter equation}, Trans. Amer. Math. Soc. {\bf 349}
(1997), {\it 9}, 3823--3836.

\item{[77]} {\sc Bellissard, J., Iochum, B.}, a) {\it Homogeneous self dual cones versus Jordan algebras.
The theory revisited}, Ann. Inst. Fourier {\bf 28} (1978), 27--67.

b) {\it L'alg\`{e}bre de Jordan d'un c\^{o}ne autopolaire facialement homog\`{e}ne}, C.R. Acad.
Sci. Paris {\bf 288} (1979), 229--232.

c) {\it Spectral theory for facially homogeneous symmetric self-dual cones}, Math. Scand.
{\bf 45} (1979), {\it 1}, 118--126.

d) {\it Homogeneous self-dual cones and Jordan algebras}, in
``Quantum fields-alge\-bras, processes"
(Bielefeld 1978), 153--165 (1980).

\item{[78]} {\sc Benito, P., Draper, C., Elduque, A.}, {\it Invariant connections on symmetric spaces},
in Nonassociative algebra and its applications (Sao Paulo, 1998), 21--33, Lecture Notes in Pures and
Appl. Math. {\bf 211}, Dekker, New York, 2000.

\item{[79]} {\sc Benkart, G.}, {\it Derivations and invariant forms of Lie algebras graded by finite
root systems}, Canad. J. Math. {\bf 50} (1998), {\it 2}, 225--241.

\item{[80]} {\sc Benkart, G., Fern\' andez Lopez, A.}, {\it The Lie inner ideal structure of associative rings revisited}, Comm. Algebra {\bf 37} (2009), {\it 11}, 3833-3850. 

\item{[81]} {\sc Benkart, G., Madariaga, S., P\' erez-Izquierdo, J.M.}, {\it Hopf algebras with triality}, (submited).

\item{[82]} {\sc Benkart, G., Perez-Izquierdo, J.M.}, {\it A quantum octonion algebra},
Trans. Amer. Math. Soc. {\bf 352} (2000), {\it 2}, 935--968.

\item{[83]} {\sc Benkart, G., Smirnov, O.}, {\it Lie algebras graded by the root system BC$_1$},
J. Lie Theory {\bf 13} (2003), {\it 1}, 91--132.

\item{[84]} {\sc Benkart, G., Zelmanov, E.}, {\it Lie algebras graded by finite root systems
and intersection matrix algebras}, Invent. Math. {\bf 126} (1996), 1--45.

\item{[85]} {\sc Benslimane, M., Boudi, N.}, {\it Noetherian Jordan Banach algebras are finite-dimensional},
in {\it Proceedings of the Internat. Conference on Jordan structures} (Malaga, 1997),
Castelon Serrano, Cuenca Mira, Fern\'andez L\'opez, Mar\-t\'\i n Gonzalez (eds.), Malaga, 1999, 53--57.

\item{[86]} {\sc Benslimane, M., Marhnine, H., Zarhouti, C.}, {\it Modular annihilator Jordan pairs},
Extracta Math. {\bf 16} (2001), {\it 1}, 63--77.

\item{[87]} {\sc Benz, W.}, a) {\it On Barbilian domains over commutative rings}, J. Geom.
{\bf 12} (1979), {\it 2}, 146--151.

b) {\it On mappings preserving a single Lorentz-Minkowski-distance in Hjelmslev's
classical plane}, Aequationes Math. {\bf 28} (1985), {\it 1-2}, 80--87.

\item{[88]} {\sc Berezin, F.A.}, a) {\it Quantization}, Math. USSR-Izv. {\bf 8} (1974), 1109--1165.

b) {\it Quantization in complex symmetric spaces}, Math. USSR-Izv. {\bf 9}
(1975), 341--379.

c) {\it Introduction to algebra and analysis with anticommuting variables} (in Russian),
Moscow Univ. Press, 1983.

\item{[89]} {\sc Bergen, J., Grzeszczuk, P.}, {\it Simple Jordan color algebras arising 
from associative graded algebras}, J. Algebra {\bf 246} (2001), {\it 2}, 915--950.

\item{[90]} {\sc Berger, C.A., Coburn, L.A.}, {\it Toeplitz operators and quantum mechanics},
J. Funct. Anal. {\bf 68} (1986), {\it 3}, 273--299.

\item{[91]} {\sc Berger, M.}, {\it Les espaces sym\'etriques non compacts}, Ann. Ec. Norm.
Sup. (3) {\bf 74} (1997), 85--177.

\item{[92]} {\sc Bergvelt, M.J., Rabin, J.M.}, {\it Supercurves, their Jacobians, and
super {\rm KP} equations}, Duke Math. J. {\bf 98} (1999), 1--57.

\item{[93]} {\sc Berline, N., Vergne, M.}, {\it Equations de Hua et noyau de Poisson}, in Lecture
Notes in Math. {\bf 880} (1981), Springer-Verlag, Berlin, 1--51.

\item{[94]} {\sc Berman, S., Gao, Y., Krylyuk, Y., Neher, E.}, {\it The alternative
torus and the structure of elliptic quasi-simple Lie algebras of type $A_2$},
Trans. Amer. Math. Soc. {\bf 347} (1995), 4315--4363.

\item{[95]} {\sc Berndt, J.}, {\it Geometry of weakly symmetric spaces}, in
{\it Proc. Internat. Workshop Diff. Geom. Appl.} (Bra\c sov, Romania, 1999), 66--73, Bra\c sov, 2000.

\item{[96]} {\sc Berndt, J., Tricerri, F., Vanhecke, L.}, {\it Generalized
Heisenberg groups and Damek-Ricci harmonic spaces}, Lecture Notes in
Math. {\bf 1598}, Springer-Verlag, Berlin--Heidelberg--New York,
1995.

\item{[97]} {\sc Berndt, J., Vanhecke, L.}, a) {\it Two natural generalizations of
locally symmetric spaces}, Differential Geometry Appl. {\bf 2} (1992),
57--80.

b) {\it Aspects of the geometry of the Jacobi operator}, Riv. Mat. Univ. Parma
(5) {\bf 3} (1994), 91--108.

\item{[98]} {\sc Bertram, W.}, a) {\it Dualit\'{e} des espaces riemannians sym\'etriques et
analyse harmonique}, Ph.D. Thesis, Universit\'{e} Paris 6, 1994.

b) {\it Un th\'{e}or\`{e}me de Liouville pour les alg\`{e}bres de Jordan}, Bull. Soc. Math.
Francaise {\bf 124} (1996), 299--327.

c) {\it On some causal and conformal groups}, J. Lie Theory {\bf 6} (1996), 215--244.

d) {\it Complexifications of symmetric spaces and Jordan theory}, Preprint, TU Clausthal,
1997, Trans. Amer. Math. Soc. {\bf 353} (2001), {\it 6}, 2531--2556.

e) {\it Algebraic structures of Makarevich spaces}, I, Transformation Groups {\bf 3}
(1998), 3--32.

f) {\it Jordan algebras and conformal geometry}, in {\it Positivity in Lie Theory: Open
Problems}, Hilgert et al. (eds.), Walter de Gruyter GmbH \& Co., Berlin--New York, 1998,
1--20.

g) {\it From vector spaces to symmetric spaces}, in {\it Lie Theory and its Applications
in Physics III}, World Scientific, Singapore, 2000, 99--109.

h) {\it The geometry of Jordan and Lie structures}, Lecture Notes in Math. {\bf 1754},
Springer-Verlag, Berlin, 2000.

i) {\it Generalized projective geometries: general theory and equivalence with Jordan structures}, Adv. in Geometry {\bf 2} (2002), 329--369.

j) {\it The geometry of null systems, Jordan algebras and von Standt's Theorem}, 
Ann. Inst. Fourier {\bf 53} (2003), {\it 1}, 193--225.

k) {\it Complex and quaternionic structures on symmetric spaces-correspondence with
Freudenthal-Kantor triple systems}, in {\it Theory of Lie groups and manifolds},
Sophia Kokyuroku in Math. {\bf 45} (2002), 57--76.

l) {\it Generalized projective geometries} in {\it Proc. Internat. Workshop Diff. Geom.
Appl.} (Timi\c soara, September, 2001), Ann. Univ. de Vest Timi\c soara, Mat.-Inf.
{\bf 39} (2001), 57--65.

m) {\it From linear algebra via affine algebra to projective algebra}, Linear Alg. Appl.
{\bf 378} (2004), 109--134.

n) {\it Differential geometry over general base fields and rings}, in {\it Modern Trends
in Geometry and Topology} (Proc. Internat. Workshop Diff. Geom. Appl., Deva, September 2005), 
Cluj University Press, 2006, 95--101.

o) {\it Homotopes and conformal deformations of symmetric spaces}, J. Lie Theory
{\bf 18} (2008), {\it 2}, 301--333.

p) {\it Differential geometry, Lie groups and symmetric spaces over general base fields
and rings}, Mem. Amer. Math. Soc. {\bf 192} (2008), {\it 900}, 186~+~V.

r) {\it Jordan structures and non-associative geometry}, in {\it Developments and trends in infinite-dimensional Lie theory}, K.-H. Neeb, A. Pianzola (Eds.), Progress in Math. {\bf 288} (2011), Birkh\"auser, New York, 221-241, arXiv: math.RA/0706.1406.

s) {\it Is there a Jordan geometry underlying quantum mechanics?} Internat. J. Theoretical Physics {\bf 47} (2008), {\it 2}, 2754-2782, arXiv: math-ph/0801.3069v1.

t) {\it On the Hermitian projective line as a home for the geometry of quantum theory}, AIP Conference Proceedings {\it 1079}, 14-25, American Institute of Physics, New York, 2008.

u) Review of the book {\it A taste of Jordan algebras} by Kevin McCrimmon, Springer-Verlag, 
New York, 2004, SIAM Review {\bf 47} (2005), {\it 1}, 172--174.

\item{[99]} {\sc Bertram, W., Bieliavsky, P.}, a) {\it Homotopes of symmetric spaces I: Construction by algebras with two involutions}, Jordan Theory preprints No. 296 (13 Nov. 2010).

b) {\it Homotopes of symmetric spaces II: Structure variety and classification}, Jordan Theory preprints No. 297 (13 Nov. 2010).

\item{[100]} {\sc Bertram, W., Didry, M.}, {\it Symmetric bundles and representations of Lie triple systems}, J. Gen. Lie Theory Appl. {\bf 3} (2009), {\it 4}, 261-284.
arXiv: math.DG/0710.1543.

\item{[101]} {\sc Bertram, W., Gl\"ockner, H., Neeb, K.-H.}, {\it Differential calculus over general
base fields and rings}, Expo. Math.
{\bf 22} (2004), 213--282.

\item{[102]} {\sc Bertram, W., Hilgert, J.}, a) {\it Hardy spaces and analytic continuation of Bergman
spaces}, Bull. Soc. Math. Francaise {\bf 126} (1998), 435--482.

b) {\it Geometry of symmetric spaces via Jordan structures}, invited lecture
at the Internat.
Workshop Diff. Geom. and its Appl. (Bra\c sov, Romania, 1999),
Bull. Transilvania Univ. Bra\c{s}ov {\bf 8(43)}, series B (2001), 7--18.

\item{[103]} {\sc Bertram, W., Kinyon, M.}, a) {\it Associative geometries I: Torsors, linear relations and Grassmannians}, Jordan Theory preprints No. 271 (01 Apr. 2009), Lie Theory {\bf 20} (2010), {\it 2}, 215-252.

b) {\it Associative geometries II: Involutions, the classical torsors, and their homotopes}, Jordan Theory preprints No. 276, J. Lie Theory {\bf 20} (2010), {\it 2}, 253-282.

\item{[104]} {\sc Bertram, W., L\"owe, H.}, {\it Inner ideals and intrinsic subspaces},
Adv. in Geo\-metry {\bf 8} (2008), 53--85.

\item{[105]} {\sc Bertram, W., Neeb, K.-H.}, a) {\it Projective completions of Jordan pairs. Part I: 
The generalized projective geometry of a Lie algebra}, J. Algebra {\bf 277} (2004), {\it 2}, 474--519.

b)  {\it Projective completions of Jordan pairs. Part II: 
Manifold structures and symmetric spaces}, Geometriae Dedicata
{\bf 112} (2005), {\it 1}, 73--113.

\item{[106]} {\sc Bhambri, S.K.}, {\it On generalized perspectivities of projective
planes}, Ph.D.~Thesis, Westfield College, Univ. of London, 1981.

\item{[107]} {\sc Biedenharn, L.C.}, a) Invited paper presented at the XVIII
Internat. Coll. on Group Theoretical Methods in Physics (Moscow, USSR, June
4--9, 1990).

b) in Internat. J. Theor. Phys. {\bf 32} (1993) {\it 10}, 1789.

\item{[108]} {\sc Bilz, H., B\" uttner, H., Fr\" ohlich, H.}, Z. Naturfarsch. {\bf 36} B (1981), 208.

\item{[109]} {\sc Bingen, F.}, {\it Geometrie projective sur un anneau semiprimaire},
Acad. Roy. Belg. Bull. Cl. Sci. (5), {\bf 52} (1966), 13--24.

\item{[110]} {\sc Binz, E., Pods, S.}, {\it The geometry of Heisenberg groups},
Amer. Math. Soc. Mathematical Surveys and Monographs {\bf 151}, 2008, 284~pp.

\item{[111]} {\sc Birkhoff, G., von Neumann, J.}, {\it The logic of quantum mechanics},
Ann. Math. (2) {\bf 37} (1936), 823--843.

\item{[112]} {\sc Bix, R.}, a) {\it Octonion plans over local rings}, Trans. Amer. Math. Soc.,
{\bf 261} (1980), {\it 2}, 417--438;

b) {\it Isomorphism theorems for octonion planes over local rings},
Trans. Amer. Math. Soc. {\bf 226} (1981), {\it 2}, 423--440.

c) {\it Octonion planes over Euclidean domains}, Tagungsbericht Jordan-Algebren
(1985), Oberwolfach.

\item{[113]} {\sc Blazic, N., Bokan, N., Gilkey, P., Rakic, Z.}, {\it Pseudo-Riemannian
Osserman manifolds}, Balkan J. Geom. Appl. {\bf 2} (1997), 1--12.

\item{[114]} {\sc Blind, B.}, {\it Distributions homogenes sur une alg\`ebre de Jordan}, Bull. Soc.
Math. France {\bf 125} (1997), {\it 4}, 493--528.

\item{[115]} {\sc Blunck, A.}, {\it Chain spaces over Jordan systems},
Abh. Math. Sem. Univ. Hamburg {\bf 64} (1994), 33--49.

\item{[116]} {\sc Bobenko, A.}, {\it All constant mean curvature tori in
${\bf R}^{3}$, ${\bf S}^{3}$, ${\bf H}^{3}$ in terms of theta-functions},
Math. Ann. {\bf 290} (1991), 209--245.

\item{[117]} {\sc Boldin, A.Yu., Safin, S.S., Sharipov, R.A.}, {\it On an old
article of Tzitzeica and the inverse scattering method}, J. Math.
Phys. {\bf 34} (1993), {\it 12}, 5801--5809.

\item{[118]} {\sc Bommireddy, M., Nagamuni Reddy, L., Sitaram, K.}, {\it Planar ternary rings
with zero of translation, Moufang and Desarguesian planes}, Arch. Math. (Basel)
{\bf 45} (1985), {\it 5}, 476--480.

\item{[119]} {\sc Bonsall, F.F.}, {\it Jordan algebras spanned by Hermitian
elements of a Banach algebra}, Math. Proc. Cambridge Philos. Soc. {\bf 81} (1977), 3--13.

\item{[120]} {\sc Bordemann, M., Walter, M.}, {\it Solutions of the quantum Yang-Baxter equation for symmetric spaces}, arXiv:math/0107018v1 [math.QA] 3 Jul. 2001.

\item{[121]} {\sc Bottacin, F.}, a) {\it Group varieties related to {\rm KP} hierarchy},
Ann. Mat. Pura Appl. (4) {\bf 162} (1992), 215--226.

b) {\it Explicit construction of group varieties related to soliton solutions of the
{\rm KP} hierarchy}, Ann. Mat. Pura Appl. (4) {\bf 164} (1993), 121--132.

\item{[122]} {\sc Boudi, N., Marhnine, H., Zarhouti, C.}, {\it Additive derivations
on Jordan Banach pairs}, Comm. Algebra 
{\bf 32} (2004), {\it 9}, 3609--3625.

\item{[123]} {\sc Boudi, N., Marhnine, H., Zarhouti, C., Fern\'andez L\'opez, A., 
Garcia Rus, E.}, {\it Noetherian Banach Jordan pairs}, Math. Proc. Cambridge Philos. Soc.
{\bf 130} (2001), {\it 1}, 25--36.

\item{[124]} {\sc Bouwknegt, P., Hannabuss, K., Mathai, V.},  {\it Nonassociative tori and applications to $T$-duality}, Commun. Math. Phys. {\bf 264} (2006), 41-69.

\item{[125]} {\sc Bowick, M.J., Rajeev, S.G.}, {\it String theory as the
K\"{a}hler geometry of loop space}, Phys. Rev. Lett. {\bf 58} (1987),
{\it 6}, 535--538.

\item{[126]} {\sc Bowling, J., McCrimmon, K.}, {\it Outer fractions in quadratic Jordan algebras}, J. Algebra {\bf 312} (2007), {\it 1}, 56-73.

\item{[127]} {\sc Bracken, P., Grundland, A.M., Martina, L.}, {\it The Weierstrass-Enneper
system for constant mean curvature surfaces and the completely integrable sigma
model}, J. Math. Phys. {\bf 40} (1999), {\it 7}, 3379--3403.

\item{[128]} {\sc Brada, C.}, {\it Elements de la g\'eom\'etrie des octaves de Cayley}, Dissertation,
Univ. d'Aix-Marseille, I, Marseille, 1984, Univ. Claude-Bernard, D\'epart\'ement de
Math\'ematique, Lyon, 1986.

\item{[129]} {\sc Brada, C., Pecaut-Tison, F.}, {\it G\'eom\'etrie du plan projectif des octaves de
Cayley}, Geom. Dedicata {\bf 23} (1987), {\it 2}, 131--154.

\item{[130]} {\sc Braun, H.}, {\it The Riccati differential equation in Jordan pairs},
Tagungsbereicht 35/1979, Jordan-Algebren (18.8-25.8, 1979), Oberwolfach.

\item{[131]} {\sc Braun, H., Koecher, M.}, {\it Jordan-Algebren}, Springer-Verlag, Berlin--Heidel\-berg--New York, 1966.

\item{[132]} {\sc Braun, R.}, {\it A survey on Gelfand-Naimark theorems}, Schriftenreiche des Math.
Inst. der Univ. M\"{u}nster, Ser. 2, Univ. M\"{u}nster, Math. Inst.,
M\"{u}nster, 1987.

\item{[133]} {\sc Braun, R., Kaup, W., Upmeier, H.}, a) {\it On the automorphisms
of circular and Reinhardt domains in complex Banach spaces}, Manuscripta Math. {\bf 25}
(1978), 97--135.

b) {\it A holomorphic characterization of Jordan
$C^{\ast}$-algebras}, Math. Z. {\bf 161} (1978), 277--290.

\item{[134]} {\sc Br\^ anzei, D.}, a) {\it Structures of generalized manifold type} (in Romanian), Ph. D. Thesis, University "Al. I. Cuza", Ia\c si, 1975.

b) {\it Structures affines et op\' erations ternaires}, An. \c Stiin\c t. Univ. "Al. I. Cuza", Ia\c si, Sect. I Mat. {\bf 23} (1977), {\it 1}, 33-38.

\item{[135]} {\sc Bremner, M.R.}, {\it Jordan algebras arising from intermolecular recombination}, Jordan Theory preprints No. 197 (10 Aug. 2005), ACM SIGSAM Bull. {\bf 39} (2005), {\it 4}, 106-117.

\item{[136]} {\sc Bremner, M.R., Peresi, L.A.}, {\it Nonhomogeneous subalgebras of Lie and special Jordan superalgebras}, J. Algebra {\bf 322} (2009), {\it 6}, 2000-2026.

\item{[137]} {\sc Bre\v sar, M., Cabrera, M., Fo\v sner, M., Villena, A.R.}, 
{\it Lie triple ideals and Lie triple epimorphisms on Jordan and Jordan-Banach algebras},
Studia Math. {\bf 169} (2005), {\it 3}, 207--228.

\item{[138]} {\sc Brezhnev, Y.V.}, {\it Darboux transformation and some multi-phase solutions
of the Dodd-Bullough-Tzitzeica equation:} $U_{xt}={\rm e}^U-{\rm e}^{-2U}$, Phys. Lett. A {\bf 211}
(1996), 94--100.

\item{[139]} {\sc Brezuleanu, A., R\u adulescu, D.C.},
a) {\it About full or injective lineations},
J. Geom. {\bf 23} (1984), 45--60.

b) {\it Characterizing lineations defined on open subsets of projective spaces over
ordered division rings}, Abh. Math. Sem. Univ. Hamburg {\bf 55} (1985), 171--181.

\item{[140]} {\sc Brody, D.C., Hughston, L.P.}, {\it Geometric quantum mechanics},
J. Geom. Phys. {\bf 38} (2001), 19--53.

\item{[141]} {\sc Bro\v{z}ikova, E.}, {\it Homomorphisms of Jordan algebras and homomorphisms
  of projective planes}, Czechoslovak Math. J. {\bf 36} (1986), {\it 1}, 48--58.

\item{[142]} {\sc Bueken, P.}, a) {\it Reflections and rotations in contact geometry}, Doctoral
Dissertation, Catholic University Leuven, 1992.

b) {\it Curvature homogeneous pseudo-Riemannian manifolds}, invited talk at
the International Workshop Diff. Geom. and its Appl. (Bra\c sov, Romania, 1999).

\item{[143]} {\sc Buliga, M.}, {\it Braided spaces with dilations and sub-Riemannian symmetric spaces}, arXiv:1005.5031v1[math.GR] 27 May 2010, in GEOMETRY, Exploratory Workshop on Differential Geometry and its Applications (Ia\c si, 2009), Cluj University Press, D. Andrica \& S. Moroianu (Eds.), 2011, 21-35.

\item{[144]} {\sc Bunce, L.J., Chu, C.-H.}, a) {\it Dual spaces of $J^{\ast}$-triple and the
Radon-Nikodym property}, Math. Z. {\bf 208} (1991), 327--334.

b) {\it Compact operations, multipliers and Radon-Nikodym property in $JB^{\ast}$-triples}, Pacific J. Math.
{\bf 153} (1992), {\it 2}, 249--265.

\item{[145]} {\sc Bunce, L.J., Chu, C.-H., Stach\'o, L., Zalar, B.}, {\it On prime $JB^{\ast}$-triples}, Quart. J. Math.
Oxford {\bf 49} (1998), 279--290.

\item{[146]} {\sc Bunce, L.J., Wright, J.D.M.}, a) {\it Quantum measures and
states on Jordan algebras}, Comm. Math. Phys. {\bf 98} (1985), {\it 2},
187--202; for a more recent contribution concerning states on Jordan algebras, see Bunce, L.J.,
Hamhalter, J., {\it Extension of Janch-Piron states on Jordan algebras},
Math. Proc. Cambridge Philos. Soc. {\bf 119} (1996), {\it 2}, 279--286.

b) {\it Quantum logics and convex geometry}, Comm. Math. Phys. {\bf 101}
(1985), {\it 1}, 87--96.

\item{[147]} {\sc Burban, I., Henrich, T.}, {\it Semi-stable vector bundles on elliptic curves and the associative Yang-Baxter equation}, arXiv:1011.4591[math.AG] 20 Nov. 2010.

\item{[148]} {\sc Burban, I., Kreussler, B.}, {\it Vector bundles on degenerations of elliptic curves and Yang-Baxter equations}, arXiv:0708.1685.

\item{[149]} {\sc Burdujan, I.}, a) {\it On homogeneous systems associated to a binary algebra}, An. St. Univ. "Ovidius" Constan\c ta {\bf 2} (1994), 31-38

b) {\it Embeddings of nonassociative algebras}, Italian J. of Pure and Appl. Math. {\bf 7} (2000), 85-94.

c) {\it Some new properties of Lie triple systems}, ROMAI J. {\bf 2} (2006), {\it 1}, 7-23.

\item{[150]} {\sc Calder\'on Martin, A.J., Mart\'\i n Gonzalez, C.}, {\it A structure
theory for Jordan $H^*$-pairs}, Boll. Unione Math. Ital. Sez. B Artic. Ric. Mat. (8)
{\bf 7} (2004), {\it 1}, 61--77.

\item{[151]} {\sc Campos, T.M.M., Holgate, P.}, {\it Algebraic isotopy in genetics}, IMA J. Math. Appl. Med. Biol. {\bf 4} (1987), {\it 3}, 215-222.

\item{[152]} {\sc Cantarini, N., Kac, V.G.}, {\it Classification of linearly compact simple Jordan and generalized Poisson superalgebras}, J. Algebra {\bf 313} (2007), {\it 1}, 100-124.

\item{[153]} {\sc Capelli, A.}, {\it \"Uber die Zur\"uckf\"uhrung der Cayley'schen Operation
$\Omega$ auf gew\"oh\-niche Polar-Operationen}, Math. Ann. {\bf 29} (1887), 331--338.

\item{[154]} {\sc Cartan, E.}, a) {\it  Sur les domaines born\'{e}s homog\`{e}nes de l'espace
de $n$ variables complexes}, Abh. Math. Sem. Univ. Hamburg
{\bf 11} (1935), 116--162.

b) {\it Sur des familles remarquables d'hypersurfaces isoparametriques dans
les espaces sph\'{e}riques}, Math. Z. {\bf 45} (1939), 335--367.

c) {\it Sur quelques familles remarquables d'hypersurfaces}, C.R. Congr\`{e}s Math.
Li\`{e}ge (1939), 30--41.

\item{[155]} {\sc Carter, D.S., Vogt, A.}, {\it Collinearity-preserving functions between
Desarguesian planes}, Mem. Amer. Math. Soc. {\bf 27} (1980), {\it 235}.

\item{[156]} {\sc Castellon, S.A., Cuenca Mira, J.A., Mart\'\i n Gonzalez, C.}, {\it
Special Jordan $H^{\ast}$-triple systems}, Comm. Algebra {\bf 20} (2000), {\it 10}, 4699--4706.

\item{[157]} {\sc Cecotti, S., Ferrara, S., Girardello, L., Porrati, M.},
{\it Super-K\"{a}hler geometry in supergravity and superstrings}, Phys.
Lett. B {\bf 185} (1987), {\it 3-4}, 345--350.

\item{[158]} {\sc Cederwall, M.}, {\it Jordan algebra dynamics}, Phys. Lett.
B {\bf 210} (1988), 169--172.

\item{[159]} {\sc Certain, J.}, {\it The ternary operation $(abc) = ab^{-1}c$ of a group}, Bull. Math. Soc. {\bf 49} (1943), 869-877.

\item{[160]} {\sc Chapline, G., G\"{u}naydin, M.}, UCRL preprint
95290 (August 1986), unpublished.

\item{[161]} {\sc Chaput, P.-E.}, a) {\it Scorza varieties and Jordan algebras},
Indag. Math. {\bf 14} (2003), 169--182.

b) {\it Geometry over composition algebras: Projective geometry}, J. Algebra {\bf 298} (2006), {\it 2}, 340-362.

\item{[162]} {\sc Chein, O., Goodaire, E.G.}, {\it Minimal nonassociative nilpotent Moufang loops}, J. Algebra {\bf 268} (2003), {\it 1}, 327-342.

\item{[163]} {\sc Chi, Q.S.}, a) {\it A curvature characterization of certain
locally rank-one symmetric spaces}, J. Diff. Geom. {\bf 28} (1988), 187--202.

b) {\it Curvature characterization and classification of rank-one symmetric spaces},
Pacific J. Math. {\bf 150} (1991), 31--41.

\item{[164]} {\sc Chin, C.H., Chen, H.H.}, {\it On spinor representation of the infinite-dimensional
grassmannian}, Nucl. Phys. B. (Proc. Suppl.) {\bf 6} (1989), 419--421.

\item{[165]} {\sc Chow, W.L.}, {\it On the geometry of algebraic homogeneous spaces}, Ann. Math. {\bf 50} (1949), 32-67.

\item{[166]} {\sc Chu, C.-H.}, a) {\it The Radon-Nikodym property in operator
algebras}, in {\it Operator algebras
and group representations}, I, Proc. Internat. Conf. (Neptun, 1980), Pitman, 1984, 65--70.

b) {\it On the Radon-Nikodym property in Jordan algebras},
Glasgow Math. J. {\bf 24} (1983),
{\it 2}, 185--189.

c) {\it Jordan structures in Banach manifolds}, AMS/IP Studies in
Advances Mathematics {\bf 20} (2001), 201--210.

d) {\it Matrix-valued harmonic functions on groups}, J. Reine Angew. Math. {\bf 552} (2002), 15--52.

e) {\it Jordan Banach algebras in harmonic analysis}, Contemp. Math. {\bf 363} (2004), 59--68. 

f) {\it Grassmann manifolds of Jordan algebras}, Arch. Math. (Basel) {\bf 87} (2006), 179--192.

g) {\it Jordan algebras and Riemannian symmetric spaces}, in {\it Modern Trends in Geometry
and Topology} (Proc. Internat. Workshop Diff. Geom. Appl., Deva, 2005),
Cluj University Press, 2006, 133--152.

h) {\it Jordan triples and Riemannian symmetric spaces}, Adv. in Math. {\bf 219} (2008), 2029-2057.

i) {\it Jordan structures in geometry and analysis}, Cambridge University Press, Series: Cambridge Tracts in Mathematics (No. 190), ISBN: 9781107016170, Not yet published - available from January 2012.

\item{[167]} {\sc Chu, C.-H., Iochum, B.}, a) {\it On the Radon-Nikodym property in Jordan triples},
Proc. Amer. Math. Soc. {\bf 99} (1987), {\it 3}, 462--464.

b) {\it Weakly compact operators on Jordan triples}, Math. Ann. {\bf 281} (1988), {\it 3}, 451--458.

\item{[168]} {\sc Chu, C.-H., Isidro, J.M.}, {\it Manifolds of tripotents in $JB^{\ast}$-triples},
Math. Z. {\bf 233} (2000), 741--754.

\item{[169]} {\sc Chu, C.-H., Lau, A.T.-M.}, a) {\it Harmonic functions on groups and Fourier algebras},
Lecture Notes in Math. {\bf 1782}, Springer-Verlag, Berlin, 2002.

b) {\it Jordan structures in harmonic functions and Fourier algebras on homogeneous spaces},
Math. Ann. {\bf 336} (2006), 803--840.

\item{[170]} {\sc Chu, C.-H., Mackey, M.}, {\it Isometries between $JB^*$-triples},
Math. Z. {\bf 251} (2005), 615--633.

\item{[171]} {\sc Chu, C.-H., Mellon, P.}, {\it Jordan structures in Banach spaces and symmetric manifolds},
Expos. Math. {\bf 16} (1998), 157--180.

\item{[172]} {\sc Chu, C.-H., Wong, N.-C.}, {\it Isometries between $C^*$-algebras},
Rev. Mat. Iberoamericana {\bf 20} (2004), 87--105.

\item{[173]} {\sc Cirelli, R., Gatti, M., Mani\`a, A.}, {\it The pure state 
space of quantum mechanics as Hermitian symmetric space}, J. Geom. Phys.
{\bf 45} (2003), 267--284.

\item{[174]} {\sc Clerc, J.-L.}, a) {\it Zeta distributions associated to a representation
of a Jordan algebra}, Math. Z. {\bf 239} (2002), {\it 2}, 263--276.

b) {\it L'indice de Maslov g\'en\'eralis\'e}, J. Math. Pure Appl. {\bf 83} (2004), 
99--114.

c) {\it An invariant for triples in the Shilov boundary of a bounded symmetric 
domain}, Comm. Analysis Geom. {\bf 15} (2007), 147--173.

d) {\it Special prehomogeneous vector spaces associated to $F_4$, $E_6$, $E_7$ and $E_8$ and
simple Jordan algebras of rank $3$}, J. Algebra {\bf 264} (2003), 98--128.

\item{[175]} {\sc Clerc, J.-L., Koufany, K.}, {\it Primitive du cocycle de Maslov g\'en\'eralis\'e}, 
Math. Ann. {\bf 337} (2007), 91--138.

\item{[176]} {\sc Clerc, J.-L., Neeb, K.-H.}, {\it Orbits of triples in the Shilov boundary of a bounded
symmetric domain}, Transformation Groups {\bf 11} (2006), 387--426.

\item{[177]} {\sc Clerc, J.-L., \O rsted, B.}, a) {\it Maslov index revisited},
Transformation Groups {\bf 6} (2001), 303--320.

b) {\it The Gromov norm of the K\"ahler class and the Maslov index}, Asian J. Math.
{\bf 7} (2003), 269--296.

c) {\it Corrigendum to the Gromov norm  of the K\"ahler class and the Maslov index}, Asian J. Math.
{\bf 8} (2004), 391--394.

\item{[178]} {\sc Cleven, J.}, {\it Coordinatization of Jordan algebras over locally ringed spaces},
in
{\it Nonassociative algebra and its applications} (Oviedo, 1993), Santos Gonzalez (ed.),
Kluwer, Dordrecht, 1994, 99--105.

\item{[179]} {\sc Clifford, A.H.}, {\it A system arising from a weakened set of a group postulates}, Ann. Math. {\bf 34} (1953), 865-871.

\item{[180]} {\sc Climescu, A.}, {\it La repr\' esentation par des matrices groupales des groupoides multiplicatif d'une alg\`ebre non-associative}, Bul. Inst. Politehn. Ia\c si (N.S.) {\bf 2} (1956), {\it 3-4}, 9-18.

\item{[181]} {\sc Cochran, W.G.}, {\it The distribution of quadratic forms in a normal system with applications to analysis of covariance}, Proc. Cambridge Philos. Soc., {\bf 30} (1934), 178-191.

\item{[182]} {\sc Connes, A.}, in Publ. Math. IHES {\bf 62} (1985), 257.

\item{[183]} {\sc Conte, E.}, {\it Biquaternion Quantum Mechanics}, Pitagora Editrice, Bologna, Italy, 2000.

\item{[184]} {\sc Conte, R., Musette, M., Grundland, A.M.}, {\it B\"{a}cklund transformation of partial
differential equations from the Painlev\'{e}-Gambier classification II.
Tzi\-tzeica equation}, J. Math. Phys. {\bf 40} (1999), 2092--2106.

\item{[185]} {\sc Cort\'es, T.}, {\it Local algebras and primitivity of Jordan algebras}, in {\it Proceedings of the
Internat. Conference on Jordan structures} (Malaga, 1997), 71--75, Castellon Serrano, Cuenca Mira,
Fern\'andez L\'opez, Mart\'\i n Gonzalez (eds.), Malaga, 1999.

\item{[186]} {\sc Corrigan, E., Hollowood, T.J.}, a) {\it A string construction of a
commutative nonassociative algebra related to the exceptional Jordan algebra}, Phys. Lett.
B {\bf 203} (1988), {\it 1--2}, 47--51.

b) {\it The exceptional Jordan algebra and the superstring},
Comm. Math. Phys. {\bf 122} (1989), {\it 3}, 393--410.

\item{[187]} {\sc Cunha, I., Elduque, A.}, {\it The extended Freudenthal magic square and Jordan
algebras}, Manuscripta Math. {\bf 123} (2007), {\it 3}, 325--351.

\item{[188]} {\sc D'andrea, F., Dabrowski, L.}, {\it Dirac operators on quantum projective spaces}, Comm. Math. Phys. {\bf 295} (2010), 731-790.

\item{[189]} {\sc D'andrea, F., Dabrowski, L., Landi, G.}, {\it The noncommutative geometry of the quantum projective plane}, Rev. Math. Phys. {\bf 20} (2008), 979-1006.

\item{[190]} {\sc D'andrea, F., Landi, G.}, {\it Anti-selfdual connections on the quantum projective plane: monopoles}, Comm. Math. Phys. {\bf 297} (2010), 841-893.

\item{[191]} {\sc Darboux, G.}, {\it Le\c{c}ons sur la th\'{e}orie g\'{e}n\'{e}rale de
surfaces et les applications g\'{e}ometriques du calcul infinitesimal}, Deuxi\`{e}me
Partie, Gauthier-Villars, Paris, 1889.

\item{[192]} {\sc Daskaloyannis, C., Kanakoglu, K.}, {\it Hopf algebraic structure of the parabosonic and
parafermionic algebras and paraparticle generalization of the Jordan-Schwin\-ger map},
J. Math. Phys. {\bf 41} (2000), {\it 2}, 652--660.

\item{[193]} {\sc Date, E., Kashiwara, M., Jimbo, M., Miwa, T.}, {\it Transformation
groups for soliton equations, non-linear integrable systems-Classical theory
and quantum theory}, in {\it Proceedings}, M. Jimbo and T. Miwa (eds.), World Scientific,
1983, 39--119.

\item{[194]} {\sc D'Atri, J.E., Dorfmeister, J., Da, Z.Y.}, {\it The isotropy representation for
homogeneous Siegel domains}, Pacific J. Math. {\bf 120} (1985), {\it 2}, 295--326.

\item{[195]} {\sc Davies, M.}, {\it Universal $G$-manifolds}, Amer. J. Math. {\bf 103} (1981), {\it 1}, 103--141.

\item{[196]} {\sc David, D., Levi, D., Winternitz, P.}, {\it Equations invariant
under the symmetry group of the KP equation}, Phys. Lett. A {\bf 129}
(1988), {\it 3}, 161--164.

\item{[197]} {\sc Davydov, A.S.}, a) {\it The theory of contraction of protein under their excitation}, J. Theor. Biol. {\bf 38} (1973), 559-569.

b) {\it Biology and quantum mechanics}, Pergamon Press, New-York, 1982.

c) {\it Solitons in molecular systems} (in Russian), Kiyev, 1984; English translation: D. Reidel, Dordrecht, 1985.

d) {\it Excitons and solitons in quasi-one-dimensional molecular structures}, Ann. Physik {\bf 43} (1986), {\it 1-2}, 93-118.

\item{[198]} {\sc De Concini, C., Fucito, F., Tirozzi, B.}, {\it Conformal theories,
grassmannians and soliton equations}, I, Nucl. Phys. B, Part. Phys.
(Netherlands) {\bf 315} (1989), {\it 3}, 681--701.

\item{[199]} {\sc Della Selva, A., Saito, S.}, {\it A simple expression for the
Sciuto three-Reggeon vertex-generating duality}, Lett. Nuovo Cimento {\bf
4} (1970), {\it 15}, 689--692.

\item{[200]} {\sc De Medts, T.}, {\it An algebraic structure for Moufang quadrangles}, Mem. Amer. Math. Soc. {\bf 173} (2005), {\it 818}, 99 pp.

\item{[201]} {\sc De Medts, T., Segev, Y.}, {\it Identities in Moufang sets}, Trans. Amer. Math. Soc. {\bf 360} (2008), {\it 11}, 5831-5852.

\item{[202]} {\sc De Medts, T., Van Maldeghem, H.}, {\it Moufang sets of type $F_4$}, Math. Z. {\bf 265} (2010), {\it 3}, 511-527.

\item{[203]} {\sc De Medts, T., Weiss, R.M.}, {\it Moufang sets and Jordan division algebras},
Math. Ann. {\bf 335} (2006), {\it 2}, 415--433.

\item{[204]} {\sc Demidov, E.E.}, {\it Noncommutative deformation of the $K -P$ hierarchy
and the universal Grassmann manifold}, in {\it Lie groups and Lie algebras},
Math. Appl. {\bf 433}, Kluwer Acad. Publ., Dordrecht, 1998, 383--391.

\item{[205]} {\sc Dickey, L.A.}, {\it Additional symmetries of {\rm KP} grassmannian, and the string equation},
Modern Phys. Lett. A {\bf 8} (1993), {\it 8}, 1259--1272.

\item{[206]} {\sc Didry, M.}, {\it Structures alg\'ebriques sur les espaces sym\'etriques},
Thesis, Nancy 2006 (available at {\tt http://www.iecn.u-nancy.fr/$^{\sim}$didrym/}).

\item{[207]} {\sc Dimakis, A., M\" uller-Hoissen, F.}, a) {\it Weakly nonassociative algebras and the Kadomtsev-Petviashvili hierarchy}, Glasgow Math. J. {\bf 51A} (2009), 49-57.

b) {\it Quasi-symmetric functions and the KP hierarchy}, J. Pure Appl. Algebra, {\bf 214} (2010), 449-460.

\item{[208]} {\sc Dirac, P.A.M.}, Invited lecture at the Southeastern Section of the American Physical Society, Columbia, SC, November 1971.

\item{[209]} {\sc Di Scala, A.}, Private communication on April 2008.

\item{[210]} {\sc Di Scala, A., Loi, A.},  {\it Symplectic duality of symmetric spaces}, Advances in
Math. {\bf 217} (2008), 2336--2352.

\item{[211]} {\sc Di Scala, A., Loi, A., Roos, G.}, {\it The bisymplectomorphism group of~a~boun\-ded symmetric domain}, Transformation Groups {\bf 13} (2008), {\it 2}, 283-304, arXiv: 0707.2125.

\item{[212]} {\sc Di Scala, A., Loi, A., Zuddas, F.}, {\it Symplectic duality between complex
domains}, Monatsh. Math. {\bf 160} (2008), 403-428, arXiv: 0803.3536.

\item{[213]} {\sc Dineen, S., Timoney, R.M.}, {\it The centroid of a $JB^{\ast}$-triple system},
Math. Scand. {\bf 62} (1988), {\it 2}, 327--342.

\item{[214]} {\sc Dobrev, V.K.}, Invited plenary lecture at the 22nd Iranian
Mathematics Conference, March 13-16, 1991, Moshhad, Iran.

\item{[215]} {\sc Doebner, H.D., Henning, J.D., L\"{u}cke, W.}, in {\it Proc.
Quantum Groups Workshop}, Clausthal 1989, H.D. Doebner and G.D.
Hennig (eds.), Lecture Notes in Physics {\bf 370}, Springer-Verlag, Berlin, 1990.

\item{[216]} {\sc Dorfmeister, J.}, a) {\it Inductive
construction of homogeneous cones}, Trans.
Amer. Math. Soc. {\bf 252} (1973), 324--349.

b) {\it Homogene Siegel-Gebiete}, Habilitationsschrift, M\"{u}nster, 1978.

c) {\it Theta functions for the special, formally real Jordan algebras $($a remark on a paper of
H. L. Resnikoff$)$}, Inventiones Math. {\bf 44} (1978), {\it 2}, 103--108.

d) {\it Quasisymmetric Siegel domains and the automorphisms of homogeneous Siegel domains},
Amer. J. Math. {\bf 102} (1980), {\it 3}, 537--563.

e) {\it Homogeneous Siegel domains}, Nagoya Math. J. {\bf 86} (1982), 39--83.

f) {\it Simply transitive groups and K\"{a}hler structures on homogeneous Siegel domains},
Trans. Amer. Math. Soc. {\bf 288} (1985), {\it 1}, 293--306.

g) {\it Soliton equations and differential
geometry}, in {\it Proc. Internat. Workshop on Diff. Geom. and its Appl.} (Bucharest, 1993),
Bull. Sc. Politehnica Univ. Bucharest, Appl. Math. and Physics
{\bf 55} (1993), {\it 3--4}, 105--145.

h) {\it Groups associated with a Grassmannian modelled on the Wiener algebra},
Integral Equations Operator Theory {\bf 17} (1993), 464--500.

i) {\it Algebraic systems in differential geometry}, in {\it Jordan Algebras}
(Oberwolfach, 1992), W. Kaup, K. McCrimmon, H.P. Petersson (eds.), 9--33,
W. de Gruyter, Berlin--New York, 1994.

j) {\it Weighted $l_1$-Grassmannians and Banach manifolds of solutions to the {\rm KP}
equation and the KdV equation}, Math. Nachr. {\bf 180} (1996), 43--73.

\item{[217]} {\sc Dorfmeister, J., Neher, E.}, a) {\it Isoparameteric triple
systems of algebra type}, Osaka J. Math. {\bf 20} (1983), {\it 1}, 145--175.

b) {\it An algebraic approach to isoparametric hypersurfaces in spheres, I, II},
T\^{o}hoku Math. J. {\bf 35} (1983), 187--224; 225--247.

c) {\it Isoparametric triple systems of $FKM$-type, I, II, III}, Abh. Math. Sem.
Univ. Hamburg {\bf 53} (1983), 191--216; Manuscripta Math.
{\bf 43} (1983) {\it 1}, 13--44; Algebra Groups Geom. {\bf 1} (1984),
305--343.

d) {\it Isoparametric triple systems with special $Z$-structure}, Algebras Groups Geom.
{\bf 7} (1990), 21--94.

\item{[218]} {\sc Dorfmeister, J., Neher E., Szmigielski, J.}, a)
{\it Automorphisms of Banach manifolds associated with the {\rm KP} equation},
Quart. J. Math. Oxford (2) {\bf 40} (1989), 161--195.

b) {\it Automorphisms of the KdV subvariety}, Integral Equations Operator Theory
{\bf 14} (1991), 132--212.

c) {\it Banach manifolds and their automorphisms associated with groups of type $C_{\infty}$ and}
$D_{\infty}$, Contemporary Mathematics {\bf 110} (1990).

\item{[219]} {\sc Dr\' apal, A.}, {\it A simplified proof of Moufang's theorem}, Proc. Amer. Math. Soc. {\bf 139} (2011), {\it 1}, 93-98.

\item{[220]} {\sc Dray, T., Janesky, J., Manogue, C.A.}, {\it Octonionic Hermitian matrices
with non-real eigenvalues}, Adv. in Appl. Clifford Alg. {\bf 10} (2000), {\it 2}, 193--216.

\item{[221]} {\sc Dray, T., Manogue, C.A.}, {\it The exceptional Jordan eigenvalue problem},
Internat. J. Theoret. Phys. {\bf 38} (1999), {\it 11}, 2901--2916.

\item{[222]} {\sc Drinfeld, V.G.}, a) in {\it Proc. Int. Congress of Math.\
Berkeley, California}, {\bf 1}, Academic Press, New York, 1986, 798.

b) in Alg. Anal. {\bf 1} (1989).

\item{[223]} {\sc Dubois-Violette, M.}, {\it Lectures on graded differential algebras and noncommutative geometry}, in {\it Noncommutative differential geometry and its applications}, Kluwer, 2001, 245-306.

\item{[224]} {\sc Dubrovin, B., Krichever, I.M.}, in Sov. Sci. Rev. {\bf 3}
(1982), 1.

\item{[225]} {\sc Dvorsky, A., Sahi, S.}, a) {\it Explicit Hilbert spaces for certain unipotent
representations}. II, Invent. Math. {\bf 138} (1999), 203--224.

b) {\it Explicit Hilbert spaces for certain unipotent
representations} III, J. Funct. Anal. {\bf 201} (2003), 430--456.

\item{[226]} {\sc Edwards, C.M.}, a) {\it On the centres of hereditary $JBW$-subalgebras of
$JBW$-algebras}, Math. Proc. Cambridge Philos. Soc. {\bf 85} (1979), {\it 2}, 317--324.

b) {\it Facial structure of a $JB$-algebra}, J. London Math. Soc. {\bf 19} (1979), {\it 2},\break 335--344.

c) {\it Multipliers of $JB$-algebras}, Math. Ann. {\bf 249} (1980), {\it 3}, 265--272.

\item{[227]} {\sc Edwards, C.M., Ferna\'andez-Polo, F.J., Hoskin, C.S., Peralta, A.M.}, {\it On the facial structure of the unit ball in a $JB^{\ast}$-triple}, J. Reine Angew. Math. {\bf 641} (2010), 123-144.

\item{[228]} {\sc Edwards, C.M., H\"ugli, R.V.}, {\it Decoherence in pre-symmetric spaces}, Revista Matematica
Complutense {\bf 21} (2008), {\it 1}, 219--249.

\item{[229]} {\sc Edwards, C.M., H\"ugli, R.V., R\"uttimann, G.T.}, 
{\it A geometric characte\-rization of structural projections on a $JB^*$-triple}, J. Funct. Anal.
{\bf 202} (2003), {\it 1}, 174--194.

\item{[230]} {\sc Edwards, C.M., R\"{u}ttimann, G.T.}, {\it Lattice of
tripotents in a $JBW^{\ast}$-triples}, Internat. J. Theoret. Phys. {\bf 34} (1995),
{\it 8}, 1347--1357.

\item{[231]} {\sc Effros, E.G., St\o rmer, E.}, a) {\it Jordan algebras of self-adjoint
operators}, Trans. Amer. Math. Soc. {\bf 127} (1967), 313--316.

b) {\it Positive projections and Jordan structure in oparator algebras}, Math. Scand. {\bf 45} (1979), 127--138.

\item{[232]} {\sc Elduque, A.}, {\it Jordan gradings on exceptional simple Lie algebras}, Proc. Amer. Math. Soc. {\bf 137} (2009), {\it 12}, 4007-4018.

\item{[233]} {\sc Elduque, A., Kamiya, N., Okubo, S.}, a) {\it Simple $(-1,-1)$ balanced Freudenthal Kantor triple systems}, Glasgow Math. J. {\bf 45} (2003), 353-372.

b) {\it $(-1,-1)$ Balanced Freudenthal Kantor triple systems and noncommutative Jordan algebras}, J. Algebra {\bf 294} (2005), 19-40.

\item{[234]} {\sc El Yacoubi, N.}, {\it Th\'{e}or\`{e}mes de structure d'alg\`{e}bres de
Jordan topologiques ou bornologiques},
Algebras Groups Geom. {\bf 16} (1999), {\it 2}, 245--267.

\item{[235]} {\sc El Yacoubi, N., Kemmoun, H.}, a) {\it $J$-diviseurs topologiques
de z\'{e}ro dans les alg\`{e}bres de Jordan
$p$-norm\'{e}es unitaires}, Period. Math. Hungar. {\bf 35} (1997), {\it 3}, 159--167.

b) $J$-{\it diviseurs topologiques de z\'{e}ro d'une alg\`ebre de Jordan localement
multiplicativemenet convexe}, Afrika Mat. (3) {\bf 9} (1998), 19--26.

\item{[236]} {\sc Emch, G.G.}, a) {\it Algebraic medhods in statistical
mechanics and quantum field theory}, John Wiley \& Sons
Inc., New York--London--Sydney--Toronto, 1972.

b) {\it Mathematical and conceptual foundations of $20$th-century physics},
North-Holland Mathematics Studies {\bf 100}, Amsterdam, 1984.

\item{[237]} {\sc Emch, G.G., King, W.P.C.}, {\it Faithful normal states on $JBW$-algebras},
in {\it Operator Algebras and Applications}, Proc. Symp. Pure Math. {\bf 38}, Part 2, Amer. Math.
Soc., 1982, 305--308.

\item{[238]} {\sc Esser, M.}, {\it Self-dual postulates for $n$-dimensional geometry}, Duke Math.
J. {\bf 18} (1951), 475--479.

\item{[239]} {\sc Etherington, I.M.H.}, a) {\it Genetic algebras}, Proc. Roy. Soc. Edinburgh {\bf 59} (1939), 242-258.

b) {\it Nonassociative algebra and the symbolism of genetics}, Proc. Roy. Soc. Edinburgh {\bf 61} (1941), 24-42.

c) {\it Duplication of linear algebra}, Proc. Edinburgh Math. Soc. {\bf 6} (1941), 222-230.

\item{[240]} {\sc Evans, J.M.}, {\it Supersymmetric Yang-Mills theories and
division algebras}, Nuclear Phys. B {\bf 298} (1988), {\it 1}, 92--108.

\item{[241]} {\sc Faddeev, L.D.}, Preprint ITP-SB-94-11 and HEP-TH 9404013
(1994).

\item{[242]} {\sc Faddeev, L.D., Reshetikin, N.Yu., Takhtajan, L.A.}, in Alg. Anal. {\bf 1} (1989),~178.

\item{[243]} {\sc Fairlie, D.B., Manogue, C.}, {\it Lorentz invariance and the
composite string}, Phys. Rev. D {\bf 34} (1986), {\it 6}, 1832--1839.

\item{[244]} {\sc Faraut, J.}, a) {\it S\'{e}ries de Taylor shperiques sur une alg\`{e}bre de Jordan}, in
{\it Repr\'{e}\-sen\-ta\-tions des groups et analyse complexe} (Luminy, 1986), 35--46, Journ\'{e}es SMF,
{\bf 24}, Univ. Poitiers, Poitiers, 1986.

b) {\it Fonctions sp\'{e}ciales sur une alg\`{e}bre de Jordan}, in {\it Jordan algebras}
(Oberwolfach, 1992), 35--54, Kaup, McCrimmon, Petersson (eds.), de Gruyter, Berlin 1994.

\item{[245]} {\sc Faraut, J., Kaneyuki, S., Kor\'anyi, A., Lu, Q.-K., Roos, G.},
{\it Analysis and geometry on complex homogeneous domains}, Progress in
Math. {\bf 185}, Birkh\"{a}user, Boston, MA, 2000.

\item{[246]} {\sc Faraut, J., Kor\'anyi, A.}, a) {\it Function spaces and reproducing kernels
on bounded symmetric domains}, J. Funct. Anal. {\bf 89} (1990), 64--89.

b) {\it Analysis on symmetric cones}, Clarendon Press, Oxford, 1994.

\item{[247]} {\sc Fastr\'{e}, J.}, {\it A Grassmannian version of the Darboux transformation},
Bull. Sci. Math. {\bf 123} (1999), 181--232.

\item{[248]} {\sc Faulkner, J.R.}, a) {\it Octonion planes defined by quadratic Jordan
algebras}, Mem. Amer. Math. Soc. {\bf 104} (1970).

b) {\it Groups with Steinberg relations and coordinatization of
polygonal geometries}, Mem. Amer. Math. Soc. {\bf 185} (1977).

c) {\it Tangent bundles for Barbilian planes}, Tagunsbericht Jordan-Algebren (1985), Oberwolfach.

d) {\it A geometric construction of Moufang planes}, Geom. Dedicata {\bf 29} (1989), {\it 2}, 133--140.

e) {\it Barbilian planes}, Geom. Dedicata {\bf 30} (1989), {\it 2}, 125--181.

f) {\it Current results on Barbilian planes}, in {\it Proc. Internat. Workshop Diff.
Geom. and its Appl.} (Bucharest, 1993), Scientific Bull. Politehnica Univ.
Bucharest, Ser. A {\bf 55} (1993), {\it 3-4}, 147--152.

g) {\it Structurable triples, Lie triples, and symmetric spaces}, Forum Math.
{\bf 6} (1994), 637--650.

h) {\it Geometry and algebraic structures}, in
{\it Jordan Algebras} (Oberwolfach, 1992), W. Kaup, K.~McCrimmon, H.P. Petersson (eds.),
Walter de Gruyter, Berlin,  1994, 55--59.

i) {\it Projective remoteness planes}, Geom. Dedicata {\bf 60} (1996), 237--275.

j) {\it Jordan pairs and Hopf algebras}, J. Algebra {\bf 232} (2000), {\it 1}, 152--196.

k) {\it Hopf duals, algebraic groups and Jordan pairs}, J. Algebra {\bf 249} (2004), {\it 1}, 91-120.

\item{[249]} {\sc Faulkner, J.R., Ferrar, J.C.}, a) {\it Exceptional Lie algebras and related
algebraic and geometric structures}, Bull. London Math. Soc. {\bf 9} (1977),
{\it 1}, 1--35.

b) {\it Anti-Jordan pairs}, Comm. Algebra {\bf 8} (1980), {\it 11}, 993.

c) {\it Homomorphisms of Moufang planes and alternative planes}, Geom. Dedicata
{\bf 14} (1983), {\it 3}, 215--223.

d) {\it Generalizing the Moufang plane}, in {\it Ring and Geometry}, R. Kaya
et al. (eds.) Reidel, Dordrecht, 1985, 235--288.

\item{[250]} {\sc Fauser, B.}, {\it Vertex functions and generalized normal-ordering by triple systems
in non-linear spinor field models}, Adv. in Appl. Clifford Alg. {\bf 10} (2000), {\it 2}, 173--192.

\item{[251]} {\sc Fava, F.}, a) {\it Variet\'{a} Riemanniane a connessione costante}, Celebrazioni Archi\-mede-Siracusa (1962).

b) {\it Ulteriori contributi allo studio delle variet\'{a} a connessione costante},
Atti. Acad. Sci. Torino {\bf 97} (1962--1963).

c) {\it Sulla metrica delle variet\'{a} Riemanniane a connessione costante di genere tre},
Rend. Sem. Mat. Univ. Pol. Torino {\bf 27} (1968).

d) {\it Sulle variet\`{a} pseudo-Riemanniane a connessione constante di genere uno e due}, in
{\it G. Titeica and D. Pompeiu Symposium} (Bucharest, 1973), 149--152, Bucharest, 1976.

\item{[252]} {\sc Fava, F., Garbiero, S.}, {\it The B\"{a}cklund problem for the equation}
$f(x,t)z_{xx}\! -\!z_{tt}\! =\!0$, Atti Accad. Sci. Torino {\bf 123} (1989), {\it 3--4}, 147--166.

\item{[253]} {\sc Fay, J.D.}, {\it Theta functions on Riemann surfaces},
Lecture Notes in Math. {\bf 352}, Springer-Verlag, Berlin, 1973.

\item{[254]} {\sc Faybusovich, L.}, a) {\it Linear systems in Jordan algebras
and primal-dual inte\-rior-point algorithms}, J. Comput. Appl. Math. {\bf 86} (1997), {\it 1}, 149--175.

b) {\it Euclidean Jordan algebras and interior-point algorithms},
Positivity {\bf 1} (1997), {\it 4}, 331--357.

c) {\it A Jordan-algebraic approach to potential-reduction algoritms}, Math. Z. {\bf 239} (2002), 117-129.

\item{[255]} {\sc Ferapontov, E.V., Schief, W.K.}, {\it Surfaces of Demoulin: differential
geometry, B\"{a}cklund transformation and integrability}, J. Geom. Phys. {\bf 30} (1999),
343--363.

\item{[256]} {\sc Fern\'andez L\'opez, A., Garc\'\i a Rus, E.},
{\it Towards a Goldie theory for Jordan pairs},
Manuscripta Math. {\bf 95} (1998), {\it 1}, 79--90.

\item{[257]} {\sc Fern\'andez L\'opez, A., Garc\'ia Rus, E., G\'omez Lozano, M.}, a) {\it Jordan socle of finitary Lie algebras}, J. Algebra {\bf 280} (2004), {\it 2}, 635-654.

b) {\it The Jordan algebra of a Lie algebra}, J. Algebra {\bf 308} (2007), 164-177.

\item{[258]} {\sc Fern\'andez L\'opez, A., G\'arcia Rus, E., G\'omez Lozano, M., Neher, E.}, {\it A construction of gradings of Lie algebras}, Int. Math. Res. Not. IMRN {\bf 2007}, {\it 16}, Art. ID rnm 051, 34 pp.

\item{[259]} {\sc Fern\'andez L\'opez, A., Garc\'\i a Rus, E., Montaner, F.},
{\it Goldie theory for Jordan algebras}, J. Algebra {\bf 248} (2002), 397--491.

\item{[260]} {\sc Fern\'andez L\'opez, A., Marhnine, H., Zarhouti, C.}, {\it Derivations on
Banach-Jordan pairs}, Quart. J. Math. {\bf 52} (2001), 269--283.

\item{[261]} {\sc Fern\'andez L\'opez, A., Toc\'on Barroso, M.I.},
{\it Strongly prime Jordan pairs with nonzero socle}, Manuscripta Math.
{\bf 111} (2003), {\it 3}, 321--340.

\item{[262]} {\sc Fern\'andez-Polo, F.J., Peralta, A.M.}, a) {\it Closed tripotents and weak
compactness in the dual space of a $JB^*$-triple}, J. London Math. Soc. (2) 
{\bf 74} (2006), {\it 1}, 75--92.

b) {\it On the facial structure of the unit ball in the dual space of a $JB^{\ast}$-triple}, Math. Ann. {\bf 348} (2010), {\it 4}, 1019-1032.

\item{[263]} {\sc Ferrar, J.C.}, a) {\it Homomorphisms of projective planes}, Tagungsbericht
33/1982, Jordan-Algebren, Oberwolfach.

b) {\it Homomorphisms of Moufang-Veldkamp planes}, Geom. Dedicata {\bf 46} (1993), 299--311.

\item{[264]} {\sc Ferrar, J.C., Veldkamp, F.D.}, {\it Neighbour-preserving homomorphisms
between projective ring planes},
Geom. Dedicata {\bf 18} (1985), 11--35.

\item{[265]} {\sc Ferreira, L.A., Gomes, J.F., Teotonio Sobrinho, P.,
Zimerman, A.H.}, a) {\it Symplectic bosons Fermi fields and super
Jordan algebras}, Phys. Lett. B {\bf 234} (1990), {\it 3}, 315--320.

b) {\it The Jordan structure of Lie and Kac-Moody algebras}, J. Phys. A
{\bf 25} (1992), {\it 19}, 5071--5088.

\item{[266]} {\sc Ferreira, L.A., Gomes, J.F., Zimerman, A.H.}, a) {\it Vertex
operator and Jordan fields}, Phys. Lett. B {\bf 214} (1988), {\it 3},
367--370.

b) {\it Jordan algebras in conformal field theories}, in {\it J.J. Giambiagi Festschrift},
125--143, World Scientific Publishing Co., River Edge, NJ, 1990.

\item{[267]} {\sc Ferus, D.}, {\it Symmetric submanifolds of Euclidian space},
Math. Ann. {\bf 247} (1980), 81--93.

\item{[268]} {\sc Ferus, D., Karcher, H., M\"{u}nzner, H.F.},
{\it Cliffordalgebren und neue isopa\-rame\-trische Hyperfl\"{a}chen}, Math. Z.
{\bf 177} (1981), {\it 4}, 479--502.

\item{[269]} {\sc Fleury, N., Rausch de Traubenberg, M., Yamaleev, R.}, in Adv. in Appl.
Clifford Alg. {\bf 3(1)} (1993), 7.

\item{[270]} {\sc Fonseca, M., Mexia, J.T., Zmy\'slony, R.}, {\it Binary operations on
Jordan algebras and orthogonal normal models},
Linear Algebra Appl. {\bf 417} (2006), {\it 1}, 75--86.

\item{[271]} {\sc Foot, R., Joshi, G.C.}, a) {\it String theories and Jordan
algebras}, Phys. Lett. B {\bf 199} (1987), {\it 2}, 203--208.

b) {\it Nonassociative formulation of bosonic strings}, Phys. Rev. D
{\bf 36} (1987), {\it 4}, 1169--1174.

c) {\it A natural framework for the minimal supersymmetric gauge
theories}, Lett. Math. Phys. {\bf 15} (1988), {\it 3}, 237--242.

d) {\it Modification of the superstring action and the exceptional
Jordan algebra}, Internat. J. Theoret. Phys. {\bf 28} (1989), {\it 3},
263--271.

e) {\it Space-time symmetries of superstring and Jordan algebras},
Internat. J. Theoret. Phys. {\bf 28} (1989), {\it 12}, 1449--1462.

f) {\it An application of the division algebras, Jordan algebras and
split composition algebras}, Internat. J. Modern. Phys. A {\bf 7}
(1992), {\it 18}, 4395--4413.

g) {\it Symmetries of certain physical theories and the Jordan
algebras}, Internat. J. Modern Phys. A {\bf 7} (1992), {\it 15},
3623--3637.

\item{[272]} {\sc Freudenthal, H.}, a) {\it Oktaven, Ausnahmengruppen und Oktaven-geometrie},
Math. Inst.\ Rijksuniversiteit Utrecht, Utrecht, 1951 (mimeographed);
new revised edition, 1960; republished in Geom. Dedicata {\bf 19} (1981), {\it 1}, 7--64.

b) {\it Lie groups in the foundations of geometry},
Adv. in Math. {\bf 1} (1964), 145--190.

\item{[273]} {\sc Friedman, Y.}, a) {\it Bounded symmetric domains and the
$JB^{\ast}$-triple structure in physics}, in {\it Jordan algebras} (Oberwolfach, 1992),
Kaup, McCrimmon,
Petersson (eds.) de Gruyter, Berlin, 1994, 61--82.

b) {\it Physical applications of homogeneous balls}, Progress in Mathematical
Physics, {\bf 40}, Birkh\"auser, Boston Inc., Boston, MA, 2005, 279~pp. (with the assistance
of Tzvi Scarr).

c) {\it Explicit solutions for relativistic acceleration and rotation},
Czechoslovak J. Phys. {\bf 55} (2005), {\it 11}, 1403--1408.

d) {\it Extended relativistic dynamics} (submitted).

\item{[274]} {\sc Friedman, Y., Gofman, Y.}, a) {\it Why does the geometric product simplify the
equations of physics?}, Internat. J. Theoret. Phys. {\bf 41} (2002), {\it 10}, 1841--1855.

b) {\it Relativistic linear spacetime transformations
based on symmetry}, Found. Phys. {\bf 32} (2002), {\it 11}, 1717--1736.

\item{[275]} {\sc Friedman, Y., Russo, B.}, {\it A new approach to  spinors and some representation of
the Lorentz group of them}, Found. Phys. {\bf 31} (2001), {\it 12}, 1733--1766.

\item{[276]} {\sc Fritsch, R., Prestel, A.}, {\it Bewertungsfortsetzungen und nichtinjektive
Kol\-lineationen}, Geom. Dedicata {\bf 13} (1982), 107--111.

\item{[277]} {\sc Fr\"ohlich, H.}, Collective Phenomena {\it 1} (1973), 101.

\item{[278]} {\sc Fucito, F.}, a) {\it A Grassmannian formulation of the $WZW$ and
minimal\break models}, Phys. Lett. B {\bf 222} (1989), {\it 3--4}, 406--410.

b) {\it Conformal theories, Grassmannians and soliton equations II}, Nucl. Phys. B, Part.
Phys. (Netherlands), {\bf 315} (1989), 702--716.

\item{[279]} {\sc Fukuma, M., Kawai, H., Nakayama, R.}, {\it Infinite-dimensional
Grassmannian structrure of two-dimensional quantum gravity},
Comm. Math. Phys. {\bf 143} (1992), 371--402.

\item{[280]} {\sc Gao, H.B.}, {\it String field theory and infinite Grassmannian},
Phys. Lett. B {\bf 206} (1988), {\it 3}, 433--436.

\item{[281]} {\sc Garc\'\i a Rus, E., G\'omez Lozano, M.}, {\it Center, centroid, extended
centroid, and quotients of Jordan systems}, Comm. Algebra  {\bf 34} (2006), {\it 12}, 4311--4326.

\item{[282]} {\sc Garc\'\i a Rus, E., Neher, E.}, a) {\it Tits-Kantor-Koecher superalgebras of Jordan
superpairs covered by grids}, Comm. Algebra {\bf 31} (2003), {\it 7}, 3335--3375.

b) {\it Semiprime, prime and simple Jordan superpairs covered by grids}, 
J. Algebra {\bf 273} (2004), {\it 1}, 1--32.

\item{[283]} {\sc Garcia-Muniz, M.A., Gonzalez-Santos, J.},
{\it Baric, Bernstein and Jordan algebras}, Comm. Algebra {\bf 26} (1998), {\it 3}, 913--930.

\item{[284]} {\sc Garcia-Rio, E., Kupeli, D.N., Vazquez-Abal, M.E., Vazquez-Lorenzo, R.}, {\it Affine Osserman
connections and their Riemann extensions}, Diff. Geom. Appl. {\bf 11} (1999), 145--153.

\item{[285]} {\sc Garner, L.E.}, {\it Fields and projective planes: a category equivalence},
Rocky Mountain J. Math. {\bf 2} (1972), 605--610.

\item{[286]} {\sc Gerold, A., Buchner, K.}, {\it Solitons and isometric immersions},
J. Math. Phys. {\bf 32} (1991), {\it 8}, 2056--2062.

\item{[287]} {\sc Gervais, J.L.}, in Comm. Math. Phys. {\bf 130} (1990), 257;
Phys. Lett. B {\bf 243} (1990), 85.

\item{[288]} {\sc Giambruno, A., Zaicev, M.}, {\it Codimension growth of special simple Jordan algebras}, Trans. Amer. Math. Soc. {\bf 362} (2010), 3107-3124.

\item{[289]} {\sc Gilbert, G.}, {\it The {\rm KP} equations and fundamental string
theory}, Comm. Math. Phys. {\bf 117} (1988), {\it 2}, 331--348.

\item{[290]} {\sc Gilkey, P., Swann, A., Vanhecke, L.}, {\it Isoparametric geodesic
spheres and a conjecture of Osserman concerning
the Jacobi operator}, Quart. J. Math. Oxford {\bf 46} (1995), 299--320.

\item{[291]} {\sc Gindikin, S., Kaneyuki, S.}, {\it On the automorphism group of the generalized
conformal structure of a symmetric $R$-space}, Diff. Geom. Appl. {\bf 8} (1998), 21--33.

\item{[292]} {\sc Gnedbaye, A.V., Wambst, M.}, {\it Jordan triples and operads},
J. Algebra {\bf 231} (2000), {\it 2}, 744--757.

\item{[293]} {\sc Goddard, P., Nahm, W., Olive, D.I., Ruegg, H., Schwimmer,
A.}, {\it Fermions and octonions}, Comm. Math. Phys. {\bf 112} (1987),
{\it 3}, 385--408.

\item{[294]} {\sc Goddard, P., Olive, D., Schwimmer, A.}, {\it The heterotic string
and a fermionic construction of $E_8$ Kac-Moody algebra}, Phys.
Lett. B {\bf 157} (1985), 393.

\item{[295]} {\sc G\'omez-Ambrosi, C., Montaner, F.}, {\it On Herstein's constructions relating Jordan and associative superalgebras}, Comm. Algebra {\bf 28} (2000), {\it 8}, 3743-3762.

\item{[296]} {\sc Gomez, C., Sierra, G.}, in Phys. Lett. B {\bf 240} (1990), 149.

\item{[297]} {\sc G\'omez Gonz\'alez, E., Mun\~oz Porras, J.M., Plaza Martin, F.J.}, {\it Prym varieties, curves with automorphisms, and the Sato Grassmannian}, Math. Ann. {\bf 327} (2003), {\it 4}, 609-639.

\item{[298]} {\sc Gonshor, H.}, a) {\it Special train algebras arising in genetics}, I, II, Proc. Edinburgh Math. Soc. {\bf 12} (1960), 41-53; {\bf 14} (1965), 333-338.

b) {\it Contributions to genetic algebras}, I, II Proc. Edinburgh Math. Soc. {\bf 17} (1971), 289-298; {\bf 18} (1973), 273-287.

c) {\it Multi-algebra duplication}, J. Math. Biol. {\bf 25} (1987), {\it 6}, 677-683.

d) {\it Derivations in genetic algebras}, Comm. Algebra {\bf 16} (1988), {\it 8}, 1525-1542.

\item{[299]} {\sc Gonz\'alez-D\'avila, J.C., Gonz\'alez-D\'avila, M.C., Vanhecke, L.},
{\it Classification of Killing-transversally  symmetric spaces}, Tsukuba J. Math. {\bf 20} (1996), 321--347.

\item{[300]} {\sc Gonzalez, S., Martinez, C.}, {\it On Bernstein algebras},
in {\it Non-Associative Algebra and its Applications} (Oviedo, July 1993), Ed.
S. Gonzalez, Kluwer Acad. Publ., Dordrecht 1994, 164--170.

\item{[301]} {\sc Goze, M., Remm, E.}, {\it Poisson algebras in terms of nonassociative algebras}, J. Algebra {\bf 320} (2008), {\it 1}, 294-317.

\item{[302]} {\sc Graham, P.J., Ledger, A.J.}, $s$-{\it Regular manifolds}, in
{\it Differential Geometry in honour of Kentaro Yano}, Tokyo, 1972, 134--144.

\item{[303]} {\sc Gray, A.}, {\it Tubes}, Addison-Wesley Publ. Co. Reading, 1990.

\item{[304]} {\sc Green, M.B., Schwarz, J.H.}, a) {\it Covariant description of
superstrings}, Phys. Lett. B {\bf 136} (1984), 367--370.

b) in Nuclear Phys. B {\bf 243} (1984), 285.

\item{[305]} {\sc Greferath, M. Schmidt, S.E.}, a) {\it A unified approach to projective
lattice geometries}, Geom. Dedicata {\bf 43} (1992), {\it 3}, 243--264.

b) {\it On Barbilian spaces in projective lattice geometries}, Geom. Dedicata {\bf 43} (1992), {\it 3},
337--349.

\item{[306]} {\sc Grgin, E., Petersen, A.}, {\it Algebraic implications of composability of physical systems},
Comm. Math. Phys. {\bf 50} (1976), 177--188.

\item{[307]} {\sc Grishkov, A.N.}, {\it On the genetic behavior of Bernstein algebras} (in Russian), Dokl. Akad. Nauk SSSR {\bf 294} (1987), {\it 1}, 27-30.

\item{[308]} {\sc Grishkov, A.N., Shestakov, I.P.}, {\it Speciality of Lie-Jordan algebras}, J. Algebra
{\bf 237} (2001), 621--636.

\item{[309]} {\sc Gross, D.J., Harvey, J.A., Martinec, E., Rohm, R.}, a) {\it Heterotic string},
Phys. Rev. Lett. {\bf 54} (1985), {\it 6}, 502--505.

b) in Nuclear Phys. B {\bf 256} (1985), 253.

\item{[310]} {\sc Grossman, N.}, {\it Hilbert manifolds without epiconjugate points}, Proc. Amer.
Math. Soc. {\bf 16} (1995), 1365--1371.

\item{[311]} {\sc Grundh\"{o}fer, T.}, a) {\it Automorphism groups of compact projective planes},
Geom. Dedicata {\bf 21} (1986), {\it 3}, 291--298.

b) {\it Compact disconnected Moufang planes are Desarguesian}, Arch.\ Math.\break (Basel) {\bf 49} (1987), {\it 2}, 124--126.

\item{[312]} {\sc Grundland, A.M., Levi, D.}, {\it On higher-order Riccati equations as
B\"{a}cklund transformations}, J. Phys. A: Math. Gen. {\bf 32} (1999), 3931--3937.

\item{[313]} {\sc Grundland, A.M., Martina, L., Rideau, G.}, {\it Partial differential equations
with differential constrants}, in {\it Advances in Mathematical Sciences-CRM's $25$ Years},
CRM Proc. Lecture Notes {\bf 11}, L. Vinet (ed.), Amer. Math. Soc., Providence, 1997, 135--154.

\item{[314]} {\sc Guadagnini, E.} et al., in Phys. Lett. B {\bf 235} (1990), 275.

\item{[315]} {\sc Gumm, H.P.}, {\it The Little Desarguesian Theorem for algebras in modular varieties},
Proc. Amer. Math. Soc. {\bf 80} (1980), {\it 3}, 393--398.

\item{[316]} {\sc G\"{u}naydin, M.}, a) {\it Exceptional supergravity
theories, Jordan algebras, and the magic square}, in {\it Proceedings}
 13$^{\rm th}$ Internat. Colloq. on Group Theoretical
Methods, W.W. Zachary (ed.), World Scientific, Singapore, 1984,
478--491.

b) {\it Vertex operator construction of non-associative algebras and
their affinizations}, J. Math. Phys. {\bf 30} (1989), {\it 4}, 933--942.

c) {\it The exceptional superspace and the quadratic Jordan
formulation of quantum mechanics}, in {\it Elementary Particles and the
Universe} (Essays in Honor of Murray Gell-Mann, Pasadena, CA, USA,
27--28 January 1989), Cambridge, UK, Cambridge Univ. Press, 1991, 99--119.

d) $N=2$ {\it superconformal algebras and Jordan triple systems},
Phys. Lett. B {\bf 255} (1991), {\it 1}, 46--50.

e) {\it On a one-parameter family of exotic superspaces in two
dimensions}, Modern Phys. Lett. A {\bf 6} (1991), {\it 35}, 3239--3250.

f) {\it Generalized conformal and superconformal group actions and
Jordan algebras}, Modern Phys. Lett. A {\bf 8} (1993), {\it 15},
1407--1416.

\item{[317]} {\sc G\"{u}naydin, M., Hyun, S.J.}, a) {\it Affine exceptional
Jordan algebra and vertex operators}, Phys. Lett. B {\bf 209} (1988),
{\it 4}, 498--502.

b) {\it Ternary algebraic construction of extended superconformal
algebras}, Modern Phys. Lett. A {\bf 6} (1991), {\it 19}, 1733--1743.

c) {\it Ternary algebraic approach to extended superconformal algebras},
Nuclear Phys. B {\bf 373} (1992), {\it 3}, 688--712.

\item{[318]} {\sc G\"{u}naydin, M., Piron, C., Ruegg, H.}, {\it Moufang plane
and octonionic quantum mechanics}, Comm. Math. Phys. {\bf 61} (1978),
{\it 1}, 69--85.

\item{[319]} {\sc G\"{u}naydin, M., Sierra, G., Townsend, P.K.}, a)
{\it Exceptional supergravity theories and the magic square}, Phys. Lett. B
{\bf 133} (1983), {\it 1--2}, 72--76.

b) In Nuclear Phys. B {\bf 242} (1984), {\it 1}, 244--267.

c) {\it Vanishing potentials in gauged $N=2$ supergravity: an
applicaton of Jordan algebras}, Phys. Lett. B {\bf 144} (1984),
{\it 1--2}, 41--45.

d) {\it Gauging the $d=5$ Maxwell/Einstein supergravity theories: more
on Jordan algebras}, Nuclear Phys. B {\bf 253} (1985), {\it 3--4},
573--608.

e) In {\it Proceedings 1985 Cambridge Workshop on supersymmetry and its
applications}, G.W. Gibbons, S.W. Hawking and P.K. Townsend (eds.),
1986, 367.

\item{[320]} {\sc G\"{u}rsey, F.}, a) in {\it Second Johns Hopkins Workshop on
``Current Problems in High Energy Particle Theory"}, edited by G. Domokos and
S. K\"{o}vosi-Domokos, Johns Hopkins Univ., Baltimore, MD, 1978, 3--15.

b) {\it Octonionic structures in particle physics}, in Lecture Notes in
Physics {\bf 94}, Springer-Verlag, Berlin--Heidelberg--New York,
1978, 508--521.

c) {\it Super Poincar\'{e} groups and division algebras}, Modern Phys. Lett.
A {\bf 2} (1988), {\it 12}, 967--976.

d) {\it Lattices generated by discrete Jordan algebras}, Modern Phys.
Lett. A {\bf 3} (1988), {\it 12}, 1155--1168.

e) {\it Discrete Jordan algebras and superstring symetries}, in {\it Knots,
topology and quantum field theories} (Florence, 1989), World
Sci. Publishing, River Edge, NJ, 1989, 579--606.

\item{[321]} {\sc G\"{u}rsey, F., Tze, C.-H.}, {\it On the role of division, Jordan
and related algebras in particle physics}, World Scientific Publishing Co., River Edge,
NJ, 1996.

\item{[322]} {\sc Guz, W.}, a) {\it On the lattice structure of quantum
logics}, Ann. Inst. H. Poincare {\bf 28} (1978), {\it 1}, 1--7.

b) {\it Conditional probability in quantum axiomatics}, Ann. Inst. H. Poincare Sect.
A {\bf 33} (1980), {\it 1}, 63--119.

\item{[323]} {\sc Haagerup, U., Hanche-Olsen, H.}, {\it Tomita-Takesaki theory for Jordan algebras}, J. Operator Theory {\bf 11} (1984), {\it 2}, 343--364.

\item{[324]} {\sc H\"{a}hl, H.}, a) {\it Sechzehndimensionale lokalkompakte Translationsebene mit Spin $(7)$ als Kollineationsgruppe}, Arch. Math. (Basel) {\bf 48} (1987), {\it 3}, 267--276.

b) {\it Eine Kennzeichung der Oktavenebene}, Nederl. Akad. Wetensch. Proc.
A {\bf 90}~=~Indag. Math. (1987), {\it 1}, 29--39.

\item{[325]} {\sc Hall, J.A.}, {\it Speciality of interconnected quadratic Jordan algebras},
Comm. Algebra {\bf 30} (2002), {\it 6}, 2593--2615.

\item{[326]} {\sc Hall, M.}, {\it The theory of groups},
The Macmillan Company, New York, 1959.

\item{[327]} {\sc Hamhalter, J.}, a) {\it Universal state space
embeddability of Jordan-Banach algebras}, Proc. Amer. Math. Soc.
{\bf 127} (1999), {\it 1}, 131--137.

b) {\it Pure states on Jordan algebras}, Math. Bohem.
{\bf 126} (2001), {\it 1}, 81--91.

\item{[328]} {\sc Han, H.K., Shin, H.J.}, {\it Grassmann manifold bosonization of QCD in two
dimensions}, J. Math. Phys. {\bf 38} (1997), 3988--3996.

\item{[329]} {\sc Hanche-Olsen, H.}, a) {\it A Tomita-Takesaki theory for
$JBW$-algebras}, Proc.
Sympos. Pure Math. {\bf 38} (1982), {\it 2}, 301--303.

b) {\it On the structure and tensor products of $JC$-algebras}, Canad. J. Math. {\bf 35}
(1983), {\it 6}, 1059--1074.

c) $JB$-{\it algebras with tensor product are $C^{\ast}$-algebras}, in Lecture Notes Math. {\bf 1132},
Springer, 1985, 223--229.

\item{[330]} {\sc Hanche-Olsen, H., St\o rmer, E.}, {\it Jordan operator algebras}, Pitman, Boston--London--Melbourne, 1984.

\item{[331]} {\sc Hanssens, G., van Maldeghem, H.}, {\it A note on
near-Barbilian planes}, Geom. Dedicata {\bf 29} (1989), {\it 2}, 233--235.

\item{[332]} {\sc Harris, L.A.}, a) {\it Bounded symmetric homogeneous
domains in infinite-dimen\-sional spaces}, in {\it Proceedings on
infinite-dimen\-sio\-nal holomorphy},
Lectures Notes in Math. {\bf 364}, 13--40, Springer-Verlag, Berlin, 1974.

b) {\it Operator Siegel domains}, Proc. Royal Soc. Edinburgh {\bf 79}
A (1977), 137--156.

c) {\it A generalization of $C^{\ast}$-algebras}, Proc. London Math.
Soc. (3) {\bf 42} (1981), {\it 2}, 331--361.

\item{[333]} {\sc Havel, V.}, {\it Ein Einbettungssatz f\"{u}r die
Homomorphismen von Moufang-Ebe\-nen}, Czechoslovak. Math. J.
{\bf 20} (1970), 340--347.

\item{[334]} {\sc Heckenberger, I. Sch\"{u}ler, A.},
{\it On FRT-Clifford algebras}, Adv. in Appl. Clifford Alg.
{\bf 10} (2000), {\it 2}, 267--296.

\item{[335]} {\sc Helgason, S.}, a) {\it Differential geometry and
symmetric spaces}, Academic Press, New York, 1962.

b) {\it Differential geometry, Lie groups and symmetric spaces}, Academic
Press, New York, 1980, a new edition in August 2001.

\item{[336]} {\sc Helmholtz, H. von}, {\it K\"urzeste Linien im Farbensystem}, Z. Psychol. Physiol. Sinnesorg {\bf 4} (1892), 108-122.

\item{[337]} {\sc Helwig, K.-H.}, a) {\it Zur Koecherschen Reduktionstheorie in Positivit\"{a}tsberei\-chen}, I, II, III,
Math. Z. {\bf 91} (1966), 152--168; 169--178; 355--362.

b) {\it Jordan-Algebren und symmetrische R\"{a}ume}, I, Math. Z.
{\bf 115} (1970), 315--349.

\item{[338]} {\sc Herstein, I.N.}, {\it On the Lie and Jordan rings of a simple associative ring}, Amer. J. Math. {\bf 77} (1955), 279-285.

\item{[339]} {\sc Herzer, A.}, {\it Affine Kettengeometrien \"{u}ber Jordanalgebren}, Geom. Dedicata
{\bf 59} (1996), {\it 2}, 181--195.

\item{[340]} {\sc Hessenberger, G.}, {\it Barnes idempotents and
capacity in Jordan pairs}, Math. Proc. Cambridge Philos. Soc.
{\bf 123} (1998), {\it 1}, 37--53.

\item{[341]} {\sc Hestenes, D.}, {\it Space-time algebra}, Gordon and Breach, New York, 1966.

\item{[342]} {\sc Heuser, A.}, {\it \"Uber den Funktionenk\"orper der Normfl\"ache einer zentral einfachen Algebra}, J. Reine Angew. Math. {\bf 301} (1978), 105-113.

\item{[343]} {\sc Hilgert, J., Neeb, K.H.}, {\it Vector valued
 Riesz distributions on Euclidean Jordan algebras},
J. Geom. Anal. {\bf 11} (2001), {\it 1}, 43--75.

\item{[344]} {\sc Hilgert, J., \'{O}lafsson, G.},
{\it Causal Symmetric Spaces-Geometry and Harmonic Analysis},
Perspectives in Mathematics, {\bf 18}, Academic Press, San Diego, 1996.

\item{[345]} {\sc Hirota, R.}, a) {\it Exact solution of the KdV equation
for multiple collisions of solitons}, Phys. Rev. Lett.
{\bf 27} (1971), {\it 18}, 1192--1194.

b) {\it Direct method of finding exact solutions of
nonlinear evolutions}, in {\it B\"{a}cklund transformations}, Lecture
Notes in Math. {\bf 515}, Springer-Verlag, Berlin, 1976, 40--68.

\item{[346]} {\sc Hirzebruch, U.}, {\it \"{U}ber Jordan-Algebren und kompakte Riemannsche
symmetri\-schen R\"{a}ume von Rang 1}, Math. Z. {\bf 90} (1965), 339--354.

\item{[347]} {\sc Ho, T., Martinez-Moreno, J., Peralta, A.M., Russo, B.},
{\it Derivations on real and complex $JB^*$-triples}, J. London Math. Soc. (2)
{\bf 65} (2002), {\it 1}, 85--102.

\item{[348]} {\sc Holgate, P.}, a) {\it Jordan algebras arising in population genetics}, Proc. Edinburgh Math. Soc. {\bf 15} (1967), 291-294; Corrigendum ibid. {\bf 17} (1970), 120.

b) {\it Genetic algebras satisfying Bernstein's stationarity principle}, J. London Math. Soc. (2) {\bf 9} (1975), 613-623.

c) {\it Free nonassociative principal train algebras}, Proc. Edinburgh Math. Soc. (2) {\bf 27} (1984), 313-319.

d) {\it The interpretation of derivations in genetics algebras}, Linear Algebra Appl. {\bf 85} (1987), 75-80.

e) {\it The entropic law in genetic algebra}, Rev. Roumaine Math. Pures Appl. {\bf 34} (1989), {\it 3}, 231-234.

\item{[349]} {\sc Hong, H.-K., Liu, C.S.}, a) {\it Prandtl-Reuss elastoplasticity: on-off switch and superposition formulae},
Internat. J. Solids Struct. {\bf 34} (1997), 4281--4304.

b) {\it Internal symmetry in the constitutive model
of perfect elastoplasticity}, Internat. J. Non-Linear Mech. {\bf 35} (2000), 447--466.

\item{[350]} {\sc Horn, G.}, a) {\it Klassifikation der $JBW^{\ast}$-Tripel vom Typ I, Dissertation}, Univ.
T\"{u}bingen, 1984.

b) {\it Coordinatization theorems for $JBW^{\ast}$-triples}, Quart. J. Math.
Oxford Ser. (2) {\bf 38} (1987), {\it 151}, 321--335.

c) {\it Classification of $JBW^{\ast}$-triples of type I}, Math. Z. {\bf 196} (1987), {\it 2}, 271--291.

d) {\it Characterization of the predual and ideal structure of a $JBW^{\ast}$-triple}, Math. Scand.
{\bf 61} (1987), {\it 1}, 117--133.

\item{[351]} {\sc Horn, G., Neher, E.}, {\it Classification of continuous
$JBW^{\ast}$-triples}, Trans. Amer. Math. Soc. {\bf 306} (1988),
{\it 2}, 553--578.

\item{[352]} {\sc Hua, L.-K.}, a) {\it Geometries of matrices I, II, III}, Trans. Amer. Math. Soc. {\bf 57} (1945), 441--481, 482-490, {\bf 61} (1947), 229-255.

b) {\it Geometry of symmetric matrices over any field of characteristic other than two}, Ann. Math. {\bf 50} (1949), 8-31.

c) {\it Harmonic analysis of functions of several complex
variables in the clasical domains}, Translations of Math.
Monographs {\bf 6}, Amer. Math. Soc., Providence, 1963.

\item{[353]} {\sc H\"{u}gli, R.V.}, a) {\it Structural projections on a $JBW^*$-triple
and GL-projections on its predual}, J. Korean Math. Soc. {\bf 41} (2004),
{\it 1}, 107--130. 

b) {\it Characterizations of tripotens in $JB^*$-triples}, Math. Scand. 
{\bf 99} (2006), 147--160.

c) {\it Collinear systems and normal contractive projections on 
$JBW^*$-triple}, Integr. Equ. Oper. Theory
{\bf 58} (2007), 315--339. 

d) {\it A commutation principle for contractive projections on 
$JBW^*$-triple}, Quart. J. Math. (submitted).

\item{[354]} {\sc H\"{u}gli, R.V., Mackey, M.}, {\it Transitivity of inner automorphisms in 
infinite-dimensional Cartan factors}, Math. Z. {\bf 262} (2009), {\it 1}, 125-141.

\item{[355]} {\sc Hulett, E.G., S\'{a}nchez, C.U.}, {\it An algebraic
characterization of $R$-spaces}, Geom. Dedicata {\bf 67} (1997), 349--365.

\item{[356]} {\sc Ibragimov, M.M., Tleumuratov, S.Zh.}, 
{\it Geometric properties of geometric tripotents in a neutral
SFS space} (Russian), Uzbek. Mat. Zh. {\bf 2004}, {\it 3}, 35--38.

\item{[357]} {\sc Ikai, H.}, {\it Jordan pairs of quadratic forms with values in invertible
modules}, Collect. Math. {\bf 58} (2007), {\it 1}, 85--100.

\item{[358]} {\sc Iochum, B.}, a) {\it C\^{o}nes autopolaires et
alg\`{e}bres de Jordan}, Th\`{e}se, Universit\'{e} de Provence,
1982; Lecture Notes in Math. {\bf 1049}, Springer-Verlag, Berlin, 1984.

b) {\it Non-associative $L^p$-spaces}, Pacific J. Math.
{\bf 122} (1986), {\it 2}, 417--434.

\item{[359]} {\sc Iochum, B., Kosaki, H.}, {\it Linear Radon-Nikodym
theorems for states on
$JBW$ and $W^{\ast}$-algebras}, Publ. Res. Inst. Math. Sci.
{\bf 21} (1985), {\it 6}, 1205--1222.

\item{[360]} {\sc Iochum, B., Loupias, G.}, {\it Banach-power-associative
algebras and Jordan-Banach algebras}, Ann. Inst. H. Poincar\'{e}, Sect. A.
 {\bf 43} (1985), {\it 2}, 211--226.

\item{[361]} {\sc Iochum, B., Loupias, G., Rodr\'\i guez, A.},
{\it Commutativity of $C^{\ast}$-algebras and associativity
of $JB^{\ast}$-algebras},
Math. Proc. Cambridge Philos. Soc. {\bf 106} (1989), 281--291.

\item{[362]} {\sc Ion, D.B.}, a) {\it On Jordan algebra of type A}
 (in Romanian), Stud. Cerc. Mat. {\bf 17} (1965), 301--310.

b) {\it On certain properties of the roots of Jordan algebras of
type A} (in Romanian), Stud.
Cerc. Mat. {\bf 18} (1966), 309--313.

c) {\it On some properties of the structure constants of Jordan
algebras of type A} (in Romanian), Stud. Cerc. Mat. {\bf 19} (1967),
1031--1037.

d) {\it Reproducing kernel Hilbert spaces and extremal problems
for scattering of particles with arbitrary spins}, Internat.
J. Theoret. Phys. {\bf 24}
(1985), {\it 12}, 1217--1231.

e) {\it Scaling and $s$-channel helicity conservation via optimal
state description of hadron-hadron scattering}, Internat.
J. Theoret. Phys. {\bf 25} (1986), {\it 12}, 1237--1279.

\item{[363]} {\sc Ion, D.B., Scutaru, H.}, {\it Reproducing kernel
Hilbert spaces and optimal state description of hadron-hadron
scattering}, Internat. J. Theoret.
Phys. {\bf 24} (1985), 355--365.

\item{[364]} {\sc Iord\u anescu, R.}, a) {\it L' \'{e}tude des
op\'{e}rateurs infinit\'{e}simaux
associ\'{e}s aux alg\`{e}bres de Jordan simples}, Ann. Univ.
Bucure\c sti {\bf 19} (1970), 83--90.

b) {\it G\'{e}om\'{e}trie diff\'{e}rentielle sur les formes
r\'{e}elles de Jordan de type} $A_I$,
Bull. Un. Mat. Ital. {\bf 4} (1970), 585--594.

c) {\it M\'{e}triques sur les formes r\'{e}elles de Jordan},
Rev. Roumaine Math. Pures Appl. {\bf 15}
(1970), 1437--1444.

d) {\it Sur les graduations sp\'{e}ciales r\'{e}elles simples},
Rev. Roumaine Math. Pures
Appl. {\bf 16} (1971), 691--700.

e) {\it Jordan algebras with applications} (mimeographed),
INCREST, Bucharest, 170 pp., 1979, ISSN 0250 3638.

f) {\it Jordan structures in geometry}, in {\it Proc. National
 Conference on Geometry and
Topology} (Piatra Neam\c t, Romania, 1983), 62--69, Ia\c si, 1984.

g) {\it Jordan structures with applications} (mimeographed),
Inst.\ of Math., Bucha\-rest, 495 pp., 1990,
ISSN 0250 3638.

h) {\it Jordan structures, Grassmann manifolds, and string theories},
in {\it Proc. National
Conf. Geometry and Topology} (Timi\c{s}oara, Romania, 1989),
101--110, Timi\c{s}oara, 1991.

i) {\it Strutture di Jordan e applicazioni},
 Quaderni dell' Universit\'{a} ``La Sapienza"
di Roma, 40 pp., 1991.

j) {\it On geometrical applications of Jordan algebras},
in {\it Proc. National Conference on Geometry
and Topology} (Bucharest, Romania, 1991), 153--164, Bucharest, 1991.

k) {\it Jordan structures -- a unifying framework for
Barbilian planes and Differential geometry}, in
{\it Proceedings Internat. Workshop Diff. Geom. and its Appl.}
(Bucharest, July 1993), 183--189 (1993).

l) {\it The geometrical Barbilian's work from a
 modern point of view}, Balkan J. of Geom. Appl. {\bf 1} (1996), {\it 1}, 31--36.

m) {\it Recent advances and new open problems in the Jordan
algebra approach to differential geometry},
Preprint Inst. Math. Bucharest, 8/1996; Revised version 36/1996, 32 pp.
(partially included in the monograph [364p]).

n) {\it Applications of Jordan structures to differential
geometry and to physics
$($Recent advances and new open problems$)$}, Preprint Inst.
Math. Bucharest 11/1998,
79 pp. (partially included in the monograph [364p]).

o) {\it Jordan algebras in ring geometries},
Preprint Inst. Math. Bucharest, 13/1999, 50 pp.
(partially included in the monograph [364p]).

p) {\it Jordan structures in geometry and physics},
Quaderni dell' Universit\`{a} ``La Sapienza" di Roma, 181 p., 2000.

r) {\it The Jordan structure method applied to
differential geometry}, ANSTI grant, Bucharest, 2000, 59 pp.
(partially included in the monograph [364u]).

s) {\it The Jordan algebra method applied to ring geometries},
ANSTI grant,\break Bucharest, 2001, 52 pp.
(partially included in the monograph [364u]).

t) {\it Open problems arised from a Vr\u{a}nceanu's topic},
Bull. Math. Soc. Sc. Math. Rouma\-nie {\bf 44 (92)}, (2001), {\it 1}, 25--41.

u) {\it Jordan structures in geometry and physics} ({\it with an Appendix on Jordan
structures in analysis}), Edit. Acad. Rom\^ane, Bucure\c sti, 2003, 201 pp., ISBN 973-27-0956-1.

v) {\it Dynamical systems and Jordan structures},
Internat. J. Pure Appl. Math. {\bf 35} (2007), {\it 1}, 127--146.

w) {\it Jordan structures in analysis, geometry and physics}, Edit. Acad. Rom\^{a}ne, Bucure\c sti, 2009, 234 pp., ISBN 978-973-27-1775-2.

\item{[365]} {\sc Iord\u anescu, R., Truini, P.}, {\it Quantum groups and
Jordan sturctures}, Sci. Bull. Univ. Politehnica, Bucharest,
Ser. A {\bf 57}--{\bf 58} (1995), {\it 1--4}, 43--60.

\item{[366]} {\sc Ishibashi, N., Matsuo, Y., Ooguri, H.}, {\it Soliton
equations and free fermions on Riemann surfaces},
 Modern Phys. Lett. A {\bf 2}
(1987), 119.

\item{[367]} {\sc Isidro, J.M., Kaup, W., Rodr\'\i guez Palacios, A.},
{\it On real forms of $JB^{\ast}$-triples},
Manuscripta Math. {\bf 86} (1995), 311--335.

\item{[368]} {\sc Isidro, J.M., Mackey, M.},
{\it The manifold of finite rank projections in
the algebra ${\cal L}(H)$ of bounded linear operators},
Expos. Math. {\bf 20} (2002), 97--116.

\item{[369]} {\sc Isidro, J.M., Stach\'o, L.L.}, a) {\it Holomorphic
automorphism groups in Banach spaces}, North Holland Math. Studies 105,
North Holland, New York-Amsterdam-Oxford, 1985.

b) {\it On the manifold of tripotents in $JB^*$-triples},
J. Math. Anal. Appl. {\bf 304} (2005), 147--157.

c) {\it On the manifold of complemented principal
inner ideals in $JB^*$-triples}, Quart. J. Math. (Oxford)
{\bf 57} (2006), 505--525.

\item{[370]} {\sc Isidro, J.-M., Vigu\'{e}, P.-V.},
{\it The group of biholomorphic automorphisms of symmetric Siegel
domains and its topology}, Ann. Scuola Norm. Sup.
Pisa {\bf 11} (1984), {\it 3}, 343--352.

\item{[371]} {\sc Jacobson, N.}, a) {\it Some groups of transformations defined by
Jordan algebras}, III, J. Reine Angew. Math. {\bf 207} (1961), 61--85.

b) {\it Structure and representations of
Jordan algebras}, Amer. Math. Soc. Colloq. Publ. {\bf 39} Amer.
Math. Soc., Providence, 1968.

c) {\it Structure theory of Jordan algebras}, Univ. of Arkansas,
Lecture Notes in Math. {\bf 5}, 1981.

d) {\it Some projective varieties defined by Jordan algebras}, J. Algebra {\bf 97} (1985), {\it 2}, 565-598.

\item{[372]} {\sc Janson, S., Peetre, J.}, {\it A new generalization of Hankel operators $($the case
of higher weights}),  Math. Nachr. {\bf 132} (1987), 313--328.

\item{[373]} {\sc Janssen, G.}, a) {\it Formal-reelle
Jordanalgebren unendlicher Dimension und\break verallgemeinerte
Positivit\"{a}tsbereiche}, J. Reine Angew. Math. {\bf 249} (1971),\break 143--200.

b) {\it Reelle Jordan-algebren mit endlicher Spur},
Manuscripta Math. {\bf 13} (1974), 237--273.

c) {\it Die Struktur endlicher schwach abgeschlossener
Jordanalgebren I, II}, Manu\-scrip\-ta Math. {\bf 16} (1975),
277--305, 307--332.

\item{[374]} {\sc Jimbo, M., Miwa, T.}, {\it Solitons and
infinite-dimensional Lie algebras}, Publ. RIMS, Kyoto Univ. {\bf 19} (1983), {\it 3}, 943--1001.

\item{[375]} {\sc Johnson, K.D., Kor\'anyi, A.}, in Ann. Math. (2)
{\bf 111} (1980), {\it 3}, 589--608.

\item{[376]} {\sc Jonas, H.}, a) {\it Sopra una classe di trasformazioni
asimtotiche, applicabili in particolare alle superficie la cui curvatura
\'{e} proporzionale alla quarta potenza della distanza del piano
tangente da un punto fisso}, Ann. Mat. (1921), 223--255.

b) {\it Die differentialgleichung der Affinsph\"{a}ren in
 einer neuen Gestalt}, Math.
Nachr. {\bf 10} (1953), {\it 5--6}, 331--352.

\item{[377]} {\sc Jones, V.F.R.}, {\it A polynomial invariant
for knots via von Neumann algebras},
Bull. Amer. Math. Soc. {\bf 12} (1985), 103--111.

\item{[378]} {\sc Jonker, P.}, {\it Restricted Lie algebras over a field of
characteristic 2}, Dissertation, Math. Inst. Rijksuniversiteit Utrecht,
Utrecht, 1968.

\item{[379]} {\sc Jordan, P.}, a) {\it \"{U}ber eine Klasse nichtassoziatiever hyperkomplexen
Algebren}, G\"{o}tt. Nachr. (1932), 569--575.

b) {\it \"{U}ber Verallgemeinerungsm\"{o}glichkeiten des Formalismus der
Quantenmecha\-nik}, G\"{o}tt. Nachr. (1933), 209--217.

c) {\it \"{U}ber die Multiplikation quantenmechanischen Gr\"{o}ssen},
Z. Phys. {\bf 80} (1933), 285--291.

d) {\it \"{U}ber eine nicht-Desarguessche Ebene projective
Geometrie}, Abh. Math. Sem. Univ. Hamburg {\bf 16} (1949), 74--76.

\item{[380]} {\sc Jordan, P., von Neumann, J., Wigner, E.},
{\it On an algebraic generalization of the quantum mechanical
formalism}, Ann. Math. {\bf 35} (1934), 29--64.

\item{[381]} {\sc Joyce, D.}, {\it A classifying invariant of knots; the knot quandle}, J. Pure Appl. Alg. {\bf 23} (1982), 37-65.

\item{[382]} {\sc Kaidi, A.M., Morales Campoy, A., Rodr\'\i gues Palacios, A.},
a) {\it A holomorphic characterization of $C^*$ and $JB^*$-algebras},
Manuscripta Math. {\bf 104} (2001), {\it 4}, 467--478.

b) {\it Geometrical properties of the product in a $C^*$-algebra},
Rocky Mountain J. Math. {\bf 31} (2001), {\it 1}, 197--213.

\item{[383]} {\sc Kamiya, N.}, a) {\it Examples of Peirce's decomposition associated with simple generalized Jordan triple systems of second order}, Bull. Polish Acad. Sciences, Math., {\bf 49} (2001), {\it 2}, 159-175.

b) {\it Examples of Peirce decomposition of generalized Jordan triple systems of second order - balanced cases}, Contemporary Math. {\bf 391} (2005), 157-165.

\item{[384]} {\sc Kamiya, N., Mondoc, D.}, {\it A new class of nonassociative algebras with involution}, Proc. Japan Acad. {\bf 84} (2008), Ser. A, {\it 5}, 68-72.

\item{[385]} {\sc Kamiya, N., Okubo, S.}, a) {\it On triple systems
and Yang-Baxter equations}, Internat. J. Math. Game Theory Algebra
{\bf 8} (1998), {\it 2--3}, 81--88.

b) {\it A construction of Jordan superalgebras from triple systems},
Internat. J. Math. Game Theory Algebra {\bf 8} (1998), {\it 2--3}, 89--92.

c) {\it A construction of simple Jordan superalgebra of $F$ type from a Jordan-Lie triple system}, Annali di Matematica {\bf 181} (2002), 339-348.

d) {\it Representations of $(\alpha, \beta, \gamma)$ triple systems}, Linear and Multilinear Algebra {\bf 58} (2010), {\it 5-6}, 617-643.

\item{[386]} {\sc Kaneyuki, S.}, {\it Determinantal varieties in Jordan triple systems and stratifications of products of symmetric $R$-spaces}, in {\it Theory of Lie groups and manifolds}, Sophia Kokyuroku in Math. {\bf 45} (2002), Reiko Miyaoka \& Hiroshi Tamaru (eds.), 1-12.

\item{[387]} {\sc  Kantor, I.L.}, a) {\it Classification of irreducible transitive
differential groups} (in Russian), Dokl. Akad. Nauk SSSR  {\bf 158} (1964), 1271--1274.

b) {\it Transitive differential groups and invariant connections on homogeneous
spaces} (in Russian), Trudy Sem. Vektor. Tenzor. Anal. {\bf 13} (1966), 310--398. 

c) {\it Nonlinear transformation groups defined by general norms of Jordan 
algebras} (in Russian), Dokl. Akad. Nauk SSSR  {\bf 172} (1967), {\it 4}, 511--514.

d) {\it Connection between Poisson brackets and Jordan
and Lie superalgebras}, in {\it Lie theory, differential equations, and representation
theory} (Montreal, 1989), Univ. Montreal, Montreal, PQ, 1990, 213--225.

e) {\it A generalization of the Jordan approach to symmetric Riemannian spaces}, in {\it The Monster and Lie
algebras} (Columbus, OH, 1996), 221--234,
Ohio State Univ. Math. Res. Inst. Publ., {\bf 7}, Walter de Gruyter, Berlin--New York, 1998.

\item{[388]} {\sc Kantor, I.L., Rowen, L.}, {\it The Peirce decomposition for generalized Jordan
triple systems of finite order}, J. Algebra {\bf 310} (2007), {\it 2}, 829--857.

\item{[389]} {\sc Kantor, I.L., Sirota, A.I., Solodovnikov, A.S.},
a) {\it A class of symmetric spaces with extended motion
group and a generalization of Poincar\'{e}'s model} (in Russian),
Dokl. Akad. Nauk SSSR {\bf 173} (1967), 511--514.

b) {\it Bisymmetric Riemannian spaces}, Izvestiya: Mathematics
{\bf 59} (1995), {\it 5},\break 963--970.

\item{[390]} {\sc Kasman, A.}, {\it Grassmannians, nonlinear wave
equations and generalized Schur functions}, in {\it Geometry and
topology in dynamics},
M. Barge and K. Kuperberg (eds.), Contemporary Math. {\bf 246},
163--174, Amer. Math. Soc., Providence, RI, 1999.

\item{[391]} {\sc Kauffman, L.}, in Internat. J. Mod. Phys. A {\bf 5} (1990), 93.

\item{[392]} {\sc Kaup, W.}, a) {\it \"{U}ber die Automorphismen
Grassmannscher Mannigfaltigkeiten unendlicher Dimension},
Math. Z. {\bf 144} (1975), 75--96.

b) {\it Algebraic characterization of symmetric complex Banach manifolds}, Math. Ann. {\bf 228} (1977), {\it 1}, 39--64.

c) {\it Jordan algebras and holomorphy}, in {Functional analysis,
holomorphy and approximation theory}, Lecture Notes in Math. {\bf 843}, 341--365, Springer-Verlag,
Berlin--Heidelberg--New York, 1981.

d) {\it \"{U}ber die Klassifikation der symmetrischen
hermiteschen Mannigfaltigkeiten unendlicher Dimension I, II,}
Math. Ann. {\bf 257} (1981),
463--486; {\bf 262} (1983), 57--75.

e) {\it A Riemann mapping theorem for bounded symmetric domains
in complex Banach spaces}, Math. Z. {\bf 183} (1983), {\it 4},
503--529.

f) {\it Hermitian Jordan triple systems and the automorphisms
of bounded symmetric domains}, in {\it Non-associative algebra
and its applications} (Oviedo, 1993), S. Gonzalez
(ed.), Math. Appl. 303, Kluwer Acad. Publ. Dordrecht, 1994, 204--214.

g) {\it On $JB^{\ast}$-triples defined by fibre bundles},
Manuscripta Math. {\bf 87} (1995), {\it 3}, 379--403.

h) {\it On spectral and singular values in $JB^*$-triples}, Proc.
Royal Irish Acad. {\bf 19A} (1996), {\it 1}, 95--103.

i) {\it On real Cartan factors}, Manuscripta Math. {\bf 92} (1997), 191--222.

j) {\it On Grassmannians associated with $JB^{\ast}$-triples},
Math. Z. {\bf 236} (2001),\break 567--584.

k) {\it Bounded symmetric domains and derived geometric structures},
Rend. Mat. Acc. Lincei {\bf 13} (2002), {\it 9}, 243--257.

l) {\it CR-manifolds and Jordan algebras}, in {\it Modern Trends in Geometry and Topo\-logy} 
(Proc. Internat. Workshop Diff. Geom. Appl., Deva, September 2005), Cluj University Press, 2006, 245--250.

\item{[393]} {\sc Kaup, W., Matsushima, Y., Ochiai, T.},
{\it On the automorphisms and
equivalences of generalized Siegel domains},
Amer. J. Math. {\bf 92} (1970), 475--497.

\item{[394]} {\sc Kaup, W., Sauter, J.}, {\it Boundary structure
of bounded symmetric domains}, Manuscripta Math. {\bf 101} (2000), 351--360.

\item{[395]} {\sc Kaup, W., Upmeier, H.}, {\it Jordan algebras and
symmetric Siegel domains in Banach spaces}, Math. Z. {\bf 157} (1977),
179--200.

\item{[396]} {\sc Kaup, W., Zaitsev, D.}, a) {\it On symmetric Cauchy-Riemann
manifolds}, Adv. in Math. {\bf 149} (2000), 145--181.

b) {\it On the CR-structure of compact group orbits associated with bounded symmetric domains}, 
Invent. Math. {\bf 153} (2003), 45--104.

c) {\it On local CR-transformations of Levi-degenerate group orbits in compact Hermitian symmetric spaces}, 
J. Eur. Math. Soc. {\bf 8} (2006), {\it 3}, 465--490.

\item{[397]} {\sc Kawamoto, N., Namikawa, Y., Tsuchiya, A., Yamada, Y.},
{\it Geometric realization of conformal field theory on Riemann
surfaces}, Comm. Math. Phys. {\bf 116} (1988), {\it 2}, 247--308.

\item{[398]} {\sc Kennedy, M.}, {\it Triangularization of Jordan algebra of Schatten operators}, Proc. Amer. Math. Soc. {\bf 136} (2008), {\it 7}, 2521-2527.

\item{[399]} {\sc Kibble, T.W.B.}, {\it Geometrization of quantum mechanics},
Comm. Math. Phys. {\bf 65} (1979), 189--201.

\item{[400]} {\sc King, D.}, {\it Quadratic Jordan superalgebras}, Comm. Algebra {\bf 29} (2001), {\it 1}, 375--401.

\item{[401]} {\sc King, W.P.C.}, a) {\it Dual structures in
$JBW$-algebras}, Ph.D. Thesis, Rochester, 1980.

b) {\it Semifinite traces on $JBW$-algebras}, Math. Proc.
Cambridge Philos. Soc. {\bf 93} (1983), 503--509.

\item{[402]} {\sc Klingenberg, W.}, a) {\it Projective Geometrien
mit Homomorphismus}, Math. Ann. {\bf 132} (1956), 180--200.

\item{[403]} {\sc Knebusch, M.}, {\it Generic splitting of quadratic forms, I, II}, Proc. London Math. Soc. (3) {\bf 33} (1976), 65-93; {\bf 34} (1977), 1-31.

b) {\it Riemannian geometry}, Walter de Gruyter, Berlin, 1982.

\item{[404]} {\sc Knizhnik, V.G., Zamolodchikov, A.B.}, in Nucl. Phys.
B{\bf 247} (1984), 83.

\item{[405]} {\sc Kn\"{u}ppel, F.}, a) {\it Regular homomorphismus
of generalized projective planes},
J. Geometry {\bf 29} (1987), {\it 2}, 170--181.

b) {\it Projective planes over rings}, Resultate Math.
{\bf 12} (1987), {\it 3--4}, 348--356.

\item{[406]} {\sc Kn\"{u}ppel, F., Salow, E.},
{\it Plane elliptic geometry over rings},
Pacific J. Math. {\bf 123} (1986), {\it 2}, 337--384.

\item{[407]} {\sc Kodama, Y.}, {\it A method for solving the
dispersionless KP equation and its exact solutions},
Phys. Lett. A {\bf 129}
(1988), {\it 4}, 223--226.

\item{[408]} {\sc Koecher, M.}, a) {\it Positivit\"{a}tsbereichen im}
${\bf R}^n$, Amer. J. Math. {\bf 79} (1957),\break 575--596.

b) {\it Analysis in reellen Jordan-Algebren}, Nachr. Akad.
 Wiss. G\"{o}ttingen, 1958, 67--74.

c) {\it Jordan algebras and their applications},
Lecture Notes of the University of Minneapolis, Minnesota, 1962;
{\it The Minnesota Notes on Jordan Algebras and Their Applications},
Lecture Notes in Math., 1710, Springer-Verlag, Berlin 1999.

d) {\it \"{U}ber eine Gruppe von rationalen Abbildungen},
Invent. Math. {\bf 3} (1967), 136--171.

e) {\it Imbedding of Jordan algebras into Lie algebras, I, II}, Amer.
J. Math. {\bf 89} (1967), 787--816; {\bf 90} (1968), 476--510.

f) {\it An elementary approach to bounded symmetric domains},
Lecture Notes, Rice Univ., Houston, Texas,
1969.

g) {\it Gruppen und Lie-Algebren von rationalen Funktionen},
Math. Z. {\bf 109} (1969), 349--392.

h) {\it Jordan algebras in differential geometry},
Actes Congr\`{e}s Internat. Math. {\bf 1} (Nice, 1970), 279--283.

i) {\it Die Riccati-sche Differenzialgleichung und
nichtassoziative Algebren}, Abh. Math. Sem. Hamburg {\bf 46} (1977),
129--141.

j) in Lectures at the Rhine-Westphalia Academy of Sciences
{\bf 307}, Westdeutscher Verlag, Opladen, 1982, 53.

\item{[409]} {\sc Kogan, I.L.}, preprint PUPT-1439 and HEP-TH 9401093 (1994).

\item{[410]} {\sc Kollross, A.}, {\it Exceptional $\mathbb{Z}_2 \times \mathbb{Z}_2$-symmetric spaces}, Pacific J. Math. {\bf 242} (2009), {\it 1}, 113-130.

\item{[411]} {\sc Kor\'{a}nyi, A., Malliavin, P.},
{\it Poisson formula and compound diffusion associated to
an overdetermined elliptic
system on the Siegel halfplane of rank two}, Acta Math.
{\bf 134} (1975), 185--209.

\item{[412]} {\sc Kor\'anyi, A., Wolf, J.}, {\it Realization of Hermitian symmetric spaces as generalized
half-planes}, Ann. Math. {\bf 81} (1965), 265--268.

\item{[413]} {\sc Koshlukov, P.}, {\it Speciality of Jordan pairs},
An. Ovidius Univ. {\bf 4} (1996), {\it 2}, 98--106.

\item{[414]} {\sc Kostant, B., Sahi, S.}, a) {\it The Capelli identity, tube domains,
and the generalized Laplace transform}, Adv. Math. {\bf 87} (1991), 71--92.

b) {\it Jordan algebras and Capelli identities},
Invent. Math. {\bf 112} (1993), 657--664.

\item{[415]} {\sc Koszul, J.-L.}, {\it Alg\`{e}bres de Jordan},
S\'{e}minaire Bourbaki, {\bf 1}, Exp. {\it 31}, 245--256, Soc.
Math. France, Paris, 1995.

\item{[416]} {\sc Koufany, K.}, a) {\it Jordan algebras, geometry of Hermitian symmetric
spaces and non-commutative Hardy spaces}, Seminar on Mathematical Sciences,
{\bf 33}, Keio University, Yokohama, 2005, 70~pp.

b) {\it Analyse et g\'eom\'etrie des domaines born\'es sym\'etriques}, Th\`ese d'Habilitation
de l'Universit\'e Nancy 1 (2006).

\item{[417]} {\sc Kovacs, A.}, {\it Generic splitting fields}, Comm. Algebra {\bf 6} (1978), 1017-1035.

\item{[418]} {\sc Kowalski, O.}, a) {\it Smooth and affine
$s$-manifolds}, Period. Math. Hungar. {\bf 8} (1977),
{\it 3--4}, 299--311.

b) {\it Generalized affine symmetric spaces}, Math. Nachr.
{\bf 80} (1977), 205--208.

c) {\it Generalized symmetric spaces}, Lecture Notes in Math. {\bf 805},
Springer-Verlag, Berlin--Heidelberg--New York, 1980.

\item{[419]} {\sc Kowalski, O., Pr\"{u}fer, F., Vanhecke, L.},
{\it D'Atri spaces}, in {\it Topics in Geometry: In Memory of Joseph D'Atri}
(Ed. S. Gindikin), {\it Progress in Nonlinear Differential Equations}, {\bf 20},
Birkh\"{a}user, Boston--Basel--Berlin, 1996, 241--284.

\item{[420]} {\sc Kowalski, O., Tricerri, F., Vanhecke, L.},
{\it Curvature homogeneous spaces with a solvable Lie group as homogeneous
model}, J. Math. Soc. Japan {\bf 44} (1992), {\it 3}, 461--484.

\item{[421]} {\sc Krutelevich, S.V.}, a) {\it The Tits-Kantor-Koecher construction and 
birepresentations of the Jordan superpair $SH(1,n)$}, Comm. Algebra
{\bf 32} (2004), {\it 6}, 2117--2148.

b) {\it Jordan algebras, exceptional groups, and Bhargava composition}, J. Algebra {\bf 314} (2007), {\it 2}, 924-977.

\item{[422]} {\sc K\"{u}hn, O.}, {\it Differentialgleichungen in Jordantripelsystemen},
Manuscripta Math. {\bf 17} (1975), {\it 4}, 363--381.

\item{[423]} {\sc Kulkarni, M.}, {\it Fundamental theorem of projective
geometry over a commutative ring}, Indian J. Pure Appl. Math.
{\bf 11} (1980), {\it 12}, 1561--1565.

\item{[424]} {\sc Kummer, H.}, {\it A constructive approach to the
formulation of quantum mecanics}, Found. Phys. {\bf 17} (1987),
{\it 1}, 1--62.

\item{[425]} {\sc Kundu, A.}, in {\it Application of solitons in Science
and Engeneering}, World Scientific, Singapore, 1994.

\item{[426]} {\sc Kuznetsova, T.A.}, {\it Duo-octavic and duoantioctavic planes}
(in Russian). Geometry, {\it 6}, 65--72, Leningrad, Gos. Ped. Inst., Leningrad, 1977.

\item{[427]} {\sc Landau, L., Lifshitz, E.M.},  
{\it On the theory of the dispersion of magnetic permeability in ferromagnetic
bodies}, Phys. Z. Sowjetunion {\bf 8} (1935), 153--169.

\item{[428]} {\sc Landi, G., Reina, C.}, {\it Symplectic dynamics on the
universal Grassmannian}, J. Geom. and Phys. {\bf 9} (1992), {\it 3}, 235--253.

\item{[429]} {\sc Landsberg, J., Manivel, L.}, {\it The projective geometry of
Freudenthal's magic square}, J. Algebra 
{\bf 239} (2001), 477--512.

\item{[430]} {\sc Landsman, N.P.}, {\it Quantization and classsicization:
from Jordan-Lie algebras of obervables to gauge fields}, Class.\ Quantum Gravity
(UK), {\bf 10} suppl. issue (1993), 101--108 (Journees Relativistes 1992
(Relativistic Days 1992), Amster\-dam, Netherlands, May 13--15, 1992).

\item{[431]} {\sc Lantz, D.C.}, {\it Uniqueness of Barbilian domains}, J. Geom.
{\bf 15} (1981), 21--27.

\item{[432]} {\sc Lassalle, M.}, a) {\it S\'{e}ries de Laurent des fonctions holomorphes dans la complexification d'un
espace symm\'{e}trique compact}, Ann. Sci. \'{E}cole Norm. Sup. (4) {\bf 11} (1978), {\it 2}, 167--180.

b) {\it Transform\'{e}es de Poisson, alg\`{e}bres de Jordan et
\'{e}quations de Hua}, C.R. Acad. Sci. Paris {\bf 294} (1982), 325--328.

c) {\it Alg\`{e}bres de Jordan, coordonn\'{e}es polaires et
equations de Hua}, C.R. Acad. Sci. Paris {\bf 294} (1982), 613--615.

d) {\it Syst\`{e}mes triple de Jordan, $R$-espaces symmetriques
et \'{e}quations de Hua}, C.R. Acad. Sci. Paris, S\'{e}r. I. Math.
 {\bf 298} (1984), {\it 20}, 501--504.

e) {\it Sur la valeur au bord du noyau de Poisson d'un domaine
born\'{e} symm\'{e}trique}, Math. Ann. {\bf 268} (1984), {\it 4},
417--423.

f) {\it Les \'{e}quations de Hua d'un domaine born\'{e}
symm\'{e}trique du type tube}, Invent. Math. {\bf 77} (1984), 129--161.

g) {\it Alg\`{e}bres de Jordan et \'{e}quations de Hua},
J. Funct. Anal. {\bf 65} (1986), {\it 2}, 243--272.

h) {\it Noyau de Szeg\"{o}, $K$-types et alg\`{e}bres de Jordan},
C.R. Acad. Sci. Paris {\bf 303} (1986), {\it 1}, 1--4.

i) {\it Alg\`{e}bres de Jordan et ensemble de Wallach},
Invent. Math. {\bf 89} (1987), {\it 2}, 375--393.

\item{[433]} {\sc Layne, S.P.}, {\it A possible mechanism for general anesthesia}, Los Alamos Science, Spring 1984, 23-26.

\item{[434]} {\sc Leasher, B.}, {\it Geometric aspects of Steinberg
groups defined
by Jordan pairs}, Comm. Algebra {\bf 23} (1995), 5355--5368.

\item{[435]} {\sc Ledger, J., Obata, M.}, {\it Affine and Riemannian
$s$-manifolds}, J. Differential Geom. {\bf 2} (1968), 451--459.

\item{[436]} {\sc Lee, S.Y., Lim, Y., Park, C.-Y.}, {\it Symmetric geodesics on conformal
compactifications of Euclidean Jordan algebras}, Bull. Austral. Math. Soc. {\bf 59}
(1999), {\it 2}, 187--201.

\item{[437]} {\sc Leissner, W.}, a) {\it Affine Barbilian-Ebenen, I, II},
J. Geom. {\bf 6} (1975), 31--57,
105--129.

b) {\it Parallelodronnic-Ebenen}, J. Geometry {\bf 8} (1976),
{\it 1--2}, 117--135.

c) {\it Barbilianbereiche} in {\it Beitr\"{a}ge zur Geometrischen Algebra},
H.J. Arnold, W.~Benz and H. Wefelscheid (eds.), 219--224,
Birkh\"{a}user, Basel, 1977.

d) {\it On classifying affine Barbilian spaces}, Resultate Math.
{\bf 12} (1987), 157--165.

e) {\it Rings of stable rank $2$ are Barbilian rings},
Results in Mathematics {\bf 20} (1991), {\it 1--2}, 530--537.

\item{[438]} {\sc Leissner, W., Severin R., Wolf, K.},
{\it Affine geometry over free unitary modules},
J. Geom. {\bf 25} (1985), {\it 2}, 101--120.

\item{[439]} {\sc Leur, J., van de}, {\it The vector $k$-constrained
{\rm KP} hierarchy and Sato's Grassmannian}, J.
Geom. and Phys. {\bf 23} (1997), 83--96.

\item{[440]} {\sc Levasseur, T., Stafford, J.T.}, {\it Invariant
differential operators on the tangent space of some symmetric spaces},
Ann. Inst. Fourier (Grenoble) {\bf 49} (1999), {\it 6}, 1711--1741.

\item{[441]} {\sc Levi, D., Winternitz, P.}, {\it The cylindrical
{\rm KP} equation; its KacMoody-Virasoro algebra and relation to
{\rm KP} equation}, Phys.\ Lett.\ A {\bf 129} (1988), {\it 3}, 165--167.

\item{[442]} {\sc Levin, J.J.}, {\it On the matrix Riccati equation},
Proc. Amer. Math. Soc. {\bf 10} (1959), {\it 4}, 519--524.

\item{[443]} {\sc Li, S.P., Peschanski, R., Savoy, C.A.},
{\it Generalized no-scale
models and classical symmetries of superstrings}, 
Phys. Letters {\bf B194} (1987), {\it 2}, 226--230.

\item{[444]} {\sc Lim, Y.}, {\it Applications of geometric means on 
symmetric cones}, Math. Ann. {\bf 319} (2001), {\it 3}, 457--468.

\item{[445]} {\sc Liu, C.-S.}, a) {\it Orientation control of particles and fibers for polymeric liquids and fiber reinforced materials},
J. Chinese Inst. Engineers {\bf 20} (1997), 443--456.

b) {\it A Jordan algebra and dynamical system with associator as vector field}, Internat. J. Non-Linear Mech. {\bf 35} (2000), {\it 3}, 421--429.

c) {\it The $g$-based Jordan algebra and Lie algebra with application to the model of visco-elastoplasticity},
J. Marine Sci. Tech. {\bf 9} (2001), 1--13.

d) {\it Applications of the Jordan and Lie algebras for some dynamical systems having internal symmetries},
Internat. J. Applied Math. {\bf 8} (2002), {\it 2}, 209--240.

e) {\it The $g$-based Jordan algebra and Lie algebra formulation of the Maxwell equations}, J. Mech.
{\bf 20} (2004), {\it 4}, 285--296.

f) {\it A new mathematical modeling of Maxwell equations: complex linear operator and complex field}, CMES, {\bf 38} (2003), {\it 1}, 25-38.

\item{[446]} {\sc L\"{o}hmus, J., Paal, E., Sorgsepp, L.}, {\it Nonassociativity in
mathematics and physics} (in Russian, English summary), in {\it Quasigroups and nonassociative
algebras in physics} (Tartu, 1989), L\"{o}hmus, J. (ed.), Eesti Teadaste Akad.
F\"{u}\"{u}ns. Inst. Uurim. {\bf 66} (1990), 8--22.

\item{[447]} {\sc Lomdhal, P.S., Layne, S.P., Bigio, I.J.}, {\it Solitons in biology}, Los Alamos Science, Spring 1984, 3-21.

\item{[448]} {\sc Loos, O.}, a) {\it Spiegelungsr\"{a}ume und homogene
symmetrische Mannigfaltigkeiten}, Dissertation, M\"{u}nchen, 1966.

b) {\it Speigelungsr\"{a}ume und homogene symmetrische R\"{a}ume}, Math. Z.
{\bf 99} (1967), 141--170.

c) {\it Symmetric spaces, I, II}, Benjamin, New York, 1969.

d) {\it Jordan triple systems, $R$-spaces and bounded symmetric domains},
Bull.
Amer. Math. Soc. {\bf 77} (1971), 558--561.

e) {\it Lectures on Jordan triples}, Lecture Notes Univ. of British
Columbia, Vancouver,
1971.

f) {\it A structure theory of Jordan pairs}, Bull. Amer. Math. Soc. USA
{\bf 80} (1974), 67--71.

g) {\it Jordan pairs}, Lecture Notes in Math., {\bf 460},
Springer-Verlag, Berlin--Heidel\-berg--New York, 1975.

h) {\it Bounded symmetric domains and Jordan pairs}, Lecture Notes,
Univ. of California at Irvine, 1977.

i) {\it On algebraic groups defined by Jordan pairs}, Nagoya Math. J.
{\bf 74} (1979), 23--66.

j) {\it Charakterisierung symmetrischer $R$-R\"{a}ume durch ihre Einheitsgitter},
Math. Z. {\bf 189} (1985), 211--226.

k) {\it Recent results on finiteness conditions in Jordan pairs},
in {\it Jordan algebras}
(Oberwolfach, 1992), 83--95, Kaup, McCrimmon, Petersson (eds.),
de Gruyter, Berlin, 1994.

l) {\it Generically algebraic Jordan algebras over commutative rings}, J. Algebra {\bf 297} (2006), {\it 2}, 474-529.

m) {\it Jordan pairs and bounded symmetric domains}, Jordan Theory preprints No. 286 (24 Feb. 2010).

\item{[449]} {\sc Loos, O., McCrimmon, K.}, {\it Speciality of Jordan triple systems},
Comm. Algebra {\bf 5} (1977), 1057--1082.

\item{[450]} {\sc Lorimer, J.W.}, {\it Dual numbers and topological Hjelmslev planes}, Canadian Math. Bull.
{\bf 26} (1983), {\it 3}, 297--302.

\item{[451]} {\sc Lutz, R.}, {\it Sur la g\'eometrie des espaces $\Gamma$-sym\'etriques}, C.R. Acad. Sci. Paris, S\'er. I Math. {\bf 293} (1981), {\it 1}, 55-58.

\item{[452]} {\sc Lyubich, Yu.I.}, a) {\it Basic concepts and theorems of evolution genetics of free populations} (in Russian), Uspekhi Mat. Nauk (5) {\bf 26} (1971), 51-116; English translation: Russian Math. Surveys (5) {\bf 26} (1971), 51-123.

b) {\it Bernstein algebras} (in Russian), Uspekhi Mat. Nauk (6) {\bf 32} (1977), 261-262.

c) {\it Classification of nonexceptional Bernstein algebras of type $(3,n-3)$} (in Russian), Vestnik Kharkov Gos. Univ. No. 254, Mekh. Mat. Protsessy Upravl. (1984), 36-42.

d) {\it Mathematical structures in population genetics}, Biomathematics {\bf 22}, Springer-Verlag, Berlin, Heidelberg, 1992.

\item{[453]} {\sc MacDonald, M.L.}, {\it Cohomological invariants of Jordan algebras with frames}, J. Algebra {\bf 323} (2010), {\it 6}, 1665-1677.

\item{[454]} {\sc Mack, G., Schomerous, V.}, a) in Comm. Math. Phys. {\bf 149}
(1992), 513.

b) Preprint HVTMP 94-B335 (1994).

\item{[455]} {\sc Mackey, M.T.}, {\it The Grassmannian manifold associated
to a bounded symmetric domain}, in {\it Finite
or infinite-dimensional complex analysis}, Kajiwara J., Li Zhong, and
Shon Kwang Ho (eds.), Marcel Dekker, Inc., New York, 2000.

\item{[456]} {\sc Mackey, M.T., Mellon, P.}, a) {\it The quasi-invertible manifold of a
$JB^{\ast}$-triple}, Extracta Math. {\bf 14} (1999), {\it 1}, 51--55.

b) {\it Compact-like manifolds associated to $JB^{\ast}$-triples}, Manuscripta Math. {\bf 106} (2001), 203-212.

\item{[457]} {\sc Madariaga, S., P\'erez-Izquierdo, J.M.}, {\it Non-existence of coassociative quantized universal enveloping algebras of the traceless octonions}, Comm. Algebra (to appear).

\item{[458]} {\sc Magnus, T.}, a) {\it Geometries over
nondivision rings}, Univ. Virginia Dissertation, Charlottesvile, 1991.

b) {\it Faulkner Geometry}, Geom. Dedicata {\bf 59} (1996), 1--28.

\item{[459]} {\sc Magri, F., Pedroni, M., Zubelli, J.P.},
{\it On the geometry of Darboux
transformations for the {\rm KP} hierarchy and its connection
with the discrete {\rm KP} hierarchy},
Comm. Math. Phys. {\bf 188} (1997), 305--325.

\item{[460]} {\sc Mahmoud, N.H.}, {\it Bessel systems for
 Jordan algebras of rank 2 and 3}, J. Math. Anal.
Appl. {\bf 234} (1999), {\it 2}, 372--390.

\item{[461]} {\sc Majid, S.}, a) in Internat. J. Mod. Phys. {\bf A5} (1990), 1.

b) {\it Quantum groups and noncommutative geometry}, J. Math. Phys. {\bf 41} (2000), {\it 6},
3892--3942.

\item{[462]} {\sc Makarevich, B.O.}, a) {\it Open symmetric orbits of
the reductive groups in symmetric $R$-spaces} (in Russian), Mat. Sb.
(N.S.) {\bf 91} (1973), {\it 3}, 390--401; English translation: Math. USSR Sbornik
{\bf 20} (1973), {\it 3}, 406--418.

b) {\it Ideal points of semisimple Jordan algebras} (in Russian),
Mat. Zametki
{\bf 15} (1974), {\it 2}, 295--305.

c) {\it Jordan algebras and orbits in symmetric $R$-spaces} (in Russian),
Trudy
Moscov. Mat. Obshch. {\bf 39} (1979), 157--179, English translation:
Trans. Moscow
Math. Soc. {\bf 39} (1979), 169--193.

\item{[463]} {\sc Makhlouf, A.}, a) {\it Hom-alternative algebras and Hom-Jordan algebra}, arXiv:\break 0909.0326v1[math.RA] 2 Sep. 2009, Internat. Electronic J. Algebra, 2010, 13 pp.

b) {\it Paradigm of non-associative Hom-algebras and Hom-superalgebras}, arXiv:\break 1001.4240v1[math.RA] 24 Jan. 2010.

\item{[464]} {\sc Makhlouf, A., Silvestrov, S.D.}, a) {\it Hom-algebra structures}, J. Gen. Lie Theory Appl. {\bf 2} (2008), {\it 2}, 51-64.

b) {\it Hom-Lie admissible Hom-coalgebras and Hom-Hopf algebras}, Published as Chapter 17, pp. 189-206, S. Silvestrov, E. Paal, V. Abramov, A. Stolin (Eds.), {\it Generalized Lie theory in Mathematics, Physics and Beyond}, Springer-Verlag, Berlin, Heidelberg, 2008.

c) {\it Notes on formal deformations of Hom-associative and Hom-Lie algebras}, Forum Mathematicum {\bf 22} (2010), {\it 4}, 715-739; arXiv: 0712.3130v1[math.RA], 2007.

d) {\it Hom-Algebras and Hom-Coalgebras}, J. Algebra and its Appl. {\bf 9} (2010), {\it 4}, 553-589; arXiv:0811.0400[math.RA], 2008.

\item{[465]} {\sc Malcev, A.I.}, {\it Analytic loops} (in Russian), Mat. Sb. {\bf 36} (1955), {\it 3}, 569-576.

\item{[466]} {\sc Malley, J.D.}, {\it Statistical Applications of
Jordan algebras}, Lecture Notes in Statistics, {\bf 91}, Springer-Verlag, New-York, 1994.

\item{[467]} {\sc Manin, Yu.I.}, a) {\it Some remarks on Koszul
algebras and quantum groups},
Ann. Inst. Fourier {\bf 37} (1987), {\it 4}, 191--205.

b) Montreal Univ. preprint, CRM-1561 (1988).

c) {\it Multiparametric quantum deformation of the general linear supergroup},
Commun. Math. Phys. {\bf 123} (1989), {\it 1}, 163--175.

d) Invited paper presented at the XVIII Int.
Coll. on Group Theoretical Methods in Physics, Moscow, USSR,
June 4--9, 1990.

\item{[468]} {\sc Manin, Yu.I., Radul, A.O.},
{\it A supersymmetric extension of the {\rm KP}
hierarchy}, Comm. Math. Phys. {\bf 98} (1985), 65--77.

\item{[469]} {\sc Marchiafava, S.}, a) {\it On manifolds with generalized quaternionic structure} (in Italien), Rend. di Matematica (3) {\bf 3}, (1970), ser. VI, 1-17.

b) {\it Locally Grassmann quaternionic manifolds}, Atti Accad. Naz. Lincei, Rend. Cl. Sci. Fis. Mat. Nat. {\bf 57} (1974), 80-89.

c) {\it A report on almost quaternionic Hermitian manifolds}, in {\it Proc. Internat. Worshop Diff. Geom. and its Appl.} (Constan\c ta, Rom\^{a}nia, 1995), An. Univ. "Ovidius" {\bf 3} (1995), 55-63.

\item{[470]} {\sc Marshall, C.D.}, {\it Calculus on subcartesian spaces}, J. Differential Geom. {\bf 10} (1975), {\it 4}, 551-574.

\item{[471]} {\sc Martinelli, E.}, a) {\it Variet\`a a struttura quaternionale generalizzata}, Atti Accad. Naz. Lincei Rend. Cl. Sci. Fis. Mat. Natur. {\bf 26} (1959), 353-362.

b) {\it Modello metrico reale dello spazio proiettivo quaternionale}, Ann. Mat. Pura Appl. {\bf 49} (1960), 73-89.

\item{[472]} {\sc Martinez, A., Perez, J.D.},
{\it Real hypersurfaces in quaternionic projective space},
Ann. Mat. Pura Appl. {\bf 145} (1986), 355--384.

\item{[473]} {\sc Martinez, C.}, a) {\it Graded simple Jordan
algebras and superalgebras}, in {\it Recent prograss in algebra}
(Taejon/Seoul, 1997), 189--197, Contemp. Math.
{\bf 234}, Amer. Math. Soc. Providence, RI, 1999.

b) {\it Jordan superalgebras}, in {\it Nonassociative
 algebra and its applications} (Sao Paolo, 1998),
211--218, Lecture Notes Pures Appl. Math.
{\bf 211}, Dekker, New York, 2000.

\item{[474]} {\sc Martinez, C., Shestakov, I., Zelmanov, E.}, a) {\it Jordan superalgebras defined by brakets}, J. London Math. Soc. (2)
{\bf 64} (2001), {\it 2}, 357--368.

b) {\it Jordan bimodules over the superalgebras $P(n)$ and $Q(n)$} (submitted).

\item{[475]} {\sc Martinez, C., Zelmanov, E.}, a) {\it Simple finite-dimensional Jordan superalgebras of prime
characteristic}, J. Algebra {\bf 236} (2001), {\it 2}, 575--629.

b) {\it Unital bimodules over the simple Jordan superalgebras $D(t)$}, Trans. Amer. Math. Soc. {\bf 358} (2006), {\it 8}, 3637-3650.

c) {\it Representation theory of Jordan superalgebras I}, Trans. Amer. Math. Soc. {\bf 362} (2010), {\it 2}, 815-846.

d) {\it Representation theory of Jordan superalgebras II}, (preprint).

\item{[476]} {\sc Martinez, J.}, a) {\it On complete normed Jordan algebras}
(in Spanish), Ph.D. Thesis, Univ. of Granada, 1977.

b) $JV$-{\it algebras}, Math. Proc. Cambridge Philos. Soc.
{\bf 87} (1980), {\it 1}, 47--50.

\item{[477]} {\sc Massam, H., Neher, E.},
{\it Estimation and testing for lattice conditional
independence models on Euclidean Jordan algebras},
Ann. Stat. {\bf 26} (1998), 1051--1082.

\item{[478]} {\sc Mathiak, K.}, {\it Valuations of skew fields
and projective Hjelmslev
spaces}, Lecture Notes in Math. {\bf 1175},
Springer, Berlin--Heidelberg--New York, 1986.

\item{[479]} {\sc McCarthy, O.D.,  Strachan, I.A.B.},
{\it Unitized Jordan algebras and dispersionless KdV equations},
Kowalevski Workshop on Mathematical Methods of Regular
Dynamics (Leeds, 2000), J. Phys. A {\bf 34} (2001), {\it 11},
2435--2442.

\item{[480]} {\sc McCrimmon, K.}, a) {\it A general theory of Jordan rings},
Proc. Nat. Acad.
Sci. USA {\bf 56} (1966), 1072--1079.

b) {\it The Freudenthal-Springer-Tits constructions of exceptional
Jordan algebras}, Trans. Amer. Math. Soc. {\bf 139} (1969), 495--510.
{\bf 148} (1970), 293--314.

c) {\it Jordan algebras and their applications},
Bull. Amer. Math. Soc. {\bf 84}
(1978), {\it 4}, 612--627.

d) {\it The Russian revolution in Jordan algebras}, Algebras
Groups Geom. {\bf 1} (1984), 1--61.

e) {\it Jordan centroids}, Comm. Algebra {\bf 27} (1999),
{\it 2}, 933--954.

f) {\it Little Jordan Clifford algebras}, Comm. Algebra
{\bf 27} (1999), {\it 6}, 2701--2732.

g) {\it Properness, strictness, and nilness in Jordan systems},
Comm. Algebra {\bf 27} (1999),
{\it 7}, 3041--3066.

h){\it A taste of Jordan algebras}, Springer-Verlag, New York, 2004.

i) {\it Grassmann speciality of Jordan supersystems}, J. Algebra {\bf 278} (2004), {\it 1}, 3-31.

\item{[481]} {\sc Mellon, P.}, {\it Dual manifolds
of $JB^{\ast}$-triples of the form} $C(X, U)$, Proc.
Roy. Irish Acad. Sect. A
{\bf 93} (1993), {\it 1}, 27--42.

\item{[482]} {\sc Merlini Giuliani, M.L., Polcino Milies, C.}, {\it The smallest simple Moufang loop}, J. Algebra {\bf 320} (2008), {\it 3}, 961-979.

\item{[483]} {\sc Meyberg, K.}, a) {\it Jordan-Tripelsystemen
und die Koecher-Konstruktion von Lie-Algebren},
Habilitationsshrift, M\"{u}nchen,
1969.

b) {\it Jordan-Tripelsystemen und die
Koecher-Konstruktion von Lie-Algebren}, Math. Z. {\bf 115} (1970), 58--78.

c) {\it Zur Konstruktion von Lie-Algebren aus Jordan-Tripelsystemen},
Manuscripta Math. {\bf 3} (1970), 115--132.

d) {\it Lectures on algebras and triple systems},
Lecture Notes, Univ. of Virginia, Charlottesville, 1972.

e) {\it Trace formulas and derivations in simple Jordan pairs},
Comm. Algebra {\bf 12} (1980), 1311--1326.

\item{[484]} {\sc Micali, A., Campo, T.M.M., Costa E Silva, M.C., Ferreira, S.M.M.}, {\it Derivations dans les alg\`ebres gametiques II}, Linear Algebra Appl. {\bf 64} (1985), 175-182 (for part I see Comm. Algebra {\bf 12} (1984), 239-293).

\item{[485]} {\sc Micali, A., Quattara, M.}, {\it Sur les alg\`ebres de Jordan g\'en\'etiques II} in "Alg\`ebres g\'en\'etiques", ed. A. Micali, Hermann, Paris, 1987.

\item{[486]} {\sc Micali, A., Revoy, P.}, {\it Sur les alg\`ebres gametiques}, Proc. Edinburgh Math. Soc. {\bf 29} (1986), {\it 2}, 187-197.

\item{[487]} {\sc Mickelsson, J.}, a) {\it String quantization on group
manifolds and the
holomorphic geometry of Diff} $S^1/S^1$, Comm. Math. Phys.
{\bf 112} (1987), 643--661.

b) {\it Kac-Moody groups, topology of the Dirac determinant bundle
and fermionization}, Comm. Math. Phys. {\bf 110} (1987), 173--183.

c) {\it Current algebra representation for the $3+1$
dimensional Dirac-Yang-Mills theory},
Comm. Math. Phys. {\bf 117} (1988), {\it 2}, 261--277.

d) {\it Sato's universal Grassmannian and group extensions},
Lett. Math. Phys. {\bf 22} (1991), 131--139.

\item{[488]} {\sc Mickelsson, J., Rajeev, S.G.}, {\it Current algebras in $d+1$ dimensions and determinant bundles over infinite-dimensional
Grassmannians}, Comm. Math. Phys. {\bf 116} (1988), {\it 3}, 365--400.

\item{[489]} {\sc Mingo, J.A.}, {\it Jordan subalgebras of Banach algebras},
J. London Math. Soc.
 (2) {\bf 21} (1980), {\it 1}, 162--166.

\item{[490]} {\sc Mizuhara, A., Shima, H.}, {\it Invariant projectively
flat connections and its
applications}, Towards 100 years after Sophus Lie (Kazan, 1998),
Lobachevskii J. Math. {\bf 4}
(1999), 99--107 (electronic).

\item{[491]} {\sc Mocanu, P.}, a) {\it On the classification of the spaces $A_3$
with constant local Euclidean connection} (in Romanian),
Com. Acad. RPR {\bf 1} (1951), {\it 3},
239--243.

b) {\it Espaces \`{a} connexion affine constante \'{e}quivalents en grand avec l'espace euclidian},
Com. Acad. RPR {\bf 2} (1952), 389--395.

\item{[492]} {\sc Montaner, F.}, a) {\it Local PI theory for Jordan systems}, I, II, J. Algebra {\bf 216} (1999), {\it 1}, 302--327; {\bf 241} (2001), {\it 2}, 473--514.

b) {\it A note on moduli of inner ideals in Jordan systems}, Comm. Algebra {\bf 30} (2002), {\it 1}, 411--429.

c) {\it Homotope polynomial identities in prime Jordan systems}, J. Pure Appl. Algebra
{\bf 208} (2007), {\it 1}, 107--116.

d) {\it Algebras of quotients of Jordan algebras}, J. Algebra {\bf 323} (2010), {\it 10}, 2638-2670.

\item{[493]} {\sc Montaner, F., Paniello, I.}, {\it Algebras of quotients of nonsingular
Jordan algebras}, J. Algebra {\bf 312} (2007), {\it 2}, 963--984.

\item{[494]} {\sc Montaner, F., Toc\'on Barroso, M.}, {\it Local Lesieur-Croisot theory of Jordan algebras}, J. Algebra {\bf 301} (2006), {\it 1}, 256-273.

\item{[495]} {\sc Moreno Galindo, A.}, {\it Extending the norm from special Jordan triple systems to their associative envelops},
in {\it Banach algebras} '97 (Blaubeuren), 363--375, Walter de Gruyter, Berlin, 1998.

\item{[496]} {\sc Moore D'Ortona, C.}, {\it Homomorphisms of projective remoteness planes}, Geom. Dedicata {\bf 72} (1998), 111--122.

\item{[497]} {\sc Moore, G., Reshetikhin, N.}, in  Nucl. Phys. B{\bf 328}
(1989), 557.

\item{[498]} {\sc Morozova, E.A., Chentsov, N.N.}, a) {\it Elementary
Jordan logics} (in Russian), Inst. Prikl. Mat. Akad. Nauk SSSR,
Moscow, Preprint, {\it 113}, 1975.

b) {\it On the theorem of Jordan-von Neumann-Wigner} (in Russian), Inst.
Prikl. Mat. Akad. Nauk SSSR, Moscow, Preprint, {\it 129}, 1975.

c) {\it The structure of the family of stationary states of a quantum
Markov chain} (in Russian), Inst. Prikl. Mat. Akad. Nauk SSSR,
Moscow, Preprint {\it 130}, 1976.

d) {\it Projective Euclidean geometry and noncommutative probability theory} (in Russian),
in {\it Discrete geometry and topology} (Russian), Trudy Mat. Inst. Steklov {\bf 196}
(1991), 105--113.

\item{[499]} {\sc Motreanu, D.}, a) {\it Methods of differential topology with applications in the cohomology of differentiable manifolds} (in Romanian), Ph. D. Thesis, University "Al. I. Cuza" Ia\c si, 1978.

b) {\it Embeddings of $C^{\infty}$-subcartesian spaces}, An. \c Stiint. Univ. "Al. I. Cuza" Ia\c si, Sec. I Math. {\bf 25} (1979), {\it 1}, 65-70.

c) {\it The category of preringed manifolds}, Memoirs of the Scientific Sections, The Romanian Academy (IV) {\bf 2} (1979), 77-85.

\item{[500]} {\sc Moufang, R.}, a) {\it Alternativek\"{o}rper und der Satz von vollst\"{a}ndigen
Vierseit} $(D_9)$, Abh. Math. Sem. Univ. Hamburg {\bf 9} (1933), 207--222.

b) {\it Zur Struktur von Alternativk\"orpern}, Math. Ann. {\bf 110} (1935), 416-430.

\item{[501]} {\sc Mulase, M.}, a) {\it Complete integrability of the {\rm KP}
equation}, Advances in Math. {\bf 54} (1984), 57--66.

b) {\it Cohomological structure in soliton equations and Jacobian
varieties}, J. Dif\-ferential Geom. {\bf 19} (1984), 403--430.

c) {\it Solvability of the super {\rm KP} equation and a generalization of the
Birkhoff decomposition}, Invent. Math. {\bf 92} (1988), {\it 1}, 1--46.

d) {\it A correspondence between an infinite Grassmannian and arbitrary
vector bundles on algebraic curves}, talk at the Amer. Math. Soc.
Meeting at Worcester, Massachusetts, April 1989.

e){\it  A new super {\rm KP} system and a characterization of the Jacobians of arbitrary
algebraic supercurves}, J. Differential Geom. {\bf 34} (1991), 651--680.

\item{[502]} {\sc M\"{u}nzner, H.F.}, {\it Isoparametrische Hyperfl\"{a}chen in Sph\"{a}ren, I, II}, Math.
Ann. {\bf 251} (1980), 57--71; {\bf 256} (1981), 215--232.

\item{[503]} {\sc Musette, M., Conte, R., Verhoeven, C.},
{\it B\"{a}cklund transformation and nonlinear superposition
formula of the Kaup-Kupershmidt and Tzitzeica equations},
in {\it B\"{a}cklund and Darboux
Transformations. The Geometry of Solitons}, A. Coley, D. Levi, R. Milson, C. Rogers, P. Winternitz (eds.), CRM Proceedings \&
Lecture Notes, {\bf 29}, Amer. Math. Soc., Providence, RI, 2001.

\item{[504]} {\sc Nagano, T.}, {\it Transformation groups on
 compact symmetric spaces}, Trans.
Amer. Math. Soc. {\bf 118} (1965), 428--453.

\item{[505]} {\sc Nagaoka, Sh.}, {\it On special functions associated with
Jordan algebras} (especially on theta functions) (Japanese),
in {\it Infinite analysis} (Kyoto, 1985), 21--30, Surikaisekikenky,
Kokyuroku, {\bf 578} (1985).

\item{[506]} {\sc Naitoh, H.}, a) {\it Pseudo-Riemannian symmetric
$R$-spaces}, Osaka J.
Math. {\bf 21} (1984), 733--764.

b) {\it Grassmann geometries on compact symmetric spaces of general type},
J. Math. Soc.
Japan {\bf 50} (1998), 557--592.
 
c) {\it Symmetric submanifolds and Jordan triple systems}, in {\it Theory of Lie groups and manifolds}, Sophia Kokyuroku in Math. {\bf 45} (2002), Reiko Miyaoka \& Hiroshi Tamaru (eds.), 21-38.

\item{[507]} {\sc Nakajima, K.}, a) {\it On half-homogeneous hyperbolic
manifolds and Siegel domains}, J. Math. Kyoto Univ.
{\bf 24} (1984), {\it 1}, 1--26.

b) {\it A note on homogeneous hyperbolic manifolds},
J. Math. Kyoto Univ. {\bf 24} (1984), {\it 1}, 189--196.

\item{[508]} {\sc Nakatani, M., Noumi, M.}, $q$-{\it hypergeometric
systems arising from quantum Grassmannians}, Funkcial. Ekvac.
{\bf 41} (1998),
{\it 3}, 363--381.

\item{[509]} {\sc Nambu, Y.}, in Phys. Rev. D {\bf 7} (8) (1973), 2405--2412.

\item{[510]} {\sc Neal, M.}, {\it Spectrum preserving linear maps on $JBW^{\ast}$-triples}, Arch. Math. {\bf 79} (2002), 258-267.

\item{[511]} {\sc Neal, M., Russo, B.}, {\it State spaces of $JB^*$-triples},
Math. Ann.  {\bf 328} (2004), {\it 4}, 585--624.

\item{[512]} {\sc Neher, E.}, a) {\it Differentialgeometrische Aspekte
von Idempotenten in reellen Jordan--Algebren}, Diplomarbeit,
M\"{u}nster, 1975.

b) {\it On the classification of Lie and Jordan triple systems},
Comm. Algebra {\bf 13} (1985), 2615--2667.

c) {\it Jordan triple systems by the grid approach},
Lecture Notes in Math. {\bf 1280},
Springer-Verlag, Berlin--Heidelberg--New York, 1987.

d) {\it Lie groups, Hermitian symmetric spaces and Jordan pairs}, in
{\it Hadronic Mechanics and Nonpotential Interactions} (Cedar Fall,
Iowa, 1990), H.C. Myung (ed.), Nova Science Publishers, Inc.
Commack, N.Y., 1992, 243--258.

e) {\it A statistical model for Euclidean Jordan algebras},
in {\it Nonassociative algebra
and its applications} (Oviedo, 1993), 145--149,
Santos Gonzalez (ed.), Kluwer, Dordrecht, 1994.

f) {\it Lie algebras graded by $3$-graded root systems and Jordan pairs covered
by a grid}, Amer. J. Math. {\bf 118} (1996), 439--491.

g) {\it Transformation groups of the Andersson-Perlman cone},
J. Lie Theory {\bf 9} (1999), 203--213.

h) {\it $3$-graded syperalgebras and Jordan superpairs},
Abstracts of the Amer. Math. Soc. {\bf 21} (2000), {\it 1}, 63.

i) {\it Quadratic Jordan superpairs covered by grids}, 
J. Algebra  {\bf 269} (2003), {\it 1},\break 28--73.

j) {\it Extended affine Lie algebras and generalizations}, in
{\it Developments and trends in infinite-dimensional Lie theory},
K.-H. Neeb, A Pianzola (Eds.), Progress in Math. {\bf 288} (2011), Birkh\"auser, Boston, 53-126.

\item{[513]} {\sc Neher, E., Yoshii, Y.}, {\it Derivations and invariant forms of Jordan and alternative tori}, Trans. Amer. Math. Soc. {\bf 355} (2003), {\it 3}, 1079-1109.

\item{[514]} {\sc Nesterov, A.I., Sabinin, L.V.}, {\it Non-associative
geometry and discrete structure of spacetime}, Comm. Math. Univ. Carolinae
{\bf 41(2)} (2000), 347--357.

\item{[515]} {\sc Neumann, J. von}, {\it On an algebraic
generalization of the quantum
mechanics formalism} (Part. I), Mat. Sb. {\bf 1} (1936), 415--484.

\item{[516]} {\sc Nichita, F.F., Popovici, B.}, {\it Yang-Baxter operators from $(\mathbb{G}, \theta)$-Lie algebras}, Romanian Reports in Physics {\bf 63} (2011), {\it 3} (to appear), arXiv:1011.2072v1[math.QA] 9 Nov 2010.

\item{[517]} {\sc Nomizu, K.}, a) {\it Elie Cartan's work on
isoparametric hypersurfaces}, Symp. Pure Math. Diff. Geom.
{\bf 27}, II, Stanford, 1973.

b) {\it Some results in E. Cartan's theory of isoparametric families of
hypersurfaces}, Bull. Amer. Math. Soc.
{\bf 79} (1973), {\it 6}, 1184--1188.

\item{[518]} {\sc Nomura, T.}, a) {\it Manifold of primitive
idempotents in a Jordan-Hilbert algebra}, J. Math. Soc. Japan
{\bf 45} (1993), {\it 1}, 37--58.

b) {\it Grassmann manifold of a $JH$-algebra}, Ann. Global Anal.
Geom. {\bf 12} (1994), 237--260.

\item{[519]} {\sc Novikov, S.P.}, {\it Solitons and Geometry}, Cambridge
Univ. Press, 1995, 58 pp.

\item{[520]} {\sc Oehmke, R.H., Sandler, R.}, {\it The collineation groups of division ring planes.
I Jordan algebras}, Bull. Amer. Math. Soc. {\bf 69} (1963), 791--793.

\item{[521]} {\sc Ogievetsky, O., Schmidke, W.B., Wess, J., Zumino, B.}, {\it $q$-Deformed Poincar\'e algebra}, Comm. Math. Phys. {\bf 150} (1992), {\it 3}, 495-518.

\item{[522]} {\sc Ohta, Y., Satsuma, J., Takahashi, D., Tokihiro, T.},
{\it An elemen\-tary introduction to Sato theory}, Progress Theor.
Phys. Suppl. {\bf 94}
(1988), 210--241.

\item{[523]} {\sc Okubo, S.}, a) in {\it Symetries in Science VI: From
the Rotation Group to Quantum Algebras}, Bregenz, Austria,
August 2--7, 1992, Plenum
Press, New York, 1993.

b) in J. Math. Phys. {\bf 34} (1993), 3273.

c) in J. Math. Phys. {\bf 34} (1993), 3292.

d) University of Rochester Report UR-1312 (1993), in {\it Proc. of the
$15$th Mon\-treal-Rochester-Syracuse-Toronto Meeting for High Energy Theories}.

e) University of Rocester Report UR-1319 (1993).

f) University of Rochester Report UR-1334 (1993).

g) {\it Super-triple systems, normal and classical Yang-Baxter
equations}, in {\it Non-Associative Algebra and Its Applications}
(Oviedo, 1993), Kluwer
Academic Publ., Santos Gonzalez (ed.), 1994, 300--308.

\item{[524]} {\sc Okubo, S., Kamiya, N.}, {\it A construction of Jordan
superalgebras from triple
systems}, in {\it Proceedings of the $20$th Symposium on semigroups,
Languages and
their Related Fields} (Mito, 1996), 15--19, Res. Rep. Inf. Sci.
Toho Univ., Chiba, 1997.

\item{[525]} {\sc \O rsted, B., Zhang, G.}, {\it Generalized principal series 
representations and tube domains}, Duke Math. J. 
{\bf 78} (1995), 335--357. 

\item{[526]} {\sc Osserman, R.}, {\it Curvature in the eighties},
Amer. Math. Monthly
{\bf 97} (1990), 731--756.

\item{[527]} {\sc Ozeki, H., Takeuchi, M.},
{\it On some types of
isoparametric hypersurfaces in spheres, I, II}, T\^{o}hoku Math.
J. {\bf 27} (1975),
515--559; {\bf 28} (1976), 7--55.

\item{[528]} {\sc Paal, E.}, {\it Note on analytic Moufang loops}, Comment. Math. Univ. Carolin. {\bf 45} (2004), {\it 2}, 349-354.

\item{[529]} {\sc Parimala, R., Sridharan, R., Thakur, M.L.},
{\it Tits' construction of Jordan
algebras and $F_4$ bundles on the plane}, Compositio Math.
{\bf 119} (1999), {\it 1}, 13--40.

\item{[530]} {\sc Parimala, R., Svresh, V., Thakur, M.L.}, {\it Jordan algebras and $F_4$ bundles
over the affine plane}, J. Algebra {\bf 198} (1997), {\it 2}, 582--607.

\item{[531]} {\sc Pedersen, G., St\o rmer, E.}, {\it Traces on Jordan algebras}, Canad. J. Math.
{\bf 34} (1982), {\it 2}, 370--373.

\item{[532]} {\sc Pedroza, A. C., Vianna, J.D.M.}, {\it On the Jordan
algebra and the symmetric formulation of classical mechanics}, J. Phys. A:
Math. Gen. {\bf 13} (1980), {\it 3},\break 825--831.

\item{[533]} {\sc Peetre, J.}, {\it Hankel forms of arbitrary weight over a symmetric domain via the
transvectant}, Rocky Mountain J. Math. {\bf 24} (1994), {\it 3}, 1065--1085.

\item{[534]} {\sc Peng, L., Zhang, G.}, {\it Tensor products of holomorphic representations and bilinear
differential operators}, J. Funct. Anal. {\bf 210} (2004), 171--192.

\item{[535]} {\sc Peralta, A.M.}, a) {\it Little Grothendieck's theorem for real $JB^*$-triples}, Math. Z.
{\bf 237} (2001), {\it 3}, 531--545.

b) {\it New advances on the Grothendieck's inequality problem for bilinear forms of
$JB^*$-triples}, Math. Inequal. Appl. {\bf 8} (2005), {\it 1}, 7--21.

\item{[536]} {\sc Peresi, L.A.}, a) {\it A note on duplication of algebras}, Linear Algebra Appl. {\bf 104} (1988), 65-69.

b) {\it On derivations of baric algebras with prescribed automorphisms}, Linear Algebra Appl. {\bf 104} (1988), 71-74.

\item{[537]} {\sc Perez-Izquierdo, J.M.}, a) private communication (Dec. 1999).

b) {\it Algebras, hyperalgebras, nonassociative bialgebras, and loops}, Adv. Math. {\bf 208} (2007), 834-876.

\item{[538]} {\sc Persits, D.B.}, a) {\it Geometries over degenerate octaves} (in Russian), Dokl.,
Akad. Nauk SSSR {\bf 173} (1967), 1010--1013.

b) {\it Geometries over degenerate octaves and degenerate antioctaves} (in Russian), Trudy
Sem. Vektor. Tenzor. Anal. {\bf 15} (1970), 165--187.

\item{[539]} {\sc Petersson, H.P.}, a) {\it Composition algebras over a field with discrete valuation}, J. Algebra {\bf 29} (1974), 414--126.

b) {\it Max Koecher's work on Jordan algebras}, in {\it Jordan algebras} (Oberwolfach, 1992), 187--195, Kaup, McCrimmon, Petersson
(eds.), de Gruyter, Berlin, 1994.

c) {\it Grids and the arithmetics of Jordan pairs}, J. Algebra {\bf 213} (1999), {\it 1}, 77--128.

\item{[540]} {\sc Petersson, H.P., Thkur, M.}, {\it The \'etale Tits process of Jordan algebras revisited}, J. Algebra {\bf 273} (2004), {\it 1}, 88-107.

\item{[541]} {\sc Pevzner, M.}, {\it Repr\'{e}sentation de Weil associ\'{e}e \`{a} une repr\'{e}sentation d'une alg\`{e}bre de Jordan},
C.R.\ Acad.\ Sci.\ Paris S\'{e}r.\ I.\ Math.\ {\bf 328} (1999), {\it 6}, 463--468.

\item{[542]} {\sc Piacentini Cattaneo, G.M.}, {\it Gametic algebras of mutation and Jordan algebras}, Rend. Mat. {\bf 13} (1980), {\it 2}, 179-186.

\item{[543]} {\sc Pickert, G.}, {\it Projective Ebenen}, Springer-Verlag, Berlin, 1955.

\item{[544]} {\sc Pirio, L., Russo, F.}, a) {\it On projective varieties $n$-covered by irreducible curves of degree $\delta$}, Comment. Math. Helvetici (to appear).

b) {\it Extremal varieties 3-rationally connected by cubics, quadro-quadratic Cremona transformations and cubic Jordan algebras} (preprint).

\item{[545]} {\sc Polishchuk, A.}, {\it Classical Yang-Baxter equation and the $A_{\infty}$-constraint}, Adv. Math. {\bf 68} (2002), {\it 1}, 56-95.

\item{[546]} {\sc Popovici, I.}, a) {\it Spaces with constant connection} (in Romanian), in
{\it Connections on differentiable manifolds}, Univ. Bucure\c{s}ti, 1979, 153--204.

b) {\it Sur les espaces \`{a} connexion constante attach\'{e}s aux alg\`{e}bres de Jordan simple
et centrales}, Ann. \c{S}tiin\c{t}. Univ. ``Al. I. Cuza" Ia\c{s}i, Sec. I Mat. (N.S.) {\bf 25}
(1979), 365--380.

\item{[547]} {\sc Popovici, I., Iord\u anescu, R., Turtoi, A.}, {\it Simple Jordan and Lie gradings considered in differential
geometry}, (in Romanian), Edit. Acad. RPR, Bucure\c{s}ti, 1971.

\item{[548]} {\sc Popovici, I., Turtoi, A.}, {\it Propri\'{e}t\'{e}s m\'{e}triques des formes r\'{e}elles
de Jordan de type} $A_{II}$, Rev. Roumaine Math. Pures Appl. {\bf 15} (1970), 303--314.

\item{[549]} {\sc Prema, G., Kiranagi, B.S.}, {\it Lie algebra bundles defined by Jordan algebra bundles}, Bull. Math. Soc.
Sci. Math. Roumanie {\bf 31} (1987), 255--264.

\item{[550]} {\sc Pressley, A., Segal, G.}, {\it Loop groups}, Clarendon Press, Oxford, 1986.

\item{[551]} {\sc Pumpl\"un, S.}, a) {\it Symmetric composition algebras over algebraic varieties}, Manuscripta Math. {\bf 132} (2010), 307-333.

b) {\it Jordan algebras over algebraic varieties}, Comm. Algebra, {\bf 38} (2010), 714-751.

\item{[552]} {\sc Putter, P.S., Yood, B.}, {\it Banach Jordan $\ast$-algebras}, Proc. London Math. Soc.
(3) {\bf 41} (1980), {\it 1}, 21--44.

\item{[553]} {\sc Racine, M.L., Zelmanov, E.I.}, {\it Simple Jordan superalgebras with semisimple even part}, J. Algebra {\bf 270} (2003), {\it 2}, 374-444.

\item{[554]} {\sc Rad\'{o}, F.}, a) {\it Darstelung nicht-injective Kollineationen eines projectiven Rau\-mes durch verllgemeinerte semilineare Abbildungen},
Math.\ Z.\ {\bf 110} (1969),\break 153--170.

b) {\it Non-injective collineations on some sets in Desarguesian projective planes and extension of non-commutative valuations}, Aequations Math. {\bf 4} (1970), 307--321.

c) {\it Affine Barbilian structures}, J. Geom. {\bf 14} (1980), 75--102.

\item{[555]} {\sc Rajeev, S.G.}, {\it K\"{a}hler geometry and string
theory}, in {\it Proc. Superstring
Workshop at Colorado}, K. Mahantopa and P.G.O.~Freund (eds.), Plenum
Publishing, 1987.

\item{[556]} {\sc Rao, C.R., Mitra, S.K.}, {\it Generalized inverse of matrices and its applications}, Wiley, New York, 1971.

\item{[557]} {\sc Resnikoff, H.L.}, a) {\it On a class of linear differential equations for
automorphic forms in several complex variables}, Amer. J. Math. {\bf 95} (1973), 321--332.

b) {\it Differential geometry and color perception}, J. Math. Biol. {\bf 1} (1974), 97-131.

c) {\it On the geometry of color perception}, Lectures on Mathematics in the Life Sciences {\bf 7} (1974), 217-232.

d) {\it Theta functions for Jordan algebras}, Invent. Math. {\bf 31} (1975), 87--104.

e) {\it Automorphic forms of singular weight are singular forms}, Math. Ann. {\bf 215} (1975), 173--193.

f) {\it Theta functions for Jordan pairs}, Tagungsbericht 35/1979, Jordan-Algebren (18.8--25.81979), Oberwolfach.

\item{[558]} {\sc Rivilis, A.A.}, {\it Homogeneous locally symmetric
domains in homogeneous spaces associated with semisimple Jordan algebras} (in
Russian), Mat. Sb. (N.S.) {\bf 82} (1970), {\it 2}, 262--290.

\item{[559]} {\sc Robertson, A.G.}, {\it Decomposition of positive projections onto
Jordan algebras}, Proc. Amer. Math. Soc. {\bf 96} (1986), {\it 3}, 478--480.

\item{[560]} {\sc Robertson, A.G., Youngson, M.A.}, {\it Positive projections with contractive complements on Jordan
algebras}, J. London Math. Soc. (2) {\bf 25} (1982), {\it 2}, 365--374.

\item{[561]} {\sc Robinson, D.W., St\o rmer, E.}, {\it Lie and Jordan structure in operator algebras}, J. Austral. Math. Soc. Ser. A {\bf 29} (1980), 129--142.

\item{[562]} {\sc Rocha, Fihlo, T.M., Vianna, J.D.M.}, a) {\it Jordan algebras and field
theory}, Internat. J. Theoret. Phys. {\bf 26} (1987), {\it 10}, 951--955.

b) {\it Identification of a commutative Jordan algebra structure in constrained classical
Hamiltonian systems}, Hadronic J. {\bf 19} (1996), {\it 5}, 473--492.

\item{[563]} {\sc Rodr\'\i guez-Palacios, A.}, {\it Jordan structures in analysis}, in {\it Jordan
Algebras} (Oberwolfach, 1992), 97--186, W. Kaup, K. McCrimmon, H.P. Petersson (eds.), W. de Gruyter, Berlin, 1994.

\item{[564]} {\sc Rogers, C., Schief, W.K.}, {\it B\"acklund and Darboux transformations:
geometry and modern applications in soliton theory}, Cambridge Texts in Appl. Math. 
{\bf 30} (2002), 413~pp.

\item{[565]} {\sc Roos, G.}, a) {\it Volume of bounded symmetric domains and compactification
of Jordan triple systems}, in {\it Lie groups and Lie algebras}, Math. Appl.
{\bf 433}, Kluwer Acad. Publ., Dordrecht, 1998, 249--259.

b) Review of the monograph [364u]: Zbl.1073.17014, Zentralblatt f\"ur Mathematik, 2005.

c) {\it Weighted Bergman kernels and virtual Bergman kernels}, Sci. China Ser. A {\bf 48} (2005), suppl., 225-237.

d) {\it Exceptional symmetric domains. Symmetries in complex analysis}, Contemp. Math. {\bf 468} (2008), Amer. Math. Soc., Providence, RI, 157-189.

\item{[566]} {\sc Roos, G.J., Vigu\'{e}, P.}, {\it Syst\`{e}mes triples de Jordan et domaines sym\'{e}triques}, Collection Travaux en cours, Hermann, Paris, 1992.

\item{[567]} {\sc Roquette, P.}, a) {\it On the Galois cohomology of the projection linear group and its applications to the construction of generic splitting fields of algebras}, Math. Ann. {\bf 150} (1963), 411-439.

b) {\it Isomorphisms of generic splitting fields of simple algebras}, J. Reine Angew. Math. {\bf 214/215} (1964), 207-226.

\item{[568]} {\sc Rosso, M.}, in C.R. Acad. Sci. Paris {\bf 305} (1987), Serie I, 587.

\item{[569]} {\sc Rothaus, O.S.}, {\it Siegel domains and representations of Jordan algebras}, Trans.
Amer. Math. Soc. {\bf 271} (1982), {\it 1}, 197--213.

\item{[570]} {\sc Rozenfeld, B. A.}, a) {\it Geometric interpretations of some Jordan algebras}, Publ. Inst. Math. (Beograd) (N.S.) {\bf 61} ({\bf 75}) (1997), 114--118.

b) {\it Geometry of planes over nonassociative algebras}, Acta Appl. Math. {\bf 50} (1998), {\it 1--2}, 103--110.

\item{[571]} {\sc Rozenfeld, B.A., Karpova, L.M.}, {\it Flag groups and reductive Lie groups}
(in Russian), Trudy Sem. Vektor. Tenzor. Anal. {\bf 13} (1966), 168--202.

\item{[572]} {\sc Rozenfeld, B.A., Zamahovski, M.P.}, {\it Simple and quasisimple Jordan algebras} (in Russian), Izv.
Vyss. Ucebn. Zaved. Matematika 1971, {\it 8} (111), 111--121; Letter to the editors: Izv.
Vyss. Ucebn. Zaved. Mathematika 1976, {\it 2} (165), 123.

\item{[573]} {\sc Ruegg, H.}, {\it Fermions and Jordan matrices}, in
{\it Group Theoretical Methods in Physics}
Proceedings, Varna 1987, Lecture Notes in Physics
{\bf 313}, H.D. Doebner and T. D. Palev (eds.), Springer-Verlag, Berlin, 1988, 575--586.

\item{[574]} {\sc Ruiz-Altaba, M.}, preprint UGVA-DPT-1993-10-838
and HEP-TH 9311069 (1993).

\item{[575]} {\sc Russo, B.}, {\it Structure of $JB^{\ast}$-triples}, in {\it Proceedings of the 1992 Oberwolfach conference on
Jordan algebras}, Walter de Gruyter, 209--280, Kaup, McCrimmon, Petersson (eds.), Walter de Gruyter,
Berlin, 1994.

\item{[576]} {\sc Russo, B., Dye, H.}, {\it A note on unitary operators in $C^{\ast}$-algebras}, Duke
Math. J. {\bf 33} (1966), 413--416.

\item{[577]} {\sc Rylands, L.J., Taylor, D.E.}, {\it Constructions for octonion and exceptional Jordan algebras}, in Special
issue dedicated to Dr. Jaap Seidel on the occasion of his 80th birthday (Oisterwijk, 1999), Des. Codes Cryptogr. {\bf 21}
(2000), {\it 1--3},\break 191--203.

\item{[578]} {\sc Sabinin, L.V.}, {\it Smooth quasigroups and loops}, Mathematics and its applications,
Kluwer, Dordrecht, {\bf 492} (1999).

\item{[579]} {\sc Sabinin, L.V., Mikheev, P.O.}, {\it Almost symmetric
and anti-symmetric spaces with affine connection} (in Russian),
Dokl. Akad. Nauk {\bf 337} (1994), {\it 4}, 454--455.

\item{[580]} {\sc Sabinin, L.V., Sabinina, L.L.}, {\it On the geometry
of trans-symmetric spaces}, in {\it Webs {\rm \&} Quasigroups}, Tver University,
1991, 117--122.

\item{[581]} {\sc Sabinina, L.L.}, {\it On geodesic loops of
trans-symmetric spaces}, Algebras Groups Geom. {\bf 12} (1995),
{\it 2}, 119--126.

\item{[582]} {\sc Safin, S.S., Sharipov, R.A.}, {\it B\"{a}cklund autotransformation for
the equation} $u_{xt}={\rm e}^u-{\rm e}^{-2u}$, Teoreticheskaya i Mat\'{e}maticheskaya Fizica
{\bf 95} (1993), 146--159; English translation: Theor. Math. Phys. {\bf 95} (1993),
462--470.

\item{[583]} {\sc  Sahi, S.}, a) {\it Explicit Hilbert spaces for certain unipotent representations},
Invent. Math. {\bf 110} (1992), 409--418.

b) {\it The Capelli identity and unitary representations},
Compositio Math. {\bf 81} (1992), 247--260.

c) {\it Unitary representations on the Shilov boundary of a symmetric tube
domain}, in {\it Representations of Groups and Algebras}, Contemp. Math. {\bf 145} (1993), 275--286,
Amer. Math. Soc., Providence.

d) {\it Jordan algebras and degenerate principal series},
J. reine angew. Math. {\bf 462} (1995), 1--18.

\item{[584]} {\sc Sa\"\i d, S.B.}, {\it The functional equation of zeta distributions associated with non-Euclidean
Jordan algebras}, Canad. J. Math. {\bf 58} (2006), 3--22.

\item{[585]} {\sc Sait\^o, S.}, a) {\it String vertex on the Grassmann
manifold}, Phys. Rev. D {\bf 37} (1988), {\it 4}, 990--995.

b) {\it Solitons and strings}, 
in {\it International Symposium on Spacetime Symmetries}, College Park, MD, USA, 24--28 May 1988, 436--455.

\item{[586]} {\sc Sakai, S.}, $C^{\ast}$-{\it algebras and $W^{\ast}$-algebras}, Ergebnisse der Mathematik
und ihrer Grenzgebiete {\bf 60}, Springer-Verlag, Berlin, 1972.

\item{[587]} {\sc Salinas, N.}, {\it The Grassmann manifold of a $C^{\ast}$-algebra and hermitian holomorphic bundles},
Oper. Theory. Adv. Appl. {\bf 28} (1988) 267--289.

\item{[588]} {\sc Saltman, D.}, {\it Norm polynomials and algebras}, J. Algebra {\bf 62} (1980), 333-345.

\item{[589]} {\sc S\'{a}nchez, C.}, a) $k$-{\it symmetric submanifolds
of} ${\bf R}^N$, Math. Ann. {\bf 270} (1985), 297--316.

b) {\it The index number of an $R$-space: an extension of a result of M. Takeuchi},
Proc. Amer. Math. Soc. {\bf 125} (1997), 893--900.

\item{[590]} {\sc S\'{a}nchez, C., Dal Lago, W., Garcia, A., Hulett, E.}, {\it On some properties
which characterize symmetric and general $R$-spaces}, Diff. Geom. Appl. {\bf 7}
(1997), 291--302.

\item{[591]} {\sc Sarath, B., Varadarajan, K.},
{\it Fundamental theorem of projective geometry}, Comm. Algebra
{\bf 12} (1984), {\it 7--8},
937--952.

\item{[592]} {\sc Satake, I.}, {\it Algebraic structures of symmetric domains}, Iwanami Shoten and Princeton
Univ. Press, 1980.

\item{[593]} {\sc Satake, I., Faraut, J.}, {\it The functional equation of zeta distributions associated with formally real
Jordan algebras}, Tohoku Math. J. {\bf 36} (1984), {\it 3}, 469--482.

\item{[594]} {\sc Sato, H.T.}, preprint OS-GE-40-93 and HEP-TH 9312174 (1993).

\item{[595]} {\sc Sato, M.}, a) {\it Soliton equations as dynamical
systems on an infinite-dimensional Grassmann manifold}, RIMS, Kokyuroku,
{\bf 439} (1981), 30--46.

b) {\it Universal Grassmann manifolds and integrable systems}, OATE Conference (Bu\c steni - Rom\^{a}nia, 1983).

\item{[596]} {\sc Sato, M., Sato. Y.}, {\it Soliton equations as
dynamical systems on an
infinite-dimensional Grassmann manifold}, in {\it Proc. US-Japan Seminar on Nonlinear
Partial Differential Equations} (Tokyo, 1982), P.D. Lax and H. Fujita (eds.),
259--271, North-Holland, Amsterdam--New York, 1983.

\item{[597]} {\sc Sauter, J.}, {\it Randstrukturen beschr\"{a}nkter symmetrischer Gebiete}, Ph.D.
Dissertation, Universit\"{a}t T\"{u}bingen, 1995.

\item{[598]} {\sc Schafer, R.D.}, 

a) {\it Structure of genetic algebras}, Amer. J. Math. {\bf 71} (1949), 121-135.

b) {\it An introduction to non-associative algebras}, Academic Press., New York, 1966; corrected reprint of the 1966 original, Dover Publications, Inc., New York, 1995.

\item{[599]} {\sc Schief, W.K.}, {\it Self-dual Einstein spaces via a permutability theorem for the Tzitzeica equation}, Phys. Lett.
A {\bf 223} (1996), 55--62.

\item{[600]} {\sc Schray, J.}, a) {\it Octonions and supersymmetry}, Ph.D. Thesis, Dept. of Physics, Oregon State University, 1994.

b) {\it The general classical solution of the syperparticle}, Class. Quant. Grav. {\bf 13} (1996), 27.

\item{[601]} {\sc Schroeck, F.E.}, {\it Quantum fields for reproducing kernel Hilbert spaces}, Rep.
Math. Phys. {\bf 26} (1988), {\it 2}, 197--210.

\item{[602]} {\sc Schwarz, A.}, a) {\it Fermionic string and universal moduli space}, Nucl.
Phys. B {\bf 317} (1989), 323--343.

b) {\it Grassmannian and string theory}, Comm. Math. Phys. {\bf 199} (1998), 1--24.

\item{[603]} {\sc Schwarz, J.H.}, in {\it Superstring and supergravity,
Scottish Universities Summer School in Physics}, SUSSP, Edinburgh, 1986.

\item{[604]} {\sc Scott, A.C}, {\it Davydov solitons in polypeptides}, Philos. Trans. Roy. Soc. London, Ser. A {\bf 315} (1985), 423-436.

\item{[605]} {\sc Segal G., Wilson, G.}, {\it Loop groups and equations
of KdV type}, Publ. IHES {\bf 61} (1985), 5--65.

\item{[606]} {\sc Seier, W.}, {\it Quasi-translations of Desarguesian affine Hjelmslev planes},
Math. Z. {\bf 177} (1981), {\it 2}, 181--186.

\item{[607]} {\sc Seligman, G.}, {\it Rational methods in Lie algebras},
Marcel Dekker, New York, 1976.

\item{[608]} {\sc Semenoff, G.W.}, {\it Nonassociative electric fields
in chiral gauge theory: an explicit construction}, Phys. Rev. Lett. {\bf 60} (1988),
{\it 8}, 680--683.

\item{[609]} {\sc Semyanistii, V.I.}, {\it Symmetric domains and Jordan algebras} (in Russian), Dokl. Akad.
Nauk SSSR {\bf 190} (1970), 788--791; Soviet. Math. Dokl. {\bf 11} (1970), 215--218.

\item{[610]} {\sc Serre, J.P.}, {\it Algebres de Lie semi-simples complexes}, W.A. Benjamin Inc., New York, 1966.

\item{[611]} {\sc Shandra, I.G.}, {\it The spaces $V_n(K)$ and Jordan
algebras} (in Russian), Dedicated to the memory of Lobachevskii,
{\it 1} (in Russian), 99--104, Kazan, Gos. Univ., Kazan, 1992.

\item{[612]} {\sc Sharipov, R.A., Yamilov, R.I.}, {\it B\"{a}cklund transformation and the
construction of the integrable boundary-value problem for the equation} $u_{xx}-
u_{tt} ={\rm e}^u -{\rm e}^{-2u}$, Trans. Inst. Math. BNC, URO AN SSSR, Ufa {\bf 66} (1991), 1--9.

\item{[613]} {\sc Sharpe, R.W.}, {\it Differential geometry -- Cartan's generalization
of Klein's Erlangen Program}, Springer-Verlag, New York, 1997.

\item{[614]} {\sc Sheppard, B.}, {\it A Stone-Weierstrass theorem for $JB^*$-triples},
J. London Math. Soc. (2) {\bf 65} (2002), {\it 2}, 381--396.

\item{[615]} {\sc Shestakov, I.}, a) {\it Quantization of Poisson superalgebras and speciality of Jordan
Poisson superalgebras}, in {\it Nonassociative algebra and its applications} (Oviedo, 1993), 372--378, Santos Gonzalez (ed.), Kluwer, Dordrecht, 1994.

b) {\it Moufang loops and alternative algebras}, Proc. Amer. Math. Soc. {\bf 132} (2004), {\it 2}, 313-316.

\item{[616]} {\sc Shiota, T.}, {\it Characterization of Jacobian
varieties in terms of soliton equations}, Invent. Math. {\bf 83} (1986), 333--382.

\item{[617]} {\sc Shkrinyar, M.J., Kapor, D.V., Stoyanovich, S.D.}, {\it Classical and quantum approach to Davydov's soliton theory}, Phys. Rev. A {\bf 38} (1988), {\it 12}, 4402-4408.

\item{[618]} {\sc Shultz, F.W.}, {\it On normed Jordan algebras which are Banach dual spaces},
J. Funct. Anal. {\bf 31} (1979), 360--376.

\item{[619]} {\sc Siddiqui, A.A.}, {\it Self-adjointness in unitary isotopes of $JB^*$-algebras},
{\it Arch. Math.} {\bf 87} (2006), 350--358.

\item{[620]} {\sc Sierra, G.}, {\it An application of the theories of
Jordan algebras and Freudenthal triple systems to particles and strings}, Class. Quant.
Gravity {\bf 4} (1987), {\it 2}, 227--236.

\item{[621]} {\sc Sikorski, R.}, {\it Differential modules}, Colloq. Math. {\bf 24} (1971-72), {\it 1}, 45-79.

\item{[622]} {\sc Simon, B.}, {\it Trace ideals and their applications},
Cambridge, Cambridge University Press, 1979.

\item{[623]} {\sc Singh, M.K.}, {\it Extending the theory of linearization of a quadratic transformation in genetic algebra}, Indian J. Pure Appl. Math. {\bf 19} (1988), {\it 6}, 530-538.

\item{[624]} {\sc Singh, M.K., Singh, D.K.}, {\it On baric algebras}, Math. Ed. (Siwan) {\bf 20} (1986), {\it 2}, 54-55.

\item{[625]} {\sc Skornyakov, L.A.}, {\it Homomorphisms of projective planes and $T$-homomor\-phisms
of ternaries} (in Russian), Mat. Sb. {\bf 43} (1957), 285--294.

\item{[626]} {\sc Smirnov, F.A.}, in {\it Introduction to Quantum Group
and Integrable Massive Models of Quantum Field Theory}, World Scientific,
Singapore, 1990, 1.

\item{[627]} {\sc Smirnov, O.}, {\it Functorial version of Tits-Kantor-Koecher construction}, paper
presented at Columbia Univ. Meeting (March 16--18, 2001), Abstracts of the Amer. Math. Soc. {\bf 22} (2001),
{\it 2}, 288.

\item{[628]} {\sc Smith, J.W.}, {\it The de Rham Theorem for general spaces}, T\^{o}hoku Math. J. {\bf 18} (1966), {\it 2}, 115-137.

\item{[629]} {\sc Smith, R.R.}, {\it On non-unital Jordan-Banach algebras}, Math. Proc. Cambridge Philos. Soc.
{\bf 82} (1977), {\it 3}, 375--380.

\item{[630]} {\sc Spallek, K.}, {\it Differenzierbare R\"aume}, Math. Ann. {\bf 180} (1969), {\it 4}, 269-296.

\item{[631]} {\sc Spanicciati, R.}, {\it Near-Barbilian planes}, Geom. Dedicata {\bf 24} (1987),
{\it 3},\break 311--318.

\item{[632]} {\sc Springer, T.A.}, a) {\it The projective octave plane},
Nederl. Akad. Wetensch. Proc. A {\bf 63}~=~Indag. Math. {\bf 22} (1960),
74--101.

b) {\it On the geometric algebra of the octave planes}, Nederl. Akad.
Wetensch. Proc. A {\bf 65}~=~Indag. Math. {\bf 24} (1962), 451--468.

c) {\it Jordan algebras and algebraic groups}, Reprint of the 1973 edition, Classic
in Mathematics, Springer-Verlag, Berlin, 1998.

\item{[633]} {\sc Springer, T.A., Veldkamp, F.D.}, a) {\it Elliptic and hyperbolic octave planes}, Nederl.
Akad. Wetensch. Proc. A {\bf 65}~=~Indag. Math. {\bf 25} (1963), 413--451.

b) {\it On Hjelmslev-Moufang planes}, Math. Z. {\bf 107} (1968), 249--263.

c) {\it Octonions, Jordan algebras and exceptional groups}, Springer Monographs in Math., Springer-Verlag,
Berlin, 2000.

\item{[634]} {\sc Stacey, P.J.}, a) {\it Type $I_2 \, JBW$-algebras}, Quart J. Math., Oxford II-nd Ser.
{\bf 33} (1982), 115--127.

b) {\it Local and global splittings in the state space of a $JB$-algebra}, Math. Ann.
{\bf 256} (1981), 497--507.

c) {\it Locally orientable $JBW$-algebras of complex type}, Quart. J. Math., Oxford II-nd Ser. {\bf 33}
(1982), 247--251.

d) {\it The structure of type I $JBW$-algebras}, Math. Proc. Cambridge Philos. Soc. {\bf 90} (1981), 477--482.

\item{[635]} {\sc Stach\'{o}, L.L.}, a) {\it On the manifold of tripotents in $JB^*$-triples},
in {\it Contemporary Geometry and Topology and Related Topics} (Proc. Internat. Workshop Diff.
Geom. Appl., Cluj-Napoca, August 2007), Cluj University Press, 2008, 251--264.

b) {\it On metric geodesics of tripotents in $JB^*$-triples}, preprint, 2007, {\tt www.math.}\break {\tt u-szeged.hu/$^\sim$stacho}.

\item{[636]} {\sc Stach\'o, L.L., Werner, W.}, a) {\it Homogeneous complex Jordan manifolds}, preprint, http://www.math.u-szeged.hu/stacho, 2009.

b) {\it Jordan manifolds}, http://www.math.u-szeged.hu/stacho, in GEOMETRY, Exploratory Worshop Diff. Geom. Appl., Ia\c si, September 2009, Cluj University Press, D. Andrica \& S. Moroianu (Eds.) 2011, 135-147.

\item{[637]} {\sc Stach\'{o}, L.L., Zalar, B.}, {\it Uniform primeness of the Jordan algebra of symmetric operators}, Proc.
Amer. Math. Soc. {\bf 126} (1998), {\it 8}, 2241--2247.

\item{[638]} {\sc Steinberg, R.}, {\it Lectures on Chevalley groups}, Yale University Notes, 1967.

\item{[639]} {\sc Steinke, C.F.}, {\it Topological affine planes composed of two Desarguesian halfplanes and projective planes with trivial collineation group},
Arch. Math. (Basel) {\bf 44} (1985), {\it 5}, 472--480.

\item{[640]} {\sc Stiles, W.S.}, {\it A modified Helmholtz line-element in brightnesscolor space}, Proc. Phys. Soc. London {\bf 18} (1946), 41-65.

\item{[641]} {\sc St\o rmer, E.}, a) {\it Positive linear maps of operator
algebras}, Acta Math. {\bf 110} (1963), 233--278.

b) {\it Jordan algebras of type I}, in Acta. Math. {\bf 115} (1965), 165--184.

c) {\it On the Jordan structure of $C^{\ast}$-algebras}, Trans. Amer. Math. Soc.
{\bf 120} (1965), 438--447.

d) {\it Irreducible Jordan algebras of self-adjoint operators}, Trans. Amer. Math. Soc. {\bf 130}
(1968), 153--166.

e) {\it Conjugacy of involutive antiautomorphisms of von Neumann algebras}, Pre\-print Univ. Oslo,
{\it 4}/1984.

\item{[642]} {\sc Strachan, I.A.B.}, {\it Jordan manifolds and dispersionless KdV equations}, J. Nonlinear Math.
Phys. {\bf 7} (2000), {\it 4}, 495--510.

\item{[643]} {\sc Sugitani, T.}, {\it Zonal spherical functions on quantum Grassmann manifolds}, J. Math.
Sci. Univ. Tokyo {\bf 6} (1999), {\it 2}, 335--369.

\item{[644]} {\sc Suh, T.}, {\it On isomorphisms of little projective groups of Cayley planes}, Nederl. Acad.
Wetensch. Proc. A {\bf 65}~=~Indag. Math. {\bf 24} (1962), 320--339.

\item{[645]} {\sc Sushkevich, A.K.}, {\it $G$-groups theory} (in Russian), Harkov-Kiev, 1937.

\item{[646]} {\sc Suzuki, N.}, a) {\it Structure of the soliton space of
Witten's gauge field equations}, Proc. Japan. Acad. {\bf 60} A (1984),
141--144.

b) {\it General solution of Witten's gauge field equations}, Proc. Japan Acad.
{\bf 60} A (1984), 252--255.

c) {\it Witten's gauge equations and an infinite-dimensional Grassamann
manifold}, Comm. Math. Phys. {\bf 113} (1987), 155--172.

\item{[647]} {\sc Sverchkov, S.R.}, {\it The Lie algebra of skew-symmetric elements and its application in the theory of Jordan algebras}, Jordan Theory preprints No. 283 (19 Feb. 2010).

\item{[648]} {\sc Svinolupov, S.I.}, a) {\it Jordan algebras and
generalized Korteweg-de Vries equations} (in Russian), Teoret. Mat. Fiz.
{\bf 87} (1991), {\it 3}, 391--403; translation in Theoret. and Math.
Phys. {\bf 87} (1991), {\it 3}, 611--620.

b) {\it Generalized Schr\"{o}dinger equations and Jordan pairs},
Comm. Math. Phys. {\bf 143} (1992), {\it 3}, 559--575.

c) {\it Jordan algebras and integrable systems} (in Russian),
Funktsional. Anal i
Prilozhen. {\bf 27} (1993), {\it 4},  40--53, 96; translation in
Functional Anal. Appl. {\bf 27} (1993), {\it 4}, 257--265 (1994).

\item{[649]} {\sc Svinolupov, S.I., Sokolov, V.V.}, a) {\it A
generalization of a Theorem of Lie, and Jordan tops} (in Russian),
Mat. Zametki {\bf 53} (1993), {\it 2}, 122--125; translation in Math.
Notes {\bf 53} (1993), {\it 1--2}, 201--203.

b) {\it Deformations of triple Jordan systems and integrable equations} (in Russian),
Teoret. Mat. Fiz. {\bf 108} (1996), {\it 3}, 388--392.

\item{[650]} {\sc Svinolupov, S.I., Yamilov, R.I.}, {\it Explicit
self-transformations for multifield Schr\"{o}dinger equations and Jordan generalizations
of the Toda chain} (in Russian), Teoret. Mat. Fiz. {\bf 98} (1994), {\it 2}, 207--219;
translation in  Theoret. Math. Phys. {\bf 98} (1994), {\it 2}, 139--146.

\item{[651]} {\sc Svirezhev, Yu.M., Pasekov, V.P.}, {\it Mathematical genetics} (in Russian), Nauka, Moscow, 1980.

\item{[652]} {\sc Takagi, H., Takahashi, T.}, {\it On the principal
curvatures of homogeneous hypersurfaces in a sphere}, in {\it Differential Geometry in
honour of K. Yano}, Kinokunija, Tokyo, 1972, 469--481.

\item{[653]} {\sc Takasaki, K.}, a) {\it On the structure of solutions to
the self-dual-Yang-Mills equations}, Proc. Japan Acad. {\bf 59} A (1983), {\it 9},
418--421.

b) {\it A new approach to the self-dual Yang-Mills equations}, Comm. Math.
Phys. {\bf 94} (1984), 35--59.

c) {\it A new approach to the self-dual Yang-Mills equations. II}, Saitama Math.
J. {\bf 3} (1985), 11--40.

d) {\it Geometry of universal Grassmann manifold from algebraic point of view},
Rev. Math. Phys. {\bf 1} (1989), {\it 1}, 1--46.

\item{[654]} {\sc Takeno, S.}, a) {\it Vibron solitons in one-dimensional molecular crystals}, Progr. Theoret. Phys. {\bf 71} (1984), {\it 2}, 395-398.

b) {\it Dynamical problems in soliton systems}, Springer, Berlin-Heidelberg-New York, 1985.

c) {\it Vibron solitons and coherent polarization in an exactly tractable oscillator-Lattice system - Applications to solitons in $\alpha$ helical proteins activity}, Progr. Theoret. Phys. {\bf 73} (1985), {\it 4}, 853-873.

d) {\it Vibron solitons and soliton-induced infrared spectra of crystaline acetanilide}, Progr. Theoret. Phys. {\bf 75} (1986), {\it 1}, 1-14.

\item{[655]} {\sc Takeuchi, M.}, {\it Cell--decompositions and Morse
equalities on certain symmetric spaces}, J. Fac. Sci. Univ. Tokyo {\bf 12} (1965),
81--192.

\item{[656]} {\sc Takhtajan, L.A.}, a) in Adv. Stud. Pure Mat. {\bf
19} (1989), 435.

b) in {\it Introduction to Quantum Group and Integrable Massive Models of
Quantum Field Theory}, World Scientific, Singapore, 1990, 69.

\item{[657]} {\sc Tam\'{a}ssy, L., Binh, T.Q.}, {\it On weakly
symmetric and weakly projective symmetric Riemannian manifolds}, in
{\it Diff. geom. and its appl.} (Eger,
1989), 663--670, Colloq. Math. Soc. Janos Bolyai, {\bf 56}, North-Holland,
Amsterdam, 1992.

\item{[658]} {\sc Tartar, L.}, {\it Sur l'\'{e}tude directe d'\'{e}quations nonlin\'{e}aires
intervenant en th\'eorie du contr\^{o}le optimal}, J. Funct. Anal. {\bf 17} (1974), {\it 1}, 1--17.

\item{[659]} {\sc Teleman, A.-M., Teleman, K.},

a) {\it A combined B\"acklund-Tzitzeica theorem}, An. Univ. Bucure\c sti {\bf 48} (1999), {\it 2}, 197-202.

b) {\it Theorems of Tzitzeica-B\"{a}cklund type}, talk at {\it Internat. Workshop Diff.\ Geom.\ and its Appl.} (Timi\c{s}oara, September 2001).

c) {\it Tzitzeica-B\"acklund theorems}, Balkan J. Geom. and its Appl. {\bf 10} (2005), {\it 1}, 108-109.

\item{[660]} {\sc Teleman, K.}, a) {\it On a theorem of Borel-Lichnerowicz} (in Russian), Rev. Rou\-maine Math. Pures Appl.
{\bf 3} (1958), 107--115.

b) {\it On the spaces with constant affine and local Euclidean connection} (in Romanian), St. Cerc. Mat. {\bf 18} (1966), {\it 6}, 783--797.

c) {\it On the mathematical work of Gheorghe Tzitzeica}, Balkan J. Geom. and its Appl. {\bf 10} (2005), {\it 1}, 59-64.

d) {\it Homage to Gheorghe Tzitzeica}, Balkan J. Geom. and its Appl. {\bf 10} (2005), {\it 1}, p. X.

\item{[661]} {\sc Theoret, J.}, {\it Geometry of a cubic Jordan algebra}, Ph.D. Thesis, University of Virginia, 2001.

\item{[662]} {\sc Tillier, A.}, {\it Sur les idempotents primitifs
d'une alg\`ebre de Jordan formelle r\'{e}elle}, C.R. Acad. Sci.
Paris {\bf 280} (1975), 767--769.

\item{[663]} {\sc Timmesfeld, F.G.}, {\it Moufang planes and the groups $E_G^K$ and $SL_2(K), K$ a Cayley division algebra}, Forum Math. {\bf 6}
(1994), 209--231.

\item{[664]} {\sc Tits, J.}, a) {\it Sur la trialit\'{e} et certains groups s'en deduisant}, IHES,
Publ. Math. {\bf 2} (1959), 37--84.

b) {\it Une classe d'alg\`ebres de Lie en relation avec les alg\`ebres de Jordan},
Indag. Math. {\bf 24} (1962), 530--535.

c) {\it Classification of buildings of spherical type and Moufang polygons:
a survey}, in {\it Colloq. Internazionale sulle Teorie Combinatorie}, Rome, September, 1973.

d) {\it Buildings of Spherical Type and Finite BN-Pairs},
Springer-Verlag, Berlin, 1974.

e) {\it On buildings and their applications}, Internat.
Congress of Math. (Vanconver, August 1974), Montreal, 1975.

f) {\it Non-existence de certains polygones g\'{e}n\'{e}ralis\'{e}s}, I,
Invent. Math. {\bf 36} (1976), 275--284.

\item{[665]} {\sc Toda, M.}, {\it Pseudospherical surfaces via moving frames and loop groups}, Ph.D. Thesis,
University of Kansas, May 2000.

\item{[666]} {\sc Topping, D.}, {\it Jordan algebras of self-adjoint operators}, Memoirs Amer.
Math. Soc. {\bf 53}, 1965.

\item{[667]} {\sc T\"{o}rner, G., Veldkamp, F.D.}, {\it Literature on geometry over rings},
J. Geom. {\bf 42} (1991), 180--200.

\item{[668]} {\sc Torrence, E.}, {\it The coordinatization of a hexagonal-Barbilian
plane by a quadratic Jordan algebra}, Univ. Virginia Dissertation, Charlottesville, 1991.

\item{[669]} {\sc Treder, H.-J.}, {\it Quantum gravity and Jordan's
nonassociative algebras}, Found. Phys. Lett. (USA) {\bf 4} (1991), {\it 6}, 601--603.

\item{[670]} {\sc Tricerri, F., Vanhecke, L.}, a) {\it Cartan
hypersurfaces and reflections}, Nihonkai Math. J. {\bf 1} (1990), 203--208.

b) {\it Geometry of a class of non-symmetric harmonic manifolds}, in
Proc. Conf. Diff. Geom. and its Appl. (Opava, CSFR, 1992), 415--426.

\item{[671]} {\sc Truini, P.}, {\it Scalar manifolds and Jordan pairs in
supergravity}, Internat. J. Theor. Phys. {\bf 25} (1986), {\it 5}, 509--525.

\item{[672]} {\sc Truini, P., Biedenharn, L.C.}, a) {\it A comment on
the dynamics of} $M_3^8$, in {\it Proceedings of the IIIth Workshop on Lie
Admissible Formulations}, Hadronic, Boston, 1981.

b) {\it An ${\cal E}_6 \otimes U(1)$ invariant quantum mechanics for a
Jordan pair}, J. Math. Phys. {\bf 23} (1982), {\it 7}, 1327--1345.

\item{[673]} {\sc Truini, P., Varadarajan, V.S.}, a) in {\it Symmetries
in Science VI: From
the Rotation Group to Quantum Algebras}, Bregenz, Austria, August 2--7,
1992, Plenum Press, New York, 1993, 731.

b) in Rev. Math. Phys. {\bf 5} (1993), {\it 2}, 363.

\item{[674]} {\sc Trushina, M.}, {\it Modular representations of the Jordan superalgebras $D(t)$ and $K_3$}, J. Algebra {\bf 320} (2008), {\it 4}, 1327-1343.

\item{[675]} {\sc Tsao, L.C.}, {\it The rationality of the Fourier coefficients of certain
Eisenstein series on tube domains}, Compositio Math. {\bf 32} (1976), {\it 3}, 225--291.

\item{[676]} {\sc Turtoi, A.}, a) {\it Espaces de Wagner et formes r\'{e}elles de Jordan de type
D}, Rev. Roumaine Math. Pures Appl. {\bf 16} (1971), 111--119.

b) {\it On Jordan algebraic fibration} (in Romanian), Stud. Cerc. Mat. {\bf 37} (1985), 484--492.

\item{[677]} {\sc Tvalavadze, M.V.}, {\it Simple decompositions of simple special Jordan
superalgebras}, Comm. Algebra {\bf 35} (2007), {\it 6}, 2008--2034.

\item{[678]} {\sc Tzitzeica, G.}, {\it Sur une nouvelle class de surfaces},
C.R. Acad. Sci. Paris {\bf 144} (1907), 1257--1259.

\item{[679]} {\sc Ueno, K.}, a) {\it Analytic and algebraic aspects of the {\rm KP} hierarchy
from the viewpoint of the universal Grassmann manifold}, in {\it Infinite-dimensional
groups with applications}, Math. Sci. Res. Inst. Publ. {\bf 4}, 
Springer-Verlag, Berlin--Heidelberg--New York--Tokyo, 1985, 335--353.

b) {\it Super {\rm KP} hierarchy and super Grassmann manifold}, in {\it Topological and geo\-metrical
methods in field theory}, World Scientific, Singapore, 1986.

\item{[680]} {\sc Ueno, K., Yamada, H.}, a) {\it A supersymmetric
extension of nonlinear integrable systems}, in {\it Proc.\ Conf.\
Topological and Geometrical Methods Fields Theory}, J.
Westerholm and Hietarinta, (eds.), World Scientific, 1986, 59--72.

b) {\it Super {\rm KP} hierarchy and super Grassmann manifold}, Lett. Math.
Phys. {\bf 13} (1987), 59--68.

c) {\it Supersymmetric extension of
the {\rm KP} hierarchy and the universal super Grassmann manifold}, in
{\it Conformal Field Theory and Solvable Lattice Model}, Advanced
Studies in Pure Mathematics {\bf 16}, Miwa et al. (eds.)
Kinokyniya, 1988.

d) {\it Soliton solutions and bilinear residue
formula for the super {\rm KP} hierarchy}, in {\it Group Theoretical Methods
in Physics Proceedings} Varna 1987, Lecture Notes in Physics
{\bf 313}, H.D. Doebner and T.D. Palev (eds.), Springer-Verlag, Berlin,
1988, 176--184.

\item{[681]} {\sc Uohashi, K., Ohara A.}, {\it Jordan algebras and dual affine 
connections on symmetric cones}, Positivity {\bf 8} (2004), {\it 4}, 369--378.

\item{[682]} {\sc Upmeier, H.}, a) {\it Derivation algebras of $JB$-algebras}, Manuscripta Math. {\bf 30} (1979), {\it 2}, 199--214; Erratum: {\bf 32} (1980), {\it 1--2}, 211.

b) {\it Derivation of Jordan $C^{\ast}$-algebras}, Math. Scand. {\bf 46} (1980), {\it 2}, 251--264.

c) {\it Automorphism groups of Jordan $C^{\ast}$-algebras}, Math. Z. {\bf 176} (1981), 21--34.

d) {\it Jordan algebras and operator theory -- a survey}, in {\it Dilatation theory, Toeplitz operators, and other topics} (Timi\c{s}oara/Herculane, 1982), 395--407, Birkh\"{a}user, Basel--Boston, 1983.

e) {\it A holomorphic characterization of $C^{\ast}$-algebras}, in {\it Functional analysis, holomorphy and approximation theory},
Proc. Semin. (Rio de Janeiro, 1981), North-Holland, Math. Stud. {\bf 86} (1984), 427--467.

f) {\it Toeplitz operators on symmetric Siegel domains}, Math. Ann. {\bf 271} (1985), {\it 3}, 401--414.

g) {\it Symmetric Banach manifolds and Jordan $C^{\ast}$-algebras}, North-Holland, Ams\-ter\-dam--New York--Oxford, 1985.

h) {\it Jordan algebras and harmonic analysis on symmetric spaces}, Amer. J. Math. {\bf 108} (1986),
{\it 1}, 1--25.

i) {\it Multivariable Toeplitz operators, Jordan algebras and index theory} in {\it Operator algebras and mathematical physics} (Iowa
City, Iowa, 1985), 461--475, Contemp. Math. {\bf 62}, Amer. Math. Soc., Providence, 1987.

j) {\it Jordan algebras in analysis, operator theory, and quantum
mechanics}, CBMS Regional Conf. Ser. Math. {\bf 67},
Amer. Math. Soc., 1987.

k) {\it Jordan algebras, complex analysis and quantization}, in {\it Jordan algebras} (Oberwolfach, 1992), 301--317, Kaup, McCrimmon, Petersson (eds.), Walter de\break Gruyter, Berlin, 1994.

l) {\it Toeplitz-Berezin quantization and non-commutative differential geometry}, in {\it Linear operators} (Warsaw, 1994), 385--400, Banach Center Publ. {\bf 38}, Polish Acad. Sci., Warsaw, 1997.

m) {\it Multivariable Toeplitz operators and index theory}, Birkh\"{a}user, Basel, 2001.

n) Review of the monograph [364u]: MR 1979748, Mathematical Reviews, 2004.

o) Review of the book {\it Physical applications of homogeneous balls} by Yaakov Friedman with the
assistance of Tzvi Scarr, Progress in Mathematical Physics, {\bf 40}, Birkh\"auser, Boston,
2005, XXIV + 279~pp., MR 2104434.

\item{[683]} {\sc Vafa, C.}, {\it Operator formulation on Riemann
surfaces}, Phys. Lett. B {\bf 190} (1987), 47--54.

\item{[684]} {\sc Vagner, V.V.}, a) {\it The theory of compound manifolds} (in Russian), Trudy Sem. Vektor. Tenzor. Anal. {\bf 8} (1950), 11-72.

b) {\it Ternary algebraic operation in the theory of coordinate structure} (in Russian), Dokl. Akad. Nauk SSSR {\bf 81} (1951), {\it 6}, 981-984.

c) {\it On the theory of partial transformations} (in Russian), Dokl. Akad. Nauk SSSR {\bf 84} (1952), 653-656.

d) {\it Foundations of differential geometry and contemporary algebras} (in Russian), in {\it Proc. Fourth Union Congress of Math.}, Leningrad, 3-12 July 1961, Akad. Nauk SSSR, 1963, 17-29.

e) {\it On the algebraic theory of coordinates atlases} (in Russian), Trudy Sem. Vektor. Tenzor. Anal. {\bf 13} (1966), 510-563; {\bf 14} (1968), 229-281.

\item{[685]} {\sc Vanhecke, L.}, {\it Geometry in normal and tubular
neighborhood}, in  {\it Proc. Workshop on Differential
Geometry and Topology} (Cala Gonone (Sardinia) 1988), Rend. Sem. Fac. Sci. Univ.
Cagliari, Supplemento al vol. {\bf 58} (1988), 73--176.

\item{[686]} {\sc Varadarajan, V.S.}, {\it Geometry of quantum theory}, Springer-Verlag,
New York, 1985.

\item{[687]} {\sc Vasiu, A.}, {\it Hjelmslev-Barbilian structures}, Anal. Numer.
Theor. Approx. {\bf 273} (50) (1985), 73--77.

\item{[688]} {\sc Vel\'asquez, R., Felipe, R.}, a) {\it Quasi-Jordan algebras}, Comm. Algebra
{\bf 36} (2008), 1580--1602.

b) {\it Split dialgebras, split quasi-Jordan algebras and regular elements},
J. Algebra Appl. (to appear).
 
\item{[689]} {\sc Veldkamp, F.D.}, a) {\it Isomorphisms of little and middle
projective groups of octave planes}, Nederl.\ Akad.\ Wetensch.\ Proc.\ A
{\bf 67}~=~Indag. Math. {\bf 27} (1964), 280--289.

b) {\it Collineation groups in Hjelmslev-Moufang planes}, Math.\ Z. {\bf 108} (1968),\break
37--52.

c) {\it Unitary groups in projective octave planes}, Compositio. Math. {\bf 19} (1967), 213--258.

d) {\it Unitary groups in Hjelmslev-Moufang planes}, Math. Z. {\bf 108} (1969), 288--312.

e) {\it Projective planes over rings of stable rank $2$}, Geom. Dedicata {\bf 11} (1981),
{\it 5}, 285--308.

f) {\it Distant-preserving homomorphisms between projective ring planes}, Nederl. Acad. Wetensch. Proc. A {\bf 88}~=~Indag. Math. {\bf 47} (1985),
{\it 4}, 443--454.

g) {\it Incidence-preserving mappings between projective ring planes}, Nederl. Akad. Wetensch. Proc. A {\bf 88}~=~Indag. Math.
{\bf 47} (1985), {\it 4}, 455--459.

h) {\it Projective ring planes and their homomorphisms}, in {\it Rings and
Geometry}, R. Kaya et al. (eds.), Reidel, Dordrecht, 1985, 289--350.

i) {\it Projective Barbilian spaces I, II}, Resultate Math. {\bf 12} (1987),
{\it 1-2}, 222--240; {\bf 12} (1987), {\it 3-4}, 434--449.

j) $n$-{\it Barbilian domains, Results in Mathematics}, {\bf 23} (1993), 177--200.

k) {\it Geometry over rings}, in {\it Handbook of incidence geometry}, F.
Buekenhout (ed.), North-Holland, Amsterdam, 1995, 1033--1084.

\item{[690]} {\sc Verder, K.}, {\it Generalized elations}, Bull. London Math.
Soc. {\bf 18} (1986), 573--579.

\item{[691]} {\sc Vey, J.}, {\it Sur la division des domaines de Siegel}, Ann. Sci.
\'{E}cole Norm. Sup. {\bf 3} (1970), 479--506.

\item{[692]} {\sc Vdovin, V.V.}, a) {\it Simple projective planes}, Arch. Math.
(Basel) {\bf 47} (1986), {\it 5}, 469--480.

b) {\it Homomorphisms of freely generated projective planes}, Comm. Algebra
{\bf 16} (1988), {\it 11}, 2209--2230.

\item{[693]} {\sc Vigu\'{e}, J.P.}, {\it Les domaines born\'{e}s sym\'{e}triques et les
syst\'{e}mes triples de Jordan}, Math. Ann. {\bf 229} (1977), 223--231.

\item{[694]} {\sc Villena, A.R.}, {\it Derivations on Jordan-Banach algebras},
Studia Math.\ {\bf 118} (1996), 205--229.

\item{[695]} {\sc Vogan, D.A.jr.}, {\it Unitary representations of reductive 
Lie groups}, Princeton Univ. Press, Princeton, NJ, 1987, 308~pp.

\item{[696]} {\sc Voyt\v echovsk\'y, P.}, {\it Generators of nonassociative simple Moufang loops over finite prime fields}, J. Algebra {\bf 241} (2001), {\it 1}, 186-192.

\item{[697]} {\sc Vr\v{a}nceanu, G.}, a) {\it Sur la repr\'{e}sentation g\'{e}od\`{e}sique des espaces
de Riemann}, Rev. Roumaine Math. Pures Appl. {\bf 1} (1956), 147--165.

b) {\it Le\c{c}ons de g\'{e}om\'{e}trie differentielle}, vol. II, Edit. RPR, Bucure\c{s}ti, 1957.

c) {\it Propri\'{e}t\'{e}s globales des espaces \`{a} connexion affine}, Bull. Math. Soc.
Sci. Math. Phys. RPR {\bf 2} (1958), 475--478.

d) {\it Sopra gli spazi di Riemann a connessione constante}, Ann. Mat. Pura Appl. (4) {\bf 58} (1962).

\item{[698]} {\sc Wadati, M., Akutsu, Y.}, {\it From solitons to knots
and links}, Progress Theor. Phys. Suppl. {\bf 94} (1988), 1--41.

\item{[699]} {\sc Walcher, S.}, a) {\it A characterization of regular Jordan pairs and its
application to Riccati differential equations}, Comm. Algebra {\bf 14} (1986), {\it 10}, 1967--1978.

b) {\it \"{U}ber polynomiale, insbesondere Riccatische, Differentialgleichungen mit Fundamentall\"{o}sungen}, Math. Ann.
{\bf 275} (1986), {\it 2}, 269--280.

c) {\it Bernstein algebras which are Jordan algebras}, Arch. Math. (Basel) {\bf 50} (1988), {\it 3}, 218-222.

d) {\it Algebras and differential equations}, Hadronic Press, Palm Harbor, 1991.

e) {\it Algebraic structures and differential equations}, in {\it Jordan algebras} (Oberwolfach,
1992), 319--326, Kaup, McCrimmon, Petersson (eds.), Walter de Gruyter, Berlin, 1994.

\item{[700]} {\sc Wang, A., Yin, W., Zhang, L., Roos, G.}, {\it The K\"ahler-Einstein metric for some Hartogs domains over symmetric domains}, Sci. China Ser. A {\bf 49} (2006), {\it 6}, 1175-1210.

\item{[701]} {\sc Watson, K.S.}, {\it Manifolds of algebraic type over Jordan pairs}, Ph. D. Thesis, Flinders University of South Australia, 1978, 219 pp.

\item{[702]} {\sc Wegrzynowsky, S.}, {\it Representation of generalized affine symmetric spaces by
$S$-structures}, Demonstratio Math. {\bf 9} (1976), {\it 4}, 707--727.

\item{[703]} {\sc Weida, R.A.}, {\it Replaceable partial spreads and the construction
of non-Desar\-guesian translation planes}, Ph.D. Thesis, Univ. of Delaware, 1986.

\item{[704]} {\sc Werner, W.}, {\it $K$-theory and symmetric spaces}, lecture at the 4$^{{\rm th}}$ European Conference on Noncommutative Geometry (Bucharest, April 25-30, 2011).

\item{[705]} {\sc Wilson, G.}, a) {\it Algebraic curves and soliton
equations}, in {\it Geometry Today}, Birkhauser, Boston, 1985, 303--329.

b) {\it Collisions of Calogero-Moser particles and adelic Grassmannian} (with an
Appendix by I.G. Macdonald), Invent. Math. {\bf 133} (1998),
{\it 1}, 1--41.

\item{[706]} {\sc Wilson, J.B.}, {\it Decomposing $p$-groups via Jordan algebras}, J. Algebra {\bf 322} (2009), {\it 8}, 2642-2679.

\item{[707]} {\sc Witten, E.}, a) {\it Non-commutative geometry
and string field theory}, Nuclear Phys. B {\bf 268} (1986), {\it 2}, 253--294.

b) {\it Quantum field theory, Grassmannians and algebraic
curves}, Comm. Math. Phys. {\bf 113} (1988), 529--600.

c) in Nucl. Phys. B{\bf 330} (1990), 285.

\item{[708]} {\sc Wittstock G.}, {\it Matrixgeordnete R\"{a}ume
und $ C^{\ast}$-algebren in der Quantenme\-cha\-nik} (unpublished).

\item{[709]} {\sc Wolf, J.}, {\it Fine structure of Hermitian symmetric spaces}, in {\it Symmetric
Spaces}, Marcel Dekker (1972), 271--357.

\item{[710]} {\sc Woronowicz, S.}, a) {\it Compact matrix pseudogroups}, Comm. Math. Phys. {\bf 111} (1987), 613--665
(see also Lett. Math. Phys. {\bf 21} (1991), 35).

b) {\it Twisted $SU(2)$ group. An example of a non-commutative differential calculus}, Publ.
RIMS {\bf 23} (1987), {\it 1}, 117--181.

c) {\it Differential calculus on compact matrix pseudogroups $($quantum groups$)$}, Comm. Math. Phys. {\bf 122} (1989), 125--170.

\item{[711]} {\sc W\"orz-Busekros, A.}, a) {\it Algebras in genetics}, Lecture Notes in Biomathematics,
{\bf 36}, Springer-Verlag, New York, 1980.

b) {\it Relationship between genetic algebras and semicommutative matrices}, Linear Algebra Appl. {\bf 39} (1981), 111-123.

c) {\it Bernstein algebras}, Arch. Math. (Basel) {\bf 48} (1987), {\it 5}, 388-398.

d) {\it Further remarks on Bernstein algebras}, Proc. London Math. Soc. (3) {\bf 58} (1989), {\it 1}, 51-68.

\item{[712]} {\sc Wright, J.D.M.}, {\it Jordan $C^{\ast}$-algebras}, Michigan J. Math.
{\bf 24} (1977), 291--302.

\item{[713]} {\sc Wright, J.D.M., Youngson, M.A.}, a) {\it A Russo Dye theorem for Jordan $C^{\ast}$-algebras},
Functional Analysis; Surveys and Recent Results, North-Holland, 1977.

b) {\it On isometries of Jordan algebras}, J. London Math. Soc. {\bf 17} (1978), 339--344.

\item{[714]} {\sc Yamada, H.}, {\it Super Grassmann hierarchies.\
A multicomponent theory}, Hiroshima Math. J. {\bf 17} (1987), 377--394.

\item{[715]} {\sc Yamagishi, K.}, {\it The KP hierarchy and
extended Virasoro algebras}, Phys. Lett. B {\bf 205} (1988), {\it 4}, 466--470.

\item{[716]} {\sc Yamaguti, K.}, {\it An algebraic system arising from the metasymplectic
geometry. I, II, III}, Bull. Fac. School Ed. Hiroshima Univ. Part
II {\bf 4} (1981), 57--63; {\bf 5} (1982), 65--74; {\bf 9} (1986), 65--73.

\item{[717]} {\sc Yamaguti, K., Atsuko, O.}, {\it On representations of
Freudenthal Kantor triple systems} $U(\varepsilon, \delta)$, Bull. Fac. School Ed. Hiroshima Univ. Part II {\bf 7} (1984), 43--51.

\item{[718]} {\sc Yamaleev, R.}, {\it Fractional power of momenta and para-Grassmann extension of Pauli equation}, in {\it Procedings of the Internat. Conference on the Theory of the Electron} (Cuantitlan, Mexico, 1995), J. Keller, Z. Oziewicz (eds.), Adv. in Appl. Clifford Alg. {\bf 7} (5) (1997), 279--288.

\item{[719]} {\sc Yin, W., Lu, K., Roos, G.}, {\it New classes of domains with explicit Bergman kernel}, Sci. China Ser. A {\bf 47} (2004), {\it 3}, 352-371.

\item{[720]} {\sc Yokota, I.}, {\it Embeddings of projective spaces into elliptic projective Lie groups},
Proc. Japan Acad. {\bf 35} (1959), 281--283.

\item{[721]} {\sc Yomosa, S.}, {\it Toda-lattice in $\alpha$-helical proteins}, J. Phys. Soc. Japan {\bf 53} (1984), {\it 10}, 3692-3698.

\item{[722]} {\sc Yoshii, Y.}, a) {\it Jordan tori}, Math. Reports Acad. Sci. Canada {\bf 18} (4) (1996), 153--158.

b) {\it Jordan analogue of Laurent polynomial algebra}, in {\it Proceedings of the Internat. Conference on Jordan
structures} (Malaga, 1997), 191--197, Castellon Serrano, Cuenca Mira, Fern\'andez L\'opez,
Mart\'\i n Gonzalez (eds.), Malaga, 1999.

c) {\it Coordinate algebras of extended affine Lie algebras of type $A_1$},
J. Algebra {\bf 234} (2000), 128--168.

d) {\it Root-graded Lie algebras with compatible grading}, Comm. Algebra {\bf 29} (2001), {\it 8},
3365--3391.

\item{[723]} {\sc Youngson, M.A.}, a) {\it A Vidav theorem for Banach Jordan algebras}, Math. Proc.
Cambridge Philos. Soc. {\bf 84} (1978), 263--272.

b) {\it Hermitian operators on Banach Jordan algebras}, Proc. Edinburgh Math. Soc. {\bf 22} (1979), 93--104.

c) {\it Equivalent norms on Banach Jordan algebras}, Math.\ Proc.\ Cambridge Philos.\ Soc. {\bf 86} (1979),
{\it 2}, 261--270.

d) {\it Non-unital Banach Jordan algebras and $C^{\ast}$-triple systems}, Proc. Edinburgh Math. Soc.
{\bf 24} (1981), 19--29.

\item{[724]} {\sc Zabrodin, A.V.}, {\it Fermions on a Riemann
surface and the KP equation} (in
Russian), Theor. Math. Phys. {\bf 78} (1989), {\it 2}, 234--247.

\item{[725]} {\sc Zamahovski, M.P.}, {\it Models of non-Euclidean and quasi-non-Euclidean spaces in the form of orbits and of quasi-simple Jordan algebras of class D} (in Russian), Geometry and Topology {\it 2} (in Russian), 94--98, Leningrad, Gos. Ped. Inst., Leningrad, 1974.

\item{[726]} {\sc Zamolodchikov, A.B., Zamolodchikov, Al.B.},
in Ann.\ Phys.\ {\bf 120} (1979), 253.

\item{[727]} {\sc Zanon, D.}, {\it Superstring effective actions
and the central charge of the Virasoro algebra on a K\"{a}hler manifold},
Phys. Lett. B {\bf 186} (1987), {\it 3--4}, 309--312.

\item{[728]} {\sc Zarikov, F.M.}, {\it The Radon-Nikodym theorem for invariant states on $JBW^{\ast}$-algebras}
(in Russian), Izv. Akad. Nauk UzSSR, Ser. Fiz.-Mat. 1985, {\it 3}, 13--17.

\item{[729]} {\sc Zelmanov, E.}, a) {\it On the theory of Jordan algebras}, in
{\it Proc.\ Internat.\ Congress Mathematicians} (August 15--24, 1983,
Warzawa), P.W.N. \& North-Holland, 1984, 455--463.

b) {\it On prime Jordan triple systems, I, II, III}, Siberian Mat. J. {\bf 24} (1983),
23--27; {\bf 25} (1984), 50--61; {\bf 26} (1985), 71--82.

\item{[730]} {\sc Zhang, G.}, {\it Jordan algebras and generalized principal series}, Math. Ann.
{\bf 302} (1995), 773--786.

\item{[731]} {\sc Zhelyabin, V.N.}, {\it Jordan $($super$)$coalgebras and Lie $($super$)$coalgebras}
(in Russian), Sibirsk. Mat. Zh. {\bf 44} (2003), {\it 1}, 87--111; translation in
Siberian Math. J. {\bf 44} (2003), {\it 1}, 73--92.

\item{[732]} {\sc Zhelyabin, V.N., Shestakov, I.P.}, {\it Simple special Jordan superalgebras with an
associative even part}, Sibirsk. Mat. Zh. {\bf 45} (2004), {\it 5}, 1046--1072; translation in
Siberian Math. J. {\bf 45} (2004), {\it 5}, 860--882.

\item{[733]} {\sc Zhevlakov K.A., Slinko, A.M., Shestakov, I.P., 
Shirshov, A.I.}, {\it Rings that are nearly associative} (in
Russian), Nauka, Moscow, 1978; English translation: Academic
Press, New York, 1982.

\end{description}

\end{document}